\documentclass[11pt
]{article}

\usepackage{amssymb}
\usepackage{amsfonts,amstext,amsmath,amsthm,
latexsym,mathrsfs,amsbsy,mparhack}

\usepackage{graphicx}
\usepackage[dvipsnames]{xcolor}

\usepackage[bbgreekl]{mathbbol}
\usepackage{bold-extra}

\usepackage{marginnote}
\usepackage{multicol}

\usepackage{makeidx}
\makeindex

\usepackage[
english,greek]{babel}

\usepackage[b]{esvect}

\usepackage{tocloft}

\usepackage{doi}

\usepackage{hyperref}

\usepackage{tikz}

\usepackage{flowchart}
\usetikzlibrary{arrows}

\usepackage{bm}

\input{74dcrop.sty}
\input{74d.sty}

\hypersetup
{%
hypertex=true,
hypertexnames=true,
hidelinks,
pagebackref=true,
citecolor = Magenta,
colorlinks=true, 
}

\usepackage{cleveref}

\begin{document}

\selectlanguage{english}

\setcounter{tocdepth}{2}


\title{On the significance of parameters and the 
complexity level in  
the Choice and Collection axioms%
\,\thanks
{The research was carried out within the framework 
of the state assignment of the Institute
for Information Transmission Problems 
of the Russian Academy of Sciences, approved 
by the Ministry
of Education and Science of the Russian 
Federation.}
}

\author{Vladimir Kanovei\thanks
{Institute for Information Transmission Problems
(Kharkevich Institute) of Russian Academy
of Sciences (IITP), Moscow, Russia, 
{\tt kanovei@iitp.ru} } 
\and 
Vassily~Lyubetsky\thanks
{Institute for Information Transmission Problems
(Kharkevich Institute) of Russian Academy
of Sciences (IITP), Moscow, Russia, 
{\tt lyubetsk@iitp.ru} } 
}

\date{\today}
\maketitle


\begin{abstract}
We make use of generalized iterations of  
Jensen forcing  
to define a cardinal-preserving generic model   
of {\bf ZF} for any $\nn\ge1$ and 
each of the following four Choice hypotheses:
\bit
\item[(1)] 
$\xDC{\fp1\nn}\land\neg\xAC{\ip1{\nn+1}}\,;$ 

\item[(2)] 
$\xAC{\od}\land\xDC{\ip1{\nn+1}}\land
\neg\xAC{\fp1{\nn+1}}\,;$ 

\item[(3)] 
$\AC\land\xDC{\fp1{\nn}}\land\neg\xDC{\ip1{\nn+1}}\,;$ 

\item[(4)] 
$\AC\land\xDC{\ip1{\nn+1}}\land\neg\xDC{\fp1{\nn+1}}\,.$ 
\eit
Thus if\/ $\zf$ is consistent and $\nn\ge1$ 
then each of these four conjunctions (1)--(4) 
is consistent with $\zf.$ 

As for the second main result, let $\pao$ be the 2nd-order 
Peano arithmetic totally without the Comprehension schema $\CA$. 
For any $\nn\ge1$, we define a cardinal-preserving generic 
model of $\zf$, and a set\/ $M\sq\pws\om$ in this 
model, such that $\stk\om M$ satisfies  
\bit
\item[(5)] 
$\pao+\xAC{\is1\iy}+\xCA{\fs1{\nn+1}}+\neg\xCA{\fs1{\nn+2}}$. 
\eit
Thus\/ $\xCA{\fs1{\nn+1}}$ does not imply $\xCA{\fs1{\nn+2}}$ 
in $\pao$ 
even in the presence of the full parameter-free  Choice 
$\xAC{\is1\iy}\,.$ 

\indent

\end{abstract}

\footnotetext{MSC 03E25, 03E35, 03F35 (Primary), 03E15 (Secondary).}


\np

\setlength{\cftbeforesecskip}{3.7pt}
\setlength{\cftbeforesubsecskip}{0.2pt}

{\def\contentsname{
Contents
%
%
\vspace*{-0.50ex}
}\small\tableofcontents}

\np



\sekt{Introduction and preliminaries}
\las{ipre}

\parf{The main results}
\las{int}

This paper studies interrelations between different 
restricted forms of the axioms of countable 
\rit{independent choice} $\AC$ and 
\rit{dependent choice} $\DC$ in set theory, and 
of the Comprehension schema $\CA$ 
in second order arithmetic. 

The forms of the axiom of choice we consider 
will be distinguished by 
limiting the set or formula that specifies the choice, 
by one or another effective (\rit{lightface}) 
or classical (\rit{boldface}) projective class, resp.\ 
$\isp1n$, $\fsp1n$. 
The following theorem, our first main result, shows that 
the following 
three factors play a role in determining the strength 
of these forms of the axiom of choice, namely, 
the variant of the axiom ($\AC$ or $\DC$), the 
projective index $n$, 
as well as the assumption (boldface classes) or 
exclusion (lightface classes) of parameters in the
definitions of choice sets. 
Note that $\od$ = 
\rit{ordinal-definable} 
sets in \ref{mt12}. 

\bte
[1st main theorem]
\lam{mt1}
Assume that\/ $\nn\ge1$.
Then there exist 
cardinal-preserving generic extensions\/ 
$\rV_1$, $\rV_2$, $\rV_3$, $\rV_4$ of\/ $\rL$, in 
each of which\/ $\zf$ holds and 
the corresponding conjunction from 
the next list holds$:$ 
\ben
\nenu 
\itlb{mt11} 
$\xDC{\fp1\nn}\land\neg\xAC{\ip1{\nn+1}}\,;$ 

\itlb{mt12} 
$\xAC{\od}\land\xDC{\ip1{\nn+1}}\land\neg\xAC{\fp1{\nn+1}}\,;$ 

\itlb{mt13} 
$\AC\land\xDC{\fp1{\nn}}\land\neg\xDC{\ip1{\nn+1}}\,;$ 

\itlb{mt14} 
$\AC\land\xDC{\ip1{\nn+1}}\land\neg\xDC{\fp1{\nn+1}}\,.$ 
\een
Thus if\/ $\zf$ is consistent then each conjunction 
does not contradict\/ $\zf.$ 
\ete

The content of Theorem~\ref{mt1} 
is graphically presented in figures 
\ref{fig1}, \ref{fig2}, \ref{fig3}. 
The figures and the theorem will be commented upon in 
Sections~\ref{comF} and \ref{com1}.

 
\begin{figure*}[p]


\begin{tikzpicture}[>=latex',font={\sf \small}]

\mtho 
\def\smbwd{0cm}

\node (dc) at (0,-1.5) [draw, process, align=left,
minimum height=0cm] {$\!\DC\!$};

\node (dcod) at (1.5,0.5) [draw, process, align=left,
minimum width=\smbwd,
minimum height=0.5cm] {$\!\xDC{\od}\!$};

\node (dcb1) at (4.2,-1.6) [draw, process, align=left,
minimum width=\smbwd,
minimum height=0cm] {$\!\xDC{\fp1{n+1}}\!$};

\node (dc1) at (5.0,0.5) [draw, process, align=center,
minimum width=\smbwd,
minimum height=0cm] {$\!\xDC{\ip1{n+1}}\!$};

\node (dcb0) at (8,-1.5) [draw, process, align=left,
minimum width=\smbwd,
minimum height=0cm] {$\!\xDC{\fp1{n}}\!$};

\node (dc0) at (9.5,0.5) [draw, process, align=center,
minimum width=\smbwd,
minimum height=0cm] {$\!\xDC{\ip1{n}}\!$};

\node (acod) at (3,-3) [draw, process, align=left,
minimum width=\smbwd,
minimum height=0cm] {$\!\xAC{\od}\!$};

\node (ac1) at (7,-3) [draw, process, align=left,
minimum width=\smbwd,
minimum height=0cm] {$\!\xAC{\ip1{n+1}}\!$};

\node (ac0) at (11,-3) [draw, process, align=center,
minimum width=\smbwd,
minimum height=0cm] {$\!\xAC{\ip1{n}}\!$};

\node (ac) at (1.0,-3.5) [draw, process, align=left,
minimum width=\smbwd,
minimum height=0cm] {$\!\AC\!$};

\node (acb1) at (5.5,-5.5) [draw, process, align=left,
minimum width=\smbwd,
minimum height=0cm] {$\!\xAC{\fp1{n+1}}\!$};

\node (acb0) at (9.5,-5.5) [draw, process, align=center,
minimum width=\smbwd,
minimum height=0cm] {$\!\xAC{\fp1{n}}\!$};

\vyk{
\node (storage1) at (0,-8.5) [draw, storage,
minimum width=\smbwd,
minimum height=1cm] {STORE};

\node (process1) at (3,-8.5) [draw, process,
minimum width=\smbwd,
minimum height=1cm] {PROCESS};
}

\draw[->,thick,double] (dc) -- node[left]{1} (dcod);
\draw[->,thick,double] (dcod) -- node[left]{5} (dcb1);
\draw[->,thick,double] (dcb1) --  node[left]{8} (dc1);
\draw[->,thick,double] (dc1) -- node[right]{12} (dcb0);
\draw[->,thick,double] (dcb0) -- node[right]{15} (dc0);

\draw[->,thick,double] (dcod) -- node[left]{4} (acod);
\draw[->,thick,double] (dc1) -- node[right]{10} (ac1);
\draw[->,thick,double] (dc0) -- node[right]{17} (ac0);

\draw[->,thick,double] (dc) -- node[left]{2} (ac);

\draw[->,thick,double] (dcb1) -- node[below]
{9\hspace*{1ex}} (acb1);

\draw[->,thick,double] (dcb0) -- node[below]
{16\hspace*{2ex}} (acb0);

\draw[->,thick,double] (ac) -- node[below]{3} (acod);
\draw[->,thick,double] (acod) -- node[above]{6} (ac1);
\draw[->,thick,double] (ac1) -- node[above]{13} (ac0);
\draw[->,thick,double] (ac) -- node[below]{7} (acb1);
\draw[->,thick,double] (acb1) -- node[above]{14} (acb0);

\draw[->,thick,double] (acb0) -- node[right]{18} (ac0);
\draw[->,thick,double] (acb1) -- node[right]{11} (ac1);


\draw[->,thick,red] (ac) -- node[red,above,thick]
{{\textcolor{red}A}} (dc1);

\draw[->,thick,red] (acod) -- node[red,left]{B} (acb1);

\draw[->,thick,red] (dc1) -- node[red,above]
{\hspace*{2ex}C} (acb1);

\draw[->,thick,red] (dcb0) -- node[red,right]{D} (ac1);


\end{tikzpicture}
\vspace*{-1ex}

\caption{Provable $\imp$ and unprovable $\kra\to$ 
implications in $\zf$}
\label{fig1}
\end{figure*}


\begin{figure*}[p]


\begin{tikzpicture}[>=latex',font={\sf \small}]

\mtho 
\def\smbwd{0cm}

\hspace*{2mm}

\node (dc) at (0,-1.5) [draw, process, align=left,
minimum height=0cm] {$\!\DC\!$};

\node (dcod) at (1.5,0.5) [draw, process, align=left,
minimum width=\smbwd,
minimum height=0.5cm] {$\!\xDC{\od}\!$};

\node (dcb1) at (4.7,-1.6) [draw, process, align=left,
minimum width=\smbwd,
minimum height=0cm] {$\!\xDC{\fp1{n+1}}\!$};

\node (dc1) at (5.0,0.5) [draw, process, align=center,
minimum width=\smbwd,
minimum height=0cm] {$\!\xDC{\ip1{n+1}}\!$};

\phantom{
\node (fdc1) at (5.4,0.5) [draw, process, align=center,
minimum width=\smbwd,
minimum height=0cm] {$\!\DC(\ip1{n+1})\!$};
}

\node (dcb0) at (8,-1.5) [draw, process, align=left,
minimum width=\smbwd,
minimum height=0cm] {$\!\xDC{\fp1{n}}\!$};

\phantom{
\node (fdcb0) at (8.4,-1.5) [draw, process, align=left,
minimum width=\smbwd,
minimum height=0cm] {$\!\xDC{\fp1{n}}\!$};
}

\node (dc0) at (9.5,0.5) [draw, process, align=center,
minimum width=\smbwd,
minimum height=0cm] {$\!\xDC{\ip1{n}}\!$};

\node (ac) at (2,-3.5) [draw, process, align=left,
minimum width=\smbwd,
minimum height=0cm] {$\!\AC\!$};

\node (acb1) at (7,-3.5) [draw, process, align=left,
minimum width=\smbwd,
minimum height=0cm] {$\!\xAC{\fp1{n+1}}\!$};

\node (acb0) at (11,-3.5) [draw, process, align=center,
minimum width=\smbwd,
minimum height=0cm] {$\!\xAC{\fp1{n}}\!$};

\draw[->,thick,double] (dc) -- node[left]{1} (dcod);
\draw[->,thick,double] (dcod) -- node[left]{4} (dcb1);
\draw[->,thick,double] (dcb1) --  node[left]{8} (dc1);
\draw[->,thick,double] (dc1) -- node[left]{12} (dcb0);
\draw[->,thick,double] (dcb0) -- node[right]{15} (dc0);

\draw[->,thick,double] (dc) -- node[left]{2} (ac);

\draw[->,thick,double] (dcb1) -- node[below]{9} (acb1);
\draw[->,thick,double] (dcb0) -- node[right]{16} (acb0);

\draw[->,thick,double] (ac) -- node[above]{7} (acb1);
\draw[->,thick,double] (acb1) -- node[above]{14} (acb0);


\draw[->,thick,red] (ac) -- node[red,left]{A}  (dc1);
\draw[->,thick,red] (dc1) -- node[red,right]{C} (acb1);
\draw[->,thick,red] (fdcb0) -- node[red,above]{D} (fdc1);

\end{tikzpicture}
\vspace*{-1ex}

\caption{Provable $\imp$ and unprovable $\kra\to$ 
implications in $\zf$ + $\xAC\od$}
\label{fig2}
\end{figure*}



\begin{figure*}[p]


\begin{tikzpicture}[>=latex',font={\sf \small}]

\mtho 
\def\smbwd{0cm}
\hspace*{2mm}

\node (dc) at (-0.5,-1.6) [draw, process, align=left,
minimum height=0cm] {$\!\DC\!$};

\node (dcod) at (1.5,0.5) [draw, process, align=left,
minimum width=\smbwd,
minimum height=0.5cm] {$\!\xDC\od\!$};

\node (dcb1) at (3.5,-1.6) [draw, process, align=left,
minimum width=\smbwd,
minimum height=0cm] {$\!\xDC{\fp1{n+1}}\!$};

\phantom{
\node (fdcb1) at (3.8,-1.6) [draw, process, align=left,
minimum width=\smbwd,
minimum height=0cm] {$\!\xDC{\fp1{n+1}}\!$};
}

\node (dc1) at (5.5,0.5) [draw, process, align=center,
minimum width=\smbwd,
minimum height=0cm] {$\!\xDC{\ip1{n+1}}\!$};

\phantom{
\node (fdc1) at (5.8,0.5) [draw, process, align=center,
minimum width=\smbwd,
minimum height=0cm] {$\!\DC(\ip1{n+1})\!$};
}

\node (dcb0) at (7.8,-1.6) [draw, process, align=left,
minimum width=\smbwd,
minimum height=0cm] {$\!\xDC{\fp1{n}}\!$};

\phantom{
\node (fdcb0) at (8.1,-1.6) [draw, process, align=left,
minimum width=\smbwd,
minimum height=0cm] {$\!\DC(\fp1{n})\!$};
}

\node (dc0) at (9.7,0.5) [draw, process, align=center,
minimum width=\smbwd,
minimum height=0cm] {$\!\DC(\ip1{n})\!$};

\draw[->,thick,double] (dc) -- node[left]{1} (dcod);
\draw[->,thick,double] (dcod) -- node[left]{4} (dcb1);
\draw[->,thick,double] (dcb1) --  node[left]{8} (dc1);
\draw[->,thick,double] (dc1) -- node[left]{12} (dcb0);
\draw[->,thick,double] (dcb0) -- node[right]{15} (dc0);


\draw[->,thick,red] (fdcb0) -- node[red,above]{D} (fdc1);
\draw[->,thick,red] (fdc1) -- node[red,below]{\kra C} (fdcb1);

\end{tikzpicture}
\vspace*{-1ex}

\caption{Provable $\imp$ and unprovable $\kra\to$ implications 
in $\zf$ + $\AC$}
\label{fig3}
\end{figure*}

Our second main result is related to the Comprehension 
schema in 2nd order Peano arithmetic $\pad$. 
For the sake of brevity, let $\pao$ be the 2nd order 
arithmetic entirely without Comprehension, and 
let $\xCA K$ be the Comprehension schema 
$\sus x\,\kaz k\,(k\in x\eqv\vpi(k))$, limited to a 
given collection $K$ of formulas $\vpi$. 
Thus $\xCA{\fs1n}$, resp.\ $\xCA{\is1n}$ is 
the Comprehension schema for $\Sigma^1_n$ formulas 
\rit{with}, resp.\ \rit{without} parameters, 
and \rit{parameters} are  
formally free variables other than $k$ above.  
Note that choice principles $\AC$ and $\DC$ 
(single axioms in set theory) 
are naturally presented as 
axiom schemata in the language of $\pad$, see 
Section~\ref{com2}.

\bte
[2nd main theorem]
\lam{mt2}
Assume that\/ $\nn\ge1$. 
Then 
\vyk{
$:$
\ben
\renu
\itlb{mt21}
}%
there is a cardinal-preserving generic extension\/ 
of\/ $\rL$, and a set\/ $M\sq\pws\om$ in this 
extension, such that\/ $\rL\cap\pws\om\sq M$ and\/  
$\stk\om M$ models\/ 
$\pao+\xAC{\is1\iy}+\xCA{\fs1{\nn+1}}+
\neg\xCA{\fs1{\nn+2}}$. 
Thus\/ $\xCA{\fs1{\nn+1}}$ does not imply\/ 
$\xCA{\fs1{\nn+2}}$ 
even in the presence of\/ 
$\xAC{\is1\iy}\,.$ 
\ete

\bcor
\lam{mc2}
The full schema\/ $\AC$ is not 
finitely axiomatizable over $\pad+\xAC{\is1\iy}$,  
and the full schema\/ $\DC$ is not 
finitely axiomatizable over $\pad+\AC$.

The full schema\/ $\CA$ is not 
finitely axiomatizable over\/ 
$\pao + \xAC{\is1\iy}$. 
\qed
\ecor

\parf{Comments on figures}
\las{comF}

All \rit{unprovable} implications on the figures are such in 
virtue of Theorem~\ref{mt1}.

All \rit{provable} implications  
are rather self-evident, 
except for $\xDC K \imp\xAC K$ for different classes $K$ 
(arrows 2, 9, 10, 16, 17) -- which are well-known 
anyway, and the implication 
$\xDC{\ip1{n+1}}\imp\xDC{\fp1n}$ 
(arrow 12) proved by Lemma~\ref{23}\ref{235}. 

We consider the \rit{Baire space} $\cN=\bn,$ whose 
\kmar{cN}%
\index{Baire space $\cN$}%
\index{space!$\cN$, Baire space}%
\index{zNc@$\cN$, Baire space}%
points are called \rit{reals} in modern set theory, as 
well as product spaces of the form $\om^m\ti\cN^k$, 
$\om$ being discrete as usual. 
Sets in these spaces are called \rit{pointsets}.
See \cite{mDST} on lightface and boldface 
\rit{projective hierarchies} of pointsets.

The next definition presents the versions of $\AC$ and 
$\DC$ used here.

\bdf
\lam{02}
Let $K$ be any poinclass 
(a collection of pointsets). 
The following \rit{axioms}, or \rit{principles} 
are considered:
\bde
\item{$\xAC K$:} if $P\sq\om\ti\cN$, $P\in K$, 
\kmar{xAC K}%
\index{axiom!$\xAC K$}%
\index{zKAC@$\xAC K$}%
and $\dom P=\om$ then there is a map $x:\om\to\cN$ 
such that $\kaz k\,P(k,x(k))$.

\item{$\xDCm K$:} if $P\sq(\cN){}^2$, $P\in K$, 
\kmar{xDCm K}%
\index{axiom!$\xDCm K$}%
\index{zKDCm@$\DCm K$}%
and $\dom P=\cN,$ then there is a map 
$x:\om\to\cN$ 
such that $\kaz k\,P(x(k),x(k{+}1))$.

\item{$\xDC K$:} if $P\sq(\cN){}^2$, $P\in K$, 
\kmar{xDC K}%
\index{axiom!$\xDC K$}%
\index{zKDC@$\xDC K$}%
and $\dom P=\cN,$ then, 
{\bf for any $a\in\cN$}, there is 
$x:\om\to\cN$ 
such that $\kaz k\,P(x(k),x(k{+}1))$ and $x(0)=a$.

\item{$\xDCs K$:} if $P\sq\cN{}^2$, $P\in K$, 
\kmar{xDCs K}%
\index{axiom!$\xDCs K$}%
\index{zKDCs@$\xDCs K$}%
$\ran P\sq\dom P,$ then, for any $a\in\dom P$ there is 
$x:\om\to\cN$ 
such that $\kaz k\,P(x(k),x(k{+}1))$ and $x(0)=a$. 
\ede
Simply $\AC$, $\DC$, $\DCm\!$, $\DCs$ mean the 
case $K=\text{all sets})$. 
\edf

This definition can be used, for instance, for 
descriptive-set-theoretic, DST for brevity, 
pointclasses $K$ of the form $\is1n$ (lightface), 
$\fs1n$ (boldface), 
$\is1\iy=\bigcup_n\is1n$, 
same for $\Pi,\Da$ --- and then the corresponding 
axiom will be called \rit{a DST form of Countable Choice}. 
Non-descriptive forms are obtained \eg\ in cases 
$K=\od$ (all ordinal-definable pointsets), or 
$K=\ROD$ (all real-ordinal-definable pointsets), or 
$K$ = all pointsets of any kind. 

The axiom of (countable) dependent choices is 
known in several slightly different versions, in particular 
$\DC,\DCm,\DCs$ as above. 
Clearly the equivalence 
$\DC\eqv\DCm\eqv\DCs$ holds in $\zf$; this is 
why $\DCm$, the minimal form so to speak, is usually 
considered (and denoted by just $\DC$) in modern 
set theory. 
However $\xDC K $ as we define it turns out to be more 
convenient in the case of DST classes $K$, 
in particular, because, as far as we know, 
claim \ref{235} 
of Lemma~\ref{23} takes the form 
$\xDCm{\ip1{n+2}}\imp\xDCm{\fp1{n}}$ 
via an argument by Guzicki~\cite{guz}. 
This leaves the 
interrelations between $\xDCm{\ip1{n+1}}$ and 
$\xDCm{\fp1{n}}$ to be an open problem. 
This is why we prefer to consider $\DC$ rather 
than $\DCm$ (in the notation of Definition~\ref{02}) 
in this paper. 
The form $\DC$ was considered, by the way, in earlier papers 
\cite{aptm,guz,Kr}.

The next lemma proves some elementary connections. 
In particular, claim \ref{233} implies that there is 
no need whatsoever to consider $\Sigma$-limited forms 
of the choice principles as they can be substituted 
by 
$\Pi$-forms.

\ble
\label{23}
\ben
\renu
\itlb{231}
$\xDCs K\imp\xDC K\imp\xDCm K$ for any\/ $K\,;$

\itlb{232}
if\/ $K$ is any boldface or lightface projective class, 
or the class\/ $\od$, or the class of all sets, then$:$ 
$\xDCm K\imp\xAC K\,;$

\itlb{233}\msur
$\xAC{\ip1n}\eqv\xAC{\is1{n+1}}$, 
$\xAC{\fp1n}\eqv\xAC{\fs1{n+1}}$, 
and the same holds for\/ $\DC$ and $\DCs\,;$

\itlb{234}\msur
$\xDC{\ip1n}\eqv\xDC{\is1{n+1}}\eqv\xDCs{\ip1{n}}\eqv\xDCs{\is1{n+1}}\,,$\\[0.8ex]
$\xDC{\fp1n}\eqv\xDC{\fs1{n+1}}\eqv\xDCs{\fp1{n}}\eqv\xDCs{\fs1{n+1}}\,,$ 

\itlb{235}
$\xDC{\ip1{n+1}}\imp\xDC{\fp1{n}}$, and hence\/ 
$\xDC{\is1\iy}\eqv\xDC{\fs1\iy}\,;$ 

\itlb{236}
$\xDC\od\eqv\xDC\ROD$.

\itlb{237}
$\xDC{\fp11}$ holds in\/ $\zf$ and in\/ $\pad$
by the $\fp11$-uniformization theorem.
\een
\ele

\bpf
\ref{231} is trivial. 
\ref{232},\ref{237} 
are standard facts, see \eg\ \cite{aptm} or \cite{guz}. 

\ref{233}
As an example, 
to prove the lightface-\DC\ claim in \ref{233} 
(also a rather known fact as a whole), 
assume that $a\in\cN,$ 
and $P\sq\cN^2$ is a $\is1{n+1}$ set with 
$\dom P=\cN.$ 
Then $P(x,y)\eqv\sus z\,Q(x,y,z)$, 
where $Q\sq\cN^3$ is $\ip1n$. 
It remains to apply $\DC(\ip1n)$ to 
any $a'\in\cN$ with $(a')_0=a$ and the $\ip1n$ set 
$P'=\ens{\ang{x,y}\in\cN^2}{Q((x)_0,(y)_0,(y)_1)}$. 
(Recall that if $x\in\cN$ then $(x)_k\in\cN$ is defined 
\kmar{$(x)_k$}
\index{z(x)k@$(x)_k$}%
by $(x)_k(j)=x(2^k(2j+1)-1$, $\kaz j$.)

\ref{234}
The scheme of the proof of the first equivalence is 
\bce
$\xDC{\ip1n}\imp\xDC{\is1{n+1}}
\imp \xDCs{\ip1n}\imp\xDCs{\is1{n+1}}\imp\xDC{\ip1n}$. 
\ece
Here the 1st and 3rd implications follow from \ref{233}, 
so it remains to establish the 2nd one. 
Assume that $P\sq\cN^2$ is a $\ip1{n}$ set with 
$\ran P\sq \dom P$, and $a\in\dom P$. 
It suffices to apply $\xDC{\ip1n}$ to 
the $\id1{n+1}$ set 
$$
Q=\ens{\ang{x,y}\in\cN^2}
{P((x)_0,(x)_1)\imp 
\big[P((y)_0,(y)_1)\land (y)_0=(x)_1\big]}
$$
and any $a'\in\cN$ with $(a')_0=a$ and 
$P((a')_0,(a')_1)$.

\ref{235} 
is a bit trickier. 
Assume that $a\in\cN,$ 
and $P\sq\cN^2$ is a $\fp1{n}$ set with 
$\dom P=\cN.$ 
Then $P(x,y)\eqv S(x,y,p)$, 
where $S\sq\cN^3$ is lightface $\ip1n$, and $p\in\cN.$  
It remains to apply $\xDC{\ip1{n+1}}$ to the $\ip1{n+1}$ 
set 
$$
Q=\ens{\ang{x,y} 
\in\cN^2}
{(y)_1=(x)_1 \land 
\big[\sus z\,S((x)_0,z,(x)_1)
\imp S((x)_0,(y)_0,(x)_1)\big]}
$$
and any $a'\in\cN$ with $(a')_0=a\yi (a')_1=p$.
Finally \ref{236} is similar to \ref{235}.
\epf

\parf{Further comments on Theorem \ref{mt1}}
\las{com1}

It is quite clear that $\mathbf{AC}\imp\DC\imp\AC$. 
Studies in the early years of modern set theory by 
G\"odel, Cohen, Levy, Jensen, 
demonstrated that neither implication is reversible 
in $\zf$, $\mathbf{AC}$ is consistent with $\zf$, 
but $\AC$ is independent of $\zf$ and   
$\DC$ is independent of $\zf+\AC$ (Jensen \cite{jDC}). 

Furthermore Levy~\cite{levy2} demonstrated that the 
generic collapse of cardinals below $\aleph_\om$ 
(called the Levy collapse, see Solovay~\cite{sol}) 
results in a generic extension of $\rL$ 
in which $\xAC{\ip12}$ fails, 
which is the strongest possible failure since 
$\xAC{\fs12}$ is a theorem of $\zf$. 

Using rather similar arguments,  
Guzicki~\cite{guz} proved that the 
Levy-style generic collapse 
below $\aleph_{\omi}$ results in a generic 
extension of $\rL$ 
in which $\xAC{\fp12}$ fails, but 
$\AC\od$ holds, so that $\AC\od$ 
(for ordinal-definable sets) does not  
imply $\xAC{\fp12}$, let alone the full $\AC$. 
This can be compared with an opposite result for the 
\rit{dependent choice} axiom $\xDC{\fs1\iy}$, 
which is equivalent to 
the \paf\ form $\xDC{\is1\iy}$ by Lemma~\ref{23}. 

Recent research has shown that similar consistency 
results can be obtained via non-collapse forcing, 
and in some cases using the consistency of 2nd order 
Peano arithmetic $\pad$ as the blanket assumption 
(see Section~\ref{70}). 

Enayat~\cite{ena} used the finite-support infinite 
product of Jensen's minimal-$\id13$-real forcing 
\cite{jenmin} to define a non-collapse 
permutation model of $\zf$ with an infinite 
Dedekind-finite $\ip12$ set of reals, which 
easily yields the refutation of $\xAC{\ip12}$. 
Friedman \ea~\cite{jml19} used another generalization of 
Jensen's forcing to get a non-collapse 
model of $\zf+\AC$ in which $\xDC{\ip12}$ fails. 
(This result by a different method was also announced 
by Simpson~\cite{simp73}, but in fact never published, 
see notes in \cite[p.~4]{jml19} and \cite[p.~5]{hlong}.) 
Our own studies \cite{kl69,kl71} provided a 
Sacks-iterated, cardinal-preserving model of 
$\zf+\xAC\od$ in which $\xAC{\fp13}$ fails, and another 
such a model of $\zf$ in which $\xAC{\fs13}$ fails --- 
which is admittedly not the 
expected optimal failure of $\xAC{\fp12}$, resp., 
$\xAC{\fs12}$ in those cases. 

%

Some results related to \paf\ versions of the Separation 
and Replacement axiom schemata in $\zfc$ 
also are known from \cite{corrada,levyp,SS}. 

Our Theorem~\ref{mt1} substantially strengthens 
the above-mentioned results and 
maintains further clarification 
of the role of the projective level and 
parameters in the descriptive-theoretic axioms $\xAC K $ 
and $\xDC K $. 
Some parts of the theorem were published, in
Russian, in a technical report \cite{vin79}.

\parf{Comments on Theorem \ref{mt2}}
\las{com2}

Following \cite{aptm,Kr,simp} we define second order 
\kmar{pad}%
\index{2nd order arithmetic!$\pad$}%
\index{theory!$\pad$}%
\index{zPA2@$\pad$}%
arithmetic $\pad$ as a theory in the language 
$\lpad$ with two sorts of variables --- 
\kmar{lpad}%
\index{2nd order Peano arithmetic!language$\lpad$}%
\index{zPA2L@$\lpad$}%
for natural numbers and for sets of them. 
We use $j,k,m,n$ for variables over $\omega$ and 
$x,y,z$ for variables over $\pws\om$, reserving capital 
letters for subsets of $\pws\om$ and other sets. 
The axioms are as follows in 
\ref{pa1}, \ref{pa2}, \ref{pa3}, \ref{pa4}:
\ben
\nenu
\itlb{pa1}%
{\bf Peano's axioms} for numbers.

\itlb{pa2}%
The {\bf Induction} schema:  
$\Phi(0) \land \kaz k\,(\Phi(k)\imp\Phi(k+1))
\limp \kaz k\,\Phi(k)$, 
for every formula $\Phi(k)$ in $\lpad$, 
and in $\Phi(k)$ we allow parameters, 
\ie, free variables other than $k$. 
(We do not formulate Induction as one sentence here
because the Comprehension schema $\CA$ will not be    
always assumed in full generality by default.)  
  
\itlb{pa3}%
{\bf Extensionality} for sets of natural numbers.
  
\itlb{pa4}%
The {\bf Comprehension} schema $\CA$:
$\sus x \,\kaz k\,(k\in x\eqv\Phi(k))$, 
for every formula $\Phi$ in which $x$ 
\kmar{CA}%
\index{axiom!Comprehension, $\CA$}%
\index{Comprehension, $\CA$}%
\index{zCA@$\CA$}%
does not occur, and in $\Phi$ we allow parameters.
\een
$\pad$ is also known as $A_2^-$ (see \eg\ 
an early survey \cite{aptm}), 
as $Z_2$ (see \eg~Simpson~\cite{simp} and 
Friedman~\cite{HFuse81}), 
az $Z_2^-$ (in \cite{schindt} or elsewhere). 
The schema of Choice (see below) is 
not included in $\pad$ in this paper. 

Let $\pao$ to be the \ref{pa1}+\ref{pa2}+\ref{pa3} 
\kmar{pao}%
\index{theory!$\pao$}%
\index{zpa20@$\pao$}%
subtheory of $\pad$ (no Comprehension).

The principles $\AC$ and $\DC$ as in Definition~\ref{02} 
can be naturally reformulated as axiom schemata in the 
context of $\pad$.

\bdf
\lam{02pa}
Let $K$ be a type of formulas of $\lpad$,  
\eg\ $\is1n$ (lightface, real parameters not allowed), 
$\fs1n$ (boldface, real parameters allowed), 
$\is1\iy=\bigcup_n\is1n$, 
same for $\Pi$. 
The next axiom schemata in $\lpad$ are considered:
\bde
\item[$\xAC K$:]
\sloppy
$\kaz k\,\sus x\,\Phi(k,x)\imp
\sus x\,\kaz k\,\Phi(k,(x)_k)$, 
for every formula $\Phi$ in $K$, 
\kmar{xAC K}%
\index{axiom!$\xAC K$ in $\pad$}%
\index{zKAC@$\xAC K$ in $\pad$}%
where as usual 
$(x)_k=\ens{j}{2^k(2j+1)-1\in x}$. 

\item[$\xDC K$:]
$\kaz x\,\sus y\,\Phi(x,y)\imp
\kaz x\,\sus z\,\kaz k\,
\big((z)_0=x\land \Phi((z)_k,(z)_{k+1})\big)$, 
for any formula $\Phi$ in $K$. 
\kmar{xDC K}%
\index{axiom!$\xDC K$ in $\pad$}%
\index{zKDC@$\xDC K$ in $\pad$}%
\item[$\xCA K$:] 
$\sus x\,\kaz k\,\big(k\in x\eqv\Phi(k,(x)_k)\big)$, 
for any formula $\Phi$ in $K$. 
\kmar{xCA K}%
\index{axiom!$\xCA K$ in $\pad$}%
\index{zKCA@$\xCA K$ in $\pad$}%
\ede
Thus for instance 
$\xCA{\fs1\iy}$ is the full Comprehension schema $\CA$ 
whereas $\xCA{\is1\iy}$ is the \paf\ subschema of 
Comprehension. 
\edf

Discussing the structure and deductive properties of 
$\pad$, Kreisel 
\cite[\S\,III, page 366]{Kr} wrote that 
the selection of subsystems ``is a central problem''. 
In particular, Kreisel notes, that 
\begin{quote}
[...] if one is convinced of the significance 
of something like a given axiom schema, 
it is natural to study details, such as the effect 
of parameters. 
\end{quote}
Recall that 
\rit{parameters} in this context are 
free variables in axiom schemata that are not 
explicitly specified --- in $\pad$, $\ZFC$, 
and other similar theories. 
Thus the most obvious way to study 
``the effect of parameters'' is to compare the strength 
of a given axiom schema with its parameter-free 
subschema, \eg\ $\xCA{\fs1n}$ \vs\ $\xCA{\is1n}$. 
Working in this direction, it is established in 
our recent papers \cite{kl69,kl71} that
1) there is a cardinal-preserving generic extension  
of\/ $\rL$, and a set\/ $M\sq\pws\om$ in this extension, 
such that\/ $\pws\om\cap\rL\sq M$ and\/ $M$ is a model of\/ 
$\pao + \xCA{\is1\iy}+ \xCA{\fs12}+\neg\,\xCA{\fs14}$, and 
2) if $\pad$ is consistent then   
$\pao + \xCA{\is1\iy}+ \xCA{\fs12}$ does not prove 
$\xCA{\fs14}$.

\parf{Brief review of the forcing notions involved}
\las{stru}

The models we built to prove Theorems \ref{mt1} and 
\ref{mt2} have their own interesting history. 
It starts with forcing by perfect sets, or 
\rit{Sacks forcing} \cite{sacks} which produces 
generic reals of minimal degree. 
Further studies discovered and studied 
countable-support \rit{iterations} of Sacks forcing 
\cite{balav,amill:map,groapp}, and then 
\rit{generalized iterations} \cite{grj}, that is, 
iterations along any partial order $I$ in the ground model 
$M$. 
In this case, a generic \dd Iarray $\w: I\to\can{}=\dn$ is 
added, so that 
the 
structure of $I$ as a poset in $M$ is reflected in the structure 
of \dd Mdegrees of reals in the extension $M[\w]$. 
This connection can be used in coding by degrees of 
constructibility, see \eg\  
\cite[p.~143]{matsur}. 

As another application of generalized Sacks forcing iterations, 
in combination with the technique of ``symmetric'' generic 
extensions, cardinal-preserving generic 
models have been constructed with 
analytically definable violations of certain forms of the axiom 
of choice in the domain of reals. 

\bpri
\lam{stru1}
Taking $\rL$ as the ground model and 
$\tup={\omi}\lom\bez\La$ in $\rL$ 
\kmar{tup}%
\index{zI@$\tup$, tuples of ctble ordinals}%
\index{tuples of ctble ordinals, $\tup$}%
(all non-$\pu$ tuples of countable ordinals) 
leads to an \dd\tup iterated Sacks generic array 
$\w\in\can\tup$ of reals as above. 
Let $\Om$ consist of all countable well-founded 
(\ie, no infinite paths) 
initial segments $\xi\sq \tup$ in $\rL$. 
Then the symmetric subclass $\rL(\W\Om[\w])\sq\rL[\w]$ 
(Definition~\ref{LWOw}), 
generated by the set 
$\ens{\bx\res\et}{\et\in\Om}$, 
is a cardinal-preserving 
model of $\zf+\AC$ in which $\DC$ fails
(Jensen~\cite{jDC}), 
and more precisely $\xDC{\ip13}$ fails. 
Some other constructions within $\rL[\w]$ lead to 
other models in \cite{kl71,kl69}, \eg\ of\vom

$\zf+\neg\xAC{\ip13}$,\vom

$\zf+\xAC{\is1\iy}+\neg\xAC{\fp13}$,\vom

$\pao+\xAC{\is1\iy}+\xCA{\fs12}+\neg\xCA{\fs14}$. 
\epri

Admittedly, counter-examples obtained this way 
are one projective level worse than can be expected 
from the known positive results. 
For instance, instead of $\zf+\AC+\neg\xDC{\ip13}$ 
in the first counter-example one may want to get a 
model for $\zf+\AC+\neg\xDC{\ip12}$, since 
$\xDC{\is12}$ 
is provable. 
This goal was achieved with the help of 
Jensen $\ip12$-real singleton forcing. 

\bpri
\lam{jfr}
Jensen's forcing \cite{jenmin} is a proper subset $J\in\rL$ 
of the Sacks forcing $S$, obtained in the form 
$J=\bigcup_{\al<\omi}J_\al$ in $\rL$, where 
$\sis{J_\al}{\al<\omi}$ is a sequence of countable 
sets $J_\al\sq S$ defined by a certain 
$\omi$-long inductive construction in $\rL$ based on the 
diamond $\Diamond_{\omi}$. 
(In fact $J$ is not a unique forcing notion
in virtue of its definition, as \eg\ the Sacks forcing, 
but rather a family of similar forcing notions obtained 
by the construction in \cite{jenmin} that has some 
freedom at each step.)

The construction is maintained, 
using the diamond $\Diamond_{\omi}$ as a
\rit{sine qua non}, so that every possible 
antichain is killed at some step --- this implies 
CCC, and every possible partial order isomorphism 
also is killed at some step --- this implies the key 
property that $J$ adjoins a single generic real. 
Then estimating the complexity of the construction 
of $J$ we get that ``being a real $J$-generic 
over $\rL$'' is a $\ip12$ property. 
Therefore $J$ adjoins a generic $\ip12$ real 
singleton.
\epri

\bpri
\lam{jfa}
Countable-support iterated Jensen forcing of length $\omb$ 
was defined and studied by Abraham~\cite{abr,abr2}. 
\epri

\bpri
\lam{stru2}
Enayat~\cite{ena} used a finite-support infinite product 
of Jensen forcing to get  
a permutation model of $\zf$ with an infinite 
Dedekind-finite $\ip12$ set of reals, which 
implies the refutation of $\xAC{\ip12}$.   
\epri

\bpri
\lam{stru3}
By \cite{kl27}, it is forced by the finite-support product of 
$\om$ copies of Jensen forcing that the set of basic 
Jensen-generic reals is a countable $\ip12$ set containing 
no $\od$ real.  
\epri

\bpri
\lam{pro3}
A somewhat modified forcing notion, say ${\dJ}'$, 
rather similar to Jensen forcing $J$, 
is used in~\cite{kl25}.  
Instead of a single generic real by $\dJ$, 
it adjoins a \dd\Eo\rit{equivalence class} 
of \dd{{\dJ}}generic reals.  
(Reals $a,b\in\dn$ are \dd\Eo equivalent if 
$a(n)=b(n)$ for all but finite $n$. 
See some generalizations in~\cite{jsl97,kl36}.)
This \dd{{\dJ}'}generic 
\dd\Eo class is a (countable) $\ip12$ set containing 
no OD elements in the extension. 

This modification is maintained so that the automorphisms 
of $\dn$ naturally associated with $\Eo$ are somehow 
excluded from the killing procedure described in 
Example~\ref{jfr}.
\epri

\bpri
\lam{stru4}
Capitalizing on Examples~\ref{stru3} and \ref{stru1}, 
a generalized finite-support \dd\tup iteration of Jensen forcing 
is defined and studied in \cite{jml19}, to prove 
(among other results) that 
$\zf+\AC+\neg\xDC{\ip12}$ holds in a model similar to $\gN$ 
of Example~\ref{stru1}.  
Some other constructions within \dd\tup iterated Jensen 
extensions of $\rL$ lead to some other 
cardinal-preserving models, \eg\ of \vom

$\zf+\neg\xAC{\ip12}$,\vom 

$\zf+\xAC{\is1\iy}+\neg\xAC{\fp12}$, \vom

$\pao+\xCA{\is1\iy}+\xCA{\fs12}+\neg\xCA{\fs13}$ \ 
(see \cite{gitpf} on the latter), \vom 

\noi
which 
suitably strengthen the results of Example~\ref{stru1}. 
\epri

Another fundamental direction in these studies was discovered by Harrington~\cite{h74}. 
This is the construction of generic models in which 
some effect is achieved at a given level $n$ of the 
projective hierarchy, 
but not at previous levels. 
The results of Theorems \ref{mt1} and \ref{mt2} belong to 
this type, of course. 

\bpri
\lam{stru5}
As a further development of Jensen forcing 
of Example~\ref{jfr}, we defined a generic extension 
$\rL[a]$ in \cite{kl75}, by a real $a$ that is $\id1{n+1}$ 
in $\rL[a]$ for a given $n\ge2$, and such that any $\is1n$ real 
in $\rL[a]$ is constructible.  
(Jensen forcing itself gives the result for $n=2$ because of 
the Shoenfiend absoluteness.)  

The technique of \cite{kl75} involves a Harrington-style 
\cite{h74} modification of the original construction of 
Jensen forcing in $\rL$ in the form 
$J=\bigcup_{\al<\omi}J_\al$, 
as in Example~\ref{jfr}. 
The modification essentially requires the sequence of $J_\al$s to 
be ``$\id1n$-generic'' itself in the sense that it
meets every $\fs1{n-1}$ set dense in the ``super-tree'' of 
all possible countable beginnings of the construction. 

The effect of such a genericity is that the resulting forcing 
relation ${J(n)}\sq S\text{ (= Sacks forcing)}$ 
turns out to be an 
elementary subforcing of $S$ with respect to the forcing 
of $\is1n$ formulas. 
This leads to all $\is1n$ reals in ${J(n)}$-generic 
extensions of $\rL$ beings constructible, since such is 
the property of Sacks-generic extensions for all $n$. 
(Harrington carried out quite a similar construction in 
\cite{h74} 
\poo\ the almost-disjoint forcing of \cite{jsad}.) 
In the rest, similar to Example~\ref{jfr}, $J(n)$ 
adjoins a single \dd{J(n)}generic real $a$, and 
``being a \dd{J(n)}generic real'' is a $\ip1n$ formula, so 
$J(n)$ adjoins a $\ip1n$ real singleton, as required.
\epri

\vyk{
$\ans a$ is $\ip1n$ and $a$ itself is $\id1{n+1}$, and 3) the 
$\id1n$-genericity of the construction obscures things enough 
for all $\is1n$ reals in $\rL[a]$ being constructible. 
}

\bpri
\lam{stru6}
As a first approximation,
the proofs of our main results can be 
seen as using suitable symmetric submodels of 
generic extensions of $\rL$ forced by 
the generalized \dd\tup iteration 
(as in Examples~\ref{stru4} and \ref{stru1}) 
of a  Harrington-style 
``$\is1{\nn-1}$-generic'' version $J(\nn)$ of  
Jensen forcing. 
Yet in fact the proof will unfold somewhat differently. 
In particular, the standard forcing-iteration technicalities, 
instrumental in \cite{jml19,gitpf,ww}, will not be 
pursued. 
We'll rather define our forcing notion as   
\begin{quote}
{a  Harrington-style ``generic'' 
(as in Example~\ref{stru5}) 
Jensen-like subforcing\/ 
$\cX$ of the generalized countable-support 
\dd\tup iteration\/ $\pei$ 
(= iterated perfect sets) 
of the Sacks forcing.} 
\end{quote}
But the ideas outlined in Examples~\ref{stru1}, \ref{stru4},  
\ref{stru5} will be implicitly included.
\epri

\vyk{
This construction realizes the idea of generalized \dd\tup 
iteration of Jensen's forcing somewhat differently than in 
\cite{jml19,gitpf,ww}, in particular, the CCC property will 
not be achieved.
}

\parf{The structure of the paper}
\las{pap}

The implementation of the plan outlined in 
Example~\ref{stru6} is organized as follows. 
It turns out that the usual approach 
to iterations of Jensen or similar 
forcing based on perfect trees, as in \cite{jml19}, 
leads to significant technical difficulties, 
which we have not been able to completely overcome, 
especially with regard to Harrington's idea 
of ``generic'' forcing constructions. 

This is why we have to turn to a purely geometric 
method of working with such iterations, 
developed in \cite{fm97,jsl99}. 
It presents the generalized \dd\tup iterated Sacks forcing as 
the set $\pei$ of {\ubf iterated perfect sets}, \ie, 
certain closed sets in spaces $\can\xi$, where 
$\can{}=\dn$ is the Cantor spase and $\xi\sq\tup$ is a 
\index{Cantor space, $\can{}=\dn$}%
\index{zD@$\can{}=\dn$, Cantor space}%
\index{Cantor space, $\can{\xi}$}%
\index{zDxi@$\can{\xi}$}%
countable initial segment in $\tup$. 
These sets are introduced and studied in Chapters \ref{ipset}
and \ref{sfc}, with the {\ubf splitting/fusion construction} 
introduced in the latter. 

Any set $\cX\in\rL$ of iterated perfect sets, satisfying some 
natural conditions, can be viewed as a forcing notion that 
adjoins a generic \dd\tup array of reals in $\cam$. 
Such forcing notions $\cX\sq\pei$, called {\ubf normal forcings},  
corresponding $\cX$-{\ubf generic arrays} $\w\in\can\tup$,  
{\ubf generic extensions} $\rL[\w]$, 
their symmetric subextensions, 
and associates forcing relations,
are studied in Chapter~\ref{rfono}. 

Chapter~\ref{subex} introduces those 
{\ubf symmetric submodels} 
of generic extensions $\rL[\w]$  which are specifically 
involved in the proof of Theorem~\ref{mt1}. 

From this moment on, 
{\ubf we fix a number $\nn\ge1$} as in 
Theorems~\ref{mt1} and \ref{mt1}. 
Four key conditions for a normal forcing $\cX\sq\pei$ 
are introduced, which guarantee  
 that those symmetric submodels bring the desired result, 
 two of them involve $\nn$:
\bde 
\item[{\ubf Fusion property},] 
which postulates for $\cX$ a well-known feature of 
the Sacks forcing and its iterations like $\pei$;  

\item[{\ubf Structure property}:]
for all\/ $\i,\j\in\tup,$ we have\/ 
$\w(\i)\in\rL[\w(\j)]$ iff\/ $\i\sq\j\,;$
\index{property!Structure property}%
\index{Structure property}%

\item[$\nn$-{\ubf Definability property},]
which claims that the binary relation\vim
\bce
$x=\w(\i)\yi y=\w(\j)$ for 
some {\bfit even\/} tuples $\i\su\j$ in $\tup$\vim
\ece
is $\ip1{\nn+1}$ 
in any suitable submodel 
of any $\cX$-generic extension $\rL[\w]$. 
\ede

A tuple of ordinals  
is \rit{even}, resp.\ \rit{odd}, if such is its last term. 

\bde
\item[{\ubf $\nn$-Odd-Expansion property:}] 
if $\xi\in\rL$, $\xi\sq\tup$ is a countable initial segment, 
$\vpi(x)$ a $\ip1\nn$ formula with 
reals in $\rL[\w\res\xi]$ as parameters, 
and $\rL[\w]\mo\sus x\,\vpi(x)$, 
then such a real $x$ exists in $\rL[\w\res\ta]$ for some $\ta\in\rL$ 
(still a countable initial segment) 
such that $\ta\bez\xi$ consists only of {\bfit odd} tuples. 
\ede

We replace the 
the $\nn$-Odd-Expansion property with a more convenient 
 property of {\ubf $\nn$-comp\-le\-teness} for $\cX$ 
in Chapter~\ref{foax}.
For this purpose, we introduce {\ubf an auxiliary forcing relation} 
$X\fo\vpi$ in $\rL$,  where $X\in\pei$ 
and $\vpi$ is a formula of a certain extension of 
the language of 2nd order arithmetic $\pad$. 
Then, a normal forcing $\cX$ is {\ubf $\nn$-complete}, if for any 
closed $\is1n$ formula $\vpi$ of the extended language, 
the set of all $X\in\cX$ satisfying $X\fo\vpi$ or $X\fo\neg\vpi$, 
is dense in $\cX$. 
This is how Harrington's idea of ``generic'' forcing notions 
(Example~\ref{stru5}) is realized within 
the background forcing notion $\pei$ in our proof.

Note that $\fo$ is connected rather with the full $\pei$ 
as the forcing notion, but if 
$\cX$ is $\nn$-complete then $\fo$ 
{\ubf coincides with the usual $\cX$-forcing relation} 
up to $\is1{n+1}$ formulas. 
This allows to show that $\nn$-Completeness implies 
$\nn$-Odd-Expansion. 
Hence the whole task related to Theorem~\ref{mt1} 
is reduced to the following: 
\ben
\fenu
\itlb{pap*}
for a given $\nn\ge1$, 
find a normal forcing $\cX$ in $\rL$, satisfying 
the Fusion, Structure, $\nn$-Definability, 
and $\nn$-Completeness properties. 
\een

The construction of such a forcing $\cX$ is carried out in 
Chapters~\ref{sek6}--\ref{VII},  
as a sort of {\ubf limit} of an $\omi$-sequence of 
countable collections of iterated perfect sets, called {\ubf rudiments}. 
Rudiments, and sequences of rudiments increasing in the 
sense of a {\ubf refinement} relation $\ssq$, are studied in 
Chapter~\ref{sek6}. 

We introduce some properties of an $\ssq$-increasing 
$\omi$-sequence of rudiments in Chapter~\ref{sek7}, which 
imply that the associated limit forcing $\cX$ satisfies \ref{pap*} 
above. 
The properties are summed up in the notion of 
{\ubf 1-5-$\nn$ extension}, 
such that \ref{pap*} is reduced to the following: 
\ben
\fenu
\atc
\itlb{pap**}
for a given $\nn\ge1$, construct an $\ssq$-increasing 
$\is\hc\nn$-definable 
$\omi$-sequence of rudiments in $\rL$, such that each term 
is a 1-5-$\nn$ extension of the subsequence of all previous terms. 
\een

We prove {\ubf the existence of 1-5-$\nn$ extensions} 
in Chapter~\ref{VI}, 
and then accomplish \ref{pap**} and the proof of 
Theorem~\ref{mt1}  in 
Chapter~\ref{VII} by the construction of a sequence required by 
taking {\ubf the $\lel$-least possible} 
1-5-$\nn$ extension at each step of the construction. 

Chapter~\ref{b} presents 
{\ubf the proof of Theorem~\ref{mt2}}. 
We use yet another symmetric submodel of an $\cX$-generic 
extension $\rL[\w]$ of $\rL$, for the same forcing $\cX$. 

The paper ends with a usual conclusion-style  
material in Chapter~\ref{frq}.
In particular, we'll touch on the evaluation 
of those proof theoretic tools used in the arguments. 
We discuss in Section~\ref{70} how the main 
consistency results of 
this paper can be obtained on the basis 
of the formal consistency of second order arithmetic $\pad$.  
This is a crucial advantage comparably to some earlier 
results, like \eg\ the above-mentioned results 
by Levy~\cite{levy2} and Guzicki~\cite{guz} which 
definitely cannot be obtained on the basis of 
the consistency of $\pad$.

\vyk{
It remains to note that topics in subsystems of 
second order arithmetic remain of big interest 
in modern studies, see \eg\ \cite{frit}, and our 
paper contributes to this research line.
}

\parf{Definability, 
constructibility, diamond prerequisites}
\las{cdia}

Recall that
$\hc=H\omi=
\ens{x}{\TC (x)\text{ is at most countable}}$, 
\kmar{hc}%
\index{hereditarily countable, $\hc$}%
\index{zHC@$\hc$, hereditarily ctble}%
the set of all \rit{hereditarily countable} sets. 
The $\in$-definability over $\hc$ is connected with the 
descriptive set theoretic definability by the following 
classical result:

\bpro
[see \eg\  25.25 in Jech \cite{jechmill}]
\lam{p60}
If\/ $n\ge1$ and\/ $X\sq\cN$ then
$$
{X\in\is1{n+1}}\eqv{X\in\is\hc{n}}
\qand 
{X\in\ip1{n+1}}\eqv{X\in\ip\hc{n}}, 
$$
and\/ ${X\in\is1{n+1}(p)}\eqv{X\in\is\hc{n}(p)}$ 
for any parameter\/ $p\in\cN,$ \etc\qed
\epro

{\ubf Assume $\rV=\rL$ in the remainder of this section.}

It is known that $\hc=\lomi$ provided $\rV=\rL$.
\index{G\"odel \weo, $\lel$}%
\index{wellodering, $\lel$}%
\index{zzz-L@$\lel$, G\"odel \weo}%
\kmar{lel}
Let $\lel$ be the G\"odel \weo\ of $\rL$.  
If $\al<\omi$ then  we let 
\kmar{lc al}%
\index{zcal@$\lc\al$}%
$\lc\al$ be the $\al$th
member of $\hc=\lomi$ in the sense of $\lel$,  
\kmar{bH al}%
\index{hereditarily countable!below $\al$, $\bH\al$}%
\index{zHCal@$\bH\al$, hered.\ ctble below $\al$}%
and $\bH\al=\ens{\lc\ga}{\ga<\al}$.   
%
The following is well-known. 

\bpro
[$\rV=\rL$]
\lam{p61}
The relation\/ ${\lel}\res\hc$
has length\/ $\omi$, therefore\/ 
$\hc=\ens{\lc\al}{\al<\omi}$ and\/ 
$\bH\al\in\hc$ for all\/ $\al<\omi$. 
In addition$:$
\ben
\renu
\itlb{p611}
${\lel}\res\hc$ is a\/ $\id\hc1$ relation, 
the set\/ $\ens{\bH\al}{\al<\omi}$ is\/  $\id\hc1$, too$;$ 

\itlb{p612}
the maps\/ $\al\mto\lc\al$ and\/ $\al\mto\bH\al$ 
are\/ $\id\hc1$ as well$;$

\itlb{p613}
the relation\/ 
${\lel}\res\hc$ is\/ {\ubf good}, in the sense 
that if\/ $p\in\hc$,
$n\ge1$, and\/ $P(\cdot,\cdot,\cdot)$ is a
ternary\/ $\id\hc n(p)$ relation on $\hc$,
then so are the binary relations\/
$\sus x\lel y\,P(x,y,z)$ 
and\/
$\kaz x\lel y\,P(x,y,z)\,.$\qed 
\een
\epro

The \rit{diamond principle} 
\index{diamond $\Diamond_{\omi}$}%
\index{zzzDiamond@$\Diamond_{\omi}$, diamond}%
$\Diamond_{\omi}$ is true in $\rL$  
by \cite[Thm 13.21]{jechmill}, 
hence there is a  $\id\hc1$ sequence of sets 
\kmar{xs al}%
\index{zSal@$\xs\al$}%
$\xs\al\sq\al$, $\al<\omi$, such that 
\ben
\Aenu
\itlb{1gret}
if\/ 
$X\sq\hc$ then the set\/ 
$\ens{\al<\omi}{\xs\al=X\cap\al}$ 
is stationary in\/ $\omi$.
\een
The $\id\hc1$-definability of the sequence is achieved 
by taking the \dd\lel least possible $\xs\al$ at each 
step $\al$ in the standard construction of $\xs\al$, 
as \eg\ in \cite{jechmill}. 
Define
$$
\baS\al=\ens{\lc\ga}{\ga\in S_\al} \text{  for }\al<\omi,
\kmar{baS al}%
\index{zSbaral@$\baS\al$}%
\quad\text{hence}\quad
\baS\al\sq
\bH\al:=\ens{\lc\ga}{\ga<\al}.
$$

We get the following  as an easy corollary of \ref{1gret} 
and Proposition~\ref{p61}. 

\bpro
[$\rV=\rL$]
\lam{p62}
The map\/ $\al\mto\baS\al$ is\/ $\id\hc1$.

If\/ $\baS{}\sq\hc$ then the set\/ 
$\ens{\al<\omi}{\baS\al=\baS{}\cap\bH\al}$ 
is stationary.\qed
\epro

\vyk{
If $\al<\omi$ then (still under $\rV=\rL$)
\kmar{jj al}
we let  $\jj\al$ be the least
ordinal $\ba<\omi$  such that $\rL_\ba$ contains
a map $f:\om\onto\al$, and put
\imaf{dpsi al}
$\dpsi\al=\ens{x\in\rL_{\jj\al}}{x\sq\al}$.
}


\sekt{Iterated perfect sets}
\las{ipset}

The proof of our main results 
involves the engine of generalized 
product-iterated Sacks forcing developed in 
\cite{fm97,jsl99} on the basis of earlier papers 
\cite{balav,grj,groapp} and others. 
We consider the constructible universe 
$\rL$ as the ground model for any forcing in the 
remainder.

\parf{Spaces and projections}
\las{prelim1}

{\ubf Arguing in $\rL$ in this section},  
we define, in $\rL$, the set 
\kmar{tup}%
$ 
\index{zI@$\tup$, tuples of ctble ordinals}%
\index{tuples of ctble ordinals, $\tup$}%
{\tup=\omi\lom\bez\ans\La\in\rL}
$ %
of all non-empty tuples 
$\i=\ang{\ga_0,\dots,\ga_{n-1}}$, $n\ge1$,
of ordinals $\ga_k<\omi$. 
The set $\tup$ is 
partially ordered by {\em the strict extension\/} 
$\su$ of tuples. 
Then $\tup$ is a tree without a root because 
$\La$, the empty tuple, is excluded. 
We put 
$$
\bay{rclcl}
\index{tuples!dyadic@$\tud$}%
\index{zI2@$\tud$, dyadic tuples}%
\tud &=& 2\lom\bez\ans\La
\kmar{tud}%
&=&\ens{\i\in\tup}{\ran \i\sq\ans{0,1}},\\[0,5ex]
\index{tuples!natural@$\tuo$}%
\index{zIom@$\tuo$, natural tuples}%
\kmar{tuo}%
\tuo &=& \om\lom\bez\ans\La &=&\ens{\i\in\tup}{\ran \i\sq\om},
\eay
$$
and generally 
\kmar{tuq al}%
\index{tuples!$\al$-bounded@$\tuq\al$}%
\index{zIal@$\tuq\al$, tuples $\al$-bounded}%
$\tuq\al=\al\lom\bez\ans\La =\ens{\i\in\tup}{\ran \i\sq\al}$,
so $\tuq\omi=\tup$. 

If $\i\in\tup$ then $\lh \i$ is the length of $\i$; 
\index{tuples!length, $\lh\i$}%
\index{length, $\lh\i$}%
\index{zlhi@$\lh\i$}%
\kmar{lh i}%
$\lh\i\ge1$ since $\La$ is excluded.

Our plan is to define a   
generic extension $\rL[\a]$ of $\rL$ 
by an array $\a=\sis{\a_\i}{\i\in\tup}$ of reals 
$\a_\i\sq\om$, in which the 
structure of iterated genericity of the reals $\a_\i$ will 
be determined by this set $\tup$. 

Let $\cpo$ 
be the set of all at most countable 
\kmar{cpo}%
\index{initial segments, $\cpo$}%
\index{zzXi@$\cpo$}%
initial segments (in the sense of $\su$) 
$\za\sq\tup$. 
If $\za\in\cpo$ then $\IS_\za$ is the set of all 
\index{zISza@$\IS_\za$}%
\index{initial segments, $\IS_\za$}%
initial segments of $\za$. 

Greek letters $\xi,\,\eta,\,\za,\,\vt\yi\ta$ will 
denote sets in $\cpo$.  

Characters $\i,\,\j$ are used to denote 
{\it elements\/} of $\tup$.

For any $\i\in\za\in\cpo,$ 
we consider initial segments
\index{initial segment!isu@$\ile\i$}%
\index{zzzisu@$\ile\i$}%
\index{initial segment!isq@$\ilq\i$}%
\index{zzzisq@$\ilq\i$}%
\index{initial segment!nsq@$\nlq\i$}%
\index{zzznsq@$\nlq\i$}%
\kmar{ile i}%
$\ile\i=\ans{\j\in\tup:\j\su\i}$,  
$\ilq\i=\ans{\j\in\tup:\j\sq\i}$,  
\kmar{ilq i\mns za nlq i}%
$\za\nlq\i=\ans{\j\in\za:\i\not\sq\j},$    
%
Clearly $\ile\i\sneq\ilq\i\sq\za$. 

\vyk{
We consider $\pws\om$ as identic to $\dn,$ 
so that both $\pws\om$ and $\pws\om^\xi$ for 
$\xi\in \cpo$ are  
Polich compact spaces, homeomorphic to each other. 
Points of $\pws\om$ will be called {\it reals\/}. 
}

Let $\can{}=\dn\sq\cN$ be the {\it Cantor space\/}.
\kmar{can}
\index{Cantor space, $\can{}=\dn$}%
\index{zD@$\can{}=\dn$, Cantor space}%
\index{Cantor space, $\can{\xi}$}%
\index{zDxi@$\can{\xi}$}%
For any set 
$\xi,$ $\can\xi$ is the product of \dd\xi many 
copies of $\can{}$ with the product topology. 
Then every $\can\xi$ is a compact space. 

\bdf
[projections]
Assume that $\et\sq\xi$ belong to $\cpo$.  

If $x\in\can\xi$ then let 
$x\dar\et=x\res\eta\in\can\eta$ denote the usual restriction.
\index{projection!dar@$\dar$}%
\index{zzzpdar@$\dar$}%
\kmar{x dar et}%
If $X\sq\can\xi$ then let 
$X\dar \et=\ens{x\dar\et}{x\in X}$. 
Moreover if $\rX$ consists of sets $X\sq\can\xi$ for 
\kmar{rX car et}%
different supersets $\xi$ of $\et$ then let 
$\rX\dar \et=\ens{X\dar\et}{X\in \rX}$. 

If $Y\sq\can\eta$ then let 
$Y\uar \xi=\ens{x\in\can\xi}{x\dar\et\in Y}$ 
\index{lifting!uar@$\uar$}%
\index{zzzpuar@$\uar$}%
\kmar{uar}%
(\rit{lifting}). 
\index{lifting}

We define $X\rsq\i=X\dar\ilq\i$, 
$X\usq\i=X\uar\ilq\i$,
\index{projection!prsq@$\rsq\i$}%
\index{zzzprsq@$\rsq\i$}%
\index{projection!pusq@$\usq\i$}%
\index{zzzprsq@$\usq\i$}%
\index{projection!prsl@$\rsl\i$}%
\index{zzzprsl@$\rsl\i$}%
\index{projection!pusl@$\usl\i$}%
\index{zzzpusl@$\usl\i$}%
%
\kmar{rsq, usq}%
and similarly 
$X\rsl\i$, $X\usl\i$, 
\kmar{rsl,usl}%
$x\rsq\i$ \etc\ for points $x$, and $\cX\rsq\i$ \etc\ for 
collections $\cX$ of sets.

Finally, we let $X\dir \i=\ens{x(\i)}{x\in X}$.
\kmar{X dir i}%
\index{zzzdir@$\dir\i$}%
(Note a different arrow.)
\edf

\parf{Iterated perfect sets and \prok}
\las{prelim2}

{\ubf We argue in $\rL$ in this section.}
To describe the key idea, 
recall that the Sacks forcing consists 
of perfect subsets of $\can{}$, 
which are exactly those of the form 
$X=\ima H\cam =\ens{H(a)}{a\in\cam}$, 
where $H:\can{}\onto X$ is a homeomorphism. 

To get a \rit{product} Sacks forcing with two factors 
(the case of a two-element unordered set 
as the generalized ``length'' of iteration), 
we have to consider sets $X\sq\can 2$ 
of the form $X= \ima H{\can 2}$ where $H$ is  
any homeomorphism defined on $\can 2$ so that 
it splits in obvious way into a pair of 
one-dimensional homeomorphisms. 

To get an \rit{iterated} Sacks forcing, 
with two stages of iteration  
(the case of a two-element ordered set 
as the ``length'' of iteration), 
we make use of sets $X\sq\can 2$ 
of the form $X=\ima H{\can 2}$, 
where $H$ is any homeomorphism defined on $\can 2$ 
such that if $H(a_1,a_2)=\ang{x_1,x_2}$ 
and $H(a'_1,a'_2)=\ang{x'_1,x'_2}$ then 
$a_1=a'_1\eqv x_1=x'_1$. 

The combined product/iteration case results 
in the following definition. 

\bdf
[\cite{fm97,jsl99}]
\lam{ips}
For any 
$\za\in\cpo,$ let $\pe\za$ 
(iterated perfect sets of dimension $\za$) 
be the collection of all
\index{iterated perfect sets, $\pe\za$}%
\index{zIPSz@$\pe\za$}%
sets $X\sq\can\za$ such that there is a homeomorphism 
\kmar{pe za}%
$H:\can\za\onto X$ satisfying
\dm
x_0\dar\xi=x_1\dar\xi\,\eqv\,H(x_0)\dar\xi=H(x_1)\dar\xi
\dm
for all $x_0,\,x_1\in\dom H$ and $\xi\in\cpo$, 
$\xi\sq\za$. 
Homeomorphisms $H$ 
\pagebreak[0] 
satisfying this requirement 
\index{projection-keeping, PKH}%
\index{homeomorphism!projection-keeping, PKH}%
will be called {\it\prok\/}, \pkh\ for brevity.
\kmar{prok}%
In other words, sets in $\pe\za$ are images 
of $\can\za$ via \pkh s. 

We put $\pei=\bigcup_{\xi\in\cpo}\pe\xi$.   
\kmar{pei}
\index{IPS, iterated perfect sets}%
\index{iterated perfect sets, $\pei$}%
\index{zIPS@$\pei$}%
Sets in $\pei$ are called {\em iterated perfect sets\/}, 
\ips\ in brief. 
If $X\in\pe\xi$ then let $\dym X=\xi$ 
(the \rit{dimension} of $X$).
\index{dimension, $\dym X$}%
\index{z11X11@$\dym X$, dimension}%
\kmar{dym X}%

We let $\pel\i=\pe{\ile\i}$, $\pele\i=\pe{\ilq\i}$ 
\kmar{pel pele}
\index{iterated perfect sets, $\pel\i$}%
\index{zIPSsi@$\pel\i$}%
\index{iterated perfect sets, $\pele\i$}%
\index{zIPSsqi@$\pele\i$}%
for the sake of brevity.
\edf

\bre
\lam{why?}
Suppose that $\za\in\cpo$ in $\rL$. 
The set $\pe\za$, defined in $\rL$, can be considered 
as a forcing notion. 
It is established in \cite[Thm 1 and Subsection 6.1]{jsl99} 
that $\pe\za$ adjoins a generic array $\w\in\can\za$ of 
reals $\w(\i)\in\can{}=\dn\yd \i\in\za$, such that each 
real $\w(\i)$ is Sacks-generic over $\rL[\w\rsl\i]$. 
Thus $\pe\za$ works as a generalized $\za$-long 
iteration of the Sacks (perfect set) forcing. 
This is why we call sets in $\pe{}$ 
\rit{iterated perfect sets}.
\ere

\bre
\lam{pu}
The empty set $\pu\in\cpo$, 
$\can\pu=\ans\pu$, $\bon=\ans\pu\in\pe\pu$.
\index{z1@$\bon$}
\ere

\ble
\lam{997}
If\/ $H$ is a PKH defined on\/ $X\in\pe\za$ 
then the image\/ $\ima H X=\ans{H(x):x\in X}$ belongs to $\pe\za$.
\ele
\bpf 
The superposition of two PKHs is a PKH.
\epf

\ble
\lam{99x}
If\/ $X\in\pe\za$, $\et\in\IS_\za$, $\i\in\za\bez\et$, 
then there exist points\/ $x,y\in X$ with\/ $x\dar\et=y\dar\et$ 
but\/ $x(\i)\ne y(\i)$.
\ele
\bpf 
There is a PKH $H:\can\za\onto X$.  
Assume \noo\ that $\et=\za\nlq\i$ 
(otherwise consider $\et'=\za\nlq\i$). 
Obviously there are points $x',y'\in \can\za$ with\/ 
$x'\dar\et=y'\dar\et$ but\/ $x'(\i)\ne y'(\i)$, hence 
$x'\rsq\i\ne y'\rsq \i$. 
Their $H$-values $x=H(x')$, $y=H(y')$ then satisfy 
$x\dar\et=y\dar\et$ but $x\rsq\i\ne y\rsq \i$. 
Yet $\ile\i\sq\et$, so that $x\rsl\i= y\rsl \i$. 
And this implies $x(\i)\ne y(\i)$.
\epf

\parf{Some basic properties of iterated perfect sets}
\las{bp}

{\ubf We argue in $\rL$ in this section.} 
Here follows a collection of some results related
to iterated perfect sets, partially taken from
\cite{fm97,jsl99}.

\ble
[Proposition 4 in \cite{jsl99}]
\lam{oldf}
Let $\za\in\cpo$.
Every set\/ $X\in\pe\za$ is closed and satisfies 
the following properties$:$ 
\ben
\renu
\itlb{perf} 
if\/ $\i\in\za$ and\/ $z\in X\rsd{\su\i}$ 
then\/ 
$\cs Xz\i=\ens{x(\i)}{x\in X\land x\rsd{\su\i}=z}$
\index{zDXzi@$\cs Xz\i$, cross-section}%
\index{cross-section, $\cs Xz\i$}%
is a perfect set in\/ $\cam$, 
\kmar{cs Xzi}%

\itlb{oz} 
if\/ $\xi\in\cpo$, $\xi\sq\za$, and a set\/ $X'\sq X$ 
is open in\/ $X$ 
(in the relative topology) 
then the projection\/ $X'\car\xi$ is 
open in $X\car\xi$ --- 
in other words, the projection from $X$ to $X\car\xi$ 
is an open map, 

\itlb{indep} 
if\/ $\xi,\eta\in\IS_\za$, 
$x\in X\car\xi$, $y\in X\car\eta$, and\/ 
$x\car(\xi\cap \eta)=y\car(\xi\cap \eta),$ then\/ 
$x\cup y\in X\car(\xi\cup \eta)$. 
\een
\ele
\bpf[sketch] 
Clearly $\can\za$ satisfies \ref{perf}, \ref{oz}, 
\ref{indep}, 
and one easily shows that \prok\  
homeomorphisms preserve the requirements.
\epf

\ble
[routine from \ref{indep}] 
\lam{no19}
Suppose that\/ $\xi,\et,\vt\in\cpo$, 
$\vt\cup\et\sq\xi$, and\/ $X\in\pe\xi$. 
Then\/ 
$X\car{(\et\cup\vt)}=(X\car\vt\uar{(\et\cup\vt)})
\cup (X\car\et\uar{(\et\cup\vt)})$.\qed 
\ele

\ble
[Lemma 5 in \cite{jsl99}] 
\lam{599}
Suppose that\/ $\xi,\et,\vt\in\cpo$, 
$\vt\cup\et\sq\xi$,  
$W\in\pe\xi$, $C\sq W\res\et$ is any set, 
and\/ $U=W\cap (C\uar\xi)$.   
Then\/ 
\ben
\renu
\itlb{599i}
$U\car\vt=(W\car\vt)\cap(C\car(\vt\cap\et)\uar\vt)\,;$  

\itlb{599ii}
if\/ $\vt=\ilq\i$, $\i\in\xi$, then\/ 
$U\rsq\i=(W\rsq\i)\cap(C\dar{\sg}\usq\i)$, 
where\/ $\sg=\et\cap\ilq\i$, in particular, if\/ 
$\i\in\et$ then\/ 
$U\rsq\i=C\dar{\sg}\usq\i$.\qed 
\een
\ele

\ble
[Lemma 6 in \cite{jsl99}] 
\lam{less}
If\/ $\xi\sq\za$ belong to\/ $\cpo$, and\/ $X\in\pe\za$, 
then\/ $X\car\xi\in\pe\xi$.\qed
\ele

\ble
[Lemma 9 in \cite{jsl99}] 
\lam{apro}
Suppose that\/ $\za\in\cpo$, $\eta\in\IS_\za$,  
$X\in\pe\za$, $Y\in\pe\eta,$ and\/ $Y\sq X\car \eta$. 
Then\/ $Z=X\cap (Y\uar \za)$ belongs to\/ $\pe\za$. 

In particular\/ $Y\uar\za\in\pe\za$ {\rm(lifting)}, as 
obviously\/ $\can\za\in\pe\za$.\qed 
\ele

\ble
[Lemma 9 in \cite{fm97}] 
\lam{fm9}
If\/ $\et\sq\xi$ belong to\/ $\cpo$,
$X,Y\in\pe\xi$, and\/ $X\res\et=Y\res\et$, 
then there is a \pkh\/  
$H:X\onto Y$ such that\/ $H(x)\car\et=x\car\et$
for all\/ $x\in X$.\qed
\ele

\ble
\lam{99}
Suppose that\/ $\et\sq\xi$ belong to\/ $\cpo$,
$X\in\pe\xi$, $Y=X\car\et\in\pe\et$, and\/
$H:\can\et\onto Y$ is a \pkh. 
Then there is a \pkh\/  
$K:\can\xi\onto X$ such that\/
$K(x)\car\et=H(x\dar\et)$
for all\/ $x\in \can\xi$. 
\ele

\bpf
The set $Y'=Y\uar\xi$ belongs to $\pe\xi$ by
Lemma~\ref{apro}.
Therefore, by Lemma \ref{fm9}, there is
a \pkh\/  
$J:Y'\onto X$ such that\/ $J(x)\car\et=x\car\et$
for all\/ $x\in Y'$.
Yet by the choice of $H$, the map
$H':\can\xi\to Y'$ defined by
$H'(x)\car\et=H(x\car\et)$ and
$H'(x)\car{(\xi\bez\et)}=x\car{(\xi\bez\et)}$
for all $x\in\can\xi$, is a
\pkh\/ $\can\xi\onto Y'$.
Thus the superposition
$K(x)=J(H'(x))$ is a \pkh\ 
$\can\xi\onto X$, and
if $x\in\can\xi$ then 
$K(x)\car\et=J(H'(x))\car\et=H'(x)\car\et
=H(x\car\et)$.
\epf

\bcor
\lam{pe2}
Let\/ $\xi,\eta\in\cpo$, $\vt=\xi\cup\eta$, 
$X\in\pe\xi$, $Y\in\pe\eta$,   
$X\car(\xi\cap\eta)=Y\car(\xi\cap\eta)$. 
Then\/ $Z=(X\uar\vt)\cap (Y\uar\vt)\in\pe\vt$,
$Z\car\xi=X$, $Z\car\et=Y$.
\ecor
\bpf
The set 
$X'=X\uar\vt$ belongs to $\pe\vt$ by 
Lemma~\ref{apro}. 
In addition, 
$X'\car{\eta}=X\car(\xi\cap\eta)\uar\eta$
by Lemma~\ref{599} (with $C=X$, $W=\can\vt$). 
Then $Y\sq X'\car{\eta}$, because 
$Y\car(\xi\cap\eta)=X\car(\xi\cap\eta)$. 
We conclude that $X'\cap (Y\uar\vt)\in\pe\vt$ 
by Lemma~\ref{apro}. 
Finally,  
$X'\cap (Y\uar\vt)=Z$ by construction.

To check that say $Z\car\xi=X$, let $x\in X$.
There is $y\in Y$ with
$x\car(\xi\cap\eta)=y\car(\xi\cap\eta)$.
Then $z=x\cup y\in Z$ by construction, and
$z\car\xi=x$.
\epf

\parf{Clopen subsets}
\las{closu}

{\ubf We argue in $\rL$ in this section.} 

The next lemma highlights the Sacks-iterated character 
of sets in $\pe\xi$ in case $\xi=\ilq\i$.
Let a {\it \pbt{}\/} be any (nonempty) tree $T\sq 2^{<\om}$ 
with no endpoints, 
such that 
$B(T)=\ans{t\in T:t\we 0\in T\cj t\we 1\in T}$, 
the set of all splitting points, is cofinal 
in $T.$ 

Let 
$\Der=\ens{T\sq\bse}{T\text{ is a perfect tree}}$, 
\kmar{Der}
\index{zPT@$\Der$, perfect trees}%
\index{perfect tree!$\Der$, perfect trees}%
a closed set in $\pws{\bse}$. 

If $T\in\Der$ then 
$[T]=\ens{x\in\dn}{\kaz k\,(x\res k\in T)}$,   
\index{zzz1T1@$[T]$, perfect set}%
\index{perfect set!$[T]$, perfect set}%
a perfect set. 

Conversely,  
$\der(X)=\ens{s\in\bse}{[s]\cap X\ne\pu}\in\Der$ 
\kmar{der}
\index{ztreeX@$\der(X)$, perfect tree}%
\index{perfect tree!$\der(X)$, perfect tree}%
for any perfect set $X\sq\dn,$  
where $[s]=\ens{x\in\dn}{s\su x}$ 
for $s\in\bse.$

\vyk{
Suppose $T$ is such a tree. 
We define:
\bit
\itemsep=1mm
\item\msur
$[T]=\ans{a\in 2^\om:\forall\,m\,(a\res m\in T)},$ a 
perfect set is $\cD= 2^\om$, 
\item 
an order isomorphism $\ba_T:2^{<\om}\hsur$ onto 
$\hsur B(T).$ We define $\ba_T(u)\in B(T)$ for every  
$u\in 2^{<\om}$ by induction on 
$\dom u,$ putting $\ba_T(u\we e)$ to be the least 
$s\in B(T)$ such that $\ba_T(u)\we e\sq s,$ for $e=0,1$,

\item 
a homeomorphism $H_T:\can{}\onto[T]$,   
$H_T(a)=\bigcup_{m}\ba_T(a\res m)$, $\kaz a\in \can{}$.
\eit
}

\ble
[Lemma 11 in \cite{jsl99}]
\lam{hh}
Assume that\/ $\i\in\tup$, 
$Y\in\pel\i,$   
$\cT$ continuously maps\/ $Y$ into\/ $\cP(2^{<\om})$ so that\/ 
$\cT(y)\in\Der$ for all\/ $y\in Y.$ 
Then\/ 
$X=\ans{x\in\can{\ilq\i}:
x\rsl\i\in Y\land x(i)\in[\cT(x\rsl\i)]}
\in\pele\i$.\qed
\ele

\vyk{
\bpf
The set $Z=Y\ares\za$ belongs to $\pe\za$ by Lemma~\ref{apro}, 
so it suffices 
to define a PKH $H:Z\onto X.$ 
Let $z\in Z.$ 
Then $y=z\res\eta\in Y$ whereas $a=z(\i)\in\cD.$ 
Define $x=H(z)\in\can\za$ so that $x\res\eta=y$ and 
$x(\i)=H_{\cT(y)}(a).$ 

Then $H$ maps $Z$ onto $X$ as every $H_{\cT(y)}$ maps $\cD$ 
onto $[\cT(y)]=\ans{x(\i):x\in X\cj x\res\eta=y}.$ 
$H$ is $1-1$ since each $H_T$ is $1-1,$ 
and $H$ is continuous because so is the map $\cT.$ 
It remains to prove that $H$ is 
projection--keeping, \ie\  the equivalence 
${z_0\res\xi=z_1\res\xi\,\eqv H(z_0)\res\xi=H(z_1)\res\xi}$ 
for all $z_0,\,z_1\in Z$ and $\xi\in\IS_\za.$ 
If $\i\not\in\xi$ 
then $\xi\sq\eta$ and $z\res\xi=H(z)\res\xi$ by definition. 
If $\i\in\xi$ then $\xi=\za,$ 
and the result is obvious as well.
\epf
}

The following is a converse to Lemma~\ref{hh}. 

Recall that perfect sets $\cs Xy \i$ are defined by 
Lemma~\ref{oldf}\ref{perf}.

\ble
\lam{hh+}
Let\/ $\i\in\tup$, $X\in\pele\i$, 
$Y=X\rsl\i\in\pel\i$,  and if\/ $y\in Y$ then\/ 
$\cT_X(y)=\der(\cs Xy \i)$.  
Then\/ $\cT_X$ continuously maps\/ $Y$ into\/ 
$\Der$.
\ele
\bpf
Let $s\in\bse$ and $Y_s=\ens{y\in Y}{s\in\cT_X(y)}$. 
Then $Y_s=X_s\rsl\i$, where 
$X_s=\ens{x\in X}{s\su x(\i)}$. 
It follows that $Y_s$ is clopen in $Y$ by 
Lemma~\ref{oldf}\ref{oz}. 
By similar reasons, the set 
$Y'_s=\ens{y\in Y}{s\nin\cT_X(y)}$ 
is clopen in $Y$ as well.
\epf

We continue with assorted results on clopen subsets 
of sets in $\pei$. 

The next lemma fails for $\pe\xi$ in case $\xi\in\cpo$ 
is not linearly ordered by $\sq$.

\ble
\lam{lin}
Let\/ $\i\in\tup$, 
$X\in\pele\i$. 
Then every set\/ $\pu\ne Y\sq X$, clopen in\/ $X$, 
belongs to\/ $\pele\i$ as well.
\ele
\bpf
We argue by induction on $\lh\i$. 
If $\lh\i=1$ then $\ilq\i=\ans\i$, and hence $\pel\i$ 
is essentially the family of all perfect sets 
$P\sq\cam$. 
Thus we can refer to the fact that a clopen 
subset of a perfect set is perfect, too. 

Now suppose that $\lh\i=\ell\ge2$, and let 
$\j=\i\res{(\ell-1)}$. 
%
%
By Lemma~\ref{997}, it suffices to consider the case   
$X=\can{\ilq\i}$, so that let $Y\sq\can{\ilq\i}$ be clopen. 
By a simple topological argument, $Y$ has the form
$Y=\bigcup_{k<n}(U_k\ti P_k)$, where all $U_k\sq\can{\ile\i}$ 
are clopen and pairwise disjoint, 
and $P_k\sq\can{}$ are clopen, so that there are perfect 
trees $T_k$ satisfying $P_k=[T_k]$. 

On the other hand, the set 
$Y'=Y\rsl\i=\bigcup_{k<n}U_k$ belongs to $\pel\i=\pele\j$ 
by the inductive hypothesis, and 
the map $\cT(y)=T_k$ in case $y\in U_k$ is continuous. 
It remains to apply Lemma~\ref{hh}.
\epf

\ble
\lam{darop}
If\/ $\et\sq\za$ belong to\/ $\cpo$, $X\in\pe\za$, 
and\/ $U\sq X$ is clopen in\/ $X$ then\/   
$U\dar\et$ is clopen in\/ $X\dar\et$.
\ele

\bpf 
By Lemma~\ref{997}, it suffices to prove the result 
for $X=\can\za$, in  
which case the result is obvious.
\epf

\ble
\lam{clop}
If\/ $\za\in\cpo$, $X\in\pe\za$, 
$U\sq X$ is open in\/ $X$, 
and\/ $x_0\in U,$ then there is a set\/  
$X'\in\pe\za$, $X'\sq U,$ clopen in\/ $X$ 
and containing\/ $x_0$.
\ele

\bpf 
By Lemma~\ref{997}, it suffices to prove the result 
for $X=\can\za.$ 
Note that if $x_0\in X'\sq\can\za$ 
and $X'$ is open in $\can\za$ then there exists a 
basic clopen set $C\sq X'$ containing $x_0.$ 
({\it Basic clopen sets\/} are those of the form
\dm
C=\ens{x\in\can\za}{u_1\subset x(\i_1)\land
\ldots\land u_m\su x(\i_m)},
\dm
where $m\in\om,$ $\i_1,...,\i_m\in\za$ are pairwise different, 
and $u_1,...,u_m\in 2^{<\om}.$) 
One easily proves that every 
set $C$ of this type actually belongs to $\pe\za$.
\epf

\ble
\lam{aclo}
Suppose that tuples\/ 
$\j\su\i$ belong to\/ $\tup$, 
$X\in\pele\i$, $Y\in\pele\j$, $Y\sq X\rsd{\sq\j}$,  
and\/ $Z=X\cap (Y\usq\i)$. 
Let\/ $\pu\ne Z'\sq Z$ be clopen in\/ $Z$. 
Then there exist sets\/ $X'\sq X$ and\/ $Y'\sq Y$, 
 clopen in resp.\ $X,\,Y$, such that\/ 
$Y'\sq X'\rsd{\sq\j}$,  
and\/ $Z'=X'\cap (Y'\usq\i)$.   
\ele

Under the conditions of the lemma, note that 
$Z\in\pe\xi$ by Lemma~\ref{apro}, whereas  
$X',Z'\in\pe\xi$, $Y'\in\pe\et$ by Lemma~\ref{lin}.

\bpf
By the compactness, there is a 
set $C\sq X$,  clopen in $X$, such that 
$Z'=Z\cap C$. 
Put $X'=C$. 
To define $Y'$, note that $C'=C\rsd{\sq\j}$ 
is clopen in $X\rsd{\sq\j}$ by Lemma~\ref{oldf}\ref{oz}. 
Therefore $Y'= Y\cap C'$ is clopen in $Y.$
\epf

\ble
\lam{99y}
If\/ $X,Y\in\pe\za$, $\et\sq\za$ belong to\/ $\cpo$, 
$\i\in\za\bez\et$, 
and\/ $X\dar\et=Y\dar\et$, 
then there exists\/ $k<\om$ and 
sets\/ $X',Y'\in\pe\za$, $X'\sq X$, $Y'\sq Y$, 
clopen in resp.\/ $X,Y$ and such that\/ 
$X'\dar\et=Y'\dar\et$, and\/ $x(\i)(k)=0$ but\/ 
$y(\i)(k)=1$ for all\/ $x\in X'$ and\/ $y\in Y'$, 
or vice versa.
\ele

\bpf 
By Lemma~\ref{99x}, there are points $x_0\in X$, $y_0\in Y$ 
with $x_0\dar\et=y_0\dar\et$ but, for some $k$, 
$x_0(\i)(k)=0$ while $y_0(\i)(k)=1$ (or vice versa). 
By Lemma~\ref{clop}, there is a set 
$A\in\pe\za$, $x_0\in A\sq X$, clopen in $X$, 
and such that $x(\i)(k)=0$ for all $x\in A$.
Then $A\dar\et$ is clopen in $X\dar\et$ by Lemma~\ref{darop}. 

Note that $x_0\dar\et\in A\dar\et$ by construction, 
therefore $y_0\dar\et\in A\dar\et$ as well. 

Furthermore, $B=\ens{y\in Y}{y\dar\et\in A\dar\et}$ 
is clopen in $Y$, and $y_0\in B$. 
Still by Lemma~\ref{clop}, there is a set 
$Y'\in\pe\za$, $y_0\in Y'\sq B$, clopen in $Y$, 
and such that $y(\i)(k)=1$ for all $y\in Y'$.

It remains to define $X'=((Y'\dar\et)\uar\za)\cap A$ and apply 
Lemma~\ref{darop} to check that $X'$ is clopen in $X$, 
and Lemma~\ref{apro} to check that $X'\in\pe\za$.
\epf

\bcor
\lam{99z1}
If\/ $X\in\pe\za$, 
and\/ $\i\ne\j$ belong to\/ $\za$, 
then there exists\/ $Z\in\pe\za$, $Z\sq X$, 
clopen in\/ $X,$ and such that\/ 
$(Z\dir\i)\cap(Z\dir\j)=\pu$. 
\ecor
\bpf
Let say $\j\not\sq\i$, so that $\i\nin\et=\za\nlq\j$. 
Lemma~\ref{99y} (with $X=Y$) yields relatively clopen 
sets $X',Y'\sq X$ in $\pe\za$ with $X'\res\et=Y'\res\et$, 
and $k<\om$, such that $x(\i)(k)=0$ for all $x\in X'$ and 
$x(\i)(k)=1$ for all $x\in Y'$. 

Now note that $U=X'\res\et=Y'\res\et\in\pe\et$ by 
Lemma~\ref{less}, and $U$ is clopen in $X\res\et$ 
by Lemma~\ref{darop}. 
Lemma~\ref{clop} implies that there is a relatively 
clopen $V\sq U$, $V\in\pe\et$, such that either 
(0) $u(\j)(k)=0$ for all $u\in V$ or  
(1) $u(\j)(k)=1$ for all $u\in U$. 
Let say (1) hold. 
Then the set $Z=X'\cap(V\uar\za)\sq X$ belongs to 
$\pe\za$ by Lemma~\ref{apro}, is clopen in $X$, and 
if $x\in Z$ then $x(\j)(k)=1$ but  $x(\i)(k)=0$ 
by construction, as required. 
\epf

We leave the proof of the following generalization 
of \ref{99y}/\ref{99z1} to the reader; 
it is rather routine and similar to the above.

\ble
\lam{99z2}
Let\/ $X,Y\in\pe\za$, 
$\et\sq\za$ belong to\/ $\cpo$, 
$X\dar\et=Y\dar\et$, 
$\i,\j\in\za$, 
and either\/ $\i\ne\j$ or\/ $\i=\j\nin\et$. 
Then there is\/ $k<\om$ and 
sets\/ $X',Y'\in\pe\za$, $X'\sq X$, $Y'\sq Y$, 
clopen in resp.\  $X,Y,$ and such that still\/ 
$X'\dar\et=Y'\dar\et,$ and\/ $x(\i)(k)=0$ but\/ 
$y(\j)(k)=1$ for all\/ $x\in X'$, $y\in Y'$, 
or vice versa.\qed
\ele

\parf{Vertical splitting}
\las{vs}

{\ubf We still argue in $\rL$.} 
Given $\i\in\za\in\cpo$, and  
a set $X\in\pe\za$, we are going to split 
$X$ into a disjoint union $U\cup V$ of sets 
in $\pe\za$ such that 
$U\dar{\za\nlq\i}=V\dar{\za\nlq\i}=X\dar{\za\nlq\i}$, 
and in the same time, if $y\in X\rsl\i$ then the 
cross-sections $\cs Uy \i$, $\cs Vy \i$ 
have strictly smaller size than 
$\cs Xy \i=\ens{x(\i)}{x\in X\land x\rsl{\i}=y}$. 

Still assuming that $\i\in\za\in\cpo$,   
$X\in\pe\za$, and $y\in X\rsl\i$, recall that   
$P=\cs Xy\i$ 
is a perfect set in $\can{}=2^\om$ by 
Lemma~\ref{oldf}\ref{perf}. 
It follows that there is a unique tuple 
$u=\cu Xy\i\in\bse$ 
\kmar{cu Xyi}%
\index{zDXyi@$\cu Xy\i$}%
\index{splitting!uXyi@$\cu Xy\i$}%
of length $m=\lh u=\cm Xy\i\in\bse,$ 
\kmar{cm Xyi}%
\index{zmXyi@$\cm Xy\i$}%
\index{splitting!mXyi@$\cu Xy\i$}%
such that $u\su p$ for all $p\in P=\cs Xy\i$, 
and in the same there exist $p_0,p_1\in P$ with 
$p_0(m)=0$ and $p_1(m)=1$.  
We let, for $e=0,1$, 
$$
\spli X\i e\,=\,
\kmar{spli Xie}%
\index{zX-ie@$\spli X\i e$, splitting}%
\index{splitting!X-ie@$\spli X\i e$}%
\ens{x\in X}{x(\i)(\cm Xy\i)=e}. 
$$ 

\ble
\lam{451}
Let\/ $\i\in\za\in\cpo$, $X\in\pe\za$, 
$X_e=\spli X\i e$, $e=0,1$. 
Then 
\ben
\renu
\itlb{451a}
the sets\/ $X_e$ belong to\/ $\pe\za$ 
and are clopen in\/ $X$, 
$X=X_0\cup X_1$, $X_0\rsq\i\cap X_1\rsq\i=\pu$, 
$X_0\dar{\za\nlq\i}=X_1\dar{\za\nlq\i}=X\dar{\za\nlq\i};$ 

\itlb{451b}
if\/ $y\in X\rsl\i$ then\/ 
$\cm{X_0,}y\i>\cm Xy\i$, 
$\cm{X_1,}y\i>\cm Xy\i$ strictly$;$

\itlb{451c}
if\/ $\ta\in\cpo$, $\i\in\ta\sq\za$, $Z=X\dar\ta$, 
$Z_e=\spli{Z}\i e$, then\/ 
$Z_e=X_e\dar\ta$ and\/ 
$X_e=X\cap(Z_e\uar\za)$. 
\een
\ele

\bpf
Claims \ref{451b}, \ref{451c} hold by construction.
 
Claim \ref{451a} {\ubf Case 1:} $\za=\ilq\i$ 
(in other words, $\i$ is the largest tuple in $\za$). 
By Lemma~\ref{oldf}\ref{oz}, if $u\in\bse$ then 
$S_u=\ens{y\in X\rsl \i}{\sus p\in\cs X y\i(u\su p)}$ 
is a set relatively clopen in $Y=X\rsl\i$. 
Therefore
\bce
$Y_u=(S_{u\we0}\cap S_{u\we1})\bez
\bigcup_{v\in2^m,v\ne u}S_v$, where $m=\lh u$
\ece
is clopen in $Y$ as well. 
Therefore, by the compactness of the spaces considered, 
the set $A=\ens{u\in\bse}{Y_u\ne\pu}$ is finite. 
It follows that, for $e=0,1$, 
\bce
$X_e:=\spli X\i e=\bigcup_{u\in A}
\ens{x\in X}{x\rsl\i\in Y_u\land u\we e\su x(\i)}$ 
\ece
is clopen in $X,$ hence $X_e\in\pele\i$ by 
Lemma~\ref{lin}. 
The rest of claims is obvious in Case 1.

\ref{451a} {\ubf Case 2:} any $\za$. 
Let $Z=X\rsq\i$, $Z_e=\spli Z\i e$. 
Then 
$X_e=X\cap(Z_e\uar\za)$ 
by \ref{451c}.   
Apply the results of Case~1 for $Z$, 
and then Lemma~\ref{apro}.
\epf

\parf{Assembling sets from projections}
\las{assem}

{\ubf We still argue in $\rL$ in this section.} 

\ble
\lam{pe1}
Assume that\/ $\xi_0,\xi_1,\xi_2,\dots\in\cpo$,
$\vt=\bigcup_n\xi_n$, and\/   
$X\in\pe\vt$.  
Then\/ $X=\bigcap_n(X\car\xi_n\uar\vt)$.
In particular,
$X=\bigcap_{\i\in\vt}(X\rsdq\i\uar\vt)$.
\ele
\bpf
The relation $X\sq X'=\bigcap_n(X\car\xi_n)\uar\vt$
is obvious. 
To prove $X'\sq X$, consider the following cases.

{\it Case 1\/}: simply $\vt=\xi_0\cup\xi_1$.
Apply Lemma~\ref{no19}.

{\it Case 2\/}: 
$\vt=\xi_0\cup\xi_1\cup\ldots\cup\xi_n$.
Argue by induction using Case 1.

{\it Case 3\/}: general case. 
By the result for Case 2, we can \noo\ assume that
$\xi_n\sq\xi_{n+1}$ for all $n$.
Then apply the compactness. 
\epf

It follows by Lemma~\ref{pe1} that each set $X\in\pe\vt$ is
fully determined by the coherent system of its projections
$X\rsd{\sq\i}=X\car{\ilq\i}\in\pele\i$, where
$\i\in\vt$ and ${\ilq\i}=\ens{\j\in\tup}{\j\sq\i}$.
The next lemma shows that conversely any coherent
system of iterated perfect sets results in a set
in $\pe\vt$.


\ble
\lam{pe4}
Let\/ $\xi_0,\xi_1,\xi_2,\dots\in\cpo$, 
$\vt=\bigcup_n\xi_n$, and  
sets\/ $X_n\in\pe{\xi_n}$ satisfy the 
{\em coherence condition} 
\ben
\fenu
\itlb{pe4*}
$X_n\car{(\xi_k\cap\xi_n)}=X_k\car{(\xi_k\cap\xi_n)}$ 
\index{coherence}%
for all $k,n$. 
\een
Then\/ $X= \bigcap_n(X_n\uar\vt)$
belongs to\/ $\pe\vt$,  and\/
$X\car\xi_n=X_n$, $\kaz n$.

In particular, if\/ $\xi_0,\xi_1,\xi_2,\dots$ 
are pairwise disjoint, then\/ \ref{pe4*} 
holds by default, hence \/ $X=\bigcap_n(X_n\uar\vt)$ 
belongs to\/ $\pe\vt$ and\/ $X\car\xi_n=X_n$, $\kaz n$.
\ele
\bpf
By Corollary~\ref{pe2}, we  
\noo\ assume that $\xi_0\sq\xi_1\sq\xi_2\sq\dots$.
Lemma~\ref{99} yields a sequence of \pkh s 
$H_n:\can{\xi_n}\onto X_n$ s.\,t.\
$H_{n+1}(x)\car\xi_n=H_n(x\dar\xi_n)$ for all $n$
and $x\in\can{\xi_{n+1}}$.
This allows us to define a \pkh\ 
$H:\can{\vt}\onto X$ by simply
$H(x)\car\xi_n=H_n(x\dar\xi_n)$ for all $n$
and $x\in\can{\vt}$.
\epf

The lemma leads to another representation of 
iterated perfect sets. 
Let $\vt\in\cpo$. 
If $X\sq\can\vt$ then the system of projections 
$X\rsq\i$, $\i\in\vt$, will be called 
\rit{the projection tree} of $X.$ 
Generally, a \rit{projection tree} is any system 
of sets $X_\i$, $\i\in\vt$, satisfying the 
{\em coherence condition} 
\index{coherence}%
\index{projection tree}%
in the form
\ben
\fenu
\atc
\itlb{cohc}
$X_\i\sq\can{\ilq\i}$, and if $\i\su\j$ belong to 
$\vt$ then $X_\i=X_\j\rsq\i$.
\een

\bcor
[of Lemma~\ref{pe4}]
\lam{pe5}
Let $\vt\in\cpo$. 
If\/ $X\in\pe\vt$ then the system of sets\/ 
$X\rsq\i$, $\i\in\vt$, satisfies\/ \ref{cohc}, and\/ 
$X= \bigcap_{\i\in\vt}(X\rsq\i\uar\vt)$.

Conversely, if sets\/   
$X_\i\in\pele\i$ satisfy\/ \ref{cohc} 
{\em(\ie, form a coherent projection tree)}, then\/
$X= \bigcap_{\i\in\vt}(X_\i\uar\vt)\in\pe\vt$ and\/
$X\rsdq\i=X_\i$ for all $\i$.\qed  
\ecor

Thus sets in $\pe\vt$ are in natural 1-1 correspondence 
with coherent projection trees of sets $X_\i\in\pele\i$.

\parf{Permutations}
\las{perm}

 
Let $\pero$ be the group of all bijections 
\index{zzPerm@$\pero$}%
\index{permutations!$\pero$}%
\kmar{pero}%
$\pi:\tup\onto\tup$, $\pi\in\rL$,
$\su$-\rit{invariant} in the sense that 
$\i\su\j\eqv \pi(\i)\su\pi(\j)$ for all $\i,\j\in\tup$. 
Thus $\pero\in\rL$. 
Bijections $\pi\in\pero$ will be called \rit{permutations}. 
Any $\pi\in\pero$ is \rit{length-preserving}, so that 
$\lh\i=\lh{\pi(\i)}$ for all $\i\in\xi$, 

The superposition 
\kmar{supo}%
\index{superposition, $\supo$}%
\index{zzzcirc@$\supo$, superposition}%
$\supo$ is the group operation: 
$(\pi\supo\ro)(\i)=\pi(\ro(\i))$. 

To define an important subgroup of $\pero$, 
recall that every ordinal $\al$ can be represented 
in the form  $\al=\la+m$, where 
$\la\in\Ord$ is a limit ordinal and $m<\om$; 
then $\al$ is called {\em odd}, resp., {\em even}, 
\index{ordinal!even, odd}%
\index{ordinal!parity}%
\index{even}%
\index{odd}%
\index{parity}%
if the number $m$ is odd, resp., even. 
A tuple $\i=\ang{\al_0,\dots,\al_k}\in\tup$ 
is {\em odd}, resp., {\em even},  
\index{tuple!even, odd}%
\index{even ordinal, tuple}%
\index{odd ordinal, tuple}%
if such is the last term $\al_k$.
\kmar{ekp}
If $\i,\j\in\tup$ then $\i\ekp\j$ will mean that
\index{zzzequiv@$\ekp$}%
$\lh\i=\lh\j$ and if $k<\lh\i$ then the ordinals 
$\i(k)$ and $\j(k)$ have the same parity. 

\bre
\lam{131}
Odd and even tuples will play different roles in the model 
construction below. 
Namely, even tuples will be involved in the coding procedures, 
whereas the role of odd tuples will be to obscure things 
enough so that the desired counterexamples will not be 
available at levels of the hierarchy lower than prescribed.
\ere

Let $\per$ be the subgroup of all
\index{zzPi@$\per$}%
\kmar{per}%
permutations $\pi\in\pero$, such that 
$\i\ekp\pi(\i)$ for every $\i\in\tup$, that is, 
\rit{parity-preserving} permutations.
\index{parity-preserving}%
\index{permutations!parity-preserving, $\per$}%
\index{permutation groups!$\per$, parity-preserving}%

\bpri
\lam{pi*ij}
Suppose that $\i,\j\in\tup$, $\lh\i=\lh\j$. 
Define a permutation $\pi=\pi_{\i\j}\in\pero$ 
\index{permutation!$\pi_{\i\j}$}%
\index{zzpiij@$\pi_{\i\j}$}%
satisfying $\pi(\i)=\j$ as follows. 
Let $\k\in\tup$. 

If $\k(0)\nin\ans{\i(0),\j(0)}$ then put $\pi(\k)=\k$. 

If $\k(0)=\i(0)$ then there is 
a largest number $1\le m\le \lh\i=\lh\j$ such that 
$\k\res m=\i\res m$. 
Then $\k=(\i\res m)\we \k'$ 
(concatenation of tuples) 
\index{tuples!concatenation}%
for some tuple $\k'\in\tup\cup\ans\La$. 
Put $\pi(\k)=(\j\res m)\we \k'$.  

Similarly, if $\k(0)=\j(0)$ then there is 
a largest number $1\le m\le \lh\i=\lh\j$ such that 
$\k\res m=\j\res m$. 
Then accordingly $\k=(\j\res m)\we \k'$  
for some $\k'\in\tup\cup\ans\La$. 
Put $\pi(\k)=(\i\res m)\we \k'$.  

Easily $\pi\in\pero$, $\pi\obr=\pi$, $\pi(\i)=\j$, 
and if $\i\ekp\j$ then $\pi\in\per$.
\epri

\vyk{
For instance if tuples $\i\ekp\j$ belong to $\tup$, 
$n=\lh\j=\lh\i$, then 
$\pi=\pi_{\i\j}\in\pep{\ilq\i\cup\ilq\j}$ is defined 
by $\pi(\i\res\ell)=\pi\obr(\i\res\ell)=\j\res\ell$ 
for all $\ell<n$. 
Accordingly, the 
{\em canonical\/} \dd\sq preserving 
and parity-preserving permutation  
$\pi'_{\i\j}\in\rer$, 
satisfying $\pi'_{\i\j}(\i)={\pi'_{\i\j}}\obr(\i)=\j$, 
is defined from $\pi_{\i\j}$ 
as in the proof of the lemma.
\index{zzpiij@$\pi_{\i\j}$}%
\index{zzpi'ij@$\pi_{\i\j}$}%
}

{\ubf Actions.} 
Any permutation $\pi\in\pero$ induces a 
transformation left-acting on several types of 
\index{permutations!action $\akt$}%
\index{action $\akt$}%
\index{zzzakt@$\;\akt\;$, action}%
objects as follows. 
\bit
\item
If $\xi\in\cpo$, or generally $\xi\sq\tup$, 
then $\pi\akt\xi:=\ima\pi\xi=\ens{\pi(\i)}{\i\in\xi}$. 
\kmar{akt}%

\item
If $\xi\sq\tup$ and $x\in\can\xi$  
then $\pi\akt x\in\can{{\pi\akt\xi}}$ is defined by 
%
$(\pi\akt x) (\pi(\i))=x(\i)$ for all $\i\in\xi$.
That is, formally $\pi\akt x=x\supo\pi\obr$, the 
superposition. 

\item
If $\xi\sq\tup$ and $X\sq\can\xi$ then 
$\pi\akt X:=\ens{\pi\akt x}{x\in X}\sq\can{{\pi\akts\xi}}$. 

\item
If $G\sq\pei$ then $\pi\akt G:=\ens{\pi\akt X}{X\in G}$.  
\eit

\vyk{
The group action property holds in all cases, 
for instance: 
\bce
$\ro\akt(\pi\akt\xi)=
\ima\ro{(\ima\pi\xi)}= 
\ima{(\ro\supo\pi)}{\xi}= 
{(\ro\supo\pi)}\akt{\xi}; 
$\\[0.7ex]
$
\ro\akt{(\pi\akt x)}=
{(x\supo\pi\obr)}\supo\ro\obr=
x\supo{(\pi\obr\supo\ro\obr)}=
x\supo{(\ro\supo\pi)\obr}=
(\ro\supo\pi)\akt x.\vhm
$
\ece
}


\ble
\lam{pi1}
Let\/ $\pi,\rho\in\pero$, $\et\in\cpo$, and\/ $v\in\can\tup$. 
Then\/ 
\ben
\renu
\itlb{pi1i}
$\pi\akt(\rho\akt v)=(\pi\supo\rho)\akt v$ --- 
the group action property$;$ 

\itlb{pi1ii}
$(\pi\akt v)\dar(\pi\akt\et)=\pi\akt(v\dar\et)$, 
equivalently, 
$(\pi\akt v)\dar\et=\pi\akt(v\dar(\pi\obr\akt\et))$. 
\een
\ele
\bpf
$\pi\akt(\rho\akt v)=(v\supo\rho\obr)\supo\pi\obr=
v\supo(\pi\supo\rho)\obr=
(\pi\supo\rho)\akt v$.
\epf

Thus in general 
$\pi\akt(v\dar\et)=(\pi\akt v)\dar(\pi\akt\et)$ 
is not equal 
to $(\pi\akt v)\dar\et$\ ! 

\ble
\lam{pi2}
If\/ $\pi\in\pero$ and\/ $X\in\pe\xi$ then 
$\pi\akt X\in\pe {\pi\akts\xi}$. 

Moreover\/ $\pi$ is an 
$\sq$-preserving automorphism of\/ $\pei$.\qed 
\ele

\sekt{Splitting/fusion construction }
\las{sfc}

{\ubf We argue in $\rL$ in this chapter.} 

We'll make use of a construction of sets in $\pe\za$ as 
$X=\bigcap_{m\in\om}\bigcup_{u\in 2^m} X_u$
where all $X_u$ belong to $\pe\za$ and
$2^m$ = all 0,1-tuples of length $m.$ 
The technique is well-known for usual perfect sets 
in Polish spaces. 
This chapter presents the modification,
originally developed for the case of iterated perfect 
sets in \cite{fm97,jsl99}. 
We outline some applications as well.

\parf{Splitting systems}
\las{spl}

First of all let us specify requirements which imply an 
appropriate behavior of a system 
of sets $X_u\in\pe\za,u\in 2^m,$ 
with respect to projections.
We need to determine, for any pair of tuples 
$u,\,v\in 2^m\msur$ $(m<\om),$
the largest  initial segment $\xi=\za[u,v]$ of $\za$
such that the projections 
$X_u\res\xi$ and $X_v\res\xi$ have to be equal,
to maintain the construction in proper way.

%

\vyk{
Let's order $\za$ by $<_\za$
so that shorter tuples \dd{<_\za}precede
longer ones, whereas within the same length the
lexicographic order acts.
Let $\phi':N\onto\za$ be the corresponding
order-preserving bijection, where $N=\card \za$,
and let $\phi=\phi_\za:\om\onto\za$ be defined
by $\phi(kN+j)=\vpi'(j)$ for all $j<N$ and $k<\om$.
}

Assume that $\za\in\cpo$ and $\phi:\om\to\tup$ is any map, 
not necessarily $\phi:\om\to\za$.
We define, for any pair of tuples 
$u,\,v\in 2^m\yi m<\om$, an initial segment
\pagebreak[0]
\dm
\left.
\bay{lcl}
\za_\phi[u,v] & = & 
\index{zzyafiuv@$\za_\phi[u,v] $}%
{\textstyle\bigcap_{\,l<m,\;u(l)\not=v(l)\,}}
\za\nlq{\phi(l)} 
\;\;\;= \\[2mm]
& = & 
\ans{\j\in\za:\neg\;\exists\,l<m 
\left(\vphantom{X^{X^x}_{X_x}}
u(l)\not=v(l)\cj \phi(l)\sq\j
\right)}
\eay
\right\}
\in\IS_\za\,.
\dm

\bdf
\lam{splis}
Let still $\za\in\cpo$ and $\phi:\om\to\tup$.
A \dd\phi{\it \cohe\ system\/} 
\index{split system}%
\index{system!\cohe\ system}%
\index{system!height}%
(rather \dd{(\phi\res m)}\cohe\  
as the notion depends only on $\phi\res m$) 
in $\pe\za$, of \rit{height} $m$,
is a  family $\sis{X_u}{u\in 2^m}$ of 
sets $X_u\in \pe\za$ satisfying,
for all $u,\,v\in 2^m:$
\ben
\def\theenumi%
{{\rm\hskip0.0ex{\rm S\arabic{enumi}}\hskip0.0ex}}
\def\labelenumi{\theenumi{:}}
\itlb{prct}\msur  
$\komp{X_u}{\za_\phi[u,v]}=\komp{X_v}{\za_\phi[u,v]}$ \ 
(\rit{projection-coherence}),\hfill {\it and\/}
\hfill$\,$

\itlb{bprct}   
if $k<m$, $\sg\in\cpo$, $\sg\sq\za$, but 
$\sg\not\sq\za_\phi[u,v]$ then 
$({X_u}\dar{\et}) \cap ({X_v}\dar{\et})=\pu$.
\een
If in addition the following strengthening of 
\ref{bprct} holds, then $\sis{X_u}{u\in 2^m}$ will be 
a \rit{strong} \dd\phi \cohe\ system:
\ben
\def\theenumi%
{{\rm\hskip0.0ex{\rm S\arabic{enumi}}\hskip0.0ex}}
\def\labelenumi{\theenumi{:}}
\atc\atc
\itlb{aprct}   
if $k<m$, $\i=\phi(k)\in\za\bez\za_\phi[u,v],$ 
then 
$({X_u}\dir{\i}) \cap ({X_v}\dir{\i})=\pu$ --- and then 
%
$X_u\rsq{\j} \cap X_v\rsq{\j}=\pu$ for all 
$\j\in\za\bez\za_\phi[u,v]$. 
\vyk{
\itlb{aprct}\  
if $\i\in\za\setminus\za_\phi[u,v]$ then 
$X_u\rsq{\i} \cap X_v\rsq{\i}=\pu$, or equivalently, 
if $k<m$,\\ 
\hspace*{0.6ex}$u(k)\ne v(k)$, and $\i=\phi(k)\in\za$ then  
$X_u\rsq{\i} \cap X_v\rsq{\i}=\pu$. 
}%
\een
We proceed with a few related definitions.
\ben
\Aenu
\itlb{splisA}
A system $\sis{X'_{u}}{u\in 2^{m}}$    
{\it \nar s\/}  
\index{split system!\nar s}%
\index{narrows@\nar}%
\kmar{\ \mns nar}%
%
$\sis{X_u}{u\in 2^m}$  
if $X'_{u}\sq X_u$ for all $u$, and 
a \rit{clopenly} \nar s, if in addition 
each $X'_{u}$ is clopen in $X_u$.

\itlb{splisB}
A system $\sis{X_{u'}}{u'\in 2^{m+1}}$ is an 
{\it \expa\/} of 
\index{split system!\spl}%
\index{expansion@\expa}%
\kmar{\ \mns expa}%
%
$\sis{X_u}{u\in 2^m}$ 
iff we have $X_{u\we e}\sq \spli{(X_u)}\i e$ 
for all $u\in 2^m$ and $e=0,1$, where $\i=\phi(m)$, 
and a \rit{clopen} \expa, if in addition 
each $X_{u\we e}$ is clopen 
in $X_u$.

\itlb{splisC}
A system $\sis{Y_{u}}{u\in 2^{m}}$ of sets 
$Y_u\in\pe\vt$, where $\za\sq\vt\in\cpo$, is a  
{\it \lif\/} of $\sis{X_u}{u\in 2^m},$ 
\index{split system!\lif}%
\index{lifting@\lif\ of split systems}%
\kmar{\ \mns lif}%
iff $Y_{u}\dar\za\sq X_u$ for all $u\in 2^m$, 
and a \rit{clopen} \lif, if in addition 
each $Y_{u}\dar\za$ is clopen 
in $X_u$.
\index{split system!\lif}%
\index{lifting@\lif}%
\qed
\een
\eDf

{\ubf A set\/ $\za\in\cpo$, and\/ $\phi:\om\to\tup$,  
remain fixed in the following lemmas.} 

\ble
\lam{disj}
Let\/ $\sis{X_u}{u\in 2^m}$ be a\/
system in\/ $\pe\za$ satisfying\/ 
\ref{prct} and\/ \ref{bprct}, 
and\/ $u,v\in 2^m.$  
Then either\/ $X_u=X_v$ 
or\/ $X_u\cap X_v=\pu$.
\ele
\bpf
If $\za_\phi[u,v]=\za$ then $X_u=X_v$ by \ref{prct}. 
If $\i\in\za\bez\za_\phi[u,v]$ then 
$({X_u}\rsq{\i}) \cap ({X_v}\rsq{\i})=\pu$ by 
\ref{bprct}, and hence    
$X_u\cap X_v=\pu$.   
\epf

The next lemma 
proves that any \cohe\ system admits a
\nar ing that honors a shrink of one of
its sets to a given smaller set in $\pe{}$. 

\ble
\lam{suz}
Let\/ $\sis{X_u}{u\in 2^m}$ be a\/
system in\/ $\pe\za$ satisfying\/ \ref{prct}, 
$u_0\in 2^m,$ $X\in\pe\za\yi X\sq X_{u_0}.$
Then the sets\/  
$Y_u={X_u\cap (X\dar\za_\phi[u,u_0]\uar\za)}$, 
$u\in 2^m,$  
belong to\/ $\pe\za$,
and the system\/ $\sis{Y_u}{u\in2^m}$ 
\nar s\/ $\sis{X_u}{u\in 2^m}$ 
and satisfies\/ \ref{prct} and\/  
$Y_{u_0}=X$ {\sf(since $\za_\phi[u_0,u_0]=\za$)}. 


If the given set\/ $X$ is clopen in\/ $X_{u_0}$, then 
each\/ $Y_u$ is clopen in\/ $X_u$.
\ele

\bpf
The sets $Y_u$ 
belong to $\pe\za$ by Lemma~\ref{apro}, because each
$X\dar\za_\phi[u,u_0]$ belongs to $\pe{\za_\phi[u,u_0]}$
by Lemma~\ref{less} (since
$X\dar\za_\phi[u,u_0]\sq X_u\dar\za_\phi[u,u_0]$). 
The clopenness claim follows from Lemma~\ref{darop}.

That the system $\sis{Y_u}{u\in2^m}$ satisfies 
\ref{prct} see the proof of Lemma 12 in  \cite{jsl99}.
\epf

There is a remarkable strengthening of the lemma.

\bcor
\lam{suz+}
Under the assumptions of Lemma~\ref{suz}, 
if in addition\/ $u_1\in2^m,$ $Y\in\pe\za$, $Y\sq X_{u_1}$,
$Y\dar\za_\phi[u_0,u_1]=X\dar\za_\phi[u_0,u_1]$, 
then the sets\/  
\bce
$Z_u={X_u\cap (X\dar\za_\phi[u,u_0]\uar\za)
\cap (Y\dar\za_\phi[u,u_1]\uar\za)}$, 
$u\in 2^m,$  
\ece
belong to\/ $\pe\za$,
and the system\/ $\sis{Z_u}{u\in2^m}$ 
\nar s\/ $\sis{X_u}{u\in 2^m}$ 
and satisfies\/ \ref{prct} and\/  
$Z_{u_0}=X$, $Z_{u_1}=Y$. 

\vyk{
then there is a\/ 
system\/ $\sis{Z_u}{u\in 2^m}$  in\/ $\pe\za,$ 
which \nar s\/ $\sis{X_u}{u\in 2^m}$   
and satisfies\/ \ref{prct} 
and\/ $Z_{u_0}=X$, $Z_{u_1}=Y$. \ 
}


If\/ $X,Y$ are clopen in resp.\ $X_{u_0},X_{u_1}$, 
then each\/ $Z_u$ is clopen in\/ $X_u$.
\ecor

\bpf
The sets
$Y_u={X_u\cap (X\dar\za_\phi[u,u_0]\uar\za)}\in\pe\za$ 
form a \ref{prct}-system $\sis{Y_u}{u\in 2^m}$, 
which \nar s $\sis{X_u}{u\in 2^m},$ with $Y_{u_0}=X$, 
by Lemma~\ref{suz}.

{\it Note that\/ $Y\sq Y_{u_1}$.}
(Indeed $Y_{u_1}={X_{u_1}\cap (X\dar\et\uar\za)}$ 
by construction, but $Y\sq  X_{u_1}$ 
and $Y\dar\et=X\dar\et$.)
It remains to apply Lemma~\ref{suz} yet again,
because 
$Z_u={Y_u\cap (Y\dar\za_\phi[u,u_1]\uar\za)}$ 
by construction.
\epf

\ble
\lam{1to2}
Let\/ $\sis{X_u}{u\in 2^m}$ be a\/
system in\/ $\pe\za$ satisfying\/ \ref{prct}. 
There is a system\/ $\sis{Y_u}{u\in2^m}$ in\/  
$\pe\za$,   
which still satisfies\/ \ref{prct},  
clopenly \nar s\/ $\sis{X_u}{u\in 2^m},$ 
and satisfies\/ \ref{aprct}  
as well.
\ele

\bpf
Pick any pair of $u_0,v_0\in 2^{m},$ and let 
$\et=\za_\phi[u_0,v_0]$, so that 
$X_{u_0}\dar\et=X_{v_0}\dar\et$ by \ref{prct}. 
Let $\i=\phi(k)\in\za\bez\et$, $k<m$. 
By Lemma~\ref{99y}, there exist 
sets\/ $U,V\in\pe\za$, $U\sq X_{u_0}$, $V\sq X_{v_0}$, 
clopen in resp.\/ $Y_{u_0},Y_{v_0}$ and such that still 
$U\dar\et=V\dar\et$,  
but\/ ${U}\dir\i\cap{V}\dir\i=\pu$. 
By Corollary~\ref{suz+}, 
there is a system\/ $\sis{X'_u}{u\in 2^{m}}$ 
of sets $X'_u\in \pe\za,$ $X'_u\sq X_u$, clopen in $X_u$, 
which 
satisfies \ref{prct} and $X'_{u_0}=U$, $X'_{u_0}=V$, 
so that ${X'_{u_0}}\dir\i\cap {X'_{v_0}}\dir\i=\pu$.

\sloppy
Thus we have succeeded to clopenly narrow 
$\sis{X_{u}}{u\in 2^{m}}$ 
to a system $\sis{X'_u}{u\in 2^{m}}$ 
still satisfying \ref{prct}, and 
also satisfying 
\ref{aprct} for a given triple of $u_0,v_0\in 2^{m}$ 
and $\i=\phi(k)\in \za\bez\za_\phi[s_0,t_0],\,k<m$.
It remains to iterate this \nar ing construction for all 
such triples.
\epf

\vyk{
\bpf 
Each set $X'_u$ belongs to $\pe\za$ by lemmas \ref{less} 
and \ref{apro}. 
We have to check only requirement~\ref{prct}. Thus let 
$u,\,v\in 2^m$ and $\xi=\za_\phi[u,v].$ We prove that 
$X'_u\res\xi=X'_v\res\xi$. 
Let in addition $\za_u=\za_\phi[u,u_0]$ 
and $\za_v=\za_\phi[v,u_0].$ Then 
\dm
X'_u\res\xi=(X_u\res\xi)\,\cap\,(X_0\res(\xi\cap\za_u)\ares\xi)
\,,\;\;
X'_v\res\xi=(X_v\res\xi)\,\cap\,(X_0\res(\xi\cap\za_v)\ares\xi)
\dm
by Lemma~\ref{pro'}. Thus it remains to prove that 
$\xi\cap\za_u=\xi\cap\za_v$ (the ``triangle'' equality).
Assume to the contrary that \eg\ $\i\in\xi\cap\za_u$
but $\i\not\in\za_v.$ 
The latter means that $\i\geq \phi(l)$ in $\za$
for some $l<m$ such that $v(l)\not=u_0(l).$
But then either $u(l)\not=u_0(l)$ -- so 
$\i\not\in\za_u,$ or $u(l)\not=v(l)$ -- so $\i\not\in\xi,$ 
contradiction.
\epf
}


\vyk{
Let $\i\in\za$ and $X\in\pe\za.$ 
Let us say that a pair of sets 
$X_0,\,X_1\in\pe\za$ 
is an {\it\dd\i splitting\/} of $X$ iff \ 
${X_0\cup X_1\sq X},\msur$ 
${X_0\rsd{\not\geq \i}=X_1\rsd{\not\geq \i}},$ and 
${X_0\rsd{\<\i}\cap X_1\rsd{\<\i}=\emps}.$ 
The splitting will be 
called {\it complete\/} if ${X_0\cup X_1= X}$ --- 
in this case we have 
${X_0\rsd{\not\geq \i}=X_1\rsd{\not\geq\i}=X\rsd{\not\geq\i}}$. 

\bass
\label{splitt}
Let\/ $\i\in\za.$ 
Every\/ $X\in\pe\za$ admits a complete\/ 
\dd\i splitting.
\eass
\bpf 
If $X=\can\za$ then we define $X_e=\ans{x\in X:x(\i)(0)=e},$ 
$e=0,1.$ 
Lemma \ref{compos} extends the result to the general 
case.
\epf
}

The next two lemmas provide \expa s and \lif s. 

\ble
\lam{pand}
Any \cohe\ system\/ 
$\sis{X_u}{u\in 2^m}$ in\/ $\pe\za$ admits a clopen \expa\/ 
by the \cohe\ system\/ $\sis{Y_{s}}{s\in 2^{m+1}}$, 
where\/ $Y_{u\we e}=\spli{(X_u)}\i e$, $\i=\phi(m)$. 
\ele

\bpf
In view of Lemma~\ref{451}\ref{451a}, 
it suffices to establish 
\ref{prct} for the new system. 
Let $s=u\we e$, $t=v\we\ve$ be tuples in $2^{m+1},$ 
$\i=\phi(m)$, $\et=\za_\phi[u,v]$, $\sg=\za\nlq\i$,
$\xi=\za_\phi[s,t]$. 
The goal is to prove (*) $X_s\dar\xi=Y_s\dar\xi$. 

{\ubf Case 1:} $\xi\sq\sg$. 
Then 
$X_s\dar\xi=X_s\dar\sg\dar\xi=X_u\dar\xi=X_u\dar\et\dar\xi$ 
(here Lemma~\ref{451} is used for the middle equality), 
and accordingly 
$X_t\dar\xi=X_v\dar\et\dar\xi$. 
Yet $X_u\dar\et=X_v\dar\et$ by \ref{prct} for 
$\sis{X_u}{u\in 2^m}$. 
This yields (*). 

{\ubf Case 2:} $\xi\not\sq\sg$. 
This means $\i\in\et$, $e=\ve$, and $\xi=\et$. 
Then 
$X_s\dar\et=\spli{(X_u\dar\et)}\i e
=\spli{(X_v\dar\et)}\i e=X_t\dar\et$ 
(by Lemma~\ref{451}\ref{451c} and \ref{prct} 
for the given system), 
which implies (*) yet again since $\xi=\et$.
\epf

\vyk{To begin with, let 
$Y_{u\we e}=X_u$ for all\/ $u\in2^m$ and\/ $e=0,1$. 
The system $\sis{Y_{s}}{s\in 2^{m+1}}$ satisfies 
\ref{prct}. 
It remains to apply Lemma~\ref{1to2}.
}

\ble
\lam{uver}
Assume that\/ $\za\sq\vt$ belong to\/ $\cpo$,   
$\sis{X_u}{u\in 2^m}$ is a\/ $\phi$-\cohe\ system\/ 
in\/ $\pe\za$, 
and\/ $Y_{u}=X_u\uar\vt$ for all\/ $u\in2^m$. 
Then\/ $\sis{Y_u}{u\in 2^m}$ is a\/ $\phi$-\cohe\ 
system\/ in\/ $\pe\vt$.
\ele

\bpf
To prove \ref{prct} for $\sis{Y_u}{u\in 2^m},$
let $u,v\in2^m.$ 
It can be the case that $\za_\phi[u,v]\sneq\vt_\phi[u,v]$, 
but definitely $\za_\phi[u,v]=\za\cap\vt_\phi[u,v]$ holds. 
Therefore 
\bce
$Y_u\dar\vt_\phi[u,v]=X_u\dar\za_\phi[u,v]\uar\vt_\phi[u,v]$,
\; 
$Y_v\dar\vt_\phi[u,v]=X_v\dar\za_\phi[u,v]\uar\vt_\phi[u,v]$.
\ece
by Lemma~\ref{599} (with $W=\can\vt$). 
However $X_u\dar\za_\phi[u,v]=X_v\dar\za_\phi[u,v]$. 
\epf

\parf{Fusion sequences}
\las{fuz}

{\ubf We argue in $\rL$ in this section.} 

Given $\za\in\cpo$, a map $\phi:\om\to\tup$ is 
\index{admissible@\dd{\za}admissible map}%
\rit{\dd{\za}admissible}, if the preimage
$\phi\obr(\i)=\ens{k}{\phi(k)=\i}$ of every
$\i\in\za$ is infinite, and in addition if 
$\j\su\i=\phi(k)$ then $\j=\phi(\ell)$ for some 
$\ell<k$.
Yet we do \rit{not} assume that $\phi(k)\in\za$, 
$\kaz k$.

\bdf
\lam{fuzD}
Suppose that $\za\in\cpo$,
and $\phi:\om\onto\za$ is \dd{\za}admissible. 

An indexed family of sets $X_u\in\pe\za,$ $u\in 2^{<\om},$
is a \dd\phi{\it fusion sequence\/} in $\pe\za$ if,
\index{fusion sequence}%
for every $m\in \om,$ 
the subfamily $\sis{X_u}{u\in 2^m}$ is a
\dd\phi split system, 
expanded by $\sis{X_u}{u\in 2^{m+1}}$  
in the sense of Definition~\ref{splis}\ref{splisB}.
\edf

\vyk{
and
\ben
\def\theenumi{{\hskip0.1ex{S-\arabic{enumi}}\hskip0.1ex}}
\setcounter{enumi}{3}
\itlb{di} 
For any $\ve>0$ there exists $m\in\om$ such that 
$\diam X_u<\ve$ for all $u\in 2^m.$
(A Polish metric on $\can\za$ is assumed to be fixed.)
\qed
\een
}%

\bte
\lam{fut}
Under the assumption of Definition~\ref{fuzD}, let\/ 
$\sis{X_u}{u\in 2^{<\om}}$ be a\/
\dd\phi fusion sequence in\/ $\pe\za.$
Then\/ 
$X=\bigcap_{m\in\om}\bigcup_{u\in 2^m} X_u$ 
belongs to\/ $\pe\za$. 
\ete

\bpf
To begin with, prove that (*) if $a\in\dn$ then the 
intersection $F(a)= \bigcap_m\bigcup X_{a\res m}$ 
is a singleton. 
Indeed if $\i\in\za$ and $m<\om$ then let $\ka^\i_m$ be 
equal to the number of all $k<m$ such that $\phi(k)=\i$. 
Thus if $u\in2^m$ and $y\in (X_u)\rsl\i$ then we have 
$\cm{X_u,}\i y\ge \ka^\i_m$ by construction. 
Now, as $\ka^\i_m\to\iy$ with $m\to\iy$, 
the set $\cs{F(a),}\i y$ is a singleton for any 
$y\in F(p)\rsl\i$. 
This implies (*) because $\za$ is well-founded. 
Thus $F(a)=\ans{f(a)}$, where $f:\dn\to\can\za$ 
is continuous, still by compactness. 

Having (*) established, we can then follow the proof 
of Theorem~14 in \cite{jsl99}. 
(Note that (*) was established in \cite{jsl99} under 
different suppositions, because the well-foundedness of 
$\za$ was not assumed there.) 
Namely we let $D=\can\za$, and then define $D_u$ by 
induction on $u\in\bse$ so that 
$D_{u\we e}=\spli{(D_u)}\i e$, where $\i=\phi(m)$ and 
$m=\lh u$. 
Then $\sis{D_u}{u\in\bse}$ is a \dd\phi fusion sequence 
in $\pe\za$ by Lemma~\ref{pand}. 

Thus similarly to (*) there is a continuous map 
$d:\dn\to\can\za$ such that for any $a\in\dn,$ 
$\bigcap_m D_{a\res m}=\ans{d(p)}$. 
Moreover, by the equality $X=X_0\cup X_1$ 
of Lemma~\ref{451}, we have $\ran d=\can\za$, 
hence $d\obr:\can\za\onto\dn$ is continuous.

If $a,b\in\dn$ then define   
$\za_\phi[a,b]=\bigcap_{m<\om}\za_\phi[a\res m,b\res m].$ 
(Note that $\za_\phi[a,b]=\za$ iff $a=b.$) 
We conclude from \ref{prct} and \ref{bprct} that\vtm

\noi
$(\dag)
\left\{
\bay{rl}
\left.
\bay{rcll}
x_a\res \za_\phi[a,b] & = & x_b\res \za_\phi[a,b] & 
\hbox{and}\\[1mm] 
d_a\res \za_\phi[a,b] & = & d_b\res \za_\phi[a,b] & 
\eay \right\} & 
\hbox{ for all }\; 
a,\,b\in \dn\\[3ex]

x_a\rsd{\<\i}\not=x_b\rsd{\<\i}\;\hbox{ and }\;
d_a\rsd{\<\i}\not=d_b\rsd{\<\i}
\hphantom{\left.
\vphantom{
\bay{rcll}
x_a\res \za_\phi[a,b] & = & x_b\res \za_\phi[a,b] & 
\hbox{and}\\[1mm] 
d_a\res \za_\phi[a,b] & = & d_b\res \za_\phi[a,b] & 
\eay
}
\right\}\;} 
& \hbox{ whenever }\;
\i\not\in\za_\phi[a,b] 
\eay
\right.
$\vspace{3mm}

\noi 
This allows to define a homeomorphism 
$H:D=\can\za\hsur$ onto 
$\hsur X$ by 
$H(d(a))=f(a)$ for all $a\in\dn.$ 
We claim that $H$ is projection-keeping --- 
which implies $X\in\pe\za$. 
Indeed let $\xi\in\cpo,\,\xi\sq\za$, and, say, 
$d(a),\,d(b)\in\can\za$, $d(a)\dar\xi=d(b)\dar\xi.$ 
Then we have 
$\xi\sq\za_\phi[a,b]$ by the second part of $(\dag),$ 
hence $f(a)\dar\xi=f(b)\dar\xi$ holds 
by the first part of $(\dag),$ as required.
\epf

\vyk
{%
\color{blue}%
\bcor
[lifting-fusion]
\lam{lfut}
Let\/ 
$\xi_1\sq\xi_2\sq\ldots\sq\xi_m\sq\ldots\sq\za\in\cpo$ 
be an increasing sequence in\/ $\cpo$, a map\/ 
$\phi:\om\to\za$ be\/ \dd{\za}admissible, 
and\/ 
\ben
\nenu
\itlb{lfut1}
$X_u\in\pe{\xi_m}$ for all\/ $m$ and\/ $u\in2^m,$

\itlb{lfut2}
$\sis{X_u}{u\in 2^{m}}$ is a\/
\dd\phi split system in\/ $\pe{\xi_m}$ for all\/ $m$, 

\itlb{lfut3}
$X_{u\we e}\dar{\xi_m}\sq X_u$ for all\/ 
$m\yi u\in2^m\yi e=0,1$.
\een
Then the set\/ 
$X=\bigcap_{m\in\om}\bigcup_{u\in 2^m} (X_u\uar\za)$ 
belongs to\/ $\pe\za$. 
\ecor
\bpf
It suffices to note that $\sis{X_u\uar\za}{u\in 2^{m}}$ 
is a\/ \dd\phi split system in\/ $\pe{\za}$, $\kaz m$, 
by Lemma~\ref{uver}, and then apply Theorem~\ref{fut}.
\epf
}%

The classical theorem, that any uncountable Borel
or $\fs11$ set includes a perfect subset,
does not generalize for 
$\pe\za:$ if $\card\za\geq 2$ then   
easily there is an uncountable closed $W\sq\can\za$
which does not include a subset in $\pe\za.$
However the following weaker claim 
(Corollary 16 in \cite{jsl99}) survives. 

\vyk{
\bcor
\lam{bor}
Assume that\/ $X\in\pe\za,$ and\/ $B\sq\can\za$ is
a Borel set. 
There exists a set\/ $Y\in\pe\za,\msur$ $Y\sq X$
such that either\/ $Y\sq B$ or\/ $Y\cap B=\emps$. 
\ecor
\bpf
Argue by induction on $\al,$ where $B\in\fs0\al.$
If $\al=1,$ 
the $B$ is open, and apply Lemma~\ref{clop}.
Otherwise 
$B=\bigcup_m B_m$ where $B_m\in\fp0{\al_m}$ 
for some $\al_m<\al.$ 
If there is $Y\in\pe\za,\msur$ $Y\sq X\cap B_m$ for some 
$m,$ then $Y\sq B.$
Otherwise, by the inductive hypothesis,
lemmas \ref{suz} and \ref{pand} yield a 
fusion sequence $\ang{X_u:u\in 2^{<\om}}$  in $\pe\za$, 
such that $X_\La\sq X$ and $X_u\cap B_m=\emps$ 
for all $m\in\om$ 
and $u\in 2^m.$ The set 
$Y=\bigcap_{m\in\om}\bigcup_{u\in 2^m}X_u$ is as required.
\epf

The result can be strengthened ! 

\bcor
\lam{ana}
Assume that\/ $X\in\pe\za,$ and\/ $A\sq\can\za$ is
a\/ $\fs11$ set. 
There exists a set\/ $Y\in\pe\za,\msur$ $Y\sq X$
such that either\/ 
$Y\sq A$ or\/ $Y\cap A=\emps$. 
\ecor
\bpf
Consider a \dd{\pe\za}generic extension $\rV^+$ of the universe 
$\rV.$ 
For a Borel or analytic set $C$ in $\rV,$ let $C^+$ denote 
the set defined in $\rV^+$ by the same construction. 
There is a condition $Y'\in\pe\za$ which 
decides $\x\in A^+,$ where $\x$ is the name for the 
\dd{\pe\za}generic element of $\can\za.$
Suppose that \eg\ $Y'$ forces $\x\in A^+.$
It is proved in \cite{jsl99} that  
$\aleph_1$ remains uncountable in $\rV^+.$
Therefore there is a 
Borel set $B\sq A$ (a {\it constituent\/} of $A$)
and a condition 
$Y\in\pe\za,\msur$ $Y\sq Y',$ which forces $\x\in B^+.$
By Corollary~\ref{bor},
we can assume that either $Y\sq B$ or $Y\cap B=\pu.$
The ``or'' case is impossible by the Shoenfield 
absoluteness, as $Y$ forces $\x\in Y^+.$
Thus $Y\sq B,$ as required.
\epf
}

\bcor
\lam{bana}
Assume that\/ $X\in\pe\za$, and a set\/ $A\sq X$ has the 
relative Baire property in\/ $X$ but not relative 
meager in\/ $X.$  
Then there is a~set\/ $Y\in\pe\za\yt Y\sq A$. 
\ecor
\bpf
It suffices to prove the result 
in case $X=\can\za.$ 
As $A$ has the BP but not meager, 
there is a basic clopen set $\pu\ne B\sq X$ 
(see the proof of Lemma~\ref{clop}) 
such that $A\cap B$ is comeager in $B$, so that there are 
relatively open dense sets $D_n\sq B$ 
satisfying $\bigcap_nD_n\sq A\cap B$. 
Now Lemmas \ref{suz} and \ref{pand} yield a 
fusion sequence $\sis{X_u}{u\in 2^{<\om}}$  in $\pe\za$, 
such that $X_\La\sq X$, and each $X_u$ is clopen and 
satisfies $X_u\sq D_m$ for all $m\in\om$ 
and $u\in 2^m.$ 
The set $Y=\bigcap_{m\in\om}\bigcup_{u\in 2^m}X_u$ 
is as required.
\epf

\parf{Uniform shrinking}
\las{174}

Say that a set $X\in\pe\za$ is \rit{uniform}, 
\index{set!uniform}%
\index{uniform}%
if for any pair of tuples $\i\su\j$ in $\za$ and any $x,y\in X$, 
we have $x(\j)=y(\j)\imp x(\i)=y(\i)$.   
As the first application of the splitting/fusion technique, we 
prove a lemma on uniform shrinking.

\ble
[in $\rL$]
\lam{174L}
If\/ $\za\in\cpo$ and\/ $X\in\pe\za$ then there is a uniform set\/ 
$Y\in\pe\za$, $Y\sq X$.
\ele

\bpf
Let  $\phi:\om\onto\za$ be \dd{\za}admissible.  
Lemmas \ref{suz}, \ref{pand}, \ref{1to2} yield a 
fusion sequence $\sis{X_u}{u\in \bse}$  in $\pe\za$, 
such that $X_\La\sq X$ and the layer $\sis{X_u}{u\in2^m}$ 
satisfies \ref{aprct} of 
Definition~\ref{splis} for all $m$. 
Then $Y=\bigcap_n\bigcup_{u\in2^n}X_u\in\pe\za$ 
by Theorem~\ref{fut}, and $Y\sq X$.  
We claim that $Y$ is uniform. 

Indeed suppose  
that $\i\su\j$ belong to $\za$, and $x, y$ in $Y$ 
satisfy $x(\i)\ne y(\i)$, say $x(\i)(k)=0$ but $y(\i)(k)=1$ 
for some $k<\om$. 
Then $x\ne y$, hence there exists $m$ and some $u\ne v$ 
in $2^m$ such 
that $x\in X_u$, $y\in X_v$. 
We can take $m$ big enough for $x'(\i)(k)=0$ for all 
$x'\in X_u$ but 
$y'(\i)(k)=1$ for all $y'\in X_v$. 
Thus $({X_u}\dir\i)\cap({X_v}\dir\i)=\pu$.

Now consider the initial segment $\et=\za_\phi[u,v]\sq\za$. 
Then $X_u\dar\et= X_v\dar\et$ by \ref{prct} of Definition~\ref{splis}. 
It follows that $\i\nin\et$ since 
$({X_u}\dir\i)\cap({X_v}\dir\i)=\pu$. 
Therefore $\j\nin\et$ either. 
But then  $({X_u}\dir\j)\cap({X_v}\dir\j)=\pu$ by 
\ref{aprct} of 
Definition~\ref{splis}.  
We conclude that $x(\j)\ne y(\j)$, as required.
\epf

\parf{Axis/avoidance shrinking}
\las{173}

{\ubf We continue to argue in $\rL$.}
Here we set up some notions related to continuous maps 
$F:\can\xi\to\cN,$ $\xi\in\cpo$. 
Recall that $\cam=\dn\sq\cN=\bn$. 
Let
$$
\bay{rcl}
\kmar{cf\_xi}
\index{functions!$\cf_\xi$}
\index{zCFxi@$\cf_\xi$, functions}
\cf_\xi &=&\ens{F:\can\xi\to\cN}
{F\text{ is continuous}},\\[0.5ex]
\kmar{cfd\_xi}
\index{functions!$\cfd_\xi$}
\index{zCFdxi@$\cfd_\xi$, functions}
\cfd_\xi &=&\ens{F:\can\xi\to\cam}
{F\text{ is continuous}}\sq\cf_\xi,
\eay
$$
and
$\cf=\bigcup_{\xi\in\cpo}\cf_\xi$, 
$\cfd=\bigcup_{\xi\in\cpo}\cfd_\xi$.
\kmar{cf, cfd}
\index{functions!$\cf$}%
\index{zCF@$\cf$, functions}%
\index{functions!$\cfd$}%
\index{zCFd@$\cfd$, functions}%
Let $\modd  f=\xi$ in case $f\in\cf_\xi$.%
\kmar{\mns modd f}%
\index{dimension, $\modd f$}%
\index{z11f11@$\modd f$, dimension}%

\bdf
[in $\rL$]
\lam{174d}
Assume that $\sg\sq\ta$ belong to $\cpo$, 
$X\in\pe\ta$, $\i\in\ta$, $S\sq\cam$, and $F\in\cf_\sg$. 

If $F(x\dar\sg)=x(\i)$ for all 
$x\in X$, then say that $F$ is 
\rit{an\/ $\i$-axis map on\/} $X$. 
\index{axis map@$\i$-axis map}%
\index{map!$\i$-axis map}%

If $F(x\dar\sg)\nin S$ 
for all $x\in X$,  then say that $F$ 
\rit{avoids\/ $S$} on $X$.
\index{avoids}%
\edf

We prove several lemmas here, related to axis maps and avoidance, 
which culminate in a dichotomy theorem (Theorem~\ref{176L}). 

\ble
\lam{1731}
If\/ $\i\in\ta\in\cpo,$  $X\in\pe\ta,$ and $F\in\cf_\ta$ 
is not an $\i$-axis~map on\/ $X$, then there is\/ 
$Y\in\pe\ta$, $Y\sq X$, such that 
$F$ avoids\/ $Y\dir\i$ on\/ $Y$.
\ele

\bpf
We have $F(x_0)\ne x(\i)$ for some $x_0\in X$, say 
$F(x_0)(k)=m$ and $x_0(\i)(k)=n$ for some $k$ and $m\ne n$. 
Thus $X'=\ens{x\in X}{F(x )(k)=m\land x (\i)(k)=n}\ne\pu$. 
But $X'$ is open in $X$. 
Take any $Y\in\pe\ta,Y\sq X'$ by Lemma~\ref{clop}.
\epf

\ble
\lam{1732}
If\/ $\et\sq\ta$ and\/ $\xi$ belong to\/ $\cpo,$  
$\i\in\ta\bez\et,$  
$X\in\pe\xi,$ $Y\in\pe\ta,$ and $F\in\cf_\xi$,  
then there exist relatively 
clopen sets\/ $X'\sq X$ and\/ $Y'\sq Y$ in resp.\/ $\pe\xi,\,\pe\ta$, 
such that\/ $Y'\dar\et=Y\dar\et$ and\/ 
$F$ avoids\/ $Y'\dir\i$ on\/ $X'.$
\ele

\bpf
Pick any $x_0\in X.$ 
Let $p_0=F(x_0)$, $Q_m=\ens{p\in\cam}{p\res m=p_0\res m}$,
\pagebreak[0]   
$$
U_m =
\ens{u\in Y\dar\et}{\sus y\in Y\,(y\dar\et=u\land y(\i)\nin Q_m)} 
$$
for all $m<\om$. 
Then $U_m\sq U_{m+1}$, $\kaz m$. 
Further, Lemma~\ref{darop} implies that each set $U_m$ is 
clopen in $Y\dar\et\in\pe\et$. 
Moreover,  we have $Y\dar\et=\bigcup_mU_m$. 
(Because if $u\in Y\dar\et$ then 
$\ens{y(\i)}{y\in Y\land y\dar\et=u}$ 
is a perfect set.) 
It follows by the compactness of $\can\et$ that 
$Y\res\et=U_m$ for some $m$.

Now let $Y'=\ens{y\in Y}{y(\i)\nin Q_m}$.  
The set $S=\ens{x\in X}{F(x)\in Q_m}$ is clopen in $X$, 
and $p_0\in S$, 
hence there exists a relatively clopen $X'\in\pe\xi,X'\sq S$. 
We claim that $X',Y'$ are as required. 
Indeed $Y'\dar\et=Y\dar\et$ holds by the choice of $m$, 
whereas $F$ avoids\/ $Y'\dir\i$ on\/ $X'$ by construction. 
It remains to show that $Y'\in\pe\ta$ and that $Y'$ 
 is relatively clopen in $Y$.  

Note that $Y'=Y\cap(V\uar\ta)$, where 
$V=\ens{v\in Y\rsq\i}{v(\i)\nin Q_m}$ is clopen in $Y\dar\et$ 
by Lemma~\ref{darop}. 
Lemma~\ref{lin} implies that $V\in\pe{\sq\i}$. 
Then $Y'\in\pe\ta$ by Lemma~\ref{apro}, as required.
\epf

\bcor
\lam{1733}
Assume that\/ $\et\sq\ta$ belong to\/ $\cpo,$ 
$X,Y\in\pe\ta,$ $X\dar\et=Y\dar\et$, $F\in\cf_\ta$,  $\i\in\ta,$   
and either\/ $\i\nin\et$ or\/ $F$ is not an\/ $\i$-axis map on\/ $X.$ 
Then there exist relatively clopen sets\/ $X'\sq X$ and\/ $Y'\sq Y$ 
in\/ $\pe\ta$, such that\/ 
$X'\dar\et=Y'\dar\et$ and\/ 
$F$ avoids\/ $Y'\dir\i$ on\/ $X'.$
\ecor
\bpf
Suppose that $\i\nin\et$. 
Then by Lemma~\ref{1732}  there exist relatively clopen sets 
$X'\sq X$ and  $Y''\sq Y$ in $\pe\ta$, such that  
$Y''\dar\et=Y\dar\et$ and $F$ avoids $Y''\dir\i$ on $X'.$ 
Take $Y'=Y''\cap(X'\dar\et\uar\ta)$, and we are done. 

Now suppose that $\i\in\et$ and $F$ is not an $\i$-axis map on $X$. 
Lemma~\ref{1731} yields a relatively clopen $X'\in\pe\ta$, 
$X'\sq X$, such that $F$ avoids\/ $X'\dir\i$ on\/ $X'.$  
Take $Y'=Y\cap(X'\dar\et\uar\ta)$, and we are done. 
\epf

\bcor
\lam{1734}
If\/ $\et\sq\ta$ belong to\/ $\cpo,$ 
$X,Y\in\pe\ta,$ $X\dar\et=Y\dar\et$, $F\in\cf_\ta$,  
$\i\in\ta\bez\et,$   
then there exist relatively clopen sets\/ $X'\sq X$ 
and\/ $Y'\sq Y$ in\/ $\pe\ta$, such that\/ 
$X'\dar\et=Y'\dar\et$ and\/ 
$(Y'\dir\i)\cap(X'\dir\i)=\pu$.
\ecor
\bpf
Use Corollary~\ref{1733} for $F(x)=x(\i)$.
\epf

\parf{Axis/avoidance dichotomy theorem}
\las{176}

And now the main result goes, a dichotomy theorem.

\bte
\lam{176L}
If\/ $\ta\in\cpo$, $X\in\pe\ta$, and\/ $F\in\cf_\ta$ 
then there is a set\/ $Y\in\pe\ta$, $Y\sq X$, 
such that 
one of the two following claims holds$:$
\ben
\renu
\itlb{176a}\msur
$F$ avoids\/ $Y\dir\i$ on\/ $Y$ 
for all\/ $\i\in\ta\,;$

\itlb{176b}
there is\/ $\j\in\ta$ such that\/ $F$ is a\/  
$\j$-axis map on\/ $Y$ and $F$ avoids\/  
$Y\dir\i$  on\/ $Y$ 
for all $\i\in\ta$, $\i\ne\j$. 
\een
\ete

\bpf
To begin with, prove that if $U\in\pe\ta$ and $\i\ne\j$ 
belong to $\ta$ 
then 
\ben
\nenu
\itlb{1761}
$F$ cannot be both $\i$-axis map on $U$ and $\j$-axis map on $U$. 
\een
Indeed suppose otherwise.  
Let say $\i\not\sq\j$, so that $\i\nin\et=\ilq\j$. 
Corollary~\ref{1734} with $X=Y=U$ (note that $\j\in\et$) 
yields sets $X',Y'\in\pe\ta$ such that $X'\cup Y'\sq U$, 
$X'\dir \j=Y'\dir\j$, but $(X'\dir \i)\cap(Y'\dir\i)=\pu$. 
Thus $X'\dir \i\ne X'\dir\j$ or $Y'\dir \i\ne Y'\dir\j$, 
both cases 
leading to a contradiction with the contrary assumption. 
This ends the proof of \ref{1761}. 

\vyk{
Then $x(\i)=x(\j)$ for all $x\in Z$. 
By Lemma~\ref{174L}, we can assume that $Z$ is uniform. 
Then we have the equivalence%
\pagebreak[1]
$$
{x\rsq\i=y\rsq\i}\leqv {x\rsq\j=y\rsq\j}\,,\quad
\text{for all }x,y\in Z\,.
$$
This equivalence is obviously preserved under any \pkh, and hence 
it must hold for $Z=\can\ta$, which is obviously false. 
}

Coming back to the theorem, we have two cases. 

{\it Case 1\/}: there exist $\j\in\ta$ and $Z\in\pe\ta$, 
$Z\sq X$, 
such that $F$ is a $\j$-axis map on $Z$. 
Let $\da=\ta\bez\ans\j$ in this case. 

{\it Case 2\/}: not case 1. 
Let $\da=\ta$ and $Z=X$ in this case. 

It follows from \ref{1761} that in both cases 
\ben
\nenu
\atc
\itlb{1762}
if $U\in\pe\ta$, $U\sq Z$, 
$\i\in\da$, then $F$ is {\ubf not} an $\i$-axis map on $U$. 
\een
\vyk{
Now we claim that 
\ben
\nenu
\atc
\atc
\itlb{1763x}
if $U,V\in\pe\ta$, $U\cup V\sq Z$, $\i\in\da$, $\et\in\cpo$, 
$\et\sq\ta$, $U\dar\et=V\dar\et$, then there exist sets 
$U',V'\in\pe\ta$, $U'\sq U$, $V'\sq V$, such that still  
$U'\dar\et=V'\dar\et$, and $(\ima F{U'})\cap({V'}\dir\i)=\pu$. 
\een

{\it Case A\/}: $\i\in\et$. 
By \ref{1762} there exists $x_0\in U$ with $F(x_0)\ne x_0\dir\i$, 
say $F(x_0)(k)=0$ but $x_0(\i)(k)=1$ for some $k<\om$. 
Then $A=\ens{x\in U}{F(x)(k)=0\land x(\i)(k)=1}\sq U$ is 
a non-empty relatively clopen subset of $U$. 
Therefore by Lemma~\ref{clop}  
there is a relatively clopen in $U$ set $A'\sq A$, $A\in\pe\ta$. 
Thus $(\ima F{A'})\cap (A'\dir\i)=\pu$ by construction. 
However the set $B'=V\cap (A'\dar\et\uar\ta)$ belongs 
to $\pe\ta$ by 
lemmas \ref{less} and \ref{apro} since $U\dar\et=V\dar\et$ 
and $A'\sq U$, and obviously $B'\dar\et=A'\dar\et$. 
Therefore, as $\i\in\et$, the sets $U'=A'$ and $V'=B'$ prove 
\ref{1763x}. 

{\it Case B\/}: $\i\in\ta\bez\et$. 
Applying Lemma~\ref{99y} twice, we get sets 
$A\sq U$ and $B,B'\sq V$ in $\pe\ta$, and some $k<\om$, 
such that $A\dar\et=B\dar\et=B'\dar\et$, 
$x(\i)(k)=0$ for all $x\in B$, 
$x(\i)(k)=1$ for all $x\in B'$, 
and either $F(x)(k)=0$ for all $x\in A$ or 
$F(x)(k)=1$ for all $x\in A$. 
Let say $F(x)(k)=0$ for all $x\in A$. 
Then the sets $U'=A$ and $V'=B'$ prove \ref{1763x}. 
}	

Now fix  any \dd{\ta}admissible map  $\phi:\om\onto\ta$. 
The next claim is a consequence of \ref{1762} and 
Corollary~\ref{1733}, by means  of 
Corollary~\ref{suz+} applied consecutively enough many times:
\ben
\nenu
\atc\atc
\itlb{1764}
If $\i\in\da$ and $m<\om$ 
then any $\phi$-split system $\sis{X_u}{u\in 2^m}$  
of sets $X_u\sq Z$ 
in $\pe\ta$ admits  a narrowing $\sis{X'_{u}}{u\in 2^{m}}$    
such that 
if $u,v\in 2^m$ then 
$F$ avoids ${X'_v}\dir\i$ on $X'_u$, and hence 
$F$ avoids ${X'_m}\dir\i$ on $X'_m=\bigcup_{u\in2^m}X'_u$.
\een
With this ``narrowing'' result, 
Lemmas \ref{suz}, \ref{pand} yield a 
fusion sequence $\sis{X_u}{u\in 2^{<\om}}$  in $\pe\ta$, 
such that $X_\La\sq Z$, and, for each $m$, 
$F$ avoids ${X_m}\dir\i$ on $X_m=\bigcup_{u\in2^m}X_u$, 
where $\i=\phi(m)\in\da$.
Then $Y=\bigcap_n\bigcup_{u\in2^n}X_u\in\pe\ta$, 
$Y\sq Z\sq X$, 
and $F$ avoids ${Y}\dir\i$ on $Y$ 
for all $\i\in\da$, 
as required.
\epf

\parf{Avoidable sets}
\las{175}

Assume that $U\in\pe{\sq\i}$, $\i\in\tup$. 
Say that a set $S\sq\cam=\dn$ is 
\rit{\dd Uavoidable on $\i$\/} if 
\index{set!avoidable@$U$-avoidable on $\i$}%
\index{Uavoidable@$U$-avoidable on $\i$}%
%
there exists a relatively clopen set 
$V\sq U$   
satisfying  $V\rsl\i=U\rsl\i$ and $S\cap(V\dir\i)=\pu$. 
Thus avoidability in this sense means that not $U$ itself 
but a certain clopen subset of $U$ with the same projection 
avoids $S$.

\bte
\lam{175t}
Suppose that\/ $\xi\in\cpo$, $X\in\pe\xi$,  $F\in\cf_{\xi}$, 
and\/ $\rU\sq\bigcup_{\i\in\tup}\pe{\sq\i}$ is a 
countable set.
Then there is a set\/ $Y\in\pe\xi$, $Y\sq X$, such that 
the image\/ $S=\ima FY$ is\/ $U$-avoidable on\/ $\i$ 
for all\/  $\i\in\tup$ and\/ 
$U\in\pe{\sq\i}\cap\rU$.
\ete

\bpf
\vyk{
For the sake of brevity, 
we write $Z\para\i U$, if $Z\in\pe\xi$, $\i\in\tup$, 
$U\in\pe{\sq\i}\cap\rU$, 
and there is $\pu\ne V\sq U$, relatively clopen in $U$, 
such that  
$V\rsl\i=U\rsl\i$, and $F$ avoids $V\dir\i$ on $Z$, \ie, 
$(\ima FZ)\cap(V\dir\i)=\pu$. 
}%
Lemma~\ref{1732} ($\ta=\ilq\i$, $\et=\ile\i$)  implies:
\ben
\nenu
\itlb{1751}
if $Z\in\pe\xi$, $\i\in\tup$, $U\in\pe{\sq\i}\cap\rU$, 
then there is 
a relatively clopen set $Z'\sq Z$, $Z'\in\pe\xi$, such that 
$\ima F{Z'}$ is $U$-avoidable on $\i$.
\een
%
Fix  any \dd{\xi}admissible map  $\phi:\om\onto\xi$. 
The next claim is a consequence of \ref{1751} and 
Corollary~\ref{1733}, by means  of 
Corollary~\ref{suz+} applied consecutively enough many times:
\ben
\nenu
\atc 
\itlb{1752}
If  
$\i\in\tup$, $U\in\pe{\sq\i}\cap \rU$,  and $m<\om$, 
then any $\phi$-split system $\sis{X_u}{u\in 2^m}$  
of sets $X_u\in\pe\xi$ 
admits  a narrowing $\sis{X'_{u}}{u\in 2^{m}}$    
in $\pe\xi$ such that $\ima F{X'_m}$ 
is ${U}$-avoidable on $\i$, 
where $X'_m=\bigcup_{u\in 2^m}X'_u$.
\een
Using this result and the countability of $\rU$, 
Lemmas \ref{suz} and \ref{pand} 
yield a fusion sequence $\sis{X_u}{u\in 2^{<\om}}$  
in $\pe\xi$, 
such that $X_\La\sq X$, and, for each $\i\in\tup$ and 
$U\in\pe{\sq\i}\cap \rU$ there is a number $m$, such that 
$\ima F{X_m}$ is \dd Uavoidable on $\i$, where 
$X_m=\bigcup_{u\in 2^m}X_u$.
Then $Y=\bigcap_m\bigcup_{u\in2^m}X_u\in\pe\xi$, 
$Y\sq X$, and $\ima F{Y}$ is \dd Uavoidable
for all $\i\in\tup$ and $U\in\pe{\sq\i}\cap \rU$. 
\epf

\bre
\lam{175r}
The theorem will be applied only in cases when the given set 
$\rU\sq\bigcup_{\i\in\tup}\pe{\sq\i}$ satisfies the property 
that if $\pu\ne V\sq U\in\rU$ is relatively clopen in $U$ 
then $V\in\rU$ as well. 
In this case, the condition of relative clopenness of $V$ 
in the definition 
of being ``\dd Uavoidable on $\i$" can be replaced by 
just $V\in\rU$, and then Theorem~\ref{175t} still holds.
\ere

\sekt{Normal forcing notions}
\las{rfono}

It will take considerable effort to actually define 
the forcing notion $\cX\sq\pei$ in the 
constructible universe $\rL$ 
for the proof of Theorem~\ref{mt1}. 
Yet we can gradually introduce some conditions on $\cX$ 
that will bring a number of useful consequences related 
to the corresponding $\cX$-generic extensions of $\rL$, 
and which will be fulfilled in the final construction 
of $\cX$. 

The first group of those conditions is wrapped up in the 
concept of {\ubf a normal forcing}, studied in Sections~\ref{rfo} 
and \ref{ker}. 
Each normal forcing $\cX$ is a forcing notion in $\rL$, 
satisfying $\cX\sq\pei$. 
It adjoins an \dd\cX{\ubf generic array} 
$\w\in\can\tup$, and we get a generic extension $\rL[\w]$ 
and various {\ubf symmetric subextensions}, introduced in 
Section~\ref{ga}. 

The associated {\ubf forcing relation} is studied in 
Section~\ref{frel}, and the effect of 
{\ubf actions by permutations $\pi\in\per$} 
in Sections~\ref{fap} and \ref{63}. 
The {\ubf Fusion property}, an important condition which 
implies continuous reading of real names, among other 
consequences, is introduced in Section~\ref{fup}. 
In particular, the background forcing $\pei$ has the 
Fusion property (section~\ref{fff}). 
We derive some consequences of the Fusion property, 
related to {\ubf various forms of $\AC$ and $\DC$}, in
Section~\ref{fupac}.

\parf{Normal forcings}
\las{rfo}

{\ubf We argue in $\rL$ in this section.} 
Any set $\cX\sq\pei$ can be viewed as a forcing 
notion, with the 
partial order $\psq$ on $\pei$ defined by: 
\index{partial order $\psq$}%
\index{zzzsqp@$\psq$}%
\imar{psq}%
$X\psq Y$ iff $\et=\dym Y\sq \dym X$ and 
${{X\dar\et}\sq Y}$. 
But we have to somehow restrict the generality, to make 
sure that $\cX$ adjoins $\tup$-arrays of reals 
(points of $\can{}$), similarly to 
$\pei$ itself. 
Recall that 
$$
\bay{rcl}
\cX\dar\et &=&
\ens{X\dar\et}{ X\in\cX \land \et\sq\dym X},\\[0.7ex]
\index{projection!dar@$\dar$}%
\index{zzzpdar@$\dar$}%
\cX\rsq\i&=&
\cX\dar\et, \text{ where }
\et=\ilq\i=\ens{\j\in\tup}{\j\sq\i},\\[0.7ex]
\index{projection!prsq@$\rsq\i$}%
\index{zzzprsq@$\rsq\i$}%
\dym X &=& \xi, \text{ in case }X\sq\can\xi, 
\index{dimension, $\dym X$}%
\index{z11X11@$\dym X$, dimension}%
\kmar{dym X}%
\eay
$$
by Section~\ref{prelim1}, for any $\cX\sq\pei$, and 
$\can{}=\dn,$ the Cantor space.

Say that a set $\cX\sq\pei$ is a \rit{normal forcing}, 
\kmar{RF}%
\index{normal forcing, $\RF$}%
\index{NF, normal forcing}%
$\cX\in\RF$ for brevity, iff the following conditions 
\ref{rfo1}--\ref{rfo5} hold:
\ben
\cenu
\itlb{rfo1}
$\cX\sq\pe{}$, and if $\ta\in\cpo$ then $\can\ta\in\cX$.

\itlb{rfo2}
If $\xi\sq\ta$ belong to $\cpo$ and 
$X\in\cX\cap\pe\ta$ then $X\dar\xi\in\cX,$ and hence 
$\cX\dar\xi=\cX\cap\pe\xi$. 
In particular the set $\bon=\ans\pu=X\dar\pu$ 
belongs to $\cX\dar\pu$, 
and $\bon\psq X$ for any $X\in\cX.$

\itlb{rfo3}
If $\xi\sq\ta$ belong to $\cpo$,   
$X\in\cX\dar\ta$, $Y\in\cX\dar\xi$, and 
$Y\sq X\dar\et$, then $X\cap(Y\uar\ta)\in\cX\dar\ta$. \ 
In particular, if $Y\in\cX\dar\xi$ then 
$Y\uar\ta\in\cX\dar\ta$.

\itlb{rfo4}
If $\ta\in\cpo$, $X\in \cX\dar\ta$, $Y\in\pe\ta$, 
$Y\sq X$ is clopen in $X$, then $Y\in \cX$.

\itlb{rfo6}
$\cX$ is $\qer$-invariant: 
if $\pi\in\qer$ and $X\in \pei$ then $X\in\cX\eqv \pi\akt X\in\cX$.

\itlb{rfo5}
If $\ta\in\cpo$, $X\in \pe\ta$, and $X\rsq\i\in \cX\rsq\i$ 
for all $\i\in\ta$, then $X\in\cX$.
\een

Quite clearly $\pei$ itself belongs to $\RF$: 
$\can\ta\in\cX$ in \ref{rfo1} holds via the identity \pkh,  
\ref{rfo2} holds by Lemma~\ref{less}, 
\ref{rfo3} holds by Lemma~\ref{apro}, 
\ref{rfo4} and \ref{rfo5} are obvious, 
\ref{rfo6} holds by Lemma~\ref{pi2}, 
so that $\pei$ is even $\pero$-invariant. 

\ble
\lam{rfoL}
Let\/ $\cX\in\RF$. 
Under the assumptions of Lemma~\ref{suz}, 
Corollary~\ref{suz+}, Lemma~\ref{1to2}, Lemma~\ref{pand},
if all the given sets\/ $X_u,X,Y$ belong to\/ $\cX$, 
then the resulting sets\/ $Y_u,Z_u,Y_s$ 
belong to\/ $\cX$ as well$.$  

Under the assumptions of Lemma~\ref{clop}, 
if\/ $X\in\cX$ then\/ $X'\in\cX,$ too.

Under the assumptions of Lemma~\ref{99y}, 
if\/ $X,Y\in\cX$ then\/ $X',Y'\in\cX$.

Under the assumptions of Lemma~\ref{pe4}, 
if\/ $X_k\in\cX$, $\kaz k$, then\/ $X\in\cX$.

Under the assumptions of Corollary~\ref{pe2}, 
if\/ $X,Y\in\cX$ then\/ $Z\in\cX$.
\ele

\bpf
Make use of \ref{rfo3} above \poo\ Lemma~\ref{suz} 
and  Corollaries~\ref{suz+} and \ref{pe2}, 
of \ref{rfo4} above \poo\ 
Lemmas~\ref{1to2}, \ref{pand}, \ref{clop}, \ref{99y}, 
and of \ref{rfo5} \poo\ Lemma~\ref{pe4}.
\epf

\bdf
\lam{noc}
If $\rP\sq\pei$ then let $\noc\rP$ 
(the \rit{normal hull} of $\rP$) 
\kmar{\mns noc rP}%
\index{normal hull!$\noc\rP$}%
\index{zNHP@$\noc\rP$, normal hull}%
be the least set $\cX\in\nf$ with $\rP\sq\cX$. 
The set $\noc\rP$ is equal to the intersection of all sets 
$\cY\in\nf$ satisfying $\rP\sq\cY$.
\edf

\parf{Kernels of normal forcings}
\las{ker}

{\ubf We still argue in $\rL$.} 
Here we show that each normal forcing $\cX$ is the  
normal hull of its smaller and simpler part called 
\rit{the kernel}. 
If $\xi\sq\tup$ then let a $\xi$-\rit{kernel} 
\index{kernel!$\xi$-kernel}%
be a system $\cK=\sis{\cK_\i}{\i\in\xi}$ of sets 
$\cK_\i\sq\pele\i$, 
satisfying \ref{ke1s}--\ref{ke6} below.
\ben
\senu
\itlb{ke1s}
If tuples $\j\su\i$ belong to $\xi$ and  
$Y\in\pele\j$ then $Y=X\rsq\j$ for some $X\in\pele\i$. 
\senu
\itlb{ke2}
If tuples $\j\su\i$ belong to $\xi$ and 
$X\in\cK_\i$ then $X\rsq\j\in\cK_\j$.

\itlb{ke3}
If tuples $\j\su\i$ belong to $\xi$,   
$X\in\cK_\i$, $Y\in\cK_\j$, and 
$Y\sq X\rsq\j$, then $Z=X\cap(Y\usq\i)\in\cK_\i$. \ 

\itlb{ke4}
If $\i\in\xi$, $X\in \cK_\i$,  
$\pu\ne Y\sq X$ is clopen in $X$, then $Y\in \cK_\i$.


\itlb{ke6}
If tuples $\j\ekp\i$ belong to $\xi$ and
$X\in \cK_i$ then $\pi_{\i\j}\akt X\in\cK_\j$. 
(See Example~\ref{pi*ij} on $\pi_{\i\j}$.)
\een
Say that $\cK$ is a \rit{strong $\xi$-kernel}, if 
\index{kernel!strong}
in addition the following \ref{ke1s} holds. 

\ben
\ssenu
\itlb{ke1}
If $\i\in\xi$ then $\can{\ilq\i}\in\cK_\i\sq\pele\i$.  
\een

\ble
\lam{***}
In the presence of\/ \ref{ke3}, 
condition\/ \ref{ke1} implies\/ \ref{ke1s}.
\ele
\bpf
As $X=\can{\ilq\i}\in\cK_\i$ by \ref{ke1}, 
the set $Z=Y\usq\i=X\cap(Y\usq\i)$ 
belongs to $\cK_\i$ by \ref{ke3},  
and obviously $Y=Z\rsq\j$. \  
\epf

\ble
\lam{r2k}
Let\/ $\cX\in\RF$. 
Then\/ $\Ker\cX=\sis{\cX\rsq\i}{\i\in\tup}$ 
\kmar{\mns Ker cX}%
\index{kernel!$\Ker\cX$}%
\index{zKerX@$\Ker\cX$, kernel}%
{\em(the \rit{kernel} of $\cX$)} 
is a strong\/ 
$\tup$-kernel.
\ele
\bpf
Infer \ref{ke1} and \ref{ke2}--\ref{ke6} from 
\ref{rfo1}--\ref{rfo6} above. 
Apply Lemma~\ref{lin} for \ref{ke4}.
\epf

Conversely, every \dd\tup kernel defines a normal forcing 
via \ref{rfo5}.

\ble
\lam{k2r}
Let\/ $\cK=\sis{\cK_\i}{\i\in\tup}$ be a strong\/ 
$\tup$-kernel.   
Then\/ $\cX=\noc\cK\in\RF$, 
$\cK=\Ker\cX$ --- so that\/ $\cX\rsq\i=\cK_\i$ 
for all\/ $\i\in\tup$, 
and if\/ $\xi\in\cpo$ then\/ 
$\cX\dar\xi$ is equal to the set\/ 
$\cY_\xi=\ens{X\in\pe\xi}{\kaz\i\in\xi\,(X\rsq\i\in\cK_\i)}$. 
\ele
\bpf
We claim that the set $\cY=\bigcup_{\xi\in\cpo}\cY_\xi$ 
belongs to $\RF$.
As \ref{rfo5} of Section \ref{rfo} obviously holds for 
$\cY$ by construction, 
we derive \ref{rfo1}--\ref{rfo6} for $\cY$  
from \ref{ke1s} and \ref{ke2}--\ref{ke6} for $\cK$. 
Here \ref{rfo1},\ref{rfo2},\ref{rfo6} are entirely obvious. 

Make use of Lemma~\ref{darop} for \ref{rfo4}. 
Now focus on \ref{rfo3}. 
Thus assume that $\xi\sq\ta$ belong to $\cpo$,   
$X\in\cY\dar\ta$, $Y\in\cY\dar\xi$, and 
$Y\sq X\dar\et$; prove that $Z=X\cap(Y\uar\ta)\in\cY\dar\ta$.   
We have to check that $Z\rsq\i\in\cK_\i$ for all $\i\in\ta$. 
If $\i\in\xi$ then $Z\rsq\i=Y\rsq\i\in\cK_\i$. 
If $\i\in\ta\bez\xi$ and $\et=\xi\cap{\ilq\i}$ then 
$Z\rsq\i=X\rsq\i\cap (Y\dar\et)\usq\i $ 
by Lemma~\ref{599}, 
hence yet again $Z\rsq\i\in\qro\cK\i$ by \ref{ke3}, 
as required. 
Thus $\cY\in\pei$, and hence $\cX\sq\cY$ by the minimality of 
$\cX$. 

Moreover $\cY\rsq\i=\cK_\i$ by construction. 
Therefore, as $\cK_\i\sq\cX$, we have $\cY\sq\cX$ by 
\ref{rfo5} of Section~\ref{rfo} for $\cX$. 
Thus $\cY=\cX$ and we are done.
\epf

We may note that in fact even dyadic 
\index{kernel!dyadic}%
\dd\tud kernels suffice to produce normal forcings. 
Recall that $\tud=2\lom\bez\ans\La$, the set of all 
non-empty dyadic tuples. 
Obviously for any $\i\in\tup$ there is a unique dyadic 
tuple $\kn\i\in\tud$ satisfying $\i\ekp\kn\i$. 
\index{zi-@$\kn \i$}%
\index{tuple!$\kn \i$}%
Indeed put $\lh{\kn\i}=\lh\i$ and 
$$
\text{for all }\,k<\lh{\kn\i}=\lh\i,\quad
\kn\i(k)=
\left\{
\bay{rcl}
0&\text{in case}& \i(k)\text{ is even}\\[0.5ex]
1&\text{in case}& \i(k)\text{ is odd}
\eay
\right.
.
\eqno(*)
$$


\ble
\lam{22w} 
Assume that\/ $2\le\al<\omi$ and\/ 
$\cK=\sis{\cK_\i}{\i\in\tuq\al}$ 
is an\/ \dd{\tuq\al}kernel. 
Put\/ $\cK\iex_\i:=\pi_{\i,\kn\i}\akt\cK_{\kn\i}$ 
for all $\i\in\tup$.
Then\/ 
$\cK\iex=\sis{\cK\iex_\i}{\i\in\tup}$   
is an \dd\tup kernel, $\cK\iex_\i=\cK_\i$ 
for all\/ $\i\in\tuq\al$, and 
if\/ $\cK$ is strong then so is\/ $\cK\iex.$\qed 
\ele

Thus to define a normal forcing $\cX$ it suffices to first 
define an auxiliary \dd{\tud}kernel $\cK$ and then let 
$\cX=\noc{\cK\iex}$ by Lemmas \ref{22w} and \ref{k2r}.

\parf
[Generic arrays, extensions, and symmetric subextensions]
{Generic arrays, extensions, and subextensions}
\las{ga}

According to the formulation of Theorem~\ref{mt2}, 
we are going to establish our main results in this paper  
by means of suitable 
generic extensions of $\rL$, the constructible universe, 
under the \rit{consistent} assumption that $\ombl<\omi$ 
in the universe, intended to imply the existence of 
generic extensions. 
The forcing notions considered in this process will be 
normal forcings as in Section~\ref{rfo} defined in $\rL$. 
As the notion of iterated perfect set and many related 
notions are definitely non-absolute, we add the following 
warning.

\bbla
\lam{inL} 
The definition of $\pei$ in Section~\ref{prelim1} and all 
other relevant definitions in 
Sections \ref{prelim1}--\ref{rfo}, 
are assumed to be relativized to $\rL$ by default, and we'll 
not bother to add the sign $^\rL$ of relativization. 
In other words, $\tup$ is $(\tup)^\rL$,  
$\cpo$ is $(\cpo)^\rL$, $\pei=(\pei)^\rL,$ 
$\per=(\per)^\rL,$ $\RF=(\RF)^\rL,$ \etc

In addition, $\ombl<\omi$ will be our blanket assumption 
in the universe.
\ebla

Under $\ombl<\omi$, if $\za\in\ \cpo$ 
(\ie, $\za\in\rL$ and $\rL\mo\za\in\cpo$) 
then every set $X\in\pe\za$ is   
a countable subset of $\can\za$ in the universe. 
However it transforms to a 
perfect set in the universe by the closure operation: 
\rit{the topological closure} $\clo X$ 
\kmar{clo X}%
\index{closure $\clo X$}%
\index{zX*@$\clo X$, closure}%
\index{operation!$\clo X$, closure}%
\index{zzz*@$\clo{\;}$, closure}%
of a set $X\in \pe\za$ 
is closed in $\can\za$ in the universe. 
(And in fact $\clo X$ satisfies the definition of 
$\pe\za$ in the universe.)

Let $\cX\sq\pei$, $\cX\in \rL$ be a normal 
forcing, that is, \ref{rfo1}--\ref{rfo5} 
of Section~\ref{rfo} hold (in $\rL$), 
and $\cX$ is ordered by $\psq$, meaning that 
\bce
\rit{if $X\psq Y$ then $X$ is a stronger condition}. 
\imat{X psq Y}
\ece
Let $G\sq\cX$ be a filter \dd\cX generic over $\rL$.  
It easily follows from Lemma~\ref{rfoL} \poo\ Lemma~\ref{clop},  
that there is a unique array 
$\w=\w[G]=\sis{\w_\i}{\i\in\tup}\in\can\tup,$ 
called \rit{\dd\cX generic array} (over $\rL$), all 
\index{array!generic array}%
\index{generic array, $\w[G]$}%
\index{generic array, $\w_\i[G]$}%
\index{zwG@$\w[G]$}%
\index{zwiG@$\w_\i[G]$}%
terms $\w_\i=\w_\i[G]=\w(\i)$ being reals 
(\ie, elements of $\cD=\dn$),  
such that the equivalence 
\bce
$\w\dar\za\in \clo X$ \ \ $\eqv$ \ \ 
$X\in G$   
\ece
holds for all $X\in\cX$ and $\za=\dym X\in\cpo.$ 
Then the model 
$\rL[G]=\rL[\w[G]]=\rL[\sis{\w_\i[G]}{\i\in\tup}]$ is an   
\index{generic extension}%
{\it \dd\cX generic extension of\/} $\rL$. 
{\ubf Equivalently,} an array\/ 
$\w\in\can\tup$ is \dd\cX generic iff the 
set $\cg\w\cap\cX$ is \dd\cX generic over $\rL$, where
$$
\cg\w=\ens{X\in\pei}
{\w\dar\za\in \clo X,\text{ where }\za=\dym X}\sq\pei 
\kmar{cg w}
\index{zGv@$\cg\w$}%
\index{generic set, $\cg\w$}%
$$
and $\clo X$ 
is the topological closure of $X\sq\can\za$ in 
$\can\za$ as above.

\ble
\lam{gras}
Assume that\/ $\cX\sq\pei$, $\cX\in \rL$ is a normal  
forcing, and\/ $\ombl<\omi$. 
If\/ $X\in\cX$ then there is an\/ $\cX$-generic 
(over\/ $\rL$) array\/ $\w\in\can\tup$
satisfying\/ $\w\dar\xi\in \clo X$, where\/ $\xi=\dym X$. 
If\/ $\w$ is such then$:$
\ben
\renu
\itlb{gras1}
if\/ $\cY\in\rL,\,\cY\sq\cX$ is pre-dense in\/ $\cX,$ 
then\/ $\cg\w\cap\cY\ne\pu\;;$

\itlb{gras2}
if\/ $\ta\in\cpo$ and some\/ 
$\cY\in\rL,\,\cY\sq(\cX\dar\ta)$ 
is pre-dense in\/ $\cX\dar\ta$ 
then\/ $\cg\w\cap\cY\ne\pu.$
\vyk{
\itlb{gras3}
if\/ $\ta\in\cpo$ and\/ $\i\in\tup\bez\ta$ then\/ 
$\x(\i)\nin\rL[\x\dar\ta]$, in particular, if\/ 
$\i\ne\j$ belong to\/ $\tup$ then\/ $\x(\i)\ne\x(\j)$.
}
\een
\ele

\bpf
\ref{gras1} is obvious. 
To prove \ref{gras2}, it suffices to show that the set 
\pagebreak[2]
$$
\cY'=\ens{X\in\cX}{\ta\sq\xi=\dym X\land 
\sus Y\in\cY(X\dar\ta\sq Y)}
$$
is dense in $\cX.$ 
{\ubf Arguing in $\rL$}, 
assume that $Z_0\in\cX,\,\et=\dym{Z_0}$. 
Let $\xi=\et\cup\ta$. 
Then $Z=Z_0\uar\xi\in\cX\dar\xi$ and 
$Z_1=Z\dar\ta\in\cX\dar\ta$ by \ref{rfo3}, \ref{rfo2}. 
By the pre-density, $Z_1$ is compatible with some 
$Y\in \cY$, so that there exists 
$U\in\cX\dar\ta,\,U\sq Y\cap Z_1$. 
Then $X=Z\cap(U\uar\xi)\in\cX\dar\xi$ by \ref{rfo3}, 
and $X\dar\ta=U\sq Y$, therefore $X\in \cY'$. 
Moreover $X\sq Z$, hence $X\psq Z_0=Z\uar\et$ by 
construction. 
This ends the density proof.
%
\epf

\bdf
[symmetric subextensions]
\lam{LWOw}
Assume that $\w\in\can\tup$ and $\Om\sq\cpo$. 
We put 
$\W\Om[\w]=
\ens{\rho\akt(\w\dar\et)}{\rho\in\per\land \et\in\Om}
$. 
\index{set!WOmw@$\W\Om[\w]$}%
\index{zWOmv@$\W\Om[\w]$}%
\kmar{W Om w}%
Note the symmetrization by $\per$!

We'll use subclasses $\rL(\W\Om[\w])$ 
\index{models!$\rL(\W\Om[\w])$ }%
\index{zLWOmw@$\rL(\W\Om[\w])$}%
of generic extensions $\rL[\w]$, $\w\in\can\tup$, 
for suitable sets $\Om\sq\cpo$ in $\rL$, as models for 
Theorem~\ref{mt1}. 
By definition, $\rL(\W\Om[\w])$ is the least transitive subclass 
of $\rL[\w]$ containing the set $\W\Om[\w]$ 
and satisfying $\zf$. 
\edf

\parf{Forcing relation}
\las{frel}

Assume that $\cX\in\RF$ is a normal forcing, \ie, 
$\cX\in\rL$ and it holds in $\rL$ that $\cX\in\RF$, 
see Blanket assumption~\ref{inL}. 
To study $\cX$-generic extensions of $\rL$, we make use 
\kmar{fla}%
\index{forcing language $\fla$}%
\index{language!forcing language $\fla$}%
\index{zLcurl@$\fla$, forcing language}%
of a \rit{forcing language} $\fla$, containing the following 
proper $\rL$-class $\nl$ of basic names: 
\kmar{nl}%
\index{basic names, $\nl$}%
\index{zNL@$\nl$, basic names}%
\bit
\item[$-$] $\namx x$ for any $x\in\rL$ --- we'll typically 
\kmar{namx x}%
\index{basic names, $\namx x$}%
\index{zxdot@$\namx x$, basic names}%
\rit{identify} $\namx x$ with $x$ itself, as usual; 
\kmar{namx x}

\item[$-$] $\npv \sg$ for any $\sg\in\per$ --- 
\kmar{npv sg}%
\index{basic names, $\npv \sg$}%
\index{zzsgv@$\npv\sg$, basic names}%
names of this form will be called \rit{unlimited}; 
\index{basic names!unlimited}%

\item[$-$] 
\rit{derived} names $\npvr \sg\et$ 
\index{basic names!derived}%
\kmar{npvr sg et}%
\index{basic names, $\npvr \sg\et$}%
\index{zzsgvet@$\npvr\sg\et$, basic names}%
for any $\sg\in\per$ and $\et\in\cpo$;

\item[$-$] 
in particular $\pv$ and $\pvr\et$ will be  
\index{basic names, $\pv$}%
\index{zv-@$\pv$, basic names}%
\kmar{pv\mns pvr et}%
\index{basic names, $\pvr\et$}%
\index{zvet@$\pvr\et$, basic names}%
shorthands for resp.\ $\npv\ve$ and $\npvr \ve\et$, 
where $\ve\in\per$ is the identity; 
\index{zzeps@$\ve$, the identity}%

\item[$\!\!-$] $\WW \Om$ for any $\Om\in\rL$, $\Om\sq\cpo$.
\kmar{WW Om}%
\index{basic names, $\WW \Om$}%
\index{zWOm@$\WW \Om$, basic names}%
\eit
All those names belong to $\rL$ as 
$\per,\cpo\in\rL$ by Blanket agreement~\ref{inL}. 

The name $\pv$ will be involved as the canonical name for 
a generic array $\w\in\can\tup$. 
Accordingly each $\npv\sg$ will work as a name for 
$\sg\akt\pv$, so in principle it is a derived name. 
Yet we'd like to have each $\npv\sg$ 
 as an independent name so to 
speak, in order to define an action of $\qer$ on 
basic names. 
Accordingly, each derived name $\npvr \sg\et$ 
will work as a name for 
$(\sg\akt\pv)\dar\et=\sg\akt(\pv\dar\et')$, where 
$\et'=\sg\obr\akt\et$ (recall Lemma~\ref{pi1}). 
Finally, $\WW \Om$ is a name for 
$\W\Om[\pv]=
\ens{\rho\akt(\pv\dar\et)}{\rho\in\per\land \et\in\Om}$.

An $\cL$-formula is {\bfit limited\/} iff it   
\index{formula!limited}%
\index{Lformula@$\fla$-formula!limited}%
\index{Lformula@$\fla$-formula!valuation, $\vpi[\w]$}%
contains unlimited names $\npv \pi$ only 
via derived names $\npvr \sg\et$, $\sg\in\qer$ 
and $\et\in\cpo$.

Given $\w\in\can\tup$ in the universe and an $\cL$-formula 
$\vpi$, we define the \rit{valuation} $\vpi[\w]$ by the 
\index{valuation, $\vpi[\w]$}%
\index{zzfi:v:@$\vpi[\w]$, valuation}%
substitution of the valuations resp.
\bce
$\namx x[\w]=x$, \ 
\index{valuation, $\namx x[\w]=x$}%
\index{zx:v:@$\namx x[\w]=x$, valuation}%
$(\npv \sg)[\w]=\sg\akt\w$, \ 
%
\index{valuation, $(\npv \sg)[\w]=\sg\akt\w$}%
\index{zzsgv:v:@$(\npv \sg)[\w]=\sg\akt\w$, valuation}%
%
$\W\Om[\w]=
\ens{\rho\akt(\w\dar\et)}{\rho\in\per\land \et\in\Om}
$
%
\index{valuation, $\W\Om[\w]$}%
\index{zWOmv@$\W\Om[\w]$, valuation}%
\kmar{W Om w}%
\ece
for any basic names resp.\ 
$\namx x$, $\npv \pi$, 
$\WW \Om$ 
in $\nl$ that occur in $\vpi$. 
All those sets belong to the extension 
$\rL[\w]=\rL[\cg\w]$, of course.

\bdf
[forcing relation]
\lam{fred}
Let\/ $\cX\in\RF$ is a normal forcing, 
in particular, $\cX\in\rL$, 
and\/ $\vpi$ be a closed\/ $\cL$-\rit{formula} 
(with names in\/ $\nl$ as parameters). 
Let\/ $X\in\cX$, $\za=\dym X$. 
We define $X\fox\cX\vpi$, iff $\vpi[\w]$ holds in 
\index{forcing!$\fox\cX$}%
\index{zzz::-X@$\fox\cX$, forcing relation}%
$\rL[\w]$ whenever $\w$ is an $\cX$-generic array 
over $\rL$, satisfying $\w\dar\za\in\clo X$.
\edf

\vyk
{\kra The next lemma 
presents some standard forcing properties. 

\ble
\lam{frp}
Under the assumptions of Definition~\ref{fred}$:$
\ben
\renu
\itlb{frp1}
The relation\/ $\fox\cX$ is\/ {\em definable} in\/ $\rL$, 
in the sense that if\/ $\vpi(x_1,\dots,x_n)$ is a 
parameter-free\/ $\in$-formula then there is another 
parameter-free\/ $\in$-formula\/ 
$\vpi'(\cdot,\cdot,x_1,\dots,x_n)$, such that if\/ 
$\cX$ is a normal forcing, $X\in\cX$, and\/ 
$t_1,\dots,t_n\in\nl$, then 
\bce
$X\fox\cX\vpi(t_1,\dots,t_n)$ \ iff\/ \ 
$\vpi'(\cX,X,t_1,\dots,t_n)$ is true in\/ $\rL$.
\ece

\itlb{frp2} 
If\/ $\w$ is an $\cX$-generic array over $\rL$, 
$\vpi$ is a closed $\cL$-\rit{formula}, 
$\w$ is an $\cX$-generic array over $\rL$, 
and $\vpi[\w]$ is true in\/ $\rL[\w]$, then there 
is a condition\/ $X\in \cg\w\cap\cX$ such that\/ 
$X\fox\cX\vpi\;.$
\vyk{\gol
\itlb{frp3} 
If\/ $X\psq Y\fox\cX\vpi$ then\/ $X\fox\cX\vpi$.

\itlb{frp4} 
If\/ $\vpi$ is closed, $X\in\cX$, and\/ 
$X\not\fox\cX\vpi$ then there is a condition\/ 
$Y\in\cX$, $Y\psq X$, $Y\fox\cX\neg\vpi$.

\itlb{frp5} 
If\/ $X\in\cX$, $\dym X\sq\xi\in\cpo$, and\/ 
$X\uar\xi\fox\cX\vpi$, then\/ $X\fox\cX\vpi$.
}%
\qed
\een
\ele
}

The next routine lemma contains an important claim;  
it involves one more definition. 
Suppose that $X\in\pei$ and $\cY\sq\pei$. 
We define 
\bde
\item[$X\sqf\bigcup\cY$,] iff there is a finite 
\index{zzzsqf@$\sqf$}%
\imar{sqf}%
set $\cY'\sq\cY$ such that 
1) $\dym Y\sq\xi=\dym X$ for all $Y\in\cY'$, 
and 2) $X\sq\bigcup_{Y\in\cY'}(Y\uar\xi)$.

\item[$X\sqd\bigcup\cY$,] iff in addition  
\index{zzzsqfd@$\sqd$}%
\imar{sqd}%
3) $(Y\uar\xi)\cap(Z\uar\xi)=\pu$ 
for all $Y\ne Z$ in $\cY'$.
\ede

\ble
\lam{553}
Under the assumptions of Definition~\ref{fred}, if\/ 
$X\in\cX$, $\cY\sq\cX$, $X\sqf\bigcup\cY$, and\/ 
$Y\fox\cX\vpi$ for all\/ $Y\in\cY$, then\/ 
$X\fox\cX\vpi$. 
\ele
\bpf
To check that every $X\in\cX$ satisfying $X\sqf\bigcup\cY$ 
is compatible with some $Y\in\cY$ use \ref{rfo4} of 
Section~\ref{rfo}, and Lemma~\ref{clop}.
\epf

\parf{Forcing and permutations}
\las{fap}

Automorphisms of forcing notions have been widely used to 
define models with various effects related to the axiom of 
choice, basically since Cohen's times. 
Define the left action of permutations $\pi\in\per$ 
on names, as follows: 
\index{permutations!action $\akt$}%
\index{action $\akt$}%
\index{zzzakt@$\;\akt\;$, action}%
\pagebreak[0]
$$
\bay{rcl}
\pi\akt\namx x &=& \namx x;\\[0.7ex]
\pi \akt \npv\sg 
&=& 
\npv{(\sg\supo\pi\obr)}, \ 
\text{ in particular, } \ 
\pi \akt \pv= \npv{(\pi\obr)}; 
\\[0.7ex]
\pi \akt \WW\Om 
&=& 
\WW{\ens{\pi\akts\xi}{\xi\in\Om}}.
\eay
$$
The group action property holds, for instance: 
$$
\ro\akt{(\pi \akt \npv\sg)}=
\ro\akt\npv{(\sg\supo\pi\obr)}=
\npv{(\sg\supo\pi\obr\supo\ro\obr)}=
\npv{(\sg\supo(\ro\supo\pi)\obr)}=
(\ro\supo\pi)\akt\npv\sg.
$$
If $\pi\in\per$ and $\vpi$ is an\/ $\cL$-formula 
then we let $\pi\vpi$ be 
obtained by the substitution of   
$\pi\akt\nu$ 
for any name $\nu$ in $\vpi$.

\vyk{
In particular, $\pv=\npv\ve$ 
(if it occurs in $\vpi$) 
is substituted in $\pi\vpi$ by 
the name $\npv{(\ve\supo\pi\obr)}=\npv{\pi\obr}$. 
}


If $\et\in\cpo$ and $\Om\sq\cpo$ then 
define the following subgroups of $\per$:
$$
\bay{rcl}
\ppi\et &=& 
\kmar{ppi et}%
\ens{\pi\in\per}
{\kaz\i\in\et\,(\i=\pi(\i))},\\[0.7ex]
\index{zzPiet@$\ppi\et$}%
\index{permutation groups!$\ppi\et$}%

\inv\Om &=& 
\kmar{inv Om}%
\index{zInvOm@$\inv\Om$}%
\index{permutation groups!$\inv\Om$}%
\ens{\pi\in\per}
{\kaz\xi\in\cpo\,(\xi\in\Om\eqv \pi\akt\xi\in\Om)}.
\vyk{
\gga(\cX)&=&
\ens{\pi\in\per}
{\kaz X\in\pei\,(X\in\cX\eqv \pi\akt X\in\cX)}. 
\index{zzGaX@$\gga(\cX)$}%
}
\eay
$$
If $\vpi$ is an $\cL$-formula,  
then let 
$$
\bay{rcl}
\inv\vpi 
\kmar{inv vpi}%
\index{zInvfi@$\inv\vpi$}%
\index{permutation groups!$\inv\vpi$}%
&=&
\bigcap\ens{\inv\Om}
{\Om=\cpo\lor \WW\Om\text{ occurs in }\vpi};\\[0.9ex]
\modd\vpi
&=&
\bigcup\ens{\sg\obr\akt\et}
{\npvr\sg\et\text{ occurs in }\vpi}, 
\ \ \text{thus $\modd\vpi\in\cpo$}.
\index{z11vpi11@$\modd\vpi$, dimension}%
\index{dimension, $\modd\vpi$}%
\kmar{modd vpi}%
\eay
$$


\ble
\lam{61}
Let\/ $\vpi$ be an\/ $\cL$-formula and\/ $\w\in\can\tup$. 
Then:
\ben
\renu
\itlb{611} 
if\/ $\pi\in\per$  then the formulas\/ 
$\vpi[\w]$ and\/ $(\pi\vpi)[\pi\akt\w]$ coincide$;$

\itlb{612} 
$\modd{\pi\vpi}=\pi\akt\modd\vpi$, 
and if\/ $\pi\in\inv\vpi$ 
then any name\/ $\WW\Om$ in\/ $\vpi$ does not change 
in\/ $\pi\vpi\;;$ 

\itlb{613} 
if\/ $\pi\in\ppi{\modd\vpi}\cap\inv\vpi$, 
and\/ $\vpi$ is a limited formula, 
then the formulas\/ 
$\vpi[\w]$, $(\pi\vpi)[\w]$ coincide$.$
\een
\ele
\bpf
\ref{611} 
Let $\npv\sg$ occur in $\vpi$. 
Then it changes to $\npv{(\sg\supo\pi\obr)}$ in $\pi\vpi$. 
It remains to note that by the group action property 
\pagebreak[0]
$$
(\sg\supo\pi\obr)\akt(\pi\akt\w)=
(\sg\supo\pi\obr\supo\pi)\akt\w=
\sg\akt\w.
$$
\vyk{
$$
\bay{ccccccccc}
\rhi\akt(\pi\akt\w)
&=&
(\rho\supo\pi\obr)\akt(\w\supo\pi\obr)
&=&(\w\supo\pi\obr)\supo(\rho\supo\pi\obr)\obr&=&\\[1ex]
&=&\w\supo\pi\obr\supo\pi\supo\rho\obr   
&=&\w\supo\rho\obr=\rho\akt\w\,,
\eay
$$
}%
%
\vyk{
The case of names $\npvr\rho\et$ is similar:
\bce
$\rhi\akt((\pi\akt\w)\dar\eti)=
\rho\akt\pi\obr\akt(\pi\akt(\w\dar\et))
=\rho\akt(\w\dar\et)$.
\ece 
}

Further, any name $\WW\Om$ in $\vpi$ changes to 
$\WW{\Om'}$, where $\Om'=\ens{\pi\akt\xi}{\xi\in\Om}$. 
Using Lemma~\ref{pi1}, we obtain: 
\pagebreak[0]
$$
\bay{rclcccccc}
\W{\Om'}[\pi\akt\w]
&=&
\ens{\rho\akt((\pi\akt\w)\dar\eti)}
{\rho\in\per\land \eti\in\Om'}
\\[1ex]
&\hspace*{-25.8ex}=&
\hspace*{-14ex}
\ens{\rho\akt((\pi\akt\w)\dar(\pi\akt\et))}
{\rho\in\per\land \et\in\Om}
=
\ens{\rho\akt(\pi\akt(\w\dar\et))}
{\rho\in\per\land \et\in\Om}\\[1ex]
&\hspace*{-25.8ex}=&
\hspace*{-14ex}
\ens{(\rho\supo\pi)\akt(\w\dar\et))}
{\rho\in\per\land \et\in\Om} 
= \ens{\rhi\akt(\w\dar\et))}
{\rhi\in\per\land \et\in\Om}, 
\eay
$$
because $\ens{\rho\supo\pi}{\rho\in\per}=\per$.

\ref{612} 
If $\WW\Om$ is a name in $\vpi$ then 
it changes to $\WW{\Om'}$ in $\pi\vpi$, 
where $\Om'=\ens{\pi\akt\et}{\et\in\Om}=\Om$ 
since $\pi\in\gga(\Om)$. 
This $\WW{\Om'}$ is identical to $\WW{\Om}$. 
Further, 
\pagebreak[0]
$$
\bay{rclccc}
\modd{\pi\vpi}
&=&
\bigcup\ens{\sgi\obr\akt\et}
{\npvr\sgi\et\text{ occurs in }\pi\vpi}\\[0.7ex]
&=&
\bigcup\ens{(\sg\supo\pi\obr)\obr\akt\et}
{\npvr\sg\et\text{ occurs in }\vpi}\\[0.7ex]
&=&
\bigcup\ens{\pi\akt(\sg\obr\akt\et)}
{\npvr\sg\et\text{ occurs in }\vpi}
&=&
\pi\akt{\modd\vpi}\,.
\eay
$$
\vyk{
Further, if $\i\in\ima\pi{\modd\vpi}$ then 
$\i=\pi(\j)$, where $\j\in\modd\vpi$, \ie, 
$\j\in\et$ and a 
name $\npvr\rho\et$ occurs in $\vpi$. 
To conclude, 
\bce 
\hfill
$\ima\pi{\modd\vpi}
=\ens{\pi(\j)}
{\j\in\et\land\npvr\rho\et\text{ occurs in }\vpi}$. 
\hfill(1)
\ece
On the other hand, if $\i\in \modd{\pi\vpi}$ then 
$\i\in\eti$, for some  
name $\npvr{\rhi}\eti$ in $\pi\vpi$. 
Then $\rhi=\rho\supo\pi\obr$ and $\eti=\pi\akt\et$, where 
$\npvr\rho\et$ occurs in $\vpi$, hence
\bce 
\hfill
$\modd{\pi\vpi}
=\ens{\pi(\j)}
{\j\in\et\land\npvr\rho\et\text{ occurs in }\vpi}$, 
\hfill(2)
\ece
where the right-hand side is obviously equal 
to the right-hand side of (1).
}

\ref{613} 
If $\npvr\sg\et$ occurs in $\vpi$ then 
it changes to $\npvr{(\sg\supo\pi\obr)}{\et}$ 
in $\pi\vpi$. 
The $\w$-valuation of 
$\npvr{(\sg\supo\pi\obr)}{\et}$ is equal 
(by Lemma~\ref{pi1}) to 
\bce
$
\sg\akt(\pi\obr \akt\w)\dar{\et} 
=
\sg\akt((\pi\obr \akt\w)\dar{(\sg\obr\akt\et)})
=
\sg\akt(\w\dar{(\sg\obr\akt\et)})
=
(\sg\akt \w)\dar\et  
$,
\ece
since $\pi$ and $\pi\obr$ are the identities on 
$\sg\obr\akt\et$ 
(because $\pi\in \ppi{\modd\vpi}$). 
But this is equal to the $\w$-valuation of 
the original derived name $\npvr\sg\et$ in $\vpi$.

If $\WW\Om$ is a name in $\vpi$ then 
it does not change in $\pi\vpi$ by \ref{612}.
\epf

\bte
\lam{621}
Assume that, in\/ $\rL$, $\cX\in\RF$ is a normal forcing, 
$\vpi$ is a closed 
$\cL$-\rit{formula}, and\/ $\pi\in\qer$. 
Let $X\in\cX$. 
Then $X\fox\cX\vpi$ iff\/ $\pi\akt X\fox\cX \pi\vpi$.
\ete

\bpf
As $\cX,\pi\in\rL$ (see Blanket agreement \ref{inL}), 
an array $\w\in\can\tup$ is $\cX$-generic over $\rL$ iff 
so is $\pi\akt\w$. 
Now the result follows from Lemma~\ref{61}\ref{611}.
\epf

\bcor
\lam{623}
Under the assumptions of Theorem~\ref{621}, suppose 
that\/ $\ta\in\cpo$, $\vpi$ is a closed limited formula, 
$\modd\vpi\sq\ta$, 
$\pi\in\ppi\ta\cap 
\inv\vpi$, 
$X\in\cX$.
Then $X\fox\cX\vpi$ iff\/ $\pi\akt X\fox\cX \vpi$.
\ecor
\bpf
The result follows from Theorem~\ref{621} and 
Lemma~\ref{61}\ref{613}.
\epf

\bcor
\lam{624}
Under the assumptions of Thm~\ref{621}, 
let\/ $\ta\sq\et$ belong to\/ $\cpo$, 
$\vpi(x)$ be a limited formula, 
$\modd\vpi\sq\ta$, 
$\pi\in\ppi\ta
\cap\inv\vpi$, 
$X\in\cX$, $\sg=\ima\pi\et$.
Then\/ $X\fox\cX(\sus x\in\rL[\pvr\et])\,\vpi(x)$ 
iff\/ 
$\pi\akt X\fox\cX (\sus x\in\rL[\pvr\sg])\,\vpi(x)$.
\ecor
\bpf
Assume that $X\fox\cX(\sus x\in\rL[\pvr\et])\,\vpi(x)$. 
Then, by Theorem~\ref{621}, 
$\pi\akt X\fox\cX
(\sus x\in\rL[\npvr{\pi\obr}{\et}])\,\pi\vpi(x)$. 
Yet if $\w\in\can\tup$ then, by Lemma~\ref{pi1}, 
$(\pi\obr\akt\w)\dar\et=\pi\obr\akt(\w\dar\sg)$,  
hence obviously 
$\rL[(\pi\obr\akt\w)\dar\et]=\rL[\w\dar\sg]$. 
We conclude that   
$\pi\akt X\fox\cX
(\sus x\in\rL[\pvr\sg])\,\pi\vpi(x)$.  
And finally, here we can replace $\pi\vpi(x)$ by 
$\vpi(x)$ by 
Lemma~\ref{61}\ref{613}.
\epf

\parf{Isolation and the narrowing theorem}
\las{63}

Suppose that $\et\in\cpo$. 
It often happens in similar cases that sentences 
relativized to $\rL[\w\dar\et]$ are decided by forcing 
conditions $X$ satisfying $\dym X\sq\et$. 
The following theorem belongs to this category. 

\bdf
[in $\rL$]
\lam{63d}
Assume that $\Ga\sq\per$ is a subgroup. 
Say that $\et\in\cpo$ is \dd\Ga\rit{isolated} 
\index{isolated!$\Ga$-isolated}%
if (*) for each $\xi\in\cpo$ with $\et\sq\xi$ 
there is a permutation $\pi\in\Ga\cap\ppi\et$ 
satisfying $\xi\cap(\pi\akt\xi)=\et$.
\edf

\ble
[in $\rL$]
\lam{iso1}
Each\/ $\et\in\cpo$ is $\qer$-isolated. 
\ele

\bpf 
Let $\et\sq\xi\in\cpo$; 
define  
$\pi\in \ppi\et$ with  
$\xi\cap(\pi\akt\xi)=\et$. 
Let $\la<\omi$ be a limit ordinal $>$ all 
ordinals $\j(k)$, where $\j\in\xi$ and $k<\lh\j$. 

Define, in $\rL$, 
$B:\omi\onto\omi$ by $B(\ga)=B\obr(\ga)=\la+\ga$ 
for all $\ga<\la$, and $B(\ga)=\ga$ for $\ga\ge\la+\la$. 
If $\i\in\tup$ then define $\ba(\i)=\i'\in\tup$ 
so that $\lh{\i'}=\lh\i$ and $\i'(\ell)=B(\i(\ell)$ for all 
$\ell<\lh\i$. 
Clearly $\ba\in\per$. 

Now let $\i\in\tup$. 
There is a largest number $m_\i\le\lh\i$ such that 
$\i\res m_\i\in\et$. 
Then $\i=(\i\res m_\i)\we \k$ for some $\k\in\tup\cup\ans\La$. 
Put $\pi(\i)=(\i\res m_\i)\we \ba(\k)$. 
\epf

\vyk{
If $\i\in\et$ then put $\pi(\i)=\i$. 

Now define $\pi(\i)$ in case $\i\nin\et$. 
There is $m_\i<\lh\i$ such that still 
$\i\res m_\i\in\et$ but $(\i\res m_\i)\we \i(m_\i)\nin\et$; 
possibly $m_i=0$. 
Consider the sets 
$H=\ens{\i'(m_\i)}{\i'\in\et\land\lh{\i'}>m_\i}$ 
and $D=\la\bez H$. 
Define $\j=\pi(\i)\in\tup$ so that  
$\lh\j=\lh\i$, $\j\res m_\i=\i\res m_\i$, 
$\j(l)=\i(l)$ in case $m_\i<l<\lh\i$, and 
{\def\theequation{\fnsymbol{equation}}
\addtocounter{equation}1
\busq{eq33}{
\j(m_\i)=
\left\{
\bay{rcl}
\la+\al,&\text{ in case }&
\al=\i(m_\i)\in D;\\[0.8ex] 
\al,&\text{ in case }&
\i(m_\i)=\la+\al, \; \al\in D;\\[0.8ex]
\i(m_\i),&\text{ in case }&
\i(m_\i)\nin D\cup\ens{\la+\al}{\al\in D}.
\eay
\right|
}}%
Easily $\pi$ is as required.
}

\bte
[the narrowing theorem, in $\rL$]
\lam{631}
Assume that\/ $X\in\cX\in\RF$, 
$\vpi$ is a closed limited\/ $\cL$-formula, 
$\et\in\cpo$ is\/ $
\inv\vpi$-isolated, 
and\/ $\modd\vpi\sq\et\sq\dym X$. 
Then $X\fox\cX\vpi$ iff\/ $X\dar\et\fox\cX \vpi$.
\ete

\bpf
Suppose to the contrary that 
$X\fox\cX\vpi$ but $X\dar\et\not\fox\cX \vpi$. 
There is a condition 
$U\in\cX$ such that $U\psq (X\dar\et)$ and 
$U\fox\cX\neg\,\vpi$. 
Let $\xi=\dym X$, $\ta=\dym U$. 
By (*) of Definition~\ref{63d}, there is a permutation 
$\pi\in 
\inv\vpi\cap\ppi\et$ 
satisfying 
$(\pi\akt{(\xi\cup\ta)})\cap(\xi\cup\ta)=\et$, in 
particular, $(\pi\akt\ta)\cap\xi=\et$. 

Let $Y=\pi\akt U$ and $\za=\modd Y=\pi\akt\ta$. 
Then $Y\in\cX$ (since $\pi\in\gga(\cX)$, 
$\xi\cap\za=\et$, and (most important!) 
$Y\fox\cX\neg\,\vpi$ by Corollary~\ref{623}. 

Furthermore, $Y\dar\et=U\dar\et$ (since $\pi\in\ppi\et$), 
in particular, $Y\dar\et\sq X\dar\et$. 
Therefore $X'=X\cap(Y\dar\et\uar\xi)\in\cX$, 
$X'\sq X$, $X'\dar\et=Y\dar\et$. 
Let $\vt=\xi\cup\za$.
It follows by 
Lemma~\ref{rfoL} \poo\  Corollary~\ref{pe2} 
that the set $Z=(X'\uar\vt)\cap(Y\uar\vt)$ belongs 
to $\cX$, and obviously $Z\psq Y$ and $Z\psq X'\sq X$. 
Thus $X$ and $Y$ are compatible in $\cX$. 
But $X,Y$ force contradictory sentences.
\epf

\bcor
\lam{632}
Assume that\/ $\cX\in\RF$, 
$\i\in\tup\bez\et$, and\/ $\w\in\can\tup$ is\/ 
$\cX$-generic. 
Then $\w(\i)\nin\rL[\w\dar\et]$.
\ecor

\bpf
Suppose towards the contrary that 
$\w(\i)\nin\rL[\w\dar\et]$.
Then there is a parameter-free $\in$-formula 
$\vpi(\cdot,\cdot,\cdot)$, and a parameter 
$p\in\rL$, such that, 
\bce
for all $k<\om$: \ \  
$\w(\i)(k)=1$ iff  
$\rL[\w\dar\et]\mo\vpi(p,\w\dar\et,k)$.
\ece
Then there exists such a condition $X\in\cX\cap \cg\w$ 
that 
\bce
\hfill
$X\fox\cX 
\kaz k\,\big(\pv(\i)(k)=1\eqv\vpi(\namx p,\pv\dar\et,k\big)$.
\hfill
(1)
\ece
Let $\xi=\dym X$. 
We may assume that $\et\sq\xi$, as otherwise replace 
$X$ by $X'=X\uar(\et\cup\xi)$, which still belongs 
to $\cX$ by \ref{rfo3} of Section \ref{rfo}. 
And we may assume that $\i\in\xi$ by the same reason.
Lemma~\ref{99y} implies that there exists $k<\om$ and 
sets $Y,Z\in\pe\xi$, clopen in $X$ and such that 
$Y\dar\et=Z\dar\et$ and $y(\i)(k)=1$ but  
$z(\i)(k)=0$ for all $y\in Y$ and $z\in Z$  
(or vice versa). 
Then $Y,Z\in\cX$ by \ref{rfo4} of Section~\ref{rfo}, 
and $Y\fox\cX \pv(\i)(k)=1$ but $Z\fox\cX \pv(\i)(k)=0$.

It follows by (1) that 
$Y\fox\cX \vpi(\namx p,\pv\dar\et,k\big)$, hence 
$Y\dar\et\fox\cX \vpi(\namx p,\pv\dar\et,k\big)$ 
by Theorem~\ref{631} 
(applicable by Lemma~\ref{iso1}). 
We have 
$Z\dar\et\fox\cX \neg\,\vpi(\namx p,\pv\dar\et,k\big)$ 
by the same reasons.
However $Y\dar\et=Z\dar\et$, which is a contradiction.
\epf

\bcor
\lam{633}
Assume that\/ $\cX\in\RF$, 
$\vpi$ is a closed limited\/ $\cL$-formula, 
$\et\in\cpo$ is\/ $
\inv\vpi$-isolated, 
and\/ $\modd\vpi\sq\et$, $\w\in\can\tup$ is\/ 
$\cX$-generic, and\/ $\rL[\w]\mo\vpi[\w]$. 
Then there is\/ $X\in\cX\dar\et\cap \cg\w$ 
such that\/ $X\fox\cX\vpi$.
\ecor
\bpf
The set 
$ 
\cD=\ens{X\in\cX\dar\et}
{X\fox\cX\vpi\,\text{ or }\,X\fox\cX\neg\,\vpi}
$ 
is pre-dense in $\cX$ by Theorem~\ref{631}.
\epf

\bcor
\lam{634}
Assume that\/ $\cX\in\RF$, 
$\vpi(\cdot)$ is a limited\/ $\cL$-formula, 
$\et\in\cpo$ is\/ $
\inv\vpi$-isolated, 
and\/ $\modd\vpi\sq\et$, $\w\in\can\tup$ is\/ 
$\cX$-generic, and\/ $A\in\rL$. 
Then the set\/ $S=\ens{a\in A}{\rL[\w]\mo\vpi(a)}$ 
belongs to\/ $\rL[\w\dar\et]$.
\ecor
\bpf
We have 
$S=\ens{a\in A}
{\sus X\in\cX\dar\et\cap \cg\w\,
(X\fox\cX\vpi(\namx a))}$. 
On the other hand, 
$\cX\dar\et\cap \cg\w=
\ens{X\in\pe\et}{\w\dar\et\in\clo X}\in\rL[\w\dar\et]$.
\epf

\bcor
\lam{641}
Assume that\/ $\cX\in\RF$, $\Om\in\rL$, $\Om\sq\cpo$ 
is\/ $\cup$-closed (under finite unions),
all\/ $\et\in\Om$ are\/ 
$
\inv\Om$-isolated, 
$\w\in\can\tup$ is\/ $\cX$-generic, 
and\/ $S\in\rL(\W\Om[\w])$, $S\sq \rL$. 
Then\/ $S\in\rL[\w\dar\et]$ for some\/ $\et\in\Om$.
\ecor
\bpf
First of all, $S\sq A$ for some $A\in\rL$. 
Then, as $S\in\rL(\W\Om[\w])$, we have 
$S=\ens{a\in A}{\rL[\w]\mo\vpi(a)}$, 
where $\vpi$ contains only 
$x\in\rL$, $\W\Om[\w]$, and some $\w\dar\et$, $\et\in\Om$, 
as parameters. 
Then\/ $S\in\rL[\w\dar\et]$ by Corollary~\ref{634}.
\epf

\bcor
\lam{642}
Assume that\/ $\cX$, $\Om$ are as in Corollary~\ref{641}, 
$\psi(\cdot)$ is a limited\/ $\cL$-formula, 
$X\in\cX$, $A\in\rL$,   
$X\fox\cX\,\sus S\in \rL(\WW\Om)\,
(S\sq A\land \psi(x))$. 
Then there exists a condition\/ $Y\in\cX,$ and\/ 
$\et\in\Om,$ such that\/ 
$Y\psq X$ and 
$Y\fox\cX\,\sus S\in \rL[\pvr\et]\,
(S\sq A\land \psi(x))$.  
\ecor
\bpf
By Lemma~\ref{gras}, there exists a \dd\cX generic 
array $\w\in\can\tup$ satisfying $X\in\cg\w$. 
There is $S\in\rL(\W\Om[\w])$ such that 
$\rL[\w]\mo\psi(S)$ and $S\sq A$. 
We have $S\in\rL[\w\dar\et]$ for some $\et\in\Om$ 
by Corollary~\ref{641}. 
Then some $Z\in\cg\w\cap\cX$ satisfies 
$Z\fox\cX\,\sus S\in \rL[\pvr\et]\,
(S\sq A\land \psi(x))$. 
But $Z$ and $X$ are compatible in $\cX$, so take any 
$Y\in\cX$ with $Y\psq X$ and $Y\psq Z$.
\epf

\vyk{
\bcor
\lam{67old}
Assume that\/ $\cX\in\RF$, 
$\vpi(\cdot)$ is a limited\/ $\cL$-formula, 
$\modd\vpi=\pu$, the set\/ $\pu$ is\/ 
$\gga(\cX)\cap\gga(\vpi)$-isolated, $\ta\in\cpo$, 
$\pi\in\gga(\cX)\cap\gga(\vpi)$, $\sg=\ima\pi\ta$. 
Suppose that\/ $\w\in\can\tup$ is\/ $\cX$-generic. 
Then, in\/ $\rL[\w]:$
$$ 
\sus x\in\cN\cap\rL[\w\dar\ta]\,\vpi(x)[\w]
\;\imp\;
\sus x\in\cN\cap\rL[\w\dar\sg]\,\vpi(x)[\w]\,.
$$
\ecor

\bpf
Recall that $\bon\in\cX\dar\pu$ 
by \ref{rfo2} of Section~\ref{rfo2}, 
and obviously $\dym \bon=\pu$. 
We conclude that 
$\bon\fox\cX \sus x\in\cN\cap\rL[\pv\dar\ta]\,\vpi(x)$ 
by Corollary~\ref{633} (with $\et=\pu$). 
It follows by \ref{624} that 
$\bon\fox\cX \sus x\in\cN\cap\rL[\pv\dar\sg]\,\vpi(x)$, 
because $\pi\akt\bon=\bon$.  
This implies the result required as 
definitely $\bon\in G_\w$.
\epf
}

\vyk{

\bcor
\lam{67}
Assume that\/ $\cX\in\RF$, 
$\vpi(\cdot)$ is a limited\/ $\cL$-formula, 
$\modd\vpi=\xi\sq\ta\in\cpo$, 
$\ta$ is\/ $\gga(\cX)\cap\gga(\vpi)$-isolated, 
$\pi\in\gga(\cX)\cap\gga(\vpi)\cap\ppi\xi$, 
$\sg=\ima\pi\ta$, $\sg\cap\ta=\xi$. 
Suppose that\/ $\w\in\can\tup$ is\/ $\cX$-generic. 
Then, in\/ $\rL[\w]:$
$$ 
\sus x\in\cN\cap\rL[\w\dar\ta]\,\vpi(x)[\w]
\;\imp\;
\sus x\in\cN\cap\rL[\w\dar\sg]\,\vpi(x)[\w]\,.
$$
\ecor

\bpf
Assume the left-hand side. 
By Corollary~\ref{633}, 
there is a condition $X\in(\cX\dar\ta)\cap \cg\w$ 
such that\/ 
$X\fox\cX\sus x\in\bn\cap\rL[\pvr\ta]\,\vpi(x)$.

Recall that $\bon\in\cX\dar\pu$ 
by \ref{rfo2} of Section~\ref{rfo2}, 
and obviously $\dym \bon=\pu$. 
We conclude that 
$\bon\fox\cX \sus x\in\bn\cap\rL[\pv\dar\ta]\,\vpi(x)$ 
by Corollary~\ref{633} (with $\et=\pu$). 
It follows by Corollary~\ref{624} that 
$\bon\fox\cX \sus x\in\bn\cap\rL[\pv\dar\sg]\,\vpi(x)$, 
because $\pi\akt\bon=\bon$.  
This implies the result required as 
definitely $\bon\in G_\w$.
\epf

}

\vyk{ 
\ble
\lam{ahu7} 
Assume that\/ $\rP\in\bfr\al$, $2\le\al<\omi$,  
sets\/ $\et\su\xi\sq\tuq\al$ in\/ $\cpo$ are finite, 
$\phi:\om\to\tuq\al$ is arbitrary,   
$\sis{X_u}{u\in 2^m}$ is a $\phi$-\cohe\ system 
in\/ $\xt\rP_\et$, 
$u_0\in2^m,$ $Y\in\xt\rP_\xi$, $Y\dar\et\sq X_{u_0}$.
Then there is a $\phi$-\cohe\ system\/  
$\sis{Y_u}{u\in 2^m}$ in\/ $\xt\rP_\xi$ 
satisfying\/ $Y_{u}\dar\et\sq X_u$, $\kaz u$, 
and\/ $Y_{u_0}=Y$. 
\ele

\bpf
First of all, we can \noo\ assume that, even more, 
$Y\dar\et=X_{u_0}$ exactly, 
just by Lemma~\ref{ahuL}\ref{ahu6} \poo\ 
Lemma~\ref{suz}.
Further, as $\xi,\et$ are finite, we can assume that 
$\xi=\et\cup\ans\i$, $\i\in\xi\bez\et$. 
Then $\ilq\i\cap\et=\ile\i$.

If $u\in2^m$ then, by \ref{fr3} above, we can define 
a set $B_u\in\pro\rP\i$ satisfying
 
(1) $B_u\rsl\i= X_u\rsl\i$, 

(2) if $u,v\in2^m$ and $X_u\rsl\i=X_v\rsl\i$ 
then $B_u=B_v$, and 

(3) $B_{u_0}=Z\rsq\i$.

\noi 
Now, each set $Y_u= (X_u\uar\xi)\cap(B_u\uar\xi)$ 
belongs to $\pe\xi$ by Corollary~\ref{pe2} 
since $B_u\rsl\i= X_u\rsl\i$ and  
$\ilq\i\cap\et=\ile\i$. 
It follows that in fact $Y_u\in\xt\rP_\xi$, because if 
$\j\in\et$ then $Y_u\rsq\j=X_u\rsq\j\in\pro\rP\j$, 
whereas if $\j\in\ilq\i$ then 
$Y_u\rsq\j=B_u\rsq\j\in\pro\rP\j$ as well. 

It remains to check 
$Y_u\dar\xi_\phi[u,v]=Y_v\dar\xi_\phi[u,v]$.\vom

{\it Case 1\/}: $\xi_\phi[u,v]\sq\et$. 
Then $\xi_\phi[u,v]=\et_\phi[u,v]$, and hence 
\bce
$Y_u\dar\xi_\phi[u,v]=X_u\dar\et_\phi[u,v]
=X_v\dar\et_\phi[u,v]=Y_v\dar\xi_\phi[u,v]$.
\ece 
(The middle equality holds 
by \ref{prct} for the system $\sis{X_u}{u\in2^m}$.)\vom

{\it Case 2\/}: $\xi_\phi[u,v]\not\sq\et$. 
Then obviously $\i\in\xi_\phi[u,v]$, hence 
$\xi_\phi[u,v]=\et_\phi[u,v]\cup\ans\i$ and 
$\ile\i\sq\et_\phi[u,v]$. 
We remind that 
$X_u\dar\et_\phi[u,v]=X_v\dar\et_\phi[u,v]$ by 
\hbox{\ref{prct},} 
and in particular $X_u\rsl\i=X_v\rsl\i$, 
and hence $B_u=B_v$ by (2).  
Therefore (*)
$Y_u\rsdq\i=B_u=B_v=Y_v\rsdq\i$. 
On the other hand,  
$Y_u\dar\et=X_u$ and $Y_v\dar\et=X_v$ 
by Corollary~\ref{pe2}, so that 
\bce
$Y_u\dar\et_\phi[u,v] =X_u\dar\et_\phi[u,v]
=X_v\dar\et_\phi[u,v]=Y_u\dar\et_\phi[u,v]$. 
\ece
We conclude that, by (*),  
\bce
$Y_u\dar(\et_\phi[u,v]\cup\ilq\i) =
Y_u\dar(\et_\phi[u,v]\cup\ilq\i)$  
\ece
by Lemma~\ref{no19}, thus 
$Y_u\dar\xi_\phi[u,v] = Y_u\dar\xi_\phi[u,v]$ 
as 
$\xi_\phi[u,v]=\et_\phi[u,v]\cup\ans\i$. 
\epf
}

\parf{Fusion property}
\las{fup}

{\ubf Arguing in $\rL$}, let $\cX\in \RF$. 
We introduce:
\bde 
\item[{\ubf Fusion property:}] 
for any sequence 
\index{property!Fusion}%
\index{Fusion property}%
$\sis{\cY_k}{k<\om}\in\rL$ of dense sets $\cY_k\sq\cX$, 
the set 
$\cY=\ens{X\in\cX}{\kaz k\,(X\sqd\bigcup\cY_k)}$ 
is dense in $\cX$ as well. 
\ede
(See before Lemma~\ref{553} on $\sqd$.)  
The Fusion property is another formalization of some 
features of the Sacks forcing. 
It somewhat differs from a more commonly used \rit{Axiom A} 
(see Jech~\cite[Definition 31.10]{jechmill}), but 
it fits better to applications in this paper.
The following theorem presents several principal 
applications. 

\bte
\lam{1022}
Assume that, in\/ $\rL$, $\cX\in \RF$ 
has the Fusion property, and\/ $\w\in\can\tup$ 
is\/ $\cX$-generic over\/ $\rL$. 
Then$:$\vim
\ben
\renu
\itlb{1022iv}
if\/ $h\in\rL[\w]$, $h:\om\to\rL$, then there is a map\/ 
$H\in\rL$ such that\/ $\dom H=\om$, and, for each\/ 
$k<\om$, $h(k)\in H(k)$ and\/ $H(k)$ is finite$;$

\itlb{1022i} 
every\/ $\rL$-cardinal remains a cardinal in\/ 
$\rL[\w]\,;$ 

\itlb{1022ii}
if\/ $x\in\bn\cap\rL[\w]$ then\/ 
$x\in\rL[\w\dar\xi]$ for some\/ $\xi\in\cpo$, 
and more general, if\/ $J\in\rL\yd J\sq\tup$ is an 
initial segment and\/ $x\in\bn\cap\rL[\w\dar J]$ then\/ 
$x\in\rL[\w\dar\xi]$ for some\/ $\xi\in\cpo\yd\xi\sq J\,;$ 

\itlb{1022iii}
if\/ $\xi\in\cpo$ 
and\/ $a\in\cN\cap\rL[\w\dar\xi]$ 
then there is a continuous map\/ $F:\can\xi\to\cN$ 
such that\/ $a=F(\w\dar\xi)$, and\/ 
$F$ is\/ {\em coded in $\rL$} 
in the sense that the restriction\/ 
$F_\rL=F\res(\rL\cap\can\xi)$ belongs to\/ $\rL$. 

\een
\ete

Note that if $F_\rL=F\res(\rL\cap\can\xi)\in\rL$ 
in \ref{1022iii} then   
$\rL\mo$ ``$F_\rL:\can\xi\to\cN$ is continuous'' and 
$F=\clo {F_\rL}$ 
(the topological closure of $F_\rL$ in $\can\xi\ti\cN$). 

\bpf
\ref{1022iv}  
There is an 
$\in$-formula $\vpi(v,k,m)$, with ordinals as 
parameters, such that 
$h=\ens{\ang{k,x}\in\om\ti\rL}{\rL[\w]\mo\vpi(\w,k,m)}$,  
and  
\ben
\nenu
\itlb{6624}
if $X\in\cX$ then 
$X\fox\cX\,
\big(
\ens{\ang{k,x}}{\vpi(\w,k,m)}\text{ is a map }\om\to\rL
\big)$.
\een
{\ubf Arguing in $\rL$}, define the sets 
$\cZ_{m}=
\ens{X\in\cX}{\sus x\,(X\fox\cX\vpi(\pv,m.\namx x)}$. 
By \ref{6624}, each $\cZ_m$ is open dense in $\cX$. 
Thus 
$\cZ=\ens{X\in\cX}{\kaz m\,(Z\sqd\bigcup\cZ_m)}$ 
is dense as well by the Fusion property. 
It follows that there exists $Z\in\cZ\cap\cg\w$, so that 
for each $m$ there exists a finite subset $\cY_m\sq\cZ_m$ 
with $Z\sq\bigcup_{Y\in \cY_m}(Y\uar\za)$, where 
$\za=\dim Z$ and $\dim Y\sq\za$ for all $Y\in\cY_m$. 

By definition and \ref{6624}, for each $m<\om$ and 
$Y\in\cY_m$ there is a unique set $x_{mY}\in\rL$ 
satisfying $Y\fox\cX\vpi(\w,m,\namx x_{mY})$. 
Let $H(m)=\ens{x_{mY}}{Y\in\cY_m}$. 
Then $H$ is as required by Lemma~\ref{553}.

\ref{1022i} is a simple corollary of \ref{1022iv},
and \ref{1022ii} is a simple corollary of \ref{1022i}. 

\vyk{
We have to prove that $\omil$ is preserved; 
the preservation of bigger cardinals holds by 
a routine cardinality argument of forcing theory. 

Suppose to the contrary that $f\in\rL[\w]$, 
$f:\om\onto\omil$. 
There is a formula $\vpi(\w,k,\al)$, 
with sets in $\rL$ as parameters, 
and such that 
$f=\ens{\ang{k,\al}\in\om\ti\omil}
{\rL[\w]\mo\vpi(\w,k,\al)}$, 
and in addition
\ben
\nenu
\itlb{6621}
if $X\in\cX$ then 
$X\fox\cX\big(\ens{\ang{k,\al}}{\vpi(\pv,k,\al)} 
\text{ is a map }\om\to\omil\big), $

\itlb{6622}
there exists $X_0\in\cX\cap G_\w$,  
$X_0\fox\cX\kaz\al<\omil\,\sus k<\om\,\vpi(\pv,k,\al)$.
\een
By \ref{6621}, if $k<\om$ then 
the set 
\bce
$\cY_k=\ens{X\in\cX}
{\sus\al<\omil\big(X\fox\cX\vpi(\pv,k,\al)\big)}$
\ece
is open dense in $\cX$. 
Therefore 
$\cY=\ens{X\in\cX}{\kaz k\,(X\sqf\bigcup\cY_k)}$ 
is dense as well, by the Fusion property. 
Take $X_0$ as in \ref{6622}, 
and let $X\in\cY$, $X\psq X_0$. 

If $k<\om$ then there is a finite $\cY'_k\sq\cY_k$ 
satisfying $\dym Y\sq\xi=\dym X$ for all $Y\in\cY'_k$, 
and $X\sq\bigcup_{Y\in\cY'_k}(Y\uar\xi)$.
Then the set 
\bce
$A_k=\ens{\al<\omil}
{\sus Y\in\cY'_k \big(Y\fox\cX\,\vpi(\pv,k,\al)\big)}$
\ece
is finite by \ref{6621}, hence $A=\bigcup_kA_k\in\rL$ 
is at most countable in $\rL$. 
Therefore there is $\al<\omil$, $\al\nin A$. 
Still by \ref{6621}, there exist some $k<\om$ 
and $Z\in\cX$, $Z\psq X$, $Z\fox\cX\vpi(k,\al)$.
But $Z$ has to be compatible with some $Y\in\cY'_k$, 
leading to a contradiction by routine arguments. 
}

\ref{1022iii}  
As $x\in\rL[\w\dar\xi]$, there is an 
$\in$-formula $\vpi(v,k,m)$, with ordinals as 
parameters, such that 
$a=\ens{\ang{k,m}\in\om\ti\om}
{\rL[\w]\mo\vpi(\w\dar\xi,k,m)}$,  
and  
\ben
\atc
\nenu
\itlb{6623}
if $X\in\cX$ then 
$X\fox\cX\,
\kaz k<\om\:\sus!\,m<\om\:\vpi(\pv\dar\xi,k,m)$.
\een
Let $\Psi(\w\dar\xi)$ be the conclusion of \ref{1022iii}   
after `then'. 
Assume towards the contrary that \ref{1022iii}  fails, 
so that there exists $X_0\in\cX\cap G_\w$,  
$X_0\fox\cX\neg\,\Psi(\pvr\xi)$. 
We may \noo\ assume by Theorem~\ref{631} that 
$\dym{X_0}=\xi$, \ie\ $X_0\in\cX\dar\xi$.

{\ubf Arguing in $\rL$}, define the sets 
$\cY_{km}=
\ens{X\in\cX\dar\xi}{X\fox\cX\vpi(\pv\res\xi,k,m)}$. 
By \ref{6623} and Theorem~\ref{631}, 
each set $\cY_k=\bigcup_m\cY_{km}$ is open dense in 
$\cX\dar\xi$. 
Therefore 
$\cY=\ens{X\in\cX\dar\xi}{\kaz k\,(X\sqd\bigcup\cY_k)}$ 
is dense as well by the Fusion property. 
It follows that there exists  
$X\in\cY$, $X\sq X_0$. 

Then for any $k<\om$ there is a finite $\cY'_k\sq\cY_k$ 
satisfying 
$X\sq\bigcup \cY'_k$, and if $Y\ne Z$ belong to $\cY'_k$ 
then $Y\cap Z=\pu$. 
Then for each $k$ we have a partition 
$\cY'_k=\bigcup_m\cY'_{km}$, where 
$\cY'_{km}=\cY_{km}\cap\cY'_k$. 
This enables us to define a continuous map $F_0:X\to\cN$ 
such that if $a \in X$ then $F_0(x)(k)=m$ iff 
$x\in\bigcup\cY'_{km}$. 
Let $F:\can\xi\to\cN$ be a continuous extension of $F_0$ 
to the whole $\can\xi$, still defined in $\rL$.
Then $X\fox\cX\Psi(\pvr\xi)$ by routine arguments, 
contrary to the choice of $X\sq X_0$. 
\epf

\parf{The case of the full forcing $\pei$}
\las{fff}

The next theorem shows that $\pei$ itself 
has the Fusion property. 

\bte
[in $\rL$]
\lam{pfp}
$\pei$ has the Fusion property. 
\ete
\bpf
Beginning the proof, we \noo\ assume that 
(*) each $\cY_k$ is open dense. \ie, if $Y\in\cY_{k}$, 
$Z\in\pei$, and $Z\psq Y$ then 
$Z\in\cY_{k}$ as well --- for if not then replace $\cY_k$ 
with $\cY'_k=\ens{Y'\in\pei}{\sus Y\in\cY_k(Y'\psq Y)}$. 

Fix some $X_0\in\pei$ and let $\et_0=\dym{X_0}$. 
Our plan is 
to define:
\ben
\nenu
\itlb{pfp1}
a sequence 
$\et_0\sq\xi_0\sq\xi_1\sq\xi_2\sq\dots$ 
of $\xi_k\in\cpo$, and $\xi=\bigcup_k\xi_k$;

\itlb{pfp2}
a $\xi$-admissible map $\phi:\om\onto\xi$, so that
\ben
\itlb{pfp2a}
if $\i\in\xi$ then the preimage 
$\phi\obr(\i)=\ens{k}{\phi(k)=\i}$ 
is infinite, 
\itlb{pfp2b} 
$\i\su\j=\phi(k)\in\xi$ implies $\i=\phi(\ell)$ 
for some $\ell<k$, 
\itlb{pfp2c} and in addition we require that 
$\phi(k)\in\xi_{k+1}$, $\kaz k$;
\een

\itlb{pfp3}
a system $\sis{X_s}{s\in\bse}$ of sets 
$X_s\in\pe{\xi_m}$ whenever $s\in 2^m$, such that 
$X_\La\sq X_0$, and 
$\sis{X_s}{s\in2^m}$ is a $\phi$-split system 
(Definition~\ref{splis}), $\kaz m$;

\itlb{pfp4}
if $s\in2^m$ and $e=0,1$ then a set 
$X_{s\we e}\sq X_s\uar{xi_{m+1}}$;

\itlb{pfp5}
finally, a set $X_s\in\cY_m$ for all  $s\in 2^{m}$.
\een

If this construction is accomplished then sets 
$Y_s=X_s\uar\xi\in\pe\xi$ form a $\phi$-fusion sequence 
by Lemma~\ref{uver}, so that 
$Y=\bigcap_m\bigcup_{s\in2^m}X_s\in\pe\xi$ by 
Theorem~\ref{fut}, and we obviously have $Y\psq X_0$, 
and $Y\sqf\bigcup\cY_m$, $\kaz m$. 

To maintain the construction, we pick any $X_\La\in\cY_0$, 
$X_\La\psq X_0$, by the density, 
let $\xi_0=\dym{X_\La}$, and let $\phi(0)$ 
be any 1-term tuple in $\xi_0$. 

Now the step $m\to m+1$, so that we assume that $\xi_m$, 
$\phi\res m$,  
and all sets 
$X_s\in\pe{\xi_m}, \,s\in2^m,$ 
are defined such that \ref{pfp1}--\ref{pfp5} hold 
wherever applicable.  

{\it Stage 1.} 
Pick any $s_0\in2^m.$ 
By the density, there is a set $Y\in\cY_{m+1}$, 
$Y\psq X_{s_0}$. 
Let $\et=\dym{Y}$; $\xi_m\sq\et$. 
Let $Y_s=X_s\uar{\et}$, so that $\sis{Y_s}{s\in2^m}$ 
is still a $\phi$-split system by Lemma~\ref{uver}, 
and $Y\sq Y_{s_0}$.
Let $Y'_s=Y_s\cap (Y\dar\et_\phi[s,s_0]\uar{\et})$
for all $s\in2^m.$ 
Then $\sis{Y'_s}{s\in2^m}$ 
is still a $\phi$-split system 
in $\pe{\et}$ by Lemma~\ref{suz}, $Y'_s\psq X_s$ 
for all $s\in2^m,$ 
and $Y'_{s_0}=Y'\in\cY_{m+1}$. 

{\it Stage 2.} 
Iterating Stage 1 (with all $s_0\in2^m$ involved one 
by one), we get a set $\za\in\cpo$ with $\xi_m\sq\za$ 
and a $\phi$-split system $\sis{Z_s}{s\in2^m}$ of sets 
$Z_s\in\pe\za,$ 
such that $Z_s\in\cY_{m+1}$ 
(here we refer to the open density assumption (*) above) 
and $Z_s\psq X_s$ for all $s\in2^m$. 
Let $\xi_{m+1}=\za$. 

{\it Stage 3.} 
We pick $\phi(m)\in\xi_{m+1}$ such that 
condition \ref{pfp2b} is preserved. 

{\it Stage 4.} 
By Lemma~\ref{pand}, there is a $\phi$-split system 
$\sis{X_u}{u\in2^{m+1}}$ is $\pe{\xi_{m+1}}$ expanding 
$\sis{Z_s}{s\in2^m}$, \ie\ $X_{s\we e}\sq Z_s$ for all 
$s\we e\in2^{m+1}$. 

As the sets $\xi_m$ obtained in the course of the 
construction are countable, we can maintain Stage 3 
at all inductive steps in such a way that 
condition \ref{pfp2a} holds. 
This ends the construction and the proof. 
\epf

\parf{Fusion property implies countable choice}
\las{fupac}

The two theorems below in this section are major 
applications of the Fusion property and 
Theorem~\ref{1022}. 
Recall that $\cN=\bn.$

\bte
\lam{65}
Assume that\/  $\cX\in \RF$ has the Fusion property, 
a set\/ $\Om\sq\cpo$, $\Om\in\rL$ is\/ 
$\cup$-closed\/ {\rm(under the finite $\cup$)}, 
%
$\ta_0\in\Om$ is\/ $\inv\Om$-isolated, and 
\ben
\fenu
\itlb{65*}
if\/ $\sis{\sg_k}{k<\om}\in\rL$ is a sequence of sets\/ 
$\sg_k\in\Om$, and\/ $\sg_k\cap\sg_\ell=\ta_0$ for all\/ 
$k\ne\ell$, then\/ $\bigcup_k\sg_k\in\Om$. 
\een
Let $\w\in\can\tup$ be\/ $\cX$-generic. 
Then\/ $\AC$ holds in\/ $\rL(\W\Om[\w])$ for all 
relations\/ $P\sq\om\ti\cN$ of class\/ 
$\od(\W\Om[\w],\w\dar\ta_0)$.  
Therefore,
\ben
\renu
\itlb{65i}
 if\/ $\ta_0=\pu$ 
then\/ $\xAC\od$ holds in\/ $\rL(\W\Om[\w])\,;$ 

\itlb{65ii}
if\/ \ref{65*} holds for all\/ $\ta_0\in\Om$  
then\/  $\AC$ holds in\/ $\rL(\W\Om[\w])$.
\een
\ete

\bpf
Fix a set $P\in\rL(\W\Om[\w])$, $P\sq\om\ti\cN,$ 
$\od(\W\Om[\w],\w\dar\ta_0)$ in $\rL(\W\Om[\w])$, 
with $\dom P=\om$. 
There is an $\in$-formula $\vpi(\cdot,\cdot,k,x)$ 
satisfying
\bce
$P=\ens{\ang{k,x}}
{\rL(\W\Om[\w])\mo\vpi(\W\Om[\w],\w\dar\ta_0,k,x)}$.
\ece 
As $\dom P=\om$, for any $k$ there is a real 
$x_k\in\cN\cap \rL(\W\Om[\w])$ with $\ang{k,x_k}\in P$, 
and then, by Corollary~\ref{641}, there is a 
set $\xi_k\in\Om$ such that $x_k\in\rL[\w\res\xi_k]$. 
Thus 
\ben
\nenu
\itlb{65*1}
$\rL(\W\Om[\w])\mo \sus x\in\rL[\w\dar\xi_k]\,
\vpi(\W\Om[\w],\w\dar\ta_0,k,x)$.
\een
Here the enumerations $k\mto x_k,\xi_k$ are maintained in 
$\rL[\w]$, not in $\rL(\W\Om[\w])$, of course. 
However Theorem~\ref{1022}\ref{1022iv} yields a map 
$H\in\rL$ such that $\dom H=\om$ and $\xi_k\in H(k)$ 
for all $k$. 
Let $\et_k=\bigcup(\Om\cap H(k))$; $\et_k\in\Om$ because 
$\Om$ is $\cup$-closed. 
Now $k\mto\et_k$ is a map in $\rL$, and $\xi_k\sq\et_k$, 
hence still $x_k\in\rL[\w\res\xi_k]$. 
We can assume that $\ta_0\sq\et_k$, $\kaz k$, of course. 
Now \ref{65*1} implies
\ben
\nenu
\atc
\itlb{65*2}
$\rL(\W\Om[\w])\mo \sus x\in\rL[\w\dar\et_k]\,
\vpi(\W\Om[\w],\w\dar\ta_0,k,x)$.
\een

Coming back to the theorem, assume to the contrary that 
\bce
$\rL(\W\Om[\w])\mo \neg \,\sus f\,\kaz k\,
\vpi(\W\Om[\w],\w\dar\ta_0,k,f(k))$. 
\ece
Putting it all together, we get 
a condition $X\in\cg\w$ which $\fox\cX$-forces 
this:
\ben
\Aenu
\itlb{65*a}
$\rL(\WW\Om)\mo\neg\, \sus f\,\kaz k\,
\vpi(\WW\Om,\pvr{\ta_0},k,f(k))$; \ and 

\itlb{65*b}
$\rL(\WW\Om)\mo \sus x\in\rL[\pvr{\et_k}]\,
\vpi(\WW\Om,\pvr{\ta_0},k,x)$, 
for each $k<\om$.
\een
We can assume that $\et=\bigcup_k\et_k\sq\modd X$. 
Then we get by Theorem~\ref{631}:  
\ben
\nenu
\atc
\atc
\itlb{65*3}
$X\dar\ta_0\fox\cX
\big(\rL(\WW\Om)\mo\neg\, \sus f\,\kaz k<\om\,
\vpi(\WW\Om,\pvr{\ta_0},k,f(k))\big)$; \ and 

\itlb{65*4}
$X\dar\et_k\fox\cX
\big(\rL(\WW\Om)\mo\sus x\in\rL[\pvr{\et_k}]\,
\vpi(\WW\Om,\pvr{\ta_0},k,x)\big)$, 
$\kaz k<\om$.
\een
This is because the formula $\big(...\big)$ in 
\ref{65*a} satisfies 
$\gga{\big(...\big)}=\Om$ and 
$\modd{\big(...\big)}=\ta_0\in\Om$, and similarly 
for \ref{65*b} with $\modd{\big(...\big)}=\et_k\in\Om$, 
and the isolation condition of the theorem is also used. 

Arguing in $\rL$ and using the 
$\inv\Om$-isolation of $\ta_0$, we get   
a sequence of permutations 
$\pi_k\in\inv\Om\cap\ppi{\ta_0}$ 
by induction, 
satisfying  $\sg_k\cap\sg_j=\ta_0$ whenever $k\ne j$,  
where $\sg_k={\pi_k}\akt\et_k\in\Om$. 
Let $Y_k=\pi_k\akt(X\dar\et_k)$.
Then
\ben
\nenu
\atc
\atc
\atc
\atc
\itlb{65*5}
$Y_k\fox\cX
\big(\rL(\WW\Om)\mo\sus x\in\rL[\pvr{\sg_k}]\,
\vpi(\WW\Om,\pvr{\ta_0},k,x)\big)$, 
$\kaz k<\om$.
\een
holds by \ref{65*4} by Corollary~\ref{624}. 
Note that $Y_k\in\cX\dar\sg_k$ by \ref{rfo6} in 
Section~\ref{rfo}. 

Note that $\sg=\bigcup_k\sg_k\in\Om$ by \ref{65*} 
of the theorem. 
The sets $Y_k$ satisfy $Y_k\dar\ta_0=X\dar\ta_0$, 
$\kaz k$, since $\pi_k\in\ppi{\ta_0}$.  
Thus $Y=\bigcap_k(Y_k\uar\sg)\in\cX\dar\sg$ 
by Lemma~\ref{rfoL} 
(\poo\ Lemma~\ref{pe4}). 
As obviously $Y\psq Y_k$, \ref{65*5} implies:

\vyk{
\ben
\nenu
\atc
\atc
\atc
\atc
\atc
\itlb{65*6}
}

\bce
$Y\fox\cX
\big(\rL(\WW\Om)\mo\kaz k\,\sus x\in\rL[\pvr{\sg}]\,
\vpi(\WW\Om,\pvr{\ta_0},k,x)\big)$, 
\ece
and hence (because any $Y$ forces that 
$\rL[\pvr{\sg}]$ is G\"odel-wellordered)
\ben
\nenu
\atc
\atc
\atc
\atc
\atc
\itlb{65*6}
$Y\fox\cX
\big(\rL(\WW\Om)\mo\sus f\,\kaz k<\om\,
\vpi(\WW\Om,\pvr{\ta_0},k,f(k))\big)$. 
\een
Now to accomplish the proof of the main claim 
of the theorem, we conclude 
that \ref{65*6} contradicts \ref{65*3} because 
$Y\dar\ta_0=X\dar\ta_0$ by construction. 

To prove \ref{65ii} of the theorem, note that every set 
$P\in\rL(\W\Om[\w])$, $P\sq\om\ti\cN$, 
belongs to $\od(\W\Om[\w],\w\dar\ta_0)$ in 
$\rL(\W\Om[\w])$ for a suitable 
$\ta_0\in\Om$.
\epf
 
\vyk{

, and the fact 
that $\Om$ is closed under finite unions, enable us 
to \noo\ assume that 
the map $k\mto\et_k$ belongs to $\rL$.

Consider the $\cL$-formula 
\bce
$\chi(k,y):=
\sus x\in\cN\cap\rL[y]\,\big(\rL(\WW\Om)\mo\vpi(k,x)\big)$; 
\ \ then:
\ece
\ben
\nenu
\itlb{651}
$\rL[\w]\mo \kaz k\,\chi(k,\WW\Om)[\w]$, 

\itlb{652}
$\gga(\chi)=\gga(\vpi)=\gga(\Om)$ and 
$\modd\chi=\modd\vpi=\ta_0$.
\een 
Moreover, we can assume that
\ben
\nenu
\atc\atc
\itlb{653}
If $X\in\cX$ then $X\fox\cX \kaz k\,\chi(k,\WW\Om)$, 
\een
for if not then replace $\vpi(k,x)$ 
from the beginning  by the formula 
\bce
$
\vpi'(k,x):= \ k<\om\land x\in\cN\land\, 
\big(\neg\,\kaz j\,\chi(j,\WW\Om)\imp\vpi(k,x)\big)
$,
\ece
which defines $P$ as well, and leads to 
\ref{651},\ref{652},\ref{653} after the according 
redefinition of $\chi$. 
Now we let
\pagebreak[0]
$$
\cY_k=\ens{X\in\cX}
{\sus \et\in\Om\,\big(\ta_0\sq\et\sq\dym X\land 
X\fox\cX\chi(k,\pvr \et)\big)}
$$
for all $k$. 
The sequence $\sis{\cY_k}{k<\om}$ belongs to $\rL$, and 
in addition
\ben
\nenu
\atc\atc
\atc
\itlb{655}
if $X\in\cX$ and $X\sqf\bigcup\cY_k$ then $X\in\cY_k$ 
--- by Lemma~\ref{553} and because $\Om$ is $\cup$-closed,  

\itlb{656}
each $\cY_k$ is open dense in $\cX$ --- 
by \ref{653} and Corollary~\ref{642},  

\itlb{657}
$\cY=\bigcap_k\cY_k$ is open dense in $\cX$ --- 
by \ref{656} and the Fusion property.  
\een

Let 
$\cZ=\ens{Z\in\cX}{\sus\da\in\Om\,\kaz k\,\big(
\ta_0\sq\da\sq\dym Z\land Z\fox\cX\chi(k,\pvr \da)
\big)}$.

\ble
\lam{658}
$\cZ$  
is an open dense set in\/ $\cX$.
\ele
\bpf[Lemma]
{\ubf Arguing in $\rL$}, let $X_0\in\cX$. 
By \ref{657}, there exists $X\in\cY$, $X\psq X_0$. 
Let $\ta=\dym X$. 
There exists (in $\rL$) a sequence of sets $\et_k\in\Om$ 
such that $\ta_0\sq\et_k\sq\ta$ and 
(\mdag) $X\fox\cX\chi(k,\pvr {\et_k})$. 

We claim that there is a sequence of permutations 
$\pi_k\in
\gga(\Om)\cap\ppi{\ta_0}$ 
satisfying $\pi_0\akt X=X$ and 
$({\pi_k}\akt\ta)\cap({\pi_j}\akt\ta)=\ta_0$ whenever 
$k\ne j$. 
Indeed, 
to define $\pi_k$ by induction, let $\pi_0$ be the 
identity. 
If $\pi_j$, $j<k$, are defined, then 
$\ta'=\ta\cup\bigcup_{j<k}({\pi_j}\akt\ta)\in\Om$, 
and hence, as $\ta_0$ is isolated, there exists 
$\pi=\pi_k\in 
\gga(\Om)\cap\ppi{\ta_0}$ 
with $({\pi_k}\akt{\ta'})\cap\ta'=\ta_0$. 

Put $Z_k=\pi_k\akt X$ and $\da_k=\pi_k\akt\ta$, 
so that $\da_k\in\Om$, $Z_k\in\cX\dar\da_k$, 
$Z_k\dar{\ta_0}=X\dar\ta_0$ for all $k$, hence 
$Z_k\dar\ta_0=Z_j\dar\ta_0$ for all $j,k<\om$. 
In addition ${k\ne j}\imp \da_k\cap\da_j=\ta_0$. 
Thus $\da=\bigcup_k\da_k\in\Om$ by \ref{65*} 
of the theorem. 
We conclude that 
$Z=\ens{x\in\can\da}{\kaz k\,(x\dar\da_k\in Z_k)}\in\pe\da$ 
by Lemma \ref{pe4}, 
and $Z\psq Z_k$, $\kaz k$, in particular, $Z\psq X\psq X_0$. 
Moreover $Z\in\cX\dar\da$ by \ref{rfo5}, Section~\ref{rfo}.

Note that $\za_k={\pi_k}\akt{\et_k}$ belongs to 
$\Om$ by the choice of $\pi_k$. 
Moreover (\mdag) implies  
$Z_k\fox\cX\chi(k,\pvr {\za_k})$ by 
Theorem \ref{621}, therefore 
$Z\fox\cX\chi(k,\pvr {\za_k})$, $\kaz k$, hence  
$Z\fox\cX\chi(k,\pvr {\da})$ as
$\za_k\sq\da_k\sq\da$. 
Thus $X\in\cZ$.
\epF{Lemma}

To accomplish the proof of the theorem, we conclude 
from the lemma that there exists $Z\in\cZ\cap G_{\w}$. 
This means, that for some $\da\in\Om$, 
\pagebreak[0]
\bce
$Z\fox\cX 
\kaz k\sus x\in\cN\cap\rL[\w\dar\da]\,P(k,x)$.
\ece
This allows us to define a choice function required 
for $P$ in $\rL(\W\Om[\w])$, because $\rL[\w\dar\da]$ 
is well-orderable.
}

\vyk{
\bcor
\lam{65c}
Under the assumptions of Theorem~\ref{65}, 
suppose that\/ \ref{65*} of the theorem 
holds for all\/ $\ta_0\in\Om$. 
Then\/  $\AC$ holds in\/ $\rL(\W\Om[\w])$.
\ecor
\bpf
Every set $P\in\rL(\W\Om[\w])$, $P\sq\om\ti\cN$, 
belongs to $\od(\W\Om[\w],\w\dar\ta_0)$ in 
$\rL(\W\Om[\w])$ for a suitable 
$\ta_0\in\Om$.
\epf
}


A somewhat simpler set of properties of $\Om$ 
leads to $\DC$ 
in classes of the form $\rL(\W\Om[\w])$, as the next 
theorem shows.

\bte
\lam{66}
Assume that\/  $\cX\in \RF$ has the Fusion property, 
a set\/ $\Om\sq\cpo$, $\Om\in\rL$ is\/ 
closed in\/ $\rL$ under countable unions, and\/ 
$\w\in\can\tup$ is\/ $\cX$-generic. 
Then\/ $\DC$ holds in\/ $\rL(\W\Om[\w])$.
\ete

\bpf
Let $\Da=\bigcup\Om$; then $\Da\in\rL$, $\Da\sq\tup$,  
$\rL(\W\Om[\w])\sq \rL[\w\dar\Da]$. 
We claim that 
$\rL(\W\Om[\w])\cap\cN= \rL[\w\dar\Da]\cap\cN$; 
this proves the theorem because the full {\bf AC} 
holds in $\rL[\w\dar\Da]$. 
In the nontrivial direction,
let $x\in\rL[\w\dar\Da]\cap\cN$. 
It follows by Theorem~\ref{1022}\ref{1022ii} 
that there is a 
($\rL$-countable!) 
$\xi\in\cpo$, $\xi\sq\Da$ 
with $x\in\rL[\w\dar\xi]$. 
Then $\xi\in\Om$, as $\Om$ is\/ 
closed in\/ $\rL$ under countable unions. 
Thus $x\in\rL(\W\Om[\w])$. 
\epf

\sekt{Choiceless generic subextensions}
\las{subex}

Thus Chapter defines and studies those generic models, 
of the form $\rL(\W\Om[\w])$, 
which will be used in the proof of Theorem~\ref{mt1}. 
The forcing notion $\cX$ is not yet defined, so our 
goal here will be to introduce some key properties 
of $\cX$ and $\cX$-generic arrays  
(the Definability, Structure, and Even Extension  
properties defined below) 
that will eventually lead to Theorem~\ref{mt1}. 

In Section~\ref{kso}, we define, in $\rL$, 
{\ubf four sets $\Om_e\sq\cpo$}, 
$e=1,2,3,4$, related to the models we'll use in the proof of 
the according {\ubf items of Theorem~\ref{mt1}}. 
We also define according {\ubf subgroups $\Ga_e\sq\per$} 
and add 
some auxiliary sets $\Om\sq\cpo$ and  according groups. 
Theorem~\ref{79} in Section~\ref{76} proves some 
{\ubf combinatorial properties} of these sets and subgroups, 
rather known in theory of symmetric generic extensions. 

The positive and negative Choice statements in items 
\ref{mt11},\ref{mt12},\ref{mt13},\ref{mt14} 
of Theorem~\ref{mt1} naturally split into 
the {\ubf three groups} defined in Section~\ref{g1}. 
Theorem~\ref{811} proves the statements of 
the {\ubf first group} in the according models  
$\rL(\W{\Ome}[\w])$, $e=2,3,4$, 
provided $\cX\in\RF$ 
has the Fusion property and $\w\in\can\tup$ is 
\dd\cX generic. 

The {\ubf second group} of Choice statements contains 
the {\rit negative} statements in items 
\ref{mt11},\ref{mt12},\ref{mt13},\ref{mt14} 
of Theorem~\ref{mt1}. 
It needs a different treatment. 
For that purpose, we introduce the {\ubf Structure} and 
{\ubf Definability} properties of a forcing $\cX$ 
in Section~\ref{82seq}, and   
derive (Theorem~\ref{82} in Section \ref{82seq,})
that they imply the negative Choice statements  
in the according models  
$\rL(\W{\Ome}[\w])$, $e=1,2,3,4$. 

Finally to provide the {\ubf third group} of statements, 
namely $\xDC{\fp1\nn}$ and $\xDC{\ip1{\nn+1}}$, 
to be true in the according models $\rL(\W{\Ome}[\w])$, 
we introduce the {\ubf Odd-Expansion} property in 
Section~\ref{91} and achieve the result required 
by Theorem~\ref{93}. 

Theorem~\ref{mt1b} summarizes the content of this chapter.

\parf{Key sets $\Ome$ and permutation groups $\Gae$}
\las{kso}

Classes of the form $\rL(\W\Om[\w])$ will serve as 
models for different parts of our main theorem. 
Here $\w\in\can\tup$ will be $\cX$-generic over $\rL$ 
for a special forcing $\cX\in\RF\cap\rL$, whereas 
$\Om\in\rL$ will be selected as special subsets of $\cpo$.

First of all, we are going to define sets 
$\Oma,\Omb,\Omc,\Omd\sq\cpo$ is $\rL$. 
This involves the notion of \rit{even} and \rit{odd} 
tuples in $\tup$ as defined in Section~\ref{perm}. 

\bdf
[in $\rL$]
\lam{71}
If $\i\sq\j$ belong to $\tup$ then $\j$ is an \rit{odd expansion} 
of $\i$, in symbol $\i\osq\j$, iff $\j(k)$ is an odd ordinal for 
all $\lh\i\le k<\lh\j$.

If $\xi,\et\in\cpo$ then $\xi$ is an \rit{odd expansion} 
of $\et$, in symbol $\et\osq\xi$, iff 
\kmar{osq}%
\index{odd expansion, $\osq$}%
\index{expansion!odd, $\osq$}%
\index{zzzosq@$\osq$}%
$\et\sq\xi$ and in addition 
all tuples $\i\in\xi\bez\et$ are odd. 
Put:\vim 
$$
\bay{rcl}
\sli\xi\al 
&\hspace*{0ex}=\hspace*{0ex}& 
\kmar{sli xi al}%
\index{slice, $\sli\xi\al$}%
\index{zzxial@$\sli\xi\al$}%
\ens{\i\in\xi}{\i(0)=\al}, \,
\text{ for any $\al<\omi$, $\xi\sq\tup$ --- 
the $\al$-\rit{slice} of $\xi$},\\[0ex]
&&\text{\,in particular } 
\kmar{axi al}%
\axi\al=\ens{\i\in\tup}{\i(0)=\al}\,;\\[1ex]
\index{slice, $\axi\al$}%
\index{zIal@$\axi\al$}%
%
\Oma 
\kmar{Oma}%
\index{key set!$\Oma$}%
\index{zzOm1@$\Oma$, key set}%
&\hspace*{-0.7ex}=\hspace*{-0.7ex}& 
\ens{\ta\in\cpo}{\sus m\,\kaz\i\in\ta\,
(\i\text{ is \ubf even }\imp \lh \i\le m)};\\[1ex]
\vyk{
\Omb &=& 
\ens{\ta\in\cpo}{\kaz\al<\omi
(\al\text{ limit }\imp \ta\cap\axi\al\in\Oma};\\[1ex]
}
\Omc 
\kmar{Omc}%
\index{key set!$\Omc$}%
\index{zzOm3@$\Omc$, key set}%
&\hspace*{-0.7ex}=\hspace*{-0.7ex}& 
\text{all $\ta\in\cpo$ which contain no infinite 
paths $\i_0\su \i_1\su \i_2\su\ldots$}\\
&&\text{of {\bf even} tuples $\i_k\in\ta$};
\eay
$$
and let 
$\Gaa=\Gac= 
\kmar{Gaa\mns Gac}%
\index{permutation groups!$\Gaa$}%
\index{zzGa1@$\Gaa$, permutation group}%
\index{permutation groups!$\Gac$}%
\index{zzGa3@$\Gac$, permutation group}%
\per$, 
\ all parity-preserving and $\su$-preserving 
$\pi:\tup\onto\tup$.
%
\vyk{
\Gab=\Gad\!&=&\!
\ens{\pi\in\qer}{\kaz\i,\j\,
\big(\pi(\i)=\j\imp
(\i(0)\text{ limit\,}\eqv\j(0)\text{ limit})\big)}.
}%
\edf

It takes more work to define $\Omb$ and $\Omd$. 
First of all, if $\al,\ba<\omi$ then define a 
\rit{shift permutation} $\pi_{\al\ba}\in\pero$ 
\index{shift permutation, $\pi_{\al\ba}\in\per$}%
\index{zzpiab@$\pi_{\al\ba}$, shift permutation}%
such that if $\i\in\tup$ then $\j=\pi_{\al\ba}(\i)$ 
satisfies $\lh\j=\lh\i$ and the following:\vom 

$-$ if $\i(0)\nin\ans{\al,\ba}$ then $\j=\i$;\vom

$-$ if $\i(0)=\al$ then $\j(0)=\ba$ and $\j(k)=\i(k)$ 
for all $0<k<\lh\i$;\vom

$-$ if $\i(0)=\ba$ then $\j(0)=\al$ and $\j(k)=\i(k)$ 
for all $0<k<\lh\i$.\vom

\noi
Note that $\pi_{\al\ba}\in\pero$, and even 
$\pi_{\al\ba}\in\per$ in case $\al,\ba$ 
have equal parity.

A routine proof of the next lemma is left to
the reader.

\ble
[in $\rL$]
\lam{72}
There is a sequence\/ 
$\sis{\aza\al}{\al<\omi\text{\,\rm succesor}}$ 
\kmar{\mns aza al}%
\index{zzyaal@$\aza\al$}%
such that$:$
\ben
\renu
\itlb{72i}
if\/ $\al<\omi$ is a successor ordinal then\/ 
$\aza\al\in\cpo$ and\/ 
$\aza\al\sq\axi\al\,;$ 

\itlb{72ii}
if $\al,\la<\omi$, $\et\in\cpo$, $\et\sq\axi\al$, 
then there is 
a~successor\/ 
$\ba>\la$ such that\/ 
${\pi_{\al\ba}}\akt\et=\aza\ba$ and the ordinals 
$\al,\ba$ 
have the same parity$.$ 
\vyk{
\itlb{72iii}
{\gol if $\la<\omi$ is limit then\/ $\aza{\la+2}=\pu$ \ 
{\rm(a technical condition)}. ????} 
}%
\qed
\een
\ele

We fix such a sequence of sets $\aza\al$ in $\rL$. 


\bdf
[in $\rL$]
\lam{73}
Put $\Omb=$, resp., $\Omd=$ all $\ta\in\cpo$ such that: 
\kmar{$\,$\mns Omb\mns Omd}%
\index{key set!$\Omb$}%
\index{zzOm2@$\Omb$, key set}%
\index{key set!$\Omd$}%
\index{zzOm4@$\Omd$, key set}%
\ben
\nenu
\itlb{73ii}
if $\al<\omi$ is a successor and $\sli\ta\al\ne\pu$ 
then $\aza\al\osq\sli\ta\al$; 

\itlb{73i}
if $\al<\omi$ is limit then $\sli\ta\al\in\Oma$, 
resp., $\sli\ta\al\in\Omc$.
\een
\vyk{
Put $\Omd=$ all $\ta\in\cpo$ satisfying \ref{73ii} and 
the following \ref{73iii} instead of \ref{73i}:
\ben
\nenu
\atc\atc
\itlb{73iii}
if $\al<\omi$ is limit then $\axi\al\cap\ta\in\Omc$.
\een
}%
In addition, put $\Gab=\Gad=$ all $\pi\in\qer$ such that 
\kmar{Gab\mns Gad}%
\index{permutation groups!$\Gab$}%
\index{zzGa2@$\Gab$, permutation group}%
\index{permutation groups!$\Gad$}%
\index{zzGa4@$\Gad$, permutation group}%
\ben
\nenu
\atc\atc 
\itlb{73iv}
if $\pi(\i)=\j$ and $\i(0)$ is limit then so is 
$\j(0)$, and 

\itlb{73v}
if $\pi(\i)=\j$ and $\al=\i(0)$ is a successor then 
$\ba=\j(0)$ is a successor either, 
and $\pi\akt{\aza\al}=\aza\ba$.
\qed
\een
\eDf
To conclude, sets $\Ome=\Oma,\Omb,\Omc,\Omd\sq\cpo$ 
and associated groups $\Gae\sq\per$ have been defined 
in $\rL$, 
mainly via conditions related to {\ubf even} tuples 
$\i\in\xi\in\Ome$, while giving {\ubf odd} tuples much 
more freedom. 
Speaking about this {\ubf distinction} between 
{\ubf  even and odd}  tuples 
in the definition of the sets $\Ome$ and their treatment,  
one may ask whether a parity-independent modifications of 
the definitions above may also work towards the proof 
of Theorem~\ref{mt1}. 
We'll explain in Section~\ref{why} that the answer is 
in the {\ubf negative}. 

Some {\ubf related sets} $\Om\sq\cpo$ 
will also be considered. 

\bdf
[in $\rL$]
\lam{74}
Let $\vt\in\cpo$. 
We first put   
\kmar{$\,$\mns Oha vt\mns Ohc vt}%
\index{key set!$\Oha\vt$}%
\index{zzOm1vt@$\Oha\vt$, key set}%
\index{key set!$\Ohc\vt$}%
\index{zzOm3vt@$\Ohc\vt$, key set}%
\bce
$\Oha\vt=\Ohc\vt=
\ens{\ta\in\cpo}{\vt\osq\ta}$ 
\kmar{Gha vt\mns Ghc vt}%
\index{permutation groups!$\Gha\vt$}%
\index{zzGa1vt@$\Gha\vt$, permutation group}%
\index{permutation groups!$\Ghc\vt$}%
\index{zzGa3vt@$\Ghc\vt$, permutation group}%
\ece
and $\Gha\vt=\Ghc\vt=\ppi\vt$ 
(all $\pi\in\per$ equal to the identity on $\vt$). 

To handle the $\ans{2,4}$-case, 
we let $\hza\al\vt=\aza\al$ if $\al<\omi$ is 
\kmar{hza al vt}%
\index{zzyaalvt@$\hza\al\vt$}%
a successor ordinal, 
and $\hza\al\vt=\sli\vt\al$ if $\al$ is limit. 
Now we define: 
$$
\bay{rcl}
\Ohb\vt=\Ohd\vt
\kmar{Ohb vt\mns Ohdvt}%
\index{key set!$\Ohb\vt$}%
\index{zzOm2vt@$\Ohb\vt$, key set}%
\index{key set!$\Ohd\vt$}%
\index{zzOm4vt@$\Ohd\vt$, key set}%
&=&
\ens{\ta\in\cpo}{\kaz\al<\omi\,\big(
\sli\ta\al\ne\pu\imp 
\hza\al\vt\osq{\sli\ta\al}\big)}
,\\[1.5ex]
%
\Ghb\vt=\Ghd\vt
\kmar{Ghb vt\mns Ghd vt}%
\index{permutation groups!$\Ghb\vt$}%
\index{zzGa2vt@$\Ghb\vt$, permutation group}%
\index{permutation groups!$\Ghd\vt$}%
\index{zzGa4vt@$\Ghd\vt$, permutation group}%
&=&
\ens{\pi\in\qer}
{\kaz\al,\ba\,\big(\pi(\ang\al)=\ang\ba \imp 
\hza\ba\vt=\pi\akt {\hza\al\vt}\big)}.
\eay
$$
Put 
$
\bay[t]{rcl}
\Oms
\kmar{Oms}%
\index{key set!$\Oms$}%
\index{zzOm*@$\Oms$, key set}%
&\!=\!&
\ens{\ta\in\Omb \text{ (equivalently, }\Omd)}
{\kaz\al\,(\al\text{ is limit}\imp\sli\ta\al=\pu)},\\[0.5ex]
\Gas
\kmar{Gas}%
\index{permutation groups!$\Gas$}%
\index{zzGa*@$\Gas$, permutation group}%
&\!=\!&
\Gab=\Gad.
\hspace*{35ex}\qed
\eay
$ 
\eDf   

\vyk{ 
$$
\hspace*{-3ex}
\bay{rcl}
\Oms
&\hspace*{-0.8ex}=\hspace*{-0.8ex}&
\ens{\ta\in\cpo}
{\kaz\i\in\ta\,\big(\al=\i(0)\text{ is a successor\,}\land 
\text{($\i$ is odd $\imp\i\in\aza\al)$}},
\hspace*{0ex}\\[1ex]
\Gas
&\hspace*{-0.8ex}=\hspace*{-0.8ex}&
\ens{\pi\in\Gab=\Gad}
{\kaz\i\in\tup\,\big(\i(0)\text{ is limit }\imp
\pi(\i)=\i\big)}. \ \ \ {\kra{\text{\bf ????}}}
\eay
\eqno
\hspace*{-6ex}
\bay{r}
\\[1ex]
\hspace*{0ex}\qed\hspace*{-1ex}
\eay
$$
\eDf
}

\parf{Invariance, isolation and other results}
\las{76}

Recall Definition~\ref{63d} on isolation. 
The next theorem contains a summary of rather simple 
properties of the sets $\Ome\sq\cpo$ and the groups 
$\Gae\sq\per$.

\bte
[in $\rL$, summary]
\lam{79}
Let\/ $e=1,2,3,4$ and\/ $\vt\in\Ome$. 
Then
\ben
\renu
\itlb{73r1}\msur
$\Oma\sq\Omc$, \ $\Omb\sq\Omd$, \ 
$\bigcup\Oma=\bigcup\Omc=\tup$, whereas\\[0.3ex] 
$\bigcup\Omb=\bigcup\Omd=
\ens{\i\in\tup} 
{\al=\i(0)\text{\rm\ is a successor}\imp 
\aza\al\osq \aza\al\cup\ilq\i};$

\itlb{73r3}
if\/ $e=1,3$ and\/ $\ta\in\cpo,$  $\ta\sq\et\in\Ome,$  
then\/ $\ta\in\Ome$ \ 
{\rm(false for $e=2,4$);}

\itlb{795}
$\vt\in\Ohe\vt\sq\Ome$,  
and if\/  $\vt\in\Omd$ then\/
$\Oms\sq\Ohb\vt=\Ohd\vt\sq\Omb\sq\Omd\,;$ 

\itlb{911}
if\/ $\xi$ and\/ $\et\osq\ta$ belong to\/ $\cpo$,  
then\/ 
$\et\in\Ome\imp \ta\in\Ome$, and\/ 
$\et\in\Ohe\xi\imp \ta\in\Ohe\xi$ --- 
{\rm take notice of this claim, it will be 
very important!}$;$ 

\itlb{794}
the sets\/ $\Ome$ are closed under 
finite unions, whereas\/ $\Oms,\,\Ohe\vt$ are closed 
under countable unions\/ {\rm(obvious)}$;$

\itlb{73r4}
$\Ome$ is\/ $\Gae$-invariant, 
$\Ohe\vt$ is $\Ghe\vt$-invariant,  
$\Oms$ is $\Gas$-invariant. 

\itlb{781} 
$\Omb$ satisfies\/ \ref{65*} of Theorem~\ref{65} 
in case\/ $\ta_0=\pu;$

\itlb{782} 
the sets\/ $\Omc,\Omd$ satisfy\/ \ref{65*} of Thm~\ref{65} 
for all\/ $\ta_0\in\Omc,$ resp.\ $\ta_0\in\Omd\,;$ 

\itlb{77}
if\/ $\xi\in\Oms$, and $\ta\in\cpo$ satisfies \ref{73ii} 
of Definition~\ref{73}, 
then there is a~permutation\/ $\pi\in\ppi{\xi}$ 
such that\/ $\pi\akt\ta\in\Oms\,;$ 

\itlb{77pu}
if\/ $e=2,4$, $\xi\in\Ome$, $\ta\in\Ohe\xi$, 
then there is a~permutation\/ $\pi\in\Ghe\xi$ 
such that\/ $\sg=\pi\akt\ta\in\Oms$ and\/ 
$\sg\cap\ta=\pu\,;$  

\itlb{791}
$\Gae\sq\inv{\Ome}$, 
each\/ $\ta\in\Ome$ is\/ $\Gae$-isolated\/$;$

\itlb{792}
$\Ghe\vt\sq\inv{\Ohe\vt}$, 
each\/ $\ta\in\Ohe\vt$ is\/ $\Ghe\vt$-isolated\/$;$

\itlb{793}
$\Gas\sq 
\Ghb\vt=\Ghd\vt$, 
each\/ $\ta\in\Oms$ is\/ $\Gas$-isolated\/$.$
\qed
\een
\ete

\bpf[in $\rL$]
Claims \ref{73r1}, 
\ref{73r3}, 
\ref{795}, 
\ref{911}, 
\ref{794},
\ref{73r4} 
are pretty routine. 

\ref{781} 
Assume that sets $\sg_k\in\Omb$ are pairwise disjoint. 
Then $\sg=\bigcup_k\sg_k\in\cpo$. 
Let $\al<\omi$ be limit. 
Then $\sli{\sg}\al=\sli{\sg_k}\al$  for some $k$ 
by the disjointness condition. 
Thus $\sli{\sg}\al\in\Oma$, as required.

\ref{782} 
Assume that $\ta_0\in\Omc$ and sets $\sg_k\in\Omb$ 
satisfy (*) $\sg_k\cap\sg_\ell=\ta_0$ for all $k\ne\ell$. 
Then any $\su$-increasing sequence in $\sg=\bigcup_k\sg_k$
entirely belongs to one of $\sg_k$, hence it cannot be 
infinite.

\ref{77}
We can \noo\ assume that $\xi\sq\ta$ 
(otherwise replace $\ta$ by $\xi\cup\ta$). 
Let  
$T=\ens{\i(0)}{\i\in\ta}$ and $\mu=\sup T$. 
If $\al\in T_0=\ens{\al'\in T}{\al'\text{ is limit}}$ 
then by \ref{72ii} of Lemma~\ref{72} there is 
a countable successor ordinal $\ba(\al)>\mu$, 
of the same parity as $\al$, such that 
${\pi_{\al{,\ba(\al)}}}\akt{\sli{\ta}\al}
=\aza{\ba(\al)}$. 
We can choose these ordinals $\ba(\al)$ so that 
$\al\ne\al'\imp \ba(\al)\ne\ba(\al')$ 
for all $\al\in T_0$. 
This allows to define $\pi\in\qer$ 
as follows:
\busq{eq11}{
\pi(\i)=
\left\{
\bay{ccl}
\i &, \text{ in case }& 
\i(0)\nin T_0\cup\ens{\ba(\al)}{\al\in T_0}\,;\\[0.8ex]
\pi_{\al{,\ba(\al)}}(\i) &, \text{ in case }& 
\i(0)\in T_0\cup\ens{\ba(\al)}{\al\in T_0}\,. 
\eay
\right.
} 
Note that $\pi\in\ppi\xi$: 
if $\i\in\xi$ then $\i(0)$ is a successor 
because $\xi\in\Oms$, and hence 
$\i(0)\nin T_0\cup\ens{\ba(\al)}{\al\in T_0}$ 
by construction, and $\pi(\i)=\i$. 

It remains to check that 
$\sg=\pi\akt\ta\in\Oms$. 
Let $\ba<\omi$ and $\sli{\sg}\ba\ne\pu$. 

{\it Case 1\/}: $\ba=\ba(\al)$ for some $\al\in T_0$. 
Then 
$\sli{\sg}\ba=
{\pi_{\al\ba}}\akt{\sli{\ta}\al}=\aza\ba=\hza\ba\pu$ 
by construction.

{\it Case 2\/}: $\ba\in T\bez T_0$, hence $\ba$ is 
a successor. 
Then $\sli{\sg}\ba=\sli{\ta}\ba$ by construction. 
Therefore $\hza\ba{}=\hza\ba\pu\esq\sli{\sg}\ba$, 
as $\ta\in\Oms=\Ohb\pu$. 

Combining the results in two cases, we get $\sg\in\Oms$. 

\ref{77pu} 
The proof is rather similar. 
Assuming that $\xi\sq\ta$ as above, we pick,  
for each $\al\in T=\ens{\i(0)}{\i\in\ta}$, 
a successor ordinal $\ba(\al)>\mu=\sup T$, 
of the same parity as $\al$, such that 
${\pi_{\al{,\ba(\al)}}}\akt{\sli{\ta}\al}
=\aza{\ba(\al)}$. 
Choose $\ba(\al)$ so that 
$\al<\al'\imp \ba(\al)<\ba(\al')$. 
Define $\pi\in\per$ as follows:
\busq{eq22}{
\pi(\i)=
\left\{
\bay{ccl}
\i &, \text{ in case }& 
\i(0)\nin T\cup\ens{\ba(\al)}{\al\in T}\,;\\[0.8ex]
\pi_{\al{,\ba(\al)}}(\i) &, \text{ in case }& 
\i(0)\in T\cup\ens{\ba(\al)}{\al\in T}\,. 
\eay
\right.
}

\ref{791}
To prove the isolation claim, 
let $\la<\omi$ be a limit ordinal $>$ all 
ordinals $\j(k)$, where $\j\in\xi$ and $k<\lh\j$. 
To handle {\ubf the case $e=1,3$},  
recall that each $\et\in\cpo$ is $\qer$-isolated 
by Lemma~\ref{iso1}. 
\vyk{
(Recall that $\Gaa=\qer$.) 
Let $\et\sq\xi\in\cpo$; 
define  
$\pi\in \ppi\xi\cap\Gaa$ satisfying 
$\xi\cap(\pi\akt\xi)=\et$. 

If $\i\in\et$ then put $\pi(\i)=\i$. 

Now suppose that $\i\nin\et$. 
Then there is $m_\i<\lh\i$ such that still 
$\i\res m_\i\in\et$ but $(\i\res m_\i)\we \i(m_\i)\nin\et$; 
possibly $m_i=0$. 
Consider the sets $H=\ens{\j(m_\i)}{\j\in\et\land\lh\j>m_\i}$ 
and $D=\la\bez H$. 
Define $\j=\pi(\i)\in\tup$ so that  
$\lh\j=\lh\i$, $\j\res m_\i=\i\res m_\i$, 
$\j(l)=\i(l)$ in case $m_\i<l<\lh\j$, and 
\busq{eq33}{
\j(m_\i)=
\left\{
\bay{rcl}
\la+\al,&\text{ in case }&
\al=\i(m_\i)\in D;\\[0.8ex] 
\al,&\text{ in case }&
\i(m_\i)=\la+\al, \; \al\in D;\\[0.8ex]
\i(m_\i),&\text{ in case }&
\i(m_\i)\nin D\cup\ens{\la+\al}{\al\in D}.
\eay
\right.
}%
}%

To handle {\ubf the case $e=2,4$},  
prove that each $\et\in\cpo$, satisfying 
\ref{73ii} of Definition \ref{73}, is $\Gab$-isolated. 
Let $\et\sq\xi\in\cpo$; let's define  
$\pi\in \ppi\xi\cap\Gab$ satisfying 
$\xi\cap(\pi\akt\xi)=\et$. 
Splitting $\tup$ into the limit and successot parts 
\bce
$\tup_0=\ens{\i\in\tup}{\i(0)\text{ is limit}}$ \ 
and \  
$\tup_1=\ens{\i\in\tup}{\i(0)\text{ is a successor}}$, 
\ece
we accordingly put  
$\et_{e}=\et\cap\tup_{e}\sq \xi_{e}=\xi\cap\tup_{e}$,  
$e=0,1$,  
define permutations $\pi_{e}$ of the domains 
$\tup_{e}$ separately, and put $\pi=\pi_0\cup\pi_1$ 
at the end. 

{\it Part 1\/}. 
We leave it to the reader to define 
$\pi_0:\tup_0\onto\tup_0$ with 
$\pi_0\res\et_0=$ the identity and 
$\xi_0\cap({\pi_0}\akt{\xi_0})=\et_0$, 
following 
the proof of Lemma~\ref{iso1}.

{\it Part 2\/}.  
We now concentrate on the construction of 
$\pi_1:\tup_1\onto\tup_1$. 

If $\i\in\et_1$ then put $\pi_1(\i)=\i$. 
{\ubf Now let $\i\in\tup_1\bez\et_1$}. 
Consider the sets
$$
A_1=\ens{\j(0)}{\j\in\et_1}
\;\sq\;
B_1=\ens{\j(0)}{\j\in\xi_1}
\;\sq\;
\ens{\al<\omi}{\al\text{ successor}}. 
$$
Following the proof of \ref{77} above, 
if $\al<\omi$ is a successor  
then by \ref{72ii} of Lemma~\ref{72} there is 
a successor $\ba(\al)>\la$, 
of the same parity as $\al$, such that 
${\pi_{\al{,\ba(\al)}}}\akt{\aza\al}
=\aza{\ba(\al)}$. 
We can choose these ordinals $\ba(\al)$ so that 
$\al<\al'\imp \ba(\al)<\ba(\al')$. 
Now, if $\i\in \tup_1$ but $\i(0)\nin A_1$ 
then put
%
\busq{eq44}{
\pi(\i)=
\left\{
\bay{ccl}
\i \!\!\!\!&, \text{ if}& 
\i(0)\nin B_1\cup\ens{\ba(\al)}{\al\in (B_1\bez A_1)}
\,;\\[0.8ex]
\pi_{\al{,\ba(\al)}}(\i) \!\!\!\!&, \text{ if}& 
\i(0)\in (B_1\bez A_1)\cup\ens{\ba(\al)}{\al\in (B_1\bez A_1)}; 
\eay
\right.
} 
following the idea of \eqref{eq11}, \eqref{eq22} above.

{\it Part 3\/}.  
We finally define $\pi_1$ on the domain 
$\tup'_1=\ens{\i\in\tup_1}{\i(0)\in A_1}$. 
Note that if $\al\in A_1$ then $\aza\al\sq\et$ since 
$\et$  satisfies \ref{73ii} of Definition \ref{73}.

If $\i\in\et$ then $\pi(\i)=\i$, see above Part 2. 
Now let $\i\in\tup'_1\bez\et$. 
Define $m_\i<\lh\i$ as in the case $e=1,3$ above and 
define $\pi(\i)$ as in the proof of Lemma~\ref{iso1}. 

{\it Finalization\/}.  
Combining the construction in Parts 1, 2, 3, we get 
the a transformation $\pi\in\ppi\xi\cup\Gab$ that 
proves the result in case $e=2,4$.

\ref{792} 
The proof is pretty similar to Part 2 in the proof of 
\ref{77} in case $e=2,4$. 

\ref{793}
The isolation claim is case $\xi=\pu$ of \ref{792}.
\epf

\parf
{First group of choice statements in the first main theorem}
\las{g1}

The content of items 
\ref{mt11},\ref{mt12},\ref{mt13},\ref{mt14} 
of Theorem~\ref{mt1} naturally splits into 
the following three groups 
of positive and negative Choice statements:
\bde
\item[weaker $\AC$ group:]
$\xAC{\od}$ in \ref{mt12} and $\AC$ in 
\ref{mt13},\ref{mt14};

\item
[negative \dd\nn group:] 
$\neg\xAC{\ip1{\nn+1}}$,  $\neg\xAC{\fp1{\nn+1}}$,  
$\neg\xDC{\ip1{\nn+1}}$, $\neg\xDC{\fp1{\nn+1}}$;

\item
[positive \dd\nn group:]
$\xDC{\fp1\nn}$, $\xDC{\ip1{\nn+1}}$, $\xDC{\fp1{\nn}}$,  
$\xDC{\ip1{\nn+1}}$.
\ede
The groups will be treated differently, and now
we are able to establish the following theorem related 
to the first group. 
The theorem also provides the full 
$\DC$ in the auxiliary models $\rL(\W{\Ohe\et}[\w])$, 
that we'll need below.

\bte
\lam{811}
Assume that\/  $\cX\in \RF$ has the Fusion property, 
and\/ $\w\in\can\tup$ is\/ $\cX$-generic. 
Then$:$\vim 
\ben
\renu
\itlb{811*1} 
$\xAC\od$ holds in\/ $\rL(\W{\Omb}[\w])\,;$ 

\itlb{811*2} 
full\/ $\AC$ holds in\/ $\rL(\W{\Omc}[\w])$ 
and in\/ $\rL(\W{\Omd}[\w])\,;$  

\itlb{811*3}
full\/ $\DC$ holds in\/ $\rL(\W{\Ohe\et}[\w])$ 
for any\/ $e=1,2,3,4$ and\/ $\et\in\Ome$.
\een
\ete

\bpf 
\ref{811*1} 
We are going to apply Theorem~\ref{65}\ref{65i}, 
therefore it suffices to check its premices for $\Omb$. 
We know that each $\et\in\Omb$ is $\Gab$-isolated by 
Theorem~\ref{79}\ref{791}. 
On the other hand, we know that 
$\Gab\sq\Gaa=\qer$,  
and we have $\Gab\sq\inv\Omb$ since 
$\Omb$ is $\Gab$-invariant by Theorem~\ref{79}\ref{73r4}. 
This proves the isolation condition of Theorem~\ref{65}. 
Moreover, $\Omb$ satisfies\/ \ref{65*} of 
Theorem~\ref{65} 
in case\/ $\ta_0=\pu$ by  Theorem~\ref{79}\ref{781}. 
It remains to apply Theorem~\ref{65}. 

\ref{811*2}
Essentially the same argument, but with  item 
\ref{65ii} of Theorem~\ref{65} instead of \ref{65i} and 
with \ref{782} of Theorem~\ref{79} instead of \ref{781}.

\ref{811*3} 
Reference to Theorem~\ref{79}\ref{794} and Theorem~\ref{66}. 
\epf

\parf
[Structure and Definability properties]
{Structure and Definability properties}
\las{82seq}

The next definition introduces conditions leading to 
level-dependent violations of some forms of countable Choice 
in the generic models considered.

\bdf
\lam{82d}
Let\/ $n<\om$ and\/ $\w\in\can\tup$. 
We define:
\bde
\item[Structure property:]
for all\/ $\i,\j\in\tup,$ we have\/ 
$\w(\i)\in\rL[\w(\j)]$ iff\/ $\i\sq\j\,;$
\index{property!Structure property}%
\index{Structure property}%

\item[$n$-{Definability property}:]
\index{property!Definability@$n$-Definability property}%
\index{Definability@$n$-Definability property}%
if\/ $\gM\sq\rL[\w]$ is a transitive class closed under 
pairs, and\/ $\rL[x]\sq \gM$ for all\/ $x\in \gM$, 
then the sets 
$\gee{}\w\cap\gM$ and\/ $\geo{}\w\cap\gM$ 
are $\ip1{n+1}$ over $\gM$, where 
$$
\hspace*{-0ex}
\bay{rcl}
\gee \gM\w
\kmar{gee gM w}%
\index{set!Eevn@$\gee \gM\w$}%
\index{zEev@$\gee \gM\w$}%
&\hspace{-1ex}=\hspace{-1ex}&
\ens{\ang{k,\w(\i)}}
{k\ge1\land \i\in\tup\text{ is even\,}\land
\lh\i=k},\\[0.8ex]
\geo \gM\w
\kmar{geo gM w}%
\index{set!Eodd@$\geo \gM\w$}%
\index{zEov@$\geo \gM\w$}%
&\hspace{-1ex}=\hspace{-1ex}&
\ens{\ang{k,\w(\i)}}
{k\ge1\land \i\in\tup\text{ is odd\,}\land
\lh\i=k}. 
\eay\vim
$$
\ede
A forcing $\cX\in \RF$  has the 
\rit{Structure} or 
$n$-\rit{Definability property}, 
if ($\cX$ forces over $\rL$ that)  
each $\cX$-generic array $\w\in\can\tup$ 
has that 
property.  
\edf

\bre
\lam{82r}
The class $\gM$ is {\ubf not} assumed to satisfy $\zf$, and the 
sets $\gee{}\w\cap\gM$ and $\geo{}\w\cap\gM$ are {\ubf not}
claimed to belong to $\gM$ in Definition~\ref{82d}. 
In fact, the proof of Theorem~\ref{mt1} below will be related 
to the case when $\gM$ satisfies $\zf$ and hence
the sets 
$\gee{}\w\cap\gM$ and $\geo{}\w\cap\gM$ do belong to $\gM$. 
However the proof of Theorem~\ref{mt2} in Chapter~\ref{b} 
involves the case when $\gM$ is not a $\zf$-class, and 
in fact the sets 
$\gee{}\w\cap\gM$ and $\geo{}\w\cap\gM$ will not belong 
to $\gM$ in that case. 
\ere

Note that, for example,   
$\perf$-generic arrays $\w$ do not have the 
$n$-Definability property for any $n$, but do have 
the Structure property.  
The construction of forcings $\cX\in \RF$ with the 
$n$-Definability 
property is quite a difficult task. 
Below, a method will be elaborated for such 
a construction. 

 \parf
{Violation of Choice}
\las{82seq,}

The next theorem shows that the properties introduced 
by Definition \ref{82d} lead to the violation of Choice in 
appropriate submodels. 
Note the difference between the lightface and boldface 
classes.

\bte
\lam{82}
Assume that\/ $n\ge1$, $\cX\in \RF$ 
has the Structure and\/ $n$-Definability properties, 
and\/ $\w\in\can\tup$ is\/ $\cX$-generic. 
Then$:$\/\vim 
\begin{multicols}{2}
\ben
\renu
\itlb{821}
$\xAC{\ip1{n+1}}$ fails in\/ $\rL(\W{\Oma}[\w])$, 

\itlb{822}
$\xAC{\fp1{n+1}}$ fails in\/ $\rL(\W{\Omb}[\w])$, 

\itlb{823}
$\xDC{\ip1{n+1}}$ fails in\/ $\rL(\W{\Omc}[\w])$, 

\itlb{824}
$\xDC{\fp1{n+1}}$ fails in\/ $\rL(\W{\Omd}[\w])$.
\een
\end{multicols} 
\ete

\bpf
We'll make use of the following 
\rit{key sets} as counterexamples:
$$
\bay{rcl}
P_1
&=&
\ens{\ang{k,\w(\i)}}
{k\ge1\land \i\in\tup\text{ is even\,}\land
\lh\i=k},\\[1ex]
P_2
&=&
\ens{\ang{k,\w(\i)}}
{k\ge1\land \i\in\tup\text{ is even\,}\land
\lh\i=k\land \i(0)=0},\\[1ex]
P_3
&=&
\ens{\ang{\w(\i),\w(\j)}}
{\i,\j\in\tup\text{ are even\,}\land
\i\su\j},\\[1ex]
P_4
&=&
\ens{\ang{\w(\i),\w(\j)}}
{\i,\j\in\tup\text{ are even\,}\land
\i\su\j\land \i(0)=0}.\text{ \kra{ or =1 ?}}
\eay
$$

\ble
\lam{83}
Let\/ $e=1,2,3,4$. 
Then\/ $P_{e}\in\rL(\W{\Ome}[\w])$ and$:$
%
\ben
\aenu
\itlb{831}  
$P_1$ is\/ 
$\ip1{n+1}$ in\/ $\rL(\W{\Oma}[\w]);$

\itlb{832}  
$P_2$ is\/ 
$\fp1{n+1}\land\fs12$ in\/ $\rL(\W{\Omb}[\w])$, 
hence just\/ $\fp1{n+1}$ in case\/ $n\ge2;$ 

\itlb{833}  
$P_3$ is\/ 
$\ip1{n+1}\land\is12$ in\/ $\rL(\W{\Omc}[\w])$, 
hence just\/ $\ip1{n+1}$ in case\/ $n\ge2;$

\itlb{834}  
$P_4$ is\/ 
$\fp1{n+1}\land\fs12$ in\/ $\rL(\W{\Omd}[\w])$, 
hence just\/ $\fp1{n+1}$ in case\/ $n\ge2$.
\een
\ele

By $\fp1{n+1}\land\fs12$ in \ref{832} and \ref{834} 
we mean the definability by a conjunction 
of a $\fp1{n+1}$ formula and a $\fs12$ formula with 
real parameters,  and $\ip1{n+1}\land\is12$ in \ref{833} 
is understood similarly (no parameters).

\bpf[Lemma]
If $e=1,2,3,4$ then define 
$\xS_e:=\gee{\rL(\W{\Ome}[\w])}\w$  
(that is, $\gee \gM\w$ as in Definition \ref{82d} 
with $\gM=\rL(\W{\Ome}[\w])$), and   
$$
\xS^0_e=
\ens{\ang{k,\w(\i)}\in\xS_e}{\i(0)=0} = 
\ens{\ang{k,\w(\i)}\in\xS_e}{\w(\ang0)\in\rL[\w(\i)]}  
$$
(the equality holds by the Structure property of $\w$). 
We may note that $\bigcup\Oma=\bigcup\Omc=\tup$, 
whereas $\axi0\sq \bigcup\Omb=\bigcup\Omd\sneq\tup$ 
by Theorem~\ref{79}\ref{73r1}. 
\pagebreak[0]
It follows that $\w(\i)\in\rL(\W{\Ome}[\w])$ for all 
$\i\in\tup$ in case $e=1,3$, whereas 
$\w(\i)\in\rL(\W{\Ome}[\w])$ for $e=2,4$ provided 
$\i(0)=0$. 
Therefore, by the $n$-Definability property of $\w$,   
$\xS_e$ is $\ip1{n+1}$ in $\rL(\W{\Ome}[\w])$ 
for $e=1,3$, 
but $\xS_e^0$ is $\fp1{n+1}$ in $\rL(\W{\Ome}[\w])$ 
(with $p=\w(\ang0)\in\can{}$ as the only parameter)
in case $n\ge2$, and is $\fp1{n+1}\land\fs12$   
in case $n=1$ because ``$x\in\rL[y]$'' is 
a $\is12$ formula. 

\ref{831} 
We immediately conclude that $P_1=\xS_1$ is $\ip1{n+1}$ 
in $\rL(\W{\Oma}[\w])$.

\ref{832} 
Similarly $P_2=\xS^0_2$ is $\fp1{n+1}\land\fs12$ 
in $\rL(\W{\Omb}[\w])$. 

\ref{833} 
Using the Structure property of $\w$, 
we observe that 
$$
P_3
=
\ens{\ang{x,y}}
{\sus k<\ell\,(\ang{k,x}\in P_1\land\ang{\ell,y}\in P_1
\land x\in\rL[y]\land y\nin\rL[x]}.
$$
Thus $P_3$ is $\ip1{n+1}\land\is12$ in $\rL(\W{\Ome}[\w])$.

\ref{834} 
follows from \ref{833} 
similarly to \ref{831}$\imp$\ref{832}.
\epF{Lemma}


\ble
[premises]
\lam{84}
The premises of the choice principles hold$:$ 
\bce
$\dom{P_1}=\om\bez\ans0$, \qquad
$\dom{P_2}=\om\bez\ans{0,1}$,  
\\[1ex]
$\ran{P_3}\sq\dom{P_3}$, \qquad
$\ran{P_4}\sq\dom{P_4}$.
\ece
\ele
\bpf[Lemma]
Assume that $k\ge1$. 
Let $\i=1^k$ ($k$ terms equal to 1). 
Then $\ang{k,\w(\i)}\in P_1$, hence $k\in\dom{P_1}$. 

If $k\ge2$ and 
$\i=0\we 1^{k-1}$, 
then $\ang{k,\w(\i)}\in P_2$, hence $k\in\dom{P_2}$. 

Similarly,  
$\ran{P_3}=\ens{\w(\i)}{\i\in\tup\land\lh\i\ge2}\sq
\dom{P_3}=\ens{\w(\i)}{\i\in\tup}$.

Finally,  we have
$\ran{P_4}=\ens{\w(\i)}
{\i\in\tup\land\lh\i\ge2\land\i(0)=0}$,  
whereas $\dom{P_4}=\ens{\w(\i)}{\i\in\tup\land\i(0)=0}$.
\epF{Lemma} 

Coming back to Theorem~\ref{82}, we finally show that   the 
choice functions required do not exist is the 
corresponding models. 

\ref{821}
We claim that there is no function $f\in\rL(\W{\Oma}[\w])$ 
such that $\ang{k,f(k)}\in P_1$ for all $k\ge1$. 
Indeed suppose to the contrary that $f$ is such a function. 
Corollary~\ref{641} implies $f\in\rL[\w\dar\et]$ for some 
$\et\in\Oma$. 
If $k\ge1$ then by definition $f(k)=\w(\i_k)$ for some 
even $\i_k\in\tup$ with $\lh{\i_k}=k$, 
and we have $\i_k\in\et$ by Corollary~\ref{632}. 
On the other hand, by definition there is $m<\om$ such that 
$\lh\i\le m$ for all even $\i\in\et$, hence
$\lh{\i_k}\le m$ for all $k$, which contradicts the above. 

To conclude, $P_1$ witnesses that 
$\xAC{\ip1{n+1}}$ fails in\/ $\rL(\W{\Oma}[\w])$, because 
$\dom{P_1}=\om\bez\ans0$ by Lemma~\ref{84}.

\ref{822}
A very similar argument shows that 
$\xAC{\fd1{n+2}}$ fails in\/ $\rL(\W{\Omb}[\w])$ 
via $P_2$. 
The failure of $\xAC{\fp1{n+1}}$ then follows by 
Lemma~\ref{23}\ref{233}. 


\ref{823}
We claim that no function $f\in\rL(\W{\Omc}[\w])$ 
satisfies $\ang{f(k),f(k+1)}\in P_3$ for all $k$. 
Indeed otherwise such a function $f$ belongs to 
$\rL[\w\dar\et]$ for some $\et\in\Omc$, by 
Corollary~\ref{641}. 
If $k<\om$ then by definition $f(k)=\w(\i_k)$ and 
$f(k+1)=\w(\i_{k+1})$ for some even  
$\i_k,\i_{k+1}\in\et$ 
with $\i_k\su\i_{k+1}$, by Corollary~\ref{632}. 
In other words, the set 
$\et'=\ens{\i\in\et}{\i\text{ is even}}\in\rL$ 
is $\su$-ill-founded in $\rL(\W{\Omc}[\w])$. 
Then $\et'$ is ill-founded in $\rL$ as well, 
which contradicts the definition of $\Omc$. 

Thus $P_3$ witnesses the failure of  
$\xDCs{\id1{n+2}}$ in\/ $\rL(\W{\Omc}[\w])$, as 
$\ran{P_3}\sq\dom{P_3}$ by Lemma~\ref{84}. 
Lemma~\ref{23}\ref{234} helps to improve this 
to the failure of  $\xDC{\ip1{n+1}}$.


\ref{824}
The same argument with $P_4$.
\epF{Theorem~\ref{82}}

\parf{Odd-Expansion property}
\las{91}

Recall 
the notion of \rit{odd expansion} $\osq$ 
of Definition~\ref{71}.

\bdf
\lam{92}
Let\/ $n<\om$ and\/ $\w\in\can\tup$.
\bde 
\item[\dd nOdd-Expansion, or \dd noe, property of $\w$:] 
\index{property!Odd-Expansion, ($n$)-oe}%
\index{Odd-Expansion, ($n$)-oe}%
\index{oe@$n$-oe}%
for every $\et\in\cpo$ and $\ip1n$ formula 
$\vpi(\cdot)$,
with parameters in $\rL[\w\dar\et]$, if 
$\rL[\w]\mo\sus x\,\vpi(x)$ then there is 
an odd expansion $\ta\in\cpo$ of $\et$ and 
$x\in\rL[\w\dar\ta]$ such that $\rL[\w]\mo\vpi(x)$.
\ede
A forcing notion $\cX\in \RF$  has the 
\rit{\dd nOdd-Expansion property}, 
if ($\cX$ forces over $\rL$ that)  each $\cX$-generic 
array $\w\in\can\tup$ 
has the $n$-oe property.  
\edf

\vyk{
\bdf
\lam{92}
Let\/ $n<\om$. 
Say that\/ $\w\in\can\tup$ has the 
$n$-\rit{odd-expansion}, or ($n$)-oe, property, 
\index{property!odd-expansion, ($n$)-oe}%
\index{odd-expansion, ($n$)-oe}%
\index{oe@$n$-oe}%
iff for every $\et\in\cpo$ and a $\ip1n$ formula 
$\vpi(\cdot)$,
with reals in $\rL[\w\dar\et]$ as parameters, if 
$\sus x\,\vpi(x)$ is true in $\rL[\w]$ then there is 
an odd expansion $\ta\in\cpo$ of $\et$ and some 
$x\in\rL[\w\dar\ta]$ such that $\rL[\w]\mo\vpi(x)$.

A forcing notion $\cX\in \RF$  has the 
$n$-\rit{oe property}, 
if ($\cX$ forces over $\rL$ that)  each $\cX$-generic 
array $\w\in\can\tup$ 
has the $n$-oe property.  
\edf
}

This property is used through the following
lemma.

\ble
\lam{922}
Suppose that\/ $n<\om$, $e=1,2,3,4$, 
and\/ $\w\in\can\tup$ has the\/ 
$n$-Odd-Expansion property.  
Then\/ 
\ben
\renu
\itlb{922i}
$\rL(\W{\Ome}[\w])$ is an elementary submodel of\/ 
$\rL[\w]$ \poo\ all\/ $\is1{n+1}$ formulas with parameters 
in\/ $\rL(\W{\Ome}[\w])$, \ and\/ 

\itlb{922ii} 
if\/ $\xi\in\Ome$ then\/ $\rL(\W{\Ohe\xi}[\w])$ is 
an elementary submodel of\/ 
$\rL[\w]$ \poo\ all\/ $\is1{n+1}$ formulas with parameters 
in\/ $\rL(\W{\Ohe\xi}[\w])$.
\een
\ele
\bpf[sketch]
For $\is12$ formulas apply the Shoenfield absoluteness. 
The step is carried out straightforwardly using 
Theorem~\ref{79}\ref{911}. 
\epf

\bre
\lam{case1}
If $n=1$ then 
$n$-Odd-Expansion property and Lemma~\ref{922} 
definitely hold for each 
$\w$ by the Shoenfield absoluteness theorem \cite{sho}.
\ere

Now let's infer come corollaries.

\bte
\lam{93}
Assume that\/ 
$\cX\in \RF$ has the Fusion property, 
$n\ge1$, 
and\/ $\w\in\can\tup$ is\/ $\cX$-generic and has 
the\/ $n$-oe property. 
Then\/\vim 
\ben
\renu
\itlb{931} 
$\DC(\fp1n)$ holds in\/ $\rL(\W{\Oma}[\w])$ and in\/ 
$\rL(\W{\Omc}[\w])$, 

\itlb{932} 
$\DC(\ip1{n+1})$ 
{\rm(lightface!)} holds in\/ 
$\rL(\W{\Omb}[\w])$ and in\/ $\rL(\W{\Omd}[\w])$.
\een
\ete

\bpf
\ref{931}
Consider a $\ip1n$ formula $\vpi(x,y)$ such that
\ben
\fenu
\itlb{931*}
$\rL(\W{\Oma}[\w])\mo\kaz x\,\sus y\,\vpi(x,y)$, 
\een
and with parameters in $\rL(\W{\Oma}[\w])$. 
Let $x_0\in\cN\cap\rL(\W{\Oma}[\w])$.
There is $\xi\in\Oma$ such that $x_0$ and 
all parameters in $\vpi$ belong to $\rL[\w\dar \et]$. 
Consider the submodel 
$\rL(\W{\Oha\xi}[\w])\sq\rL(\W{\Oma}[\w])$. 
Thus $\xi\in\Oha\xi$, and hence 
$x_0$ and 
all parameters in $\vpi$ belong to $\rL(\W{\Oha\xi}[\w])$. 
However 
\ben
\fenu
\atc
\itlb{931**}
$\rL(\W{\Oha\xi}[\w])$ is an 
elementary submodel of $\rL(\W{\Oma}[\w])$ \poo\ 
all $\is1{n+1}$ formulas with reals in 
$\rL(\W{\Oha\xi}[\w])$ as parameters, 
by Lemma~\ref{922}. 
\een
Therefore 
$\rL(\W{\Oha\xi}[\w])\mo\kaz x\,\sus y\,\vpi(x,y)$ 
by \ref{931*}. 
Moreover, $\rL(\W{\Oha\xi}[\w])\mo\DC$ by 
Theorem~\ref{811}\ref{811*3}. 
This allows to define a sequence 
$\sis{x_k}{k<\om}\in \rL(\W{\Oha\xi}[\w])$ of reals, 
beginning with the $x_0$ given above, and satisfying 
$\rL(\W{\Oha\xi}[\w])\mo\vpi(x_k,x_{k+1})$, $\kaz k$. 
It remains to refer to \ref{931**} in order to return 
to $\rL(\W{\Oma}[\w])$.  

The proof for $\rL(\W{\Omc}[\w])$ is pretty similar.\vtm

\vyk{
We claim that (\mdag) $\rL(\W{\Oha\xi}[\w])$ is an 
elementary submodel of $\rL(\W{\Oma}[\w])$ \poo\ 
all $\is1{n+1}$ formulas with reals in 
$\rL(\W{\Oha\xi}[\w])$ as parameters. 

Let $\psi(\cdot)$ be a $\ip1n$ 
formula with parameters in $\rL(\W{\Oha\xi}[\w])$, 
such that $\rL(\W{\Oma}[\w])\mo\sus x\,\psi(x)$. 
All parameters in $\psi$ belong to $\rL[\w\dar \vt]$ 
for some $\vt\in\Oha\xi$. 
It follows, by ($n$)-oe, that there is 
an odd expansion $\ta\in\cpo$ of $\vt$, and some 
$x\in\rL[\w\dar\ta]$, such that $\rL[\w]\mo\psi(x)$. 
However $\ta\in\Oha\xi$ by Lemma~\ref{79}\ref{911}. 
Thus $x\in \rL(\W{\Oha\xi}[\w])$, and this completes 
the proof of (\mdag).
 
It is an immediate corollary of (*) that
$\rL(\W{\Oha\xi}[\w])\mo\kaz x\,\sus y\,\vpi(x,y)$. 
Moreover, $\rL(\W{\Oha\xi}[\w])\mo\DC$ by 
Theorem~\ref{811}\ref{811*3}. 
This allows to define a sequence 
$\sis{x_k}{k<\om}\in \rL(\W{\Oha\xi}[\w])$ of reals, 
beginning with the $x_0$ given above, and satisfying 
$\rL(\W{\Oha\xi}[\w])\mo\vpi(x_k,x_{k+1})$, $\kaz k$. 
It remains to refer to (\mdag).
}

\ref{932} 
This part involves trickier arguments
contained in two lemmas.

\ble
\lam{95}
Assume that\/ $\xi\in\Omb$, $\vpi(y)$ is a 
parameter-free\/ $\is1\iy$ formula, and\/ 
$\rL(\W{\Ohb\xi}[\w])\mo\sus y\,\vpi(y)$. 
Then there is\/ $y\in\rL(\W{\Oms}[\w])$ such that\/ 
$\rL(\W{\Ohb\xi}[\w])\mo\vpi(y)$. 
The same for\/ $\Omd$ and\/ $\Ohd\xi$.
\ele

\bpf[Lemma]
The $\cL$-formula
$$
\chi(U)\,:=\,
\sus y\in\cN\cap\rL(U)\,
\big(
y\in\WW{\Ohb\xi}\land \rL(\WW{\Ohb\xi})\mo\vpi(y)
\big)
$$
satisfies 
$\modd\chi=\pu$ and $\inv\chi=\inv{\Ohb\xi}$. 
Under the assumptions of the lemma, 
$\rL[\w]\mo\chi(\WW{\Ohb\xi})[\w]$, 
via some $y\in\W{\Ohb\xi}[\w]$.
Then $y\in\rL[\w\dar\ta]$, $\ta\in\Ohb\xi$, 
by Corollary~\ref{641} 
(in which the isolation condition follows from 
Theorem~\ref{97}\ref{792}). 
Thus $\rL[\w]\mo\chi(\pvr\ta)[\w]$. 
Corollary~\ref{633} yields 
a condition $X\in(\cX\dar\ta)\cap \cg\w$ 
such that (1) $X\fox\cX\chi(\pvr\ta)$. 
We claim that
\ben
\nenu
\atc
\itlb{952}
$X\fox\cX\chi(\WW\Oms)$ --- 
that obviously implies the lemma.
\een
Suppose towards the contrary that \ref{952} fails. 
Then (3) $Y\fox\cX\neg\,\chi(\WW\Oms)$ holds for some 
$Y\in\cX$, $Y\psq X$, but still 
(4) $Y\fox\cX\chi(\pvr\ta)$ by (1). 
We may assume that $\xi\sq\ta$, and that $\dym Y=\ta$ 
by Theorem~\ref{631}. 

By Theorem~\ref{79}\ref{77pu},   
there is a permutation 
$\pi\in\Ghb\xi$ satisfying $\sg=\pi\akt\ta\in\Oms$ 
and $\sg\cap\ta=\pu$. 
We may note that $\Ghb\xi\sq\inv{\Ohb\xi}\sq\per$, 
so that $\pi\in
\inv\chi$. 
Then we have from (4) by Corollary~\ref{624} 
that $S\fox\cX\chi(\pvr\sg)$, where $S=\pi\akt Y$, and 
further (5) $S\fox\cX\chi(\WW\Oms)$ as $\sg\in\Oms$.

However conditions $S$ and $Y$ are compatible because 
$\ta\cap\sg=\pu$. 
Thus (5) contradicts to (3), which proves \ref{952} and 
the lemma.
\epF{Lemma}

\ble
\lam{96}
Assume that\/ $\vpi(y)$ is a\/ $\is1\iy$ formula 
with parameters in\/ $\cN\cap\rL(\W{\Oms}[\w])$, 
and there is\/ $y\in\rL(\W{\Omb}[\w])$ such that\/ 
$\rL[\w]\mo\vpi(y)$. 
Then there is\/ $x\in\rL(\W{\Oms}[\w])$ such that\/ 
$\rL[\w]\mo\vpi(x)$.
\ele

\bpf[Lemma]
By Corollary~\ref{641}, there is $\xi\in\Oms$ such that 
all parameters in $\vpi(\cdot)$ belong to $\rL[\w\dar\xi]$. 
Then there is an $\cL$-formula $\psi(\cdot)$ that contains 
only $\pvr\xi$ and some $\namx z$, $z\in\rL$, 
as the only $\cL$-names, 
and such that $\psi(\cdot)[\w]$ is identic to $\vpi(\cdot)$.
Let $\chi(U)$ be the $\cL$-formula: 
$
(\sus x\in\cN\cap\rL(U))\,\psi(x) 
$. 
Then obviously $\modd\chi=\modd\psi=\xi$ and 
$\inv\chi=\per$.

By definition we have $\rL[\w]\mo\chi(\WW\Omb)[\w]$, 
where both $\modd{\chi(\WW\Omb)}=\xi$ and 
$\inv{\chi(\WW\Omb)}=\Gab$ by the above. 
It follows by Corollary~\ref{633} that  
there is a condition $X\in(\cX\dar\xi)\cap \cg\w$ 
such that (1) $X\fox\cX\chi(\WW\Omb)$. 
We claim that
\ben
\nenu
\atc
\itlb{962}
$X\fox\cX\chi(\WW\Oms)$ --- 
which obviously proves the lemma.
\een
Suppose towards the contrary that \ref{962} fails. 
Then (3) $Y\fox\cX\neg\,\chi(\WW\Oms)$ holds for some 
$Y\in\cX$, $Y\psq X$. 
We may assume that $\dym Y=\xi$ by Theorem~\ref{631}. 
Then $Y\sq X$ and $Y\fox\cX\chi(\WW\Omb)$ by (1). 
We conclude by Corollary~\ref{642} that there is a 
condition $Z\in\cX$, $Z\psq Y$, and $\ta\in\Omb$, 
such that (4) $Z\fox\cX\chi(\pvr\ta)$. 
We can \noo\ assume that $\xi\sq\ta=\dym Z$.

Theorem~\ref{79}\ref{77} yields a permutation 
$\pi\in\ppi{\xi}$ 
with $\sg=\pi\akt\ta\in\Oms$. 
Then we have $S\fox\cX\chi(\pvr\sg)$ from (4) 
by Corollary~\ref{624}, where $S=\pi\akt Z$.
We further conclude that (5) $S\fox\cX\chi(\WW\Oms)$ 
since $\sg\in\Oms$. 

On the other hand, $S\dar\xi=Z\dar\xi$ holds because 
$\pi\in\ppi\xi$. 
Therefore $S\psq Y$ (since $\dym Y=\xi$ and $Z\psq Y$). 
It follows that (3) and (5) are contradictory. 
The contradiction proves \ref{962} and the lemma.
\epF{Lemma}

We continue 
{\ubf the proof of Theorem~\ref{93}\ref{932}}. 
Consider a \paf\ $\ip1{n+1}$ formula $\vpi(\cdot,\cdot)$, 
satisfying $\rL(\W{\Omb}[\w])\mo\kaz x\,\sus y\,\vpi(x,y)$. 
Let $u\in\cN\cap\rL(\W{\Omb}[\w])$.
Corollary~\ref{641} implies $u\in\rL[\w\dar\xi],$ 
some $\xi\in\Omb$. 
Then $\xi\in\Ohb\xi$, $u\in\rL(\W{\Ohb\xi}[\w])$. 

\ble
\lam{97}
$\rL(\W{\Ohb\xi}[\w])\mo\kaz x\,\sus y\,\vpi(x,y)$.
\ele

\bpf[Lemma]
Suppose otherwise. 
Then by Lemma~\ref{95} there is 
$p\in\cN\cap\rL(\W{\Oms}[\w])$ such that 
(\mast) 
$\rL(\W{\Ohb\xi}[\w])\mo \kaz y\,\vpi^-(p,y)$, where 
$\vpi^-(x,y)$ is the canonical $\is1{n+1}$ formula 
equivalent to $\neg\,\vpi(x,y)$. 

However $p\in\rL(\W{\Oms}[\w])$, and hence, 
we have 
$\rL(\W{\Omb}[\w])\mo\sus y\,\vpi(p,y)$ 
in our assumptions. 
Then $\rL[\w]\mo\sus y\,\vpi(p,y)$ by Lemma~\ref{922}. 
Furthermore, by Lemma~\ref{96}, there is 
$q\in\cN\cap\rL(\W{\Oms}[\w])$ satisfying  
$\rL[v]\mo\vpi(p,q)$. 
Now $p,q\in\rL(\W{\Ohb\xi}[\w])$ by 
Theorem~\ref{79}\ref{795}, and we have 
$\rL(\W{\Ohb\xi}[\w])\mo \vpi(p,q)$ 
still by Lemma~\ref{922}. 
But this contradicts (\mast).
\epF{Lemma}

Now let us accomplish the proof of 
Theorem~\ref{93}\ref{932}. 
By the last lemma, and since 
$\rL(\W{\Ohb\xi}[\w])\mo \DC$ 
(by Theorem~\ref{811}\ref{811*3}), there is a 
sequence $\sis{x_k}{k<\om}\in\rL(\W{\Ohb\xi}[\w])$ 
of reals $x_k$ satisfying $x_0=u$ and 
$\rL(\W{\Ohb\xi}[\w])\mo\vpi(x_k,x_{k+1})$, $\kaz k$. 
Then Lemma~\ref{922} implies 
$\vpi(x_k,x_{k+1})$, $\kaz k$, in $\rL(\W{\Omb}[\w])$ 
as well, as required. 
\epf

\parf{Second form of the main theorem}
\las{98}

To summarize the results achieved above, we approach  
our first main result (Theorem~\ref{mt1} in the introduction) 
via the following theorem. 

\bte
[in $\rL$]
\lam{mt1b}
If\/ $\nn\ge1$  
then there is a forcing\/ 
$\cX\in\RF$ with the Fusion, Structure, 
$\nn$-Definability, and\/ $\nn$-Odd-Expansion 
properties.
\ete

\bpf[Theorem~\ref{mt1} from Theorem~\ref{mt1b}]
Let $\w\in\can\tup$ be $\cX$-generic over $\rL$. 
\sloppy
Then:
\bit
\item[$-$]\msur 
$\xAC \od $ holds in $\rL(\W{\Omb}[\w])$, 
whereas 
the full\/ $\AC$ holds in\/ $\rL(\W{\Omc}[\w])$ 
and in\/ $\rL(\W{\Omd}[\w])$ --- by Theorem~\ref{811};

\item[$-$]\msur 
$\xAC{\ip1{\nn+1}}$, 
$\xAC{\fp1{\nn+1}}$, 
$\xDC{\ip1{\nn+1}}$, 
$\xDC{\fp1{\nn+1}}$ fail in resp.\ models
$\rL(\W{\Oma}[\w])$, 
$\rL(\W{\Omb}[\w])$, 
$\rL(\W{\Omc}[\w])$, 
$\rL(\W{\Omd}[\w])$ --- 
by Theorem~\ref{82};

\item[$-$]\msur 
$\xDC{\fp1\nn}$ holds in $\rL(\W{\Oma}[\w])$ and in  
$\rL(\W{\Omc}[\w])$, 
whereas 
$\xDC{\ip1{\nn+1}}$ holds in\/ 
$\rL(\W{\Omb}[\w])$ and in\/ $\rL(\W{\Omd}[\w])$ 
--- by Theorem~\ref{93}.
\eit

\noi 
Thus $\rL(\W{\Oma}[\w])$, 
$\rL(\W{\Omb}[\w])$, 
$\rL(\W{\Omc}[\w])$, 
$\rL(\W{\Omd}[\w])$ 
are models of $\zf$ in which conjunctions resp.\ 
\ref{mt11}, \ref{mt12}, \ref{mt13}, \ref{mt14} of 
Theorem~\ref{mt1} hold,  as required. 
\epf

Thus Theorem~\ref{mt1b} (to be proved below)  
implies Theorem~\ref{mt1}, 
the first main result of this paper.  
The next approximation will be Theorem~\ref{mt1c}.

\sekt
{Reduction of Odd-Expansion to Completeness}
\las{foax}

The goal of this Chapter is to reduce    
$n$-Odd-Expansion property of generic arrays, as 
in Definition~\ref{92}, to a property of a given 
normal forcing notion $\cX\sq\pei$, called 
$n$-{\ubf Completeness} (Section~\ref{12}). 
This property will essentially say that $\cX$ is 
an elementary substructure of $\pei$ \poo\ the 
forcing relation for $\us1n$ formulas. 
We'll begin with some technicalities, which 
include the following.
\bit
\item[$-$] 
Coding of continuous maps, with applicatios to the 
property of {\ubf continuous reading of names}
under Fusion (Section~\ref{101}). 

\item[$-$] 
The {\ubf extension} of the language of
second order arithmetic by codes of maps, 
and a forcing-type {\ubf relation $\fo$} for the 
{\ubf extended language} (Section~\ref{103}). 
It occurs that $\fo$ is invariant \poo\ all order-preserving 
transformations of $\tup$, not necessarily those 
parity-preserving.

\item[$-$] 
The invariance mentioned is one of the two main 
ingredients in the proof of the 
{\ubf narrowing} and {\ubf odd expansion} 
theorems for $\fo$ (Sections~\ref{112} and \ref{114}). 

\item[$-$] 
The other ingredient is the {\ubf action} of projection-keeping 
homomorphisms on $\fo$ defined and studied in 
Section~\ref{114}.
\eit
Then we introduce the concept of an \dd n{\ubf complete} 
normal forcing notion in Section~\ref{12}, prove the truth 
theorem for such forcing notions and according generic 
extensions (Theorem~\ref{122}), and then 
Theorem~\ref{123} shows that $n$-Completeness 
implies $n$-Odd-Expansion. 

Theorem~\ref{mt1c} summarizes the content of this chapter.

Note that this content has no relation to the 
case $\nn=1$ of Theorems \ref{mt1} and \ref{mt1b} because 
the $n$-Odd-Expansion property holds for $n=1$ anyway.

\parf{Coding of continuous maps and 
continuous reading of names}
\las{101}

The Baire space $\cN=\bn$ is a separable Polish space, 
and such is the Cantor space $\cam=\dn\sq\cN$, 
as well as every space 
of the form $\can\xi$ 
and every closed subset in such a space. 
In addition, the spaces $\cam$ and $\can\xi$ are 
compact. 
It follows from the compactness that a function 
$F:\can\xi\to\cN$ 
is continuous ($F\in\cf_\xi$, Section~\ref{173}), 
iff its \rit{graph} 
$\ens{\ang{x,F(x)}}{x\in\can\xi}$ 
(identified with $F$) 
is a closed set in $\can\xi\ti\cN.$ 
Therefore, if $F:\can\xi\to\cN$ 
is in fact continuous, and a set $X\sq\can\xi$ is 
topologically dense in $\can\xi$ then (the graph of) 
$F$ coincides with the closure $\clo{(F\res X)}$ of 
the restricted map $F\res X$ in $\can\xi\ti\cN.$ 
We take 
\pagebreak[0]
$$
\rat\xi=\ens{x\in\can\xi}
{x(\i)(k)=0\text{ for all but finite pairs }
\ang{\i,k}\in\xi\ti\om}
\kmar{rat xi}%
\index{zRat@$\rat\xi$, rationals in $\can\xi$}%
\index{set!Rat@$\rat\xi$, rationals in $\can\xi$}%
$$
($\can\xi$-rationals) 
as a \rit{canonical} countable dense set in 
$\can\xi$. 
Accordingly let 
$$
\bay{rcl}
\kmar{kont xi}%
\index{codes of functions!$\kont\xi$}
\index{functions!codes $\kont\xi$}
\index{zcCFxi@$\kont\xi$, codes}
\kont\xi
&\!=\!&
\ens{f\in\rL}
{f:\rat\xi\to\cN\land \clo f \text{ is a continuous map }
\can\xi\to\cN};\\[0.5ex]
\kmar{kond xi}%
\index{codes of functions!$\kond\xi$}
\index{functions!codes $\kond\xi$}
\index{zcCFdxi@$\kond\xi$, codes}
\kond\xi
&\!=\!&
\ens{f\in\kont\xi}
{f:\rat\xi\to\can{}, 
\text{ so $\clo f:\can\xi\to\can{}$ is continuous}}. 
\eay
$$
If $f\in\kont\xi$ then let $\modd f=\xi$. 
\kmar{modd f}%
\index{dimension, $\modd f$}%
\index{z11f11@$\modd f$, dimension}%

We further define $\komt=\bigcup_{\xi\in\cpo}\kont\xi$ 
\kmar{komt\mns komd}%
and $\komd=\bigcup_{\xi\in\cpo}\kond\xi$; 
thus $\komt,\komd\in\rL$. 
\index{codes of functions!$\komt$}%
\index{functions!codes $\komt$}%
\index{zcCF@$\komt$, codes}%
\index{codes of functions!$\komd$}%
\index{functions!codes $\komd$}%
\index{zcCFd@$\komd$, codes}%
Each $f\in\komt$ is viewed as a \rit{code} of 
the continuous map $\clo f\in\cf$, 
and each $f\in\komd$ as a \rit{code} of 
the continuous map $\clo f\in\cfd.$ 

In the particular case $\ta=\pu$ we have 
$\can\pu=\rat\pu=\ans\pu$, accordingly $\kont\pu$ 
consists of all functions $h_x(\pu)=x$, $x\in\cN,$ 
defined on $\ans\pu$.%
\vom

We would prefer to deal with continuous functions 
$F:\can\ta\to\cN$ themselves rather than their countable 
codes. 
But as any such $F$ is an uncountable set, 
this would make hardly possible to treat 
definability questions on the basis of definability 
over $\hc=\ans{\text{all hereditarily countable sets}}$. 
Anyway the following corollary is a paraphrase 
of Theorem~\ref{1022}\ref{1022iii} reflecting the 
property of continuous reading of names (of reals) 
under Fusion.

\bcor
[of Theorem~\ref{1022}\ref{1022iii}]
\lam{1011}
Assume that, in\/ $\rL$, $\cX\in \RF$   
has the Fusion property, $\w\in\can\tup$ 
is\/ $\cX$-generic over\/ $\rL$, $\ta\in\cpo$,  
and\/ $a\in\cN\cap\rL[\w\dar\ta]$.  
Then there is\/ $f\in\kont\ta$ 
such that\/ $a=\clo f(\w\dar\ta)$.\qed 
\ecor

\vyk{
The next lemma treats the question of definability.

\ble
\lam{1012}
If\/ $\ta\in\cpo$ then\/ 
$\kont\ta$ is a\/ $\is0\iy$ subset of\/ 
$(\bn)^{\rat\ta}.$ 
\ele
\bpf
Let $f:\rat\ta\to\bn$. 
Then $f\in\kont\ta$ iff $\Phi_\ta(f)$, where 
\bde
\item[$\Phi_\ta(f):=$]
$\kaz k,\ell<\om\,$\big(%
there is $m<\om$, a finite set $\xi\sq\ta$, and 
a partition $(2^m)^\xi=A\cup B$ of the (finite) set 
$(2^m)^\xi$ of dyadic multituples, such that for all 
 $x\in\can\ta$:

--- if $\sis{x(\i)\res m}{\i\in\xi}\in A$ then $f(x)(k)=\ell$, 

--- if $\sis{x(\i)\res m}{\i\in\xi}\in B$ then $f(x)(k)\ne\ell$  
\big).
\ede
This implies the lemma.
\epf

}

\parf{Forcing approximation}
\las{103}

Corollary~\ref{1011} enables us to introduce a special 
language for describing elements of $\cN=\bn$ in generic
extensions, using function codes in $\komt$ 
to be names of elements of type 1 (i.e., taking
values in $\cN$ when interpreted).

Consider the language of 2nd order Peano arithmetic 
with type-0 variables $k,l,m,n$ over $\om$ and 
type-1 variables $x,y,z,\dots$ over $\cN.$ 
The following are standard classes of formulas:
\bde
\item[$\us0\iy$]= arithmetic formulas, \ie, no type-1 
quantifiers;

\vyk{
\item[$\us10$]= formulas of the form $\sus k\,\psi(k)$, 
$\psi$ being $\us00$;

\item[$\up10$]= formulas of the form $\kaz k\,\psi(k)$, 
$\psi$ being $\us00$;
}

\item[$\us1{n+1}$]= formulas of the form $\sus x\,\psi(x)$, 
$\psi$ being $\up1n$ (or $\us0\iy$ in case $n=0$);

\item[$\up1{n+1}$]= formulas of the form $\kaz x\,\psi(x)$, 
$\psi$ being $\us1n$ (or $\us0\iy$ in case $n=0$). 
\ede

Let $\xla$ be the extension of this language 
\kmar{xla}%
\index{language!$\xla$}%
\index{Llanguage@$\xla$, language}%
\index{zL@$\xla$, language}%
by using natural numbers as type-0 
parameters and function codes $f\in \komt$ --- as 
type-1 parameters. 
Let $\ls0\iy$, $\ls1n$, $\lp1n$ be the according classes 
\kmar{ls 1n, lp1n}%
\index{formulas!$\ls0\iy$, $\ls1n$, $\lp1n$}%
\index{zlaSP@$\ls0\iy$, $\ls1n$, $\lp1n$, formulas}%
of $\xla$-formulas. 

If $\vpi$ is an $\xla$-formula then let 
$\modd\vpi=\bigcup\ens{\modd f}{f\text{ occurs in }\vpi}$. 
If $\vpi$ is $\ls1n$, then $\otr\vpi$ 
denotes the result of the canonical reduction 
of $\neg\,\vpi$ to $\lp1n$-form; 
\kmar{otr vpi}%
\index{formulas!canonical negation $\otr\vpi$}%
\index{negation $\otr\vpi$}%
\index{zzfi-@$\otr \vpi$, negation}%
similarly for $\vpi$ in $\lp1n$.  
If $\vpi$ is $\ls0\iy$ then $\otr\vpi$ is just $\neg\,\vpi$. 

If $\vpi$ is an $\xla$-formula, $\modd\vpi\sq\et\sq\tup$ 
and $v\in\can\et$, 
then the \rit{valuation} $\vva\vpi v$ 
\index{valuation, $\vva\vpi v$}%
\index{zzfi,v,@$\vva\vpi v$, valuation}%
\kmar{vva vpi v}%
is a usual 2nd order arithmetic formula with 
type-1 parameters in $\cN\cap\rL[v\dar\modd\vpi]$, 
obtained by the substitution of the valuation 
$\vva fv:=\clo f(v\dar \modd f)\in\cN$ for every code 
\index{valuation, $\vva f v$}%
\index{zzf,v,@$\vva f v$, valuation}%
\kmar{vva f v}%
$f\in\komt$ in $\vpi$. 

\bdf
[in $\rL$]
\lam{fod}
Define a relation $X\fo\vpi$, where $X\in\pei$ and 
\kmar{fo}%
\index{forcing!approximation $\fo$}%
\index{zforc@$\fo$, forcing approximation}%
$\vpi$ is a closed $\xla$-formula in 
${\ls0\iy}\cup{\bigcup_{k\ge1}(\ls1k\cup\lp1k)}$, 
by induction. 
\ben
\cenu
\itlb{fod1}
If $\vpi$ is a closed formula in 
$\ls0\iy\cup\ls11\cup\lp11$, and $X\in\pei$, 
then $X\fo\vpi$ iff $\vva\vpi x$ holds for all 
$x\in X\uar\ta$, where $\ta=\modd\vpi\cup\dym X$.

\vyk{
\itlb{fod2}
$X\fo\sus m\,\vpi(m)$ iff $X\fo\vpi(m)$ for some  
$m<\om$.

\itlb{fod3}
$X\fo\kaz m\,\vpi(m)$ iff $X\fo\vpi(m)$ for all  
$m<\om$.
}

\itlb{fod4}
If $\vpi(x)$ is a $\lp1k$ formula, $k\ge1$,
then $X\fo\sus x\,\vpi(x)$ iff $X\fo\vpi(f)$ for some  
$f\in\komt$.

\itlb{fod5}
If $\vpi$ is a closed $\lp1k$ formula, $k\ge2$, 
$X\in\pei$, then $X\fo\vpi$ iff there exists no 
$Y\in\pei$, $Y\psq X$, such that $Y\fo\otr\vpi$.
\qed
\een
\eDf

\ble
\label{111}
\ben
\renu
\itlb{111a}
If\/ $X\fo\vpi$, $Y\in\pei$, $Y\psq X$, then\/ 
$Y\fo\vpi\,;$

\itlb{111b}
$X\fo\vpi$ and\/ $X\fo\otr\vpi$ cannot hold together$;$

\itlb{111c1}
if\/ $X\in\pei$, 
$\vpi$ is a closed\/ $\ls11$ formula, 
then there exists\/ 
$Y\in\pei$, $Y\psq X$ such that\/ 
$Y\fo\vpi$ or\/ 
$Y\fo\otr\vpi\;;$

\itlb{111c}
if\/ $X\in\pei$, $k\ge2$, 
$\vpi$ is a closed\/ $\lp1k$ formula, 
and\/ $\neg\:X\fo\vpi$ then there exists\/ 
$Y\in\pei$, $Y\psq X$ such that\/ $Y\fo\otr\vpi\;;$

\itlb{111d}
if\/ $X\in\pei$, $\et=\dym X\sq\ta\in\cpo$, and\/ 
$X\uar\ta\fo\vpi$ then\/ $X\fo\vpi\,.$
\een
\ele
\bpf
Here \ref{111b},\ref{111c} hold by definition, 
\ref{111a} 
is verified by routine induction.  

To check \ref{111c1}, note that  
the set 
$U = \ens{v \in X\uar\ta}{\vva\vpi v}$ is $\fs11$,  
where $\ta=\dym X\cup\modd\psi$, hence it has the Baire  
property in $X\uar\ta$.
It follows by Corollary~\ref{bana} that there exists a set 
$Y\in\pe\ta$ such that either $Y\sq U$, 
or $U\sq(X\uar\ta)\bez U$. 
Then accordingly $Y\fo \vpi$ or $Y\fo \otr\vpi$, 
as required. 

\ref{111d} 
Lemma~\ref{apro} makes sure that 
$X\uar\ta\in\pei$. 
The proof goes by induction, and \ref{fod5} is the 
only nontrivial step.  
Suppose to the contrary that $\psi$ is $\ls1k$, 
$X\uar\ta\fo\otr\psi$, but $\neg\:X\fo\otr\psi$. 
There is $Y\in\pei$, $Y\psq X$, $Y\fo\psi$. 
Let $\xi=\dym Y$, $\za=\xi\cup\ta$, $Z=Y\uar\za$, 
$\et'=\tau\cap\xi$. 
Then $Z\psq Y$, hence $Z\fo\psi$ by \ref{111a}. 
However $Z\dar\ta=(Y\dar\et')\uar\ta$ by 
Lemma~\ref{599}. 
Here $Y\dar\et'\sq X\uar\et'$ since $Y\psq X$, 
as clearly $\et\sq\et'$. 
Therefore $Z\dar\ta\sq X\uar\et'\uar\ta=X\uar\ta$. 
Thus $Z\psq X\uar\ta$. 
We conclude that $Z\fo\otr\psi$ by \ref{111a}. 
Yet $Z\fo\psi$ as well, see above. 
This contradicts \ref{111b}.
\epf

Assume that $\et,\sg,\ta\in\cpo$,  $\xi=\et\cup\sg\cup\ta$, 
$f\in\kont\sg$, $g\in\kont\et$, $X\in\pe\ta$. 
Say that $f,g$ are \rit{valuation-equivalent}, or simply 
\index{formulas!valuation-equivalent}%
\index{valuation-equivalent}%
\index{v-equivalent}%
\rit{v-equivalent} on $X$, iff 
$\clo f(x\dar\sg)=\clo g(x\dar\et)$ for all $x\in X\uar\xi$. 
Then, $\xla$-formulas $\vpi,\psi$ are v-equivalent on $X$ 
if $\psi$ is obtained from $\vpi$ by a substitution 
of all codes $f\in\komt$ occurring in $\vpi$ with codes $g$ 
v-equivalent to $f$ on $X$. 

\ble
[in $\rL$, routine by induction]
\lam{115}
If\/ $X\fo\vpi$, and\/ $\xla$-formulas\/ $\vpi,\psi$ 
are v-equivalent on\/ $X$ then\/ $X\fo\vpi$ iff\/ 
$X\fo\psi$.\qed
\ele


\ble
[in $\rL$]
\lam{116}
Assume that\/ $X\in\pei$, 
$\vpi(x)$ is a\/ $\lp1k$-formula, $k\ge1$, 
$\ta=\dym X\cup\modd\vpi$, and\/ $X\fo\sus x\,\vpi(x)$. 
Then 
\vyk{$:$ 
\ben
\renu
\itlb{116a}
if\/ $k=0$ then there is a code\/ $g\in\kont\ta$ 
such that\/ $X\fo\vpi(g)\;;$ 

\itlb{116b}
if\/ $k\ge1$ then 
}%
there is a code\/ $g\in\kont{\xi}$ 
for some\/ $\xi\in\cpo,\,\ta\sq\xi$, 
such that\/ $X\fo\vpi(g)\;.$ 
\ele
\bpf
By definition we have $X\fo\vpi(f)$ for a code $f\in\komt$. 
Let $\sg=\modd f$ and $\xi=\sg\cup\ta$. 
\vyk{
\ref{116a}
For each $x\in\can\ta$ define $x'\in\can\xi$ by 
$x'\dar\ta=x$ and $x'(\i)=\om\ti\ans0$ for all 
$\i\in\xi\bez\ta$. 
Define $g\in\kont\ta$ by $g(x)=f(x'\dar\sg)$ for each 
$x\in\rat\ta$. 
To prove $X\fo\vpi(g)$ 
we have to check $\vva{\vpi(g)}x$ for all $x\in\can\ta$. 
However $\vva{\vpi(g)}x$ is identic to $\vva{\vpi(f)}{x'}$ 
by construction. 
On the other hand, $x'\in X\uar\xi$, which implies  
$\vva{\vpi(f)}{x'}$ because $X\fo\vpi(f)$, as required. 

\ref{116b}
}
Define $g\in\kont\xi$ by $g(x)=f(x\dar\sg)$ for each 
$x\in\rat\xi$, and use Lemma~\ref{115}. 
\epf

\parf{The narrowing and odd expansion theorems}
\las{112}

Corollary~\ref{1011} allows to view $\fo$ as a forcing-type 
relation compatible with $\pei$ as the forcing notion. 
Yet unlike the ordinary forcing $\fox\pei$, 
$\fo$ treats the $\exists$ quantifier 
over $\cN=\bn$ in the sense of continuous reading of names. 
This adds difficulty and extra work to the proof 
of the next theorem.

\bte
[the narrowing theorem, in $\rL$]
\lam{1152}
Suppose that\/ $\vpi$ is a closed\/ $\xla$-formula, 
$\modd\vpi=\et\sq\ta\in\cpo$, $X\in\pe\ta$,   
$X\fo\vpi$.   
Then\/ $X\dar\et\fo\vpi$. 
\ete

This is quite similar to Theorem~\ref{631}, 
but the proof in Section~\ref{114}  
will be somewhat more difficult because of the 
mentioned difference in the treatment of $\sus$. 
Meanwhile, here we apply Theorem~\ref{1152} in the proof 
of the following result. 
Recall Definition~\ref{71} on odd expansions.

\bte
[the odd expansion theorem, in $\rL$]
\lam{1121}
Let\/ $k<\om$, 
$\vpi(x)$ be an\/ $\lp1k$-formula, 
$\modd{\vpi(x)}=\ta_0$, $X\in\pei$,   
$X\fo\sus x\,\vpi(x)$.   
Then there is an odd expansion\/ $\ta\in\cpo$ 
of\/ $\ta_0$, and\/ $g\in\kont\ta$, such 
that\/ $X\fo\vpi(g)$. 
\ete

The next lemma will be used in the proof of 
Theorem~\ref{1121} as well.

\ble
[in $\rL$]
\lam{1122}
Assume that\/ $\ta_0\sq\sg$ belong to\/ $\cpo$. 
Then there is\/ $\pi\in\ppio{\ta_0}$ such that\/ 
$\ta=\pi\akt\sg$ is an odd expansion of\/ $\ta_0$ 
and\/ $\ta\cap\sg=\pu$. 
\ele

Recall that $\pero$ consists of all, not necessarily 
parity-preserving, permutations of $\tup$, and 
$\ppio{\ta_0}$ contains all $\pi\in\pero$ such that 
$\pi\res\ta_0$ is the identity.

\bpf
Emulating the proof of Lemma~\ref{iso1}, we 
let $\la<\omi$ be a limit ordinal 
bigger than $\sup\ens{\i(k)}{\i\in\sg\land k<\lh\i}$. 
For any $\al<\la$,   
\kmar{?? odd\mns or even}
pick an {\ubf odd} 
ordinal $\la\le\ba(\al)\le\la+\la$ such that 
$\al<\al'\imp\ba(\al)<\ba(\al')$.  
If $\al<\omi$, let $B(\al)=B\obr(\al)=\ba(\al)$, 
whereas $B(\al)=\al$ in case 
$\al\nin\la\cup\ens{\ba(\al')}{\al'<\la}$. 
Thus $B$ is a bijection of $\omi$. 

If $\i\in\tup$ then define $\j=\ro(\i)\in\tup$ such that 
$\lh\j=\lh\i$ and $\j(\ell)=\i(B(\ell))$ for all 
$\ell<\lh\j=\lh\i$, thus $\ro$ is a 
permutation in $\pero$. 

Now let $\i\in\tup$. 
Take a largest number $m_\i\le\lh\i$ such that 
$\i\res m_\i\in\ta_0$. 
Then $\i=(\i\res m_\i)\we \k$ for some $\k\in\tup\cup\ans\La$. 
Put $\pi(\i)=(\i\res m_\i)\we B(\k)$.
\epf

\bpf[Theorem~\ref{1121} from Theorem~\ref{1152}, in $\rL$]
By Lemma~\ref{111}\ref{111d}, we can assume that 
$\tau_0\sq\dym X$. 
Then by Theorem~\ref{1152}, we assume that 
$\tau_0=\dym X$ exactly. 
Now, as $X\fo\sus x\,\vpi(x)$, we have 
(*) $X\fo\vpi(f)$ 
for some $f\in\kont\sg$, $\sg\in\cpo$.   
We can \noo\ assume that $\ta_0\sq\sg$ 
(by Lemma~\ref{115}). 

Lemma~\ref{1122} yields a permutation 
$\pi\in\ppio{\ta_0}$ such that 
$\ta=\pi\akt\sg$ is an odd expansion of  $\ta_0$ 
and  $\ta\cap\sg=\pu$. 
Note that $\pi\akt X=X$ as $\tau_0=\dym X$. 

It does not take much effort to define the action of 
$\pi$ on $\komt$. 
Namely if $\xi\in\cpo$ and $\et=\pi\akt\xi$ 
then clearly $\rat\et=\pi\akt\rat\xi$ 
in the sense of Section~\ref{perm}. 
(Note that $\rat\xi\sq\can\xi$.)  
Therefore if $f\in\kont\xi$ then we naturally define 
$g=\pi\akt f\in\kont\et$ by $g(\pi\akt x)=f(x)$ for 
all $x\in\can\xi$. 

Furthermore if $\psi$ is an $\xla$-formula then 
we let $\pi\psi$ be obtained by the
substitution of $\pi\akt f$ for any code $f\in\komt$ 
in $\psi$.
As far as the given formula $\vpi(x)$ is concerned,  
note that $\pi\vpi(x)$ is 
identic to $\vpi(x)$ since $\tau_0=\modd{\vpi(x)}$. 

\ble
[routine by induction on the complexity]
\lam{1172}
If\/ $X\in\pei$ and\/ $\psi$ is an $\xla$-formula then\/ 
$X\fo\vpi$ iff\/ $\pi\akt X\fo (\pi\vpi)$.\qed
\ele

Applying the lemma to (*), we get 
$\pi\akt X\fo\pi\vpi(g)$, where $g=\pi\akt f\in\kont\ta$. 
However $\pi\akt X=X$ and $\pi\vpi(x)$ is 
identic with $\vpi(x)$, see above. 
Thus $X\fo \vpi(g)$, as required.
\epF{Thm~\ref{1121} mod Thm~\ref{1152}}

\parf{Proof of the narrowing theorem}
\las{114}

\bpf[Theorem~\ref{1152}, in $\rL$] 
Let $Y=(X\dar\et)\uar\ta$; clearly $X\dar\et=Y\dar\et$. 
Recall that the notion of 
projection-keeping homeomorphisms, or \pkh s for 
brevity, was introduced by Definition~\ref{prelim2}. 
This will be our tool for the proof of Theorem~\ref{1152}. 
In particular, Lemma~\ref{fm9} implies the existence of 
a \pkh\ $H:X\onto Y$ such that $H(x)\res\et=x\res\et$ 
for all $x\in X.$ 
{\ubf Fix such an $H$}. 

As {\ubf the first step} of the proof, we extend the 
action of $H$ as follows.
\ben
\senu
\itlb{1**}
If $\xi\sq\ta,\,\xi\in\cpo$, then a \pkh\ 
$H_\xi:X\dar\xi\onto Y\dar\xi$ is defined by 
$H_\xi(x\dar\xi)=H(x)$ for any $x\in X$. 

\itlb{2**}
Let $\za\in\cpo$ satisfy $\ta\sq\za$. 
If $x\in X'=X\uar\za$ then 
$y=H_\za(x)\in Y'=Y\uar\za$ 
is defined by $y\dar\ta=H(x\dar\ta)$ 
(thus $y\dar\ta\in Y$) and 
$y(\i)=x(\i)$ for all $\i\in\za\bez\ta$. 
\rit{We assert that $H_\za:X'\onto Y'$ is a \pkh.} 
\een

Indeed let $\sg\in\cpo$, $\sg\sq\za$, and $u,v\in X'$ 
satisfy $u\dar\sg=v\dar\sg$. 
Then in particular $u\dar\xi=v\dar\xi$, 
where $\xi=\sg\cap\ta$, and hence, by \ref{1**},  
\bce
$H_\za(u)\dar\xi=H_{\xi}(u\dar\xi)
=H_{\xi}(v\dar\xi)=H_\za(v)\dar\xi$.
\ece
But if $\i\in\sg\bez\xi$ then 
$\i\in\za\bez\ta$, so 
$H_\za(u)(\i)=u(\i)=v(\i)=H_\za(v)(\i)$. 
Overall,  
$H_\za(u)\dar\sg=H_\za(v)\dar\sg$, as required. 

We may note that $H_\za(x)\dar\et=x\dar\et$ 
since $H$ itself has this property. 

\bdf
\lam{114s} 
If still $\ta\sq\za\in\cpo$ and $x\in X\uar\za$,  
then put $H\aqt{x}=H_\za(x)$, and define 
\kmar{aqt}%
\index{action, $H\aqt Z$}%
\index{zHtZ@$H\aqt Z$, action}%
$H\aqt Z=\ens{H\aqt x}{x\in Z}$ for any 
$Z\in\pe\za$, $Z\sq X\uar\za$. 
\ben
\nenu
\itlb{114s1}
By \ref{2**} and Lemma~\ref{997}   
the map $Z\mto H\aqt Z$ is a $\psq$-preserving and 
$\dym{...}$-preserving bijection from 
$\pe{\psq X}=\ens{Z\in\pei}{Z\psq X}$ onto 
$\pe{\psq Y}=\ens{Z\in\pei}{Z\psq Y}$. 

\itlb{114s2}
 $(H\aqt Z)\dar\et=Z\dar\et$ for 
all $Z\in \pe{\psq X}$ by the above.\qed
\een
\eDf

\vyk{
Let $\za\in\cpo$ be arbitrary and $\xi=\za\cap\ta$. 
If $x\in X'=X\dar\xi\uar\za$ then 
$y=H_\za(x)\in Y'=Y\dar\xi\uar\za$ 
is defined by $y\dar\xi=H_\xi(x\dar\xi)$ 
(thus $y\dar\xi\in Y\dar\xi$) and 
$y(\i)=x(\i)$ for all $\i\in\za\bez\xi$. 
\rit{We assert that $H_\za:X'\onto Y'$ is a \pkh.} 

Indeed let $\sg\in\cpo$, $\sg\sq\za$, and $u,v\in X'$ 
satisfy $u\dar\sg=v\dar\sg$. 
Then in particular $u\dar\sg'=v\dar\sg'$, 
where $\sg'=\sg\cap\ta=\sg\cap\xi$, and hence 
$H_\za(u)\dar\sg'=H_{\sg'}(u\dar\sg')
=H_{\sg'}(v\dar\sg')=H_\za(v)\dar\sg'$. 
On the other hand, if $\i\in\sg\bez\sg'$ then 
$\i\in\za\bez\xi$, and hence 
$H_\za(u)(\i)=u(\i)=v(\i)=H_\za(v)(\i)$. 
Thus overall 
$H_\za(u)\dar\sg=H_\za(v)\dar\sg$, as required. 

As in $1^\circ,$ if $x\in X\dar\xi\uar\za$ 
then put $H\aqt{x}=H_\za(x)$.
}

The action of $H$ on $\komt$ is somewhat less 
natural because the domain of the given $H$ is a set 
$X\in\pe\ta$, perhaps a proper subset of $\can\ta.$ 

\ble
\lam{L114}
Under the assumptions above, let\/ $\sg\in\cpo$, 
and\/ $\sg\sq\et$ or\/ $\ta\sq\sg$. 
Then for any code\/ $f\in\kont\sg$ there 
is\/ $g=H\aqt f\in\kont\sg$ satisfying$:$ 
\ben
\renu
\itlb{114a}
$g=f$ and\/ 
$\clo g(H_\sg(x))=\clo f(x)$ 
for all\/ $x\in X\dar\sg$ ---  
in case\/ $\sg\sq\et\,;$ 

\itlb{114b}
$\clo g(H_\sg(x))=\clo f(x)$ for all\/ $x\in X\uar\sg$, 
in case\/ $\ta\sq\sg$.
\een
Moreover, if\/ $h\in\kont\sg$ then there exists\/ 
$f\in\kont\sg$ such that\/ $h$ is v-equivalent to\/ 
$g=H\aqt f$ 
on\/ $Y$, that is, $\clo g(y)=\clo h(y)$ for all\/ 
$y\in Y\uar\sg$. 
\ele

\bpf
\ref{114a} 
The code $g=f$ satisfies  
$\clo g(H_\sg(x))=\clo f(x)$ for all\/ $x\in X\dar\sg$, 
because $\sg\sq\et$ and 
$H(x)\res\et=x\res\et$ for all $x\in X$.  

\ref{114b}
As  $\ta\sq\sg$, $H_\sg:X\uar\sg\onto Y\uar\sg$ 
is a \pkh, see \ref{2**} above, in particular, a 
homeomorphism.	
If $y\in Y\uar\sg$ then let 
$G'(y)=\clo f(H_\sg\obr(y))$, thus $G':Y\uar\sg\to\cN$ 
is continuous. 
It has a continuous extension $G:\can\sg\to\cN$. 
Let $g=G\res\rat\sg$, so that $G=\clo g$ and $g\in\kont\sg$. 
Thus $\clo g(H_\sg(x))=\clo f(x)$ holds
for all $x\in X\uar\sg$. 
To be more specific, we let $H\aqt f$ to be the 
G\"odel-least one 
of all $g\in\kont\sg$ with this property. 
Thus $g=H\aqt f\in\kont\sg$ is defined, 
satisfying $\clo g(H_\sg(x))=\clo f(x)$ 
for all $x\in X\uar\sg$. 

Finally to prove the `moreover' claim, note that 
$F'(x)=\clo h(H_\sg(x))$ 
is a continuous map $X\uar\sg\to\cN,$ 
extend it to a continuous $F=\clo f:\can\sg\to\cN,$ 
where $f\in\kont\sg$, and let $g=H\aqt f$.
\epF{Lemma} 

\vyk{$7^\circ.$  
To conclude, given any $\sg\in\cpo$ with 
$\sg\sq\ta\lor\ta\sq\sg$, and any $f\in\kont\sg$, 
a code $g=H\aqt f\in\kont\sg$ is defined by 
$4^\circ,5^\circ,6^\circ$ such that 
$\clo g(H_\sg(x))=\clo f(x)$ holds 
for all $x\in X\dar\sg$ (in case $\sg\sq\ta$) or 
all $x\in X\uar\sg$ (in case $\ta\sq\sg$),  
and in addition $g=f$ in case $\sg\sq\et$.
}

The next definition and lemma 
continue the proof of Theorem~\ref{1152}. 

\bdf
\lam{114f} 
If $\Phi$ is a $\xla$-formula such that any $f\in\komt$ in 
$\Phi$ satisfies $\modd f\sq\et$ or $\ta\sq\modd f$, then 
$H\Phi$ is the result of substitution of $H\aqt f$  
for any $f\in\komt$ occurring in $\Phi$. 
\edf

\ble
\lam{114L}
Let\/ $\Phi$ be a closed $\xla$-formula as in 
Definition~\ref{114f}, and\/ $Z\in\pe{\psq X}$. 
Then\/ $Z\fo\Phi$ iff $H\aqt Z\fo H\Phi$.
\ele
\bpf
The case of $\Phi$ as in 
\ref{fod1} of Definition~\ref{fod}, 
as the basis of induction, routinely follow from 
the equality $\clo g(H_\sg(x))=\clo f(x)$ of 
Lemma~\ref{L114} because $Z\psq X$.
It remains to take care of the steps \ref{fod4},\ref{fod5}.

\ref{fod4}. 
Let $\Phi$ be $\sus x\,\psi(x)$. 
Assume $Z\fo\sus x\,\psi(x)$, so that 
$Z\fo\psi(f)$ for some $f\in\kont\sg$, $\sg\in\cpo$. 
By Lemma~\ref{116}, we can assume that $\ta\sq\modd f$, 
so $\psi(f)$ is still of the form as in 
Definition~\ref{114f}. 
Then $H\aqt Z\fo H(\psi(f))$ by the inductive hypothesis, 
meaning that $H\aqt Z\fo (H\psi)(g)$, where $g=H\aqt f$, 
and hence $H\aqt Z\fo \sus x\,(H\psi)(x)$, and 
$H\aqt Z\fo H\Phi$. 

To prove the inverse, we suppose that 
$H\aqt Z\fo \sus x\,(H\psi)(x)$, that is, 
${H\aqt Z}\fo (H\psi)(h)$,  
for some $h\in\kont\sg$, $\ta\sq\sg\in\cpo$. 
By Lemma~\ref{L114}, there exists $f\in\kont\sg$ 
such that\/ $h$ is v-equivalent to $g=H\aqt f$ 
on $Y$, and hence on $H\aqt Z\psq Y$ as well. 
Then $H\aqt Z\fo (H\psi)(g)$ by Lemma~\ref{115}, 
and hence $Z\fo\psi(f)$ by the inductive hypothesis, 
and 
$Z\fo\Phi$, as required.

\ref{fod5}. 
Let $\Phi$ be $\otr\psi$, where $\psi$ is a $\ls1n$ formula. 
Assume that  $Z\fo\Phi$ {\ubf fails}. 
By definition there is a condition $Z'\psq Z$, $Z'\fo\psi$. 
The inductive hypothesis implies $H\aqt Z'\fo H\psi$. 
However $H\aqt Z'\psq H\aqt Z'$, hence we conclude 
that $H\aqt Z\fo\Phi$ {\ubf fails}. 
The converse is similar.
\epF{Lemma}

Now we 
return to the formula $\vpi$ of Theorem~\ref{1152}. 
It satisfies $\modd\vpi\sq\et$, and $X\fo\vpi$. 
Lemma~\ref{114L} is applicable, so that 
$Y\fo \vpi$, because $H\vpi$ is identic to $\vpi$ 
since $\modd\vpi\sq\et$. 
This implies $X\dar\et\fo \vpi$ by 
Lemma~\ref{111}\ref{111d}.\vom

\epF{Theorems~\ref{1152} and \ref{1121}}

\bcor
\lam{1153}
Let\/ $X\in\pei$, $k<\om$, $\vpi$ is a closed\/ $\xla$ 
formula, $\et=\dym X\cup\modd\vpi$,   
$\neg\:X\fo\vpi$. 
Then there is\/ $Z\in\pe\et$, $Z\psq X$, 
$Z\fo\otr\vpi$.
\ecor
\bpf
Lemma~\ref{111}\ref{111c},\ref{111c} yields $Y\in\pei$ 
such that $\et\sq\dym Y$, $Y\psq X$, and $Y\fo \otr\vpi$.
Now let $Z=Y\dar\et$ and apply Theorem~\ref{1152}.
\epf

\parf{Complete forcing notions and 
3rd form of the main theorem}
\las{12}

After working out some technical issues with $\fo$, 
we'll prove the 
truth theorem for this forcing-type relation. 
It is based on the next definition.

\bdf
[in $\rL$]
\lam{121}
A normal forcing notion $\cX\sq\pei$ is $n$-\rit{complete} 
if for any closed formula $\vpi$ in $\bigcup_{k\le n}\ls1k$ 
the set 
$$
\Fo\vpi=\ens{X\in\cX}{X\fo\vpi\text{ or }X\fo\otr\vpi}
$$
is dense in $\cX$.
\edf

For example, the set $\pei$ is $n$-complete for each $n$ 
by Lemma~\ref{111}, \ref{111c1} and \ref{111c}. 
We will not use this fact, 
but it is useful to keep it in mind. 
In its light, $n$-complete normal forcing 
notions $\cX\sq\pei$ 
can be viewed as 
``similar to $\pei$ up to level $n$ in the sense of $\fo$''.
Let us now prove the theorem connecting $\fo$  and truth
in generic extensions.

\bte
[truth theorem]
\lam{122}
Assume that $n\ge1$ and, in\/ $\rL$, a normal forcing\/ 
$\cX\sq\pei$ is\/ $n$-complete and has the Fusion property. 
Let\/ $\w$ be an\/ $\cX$-generic array over\/ $\rL$. 
Let\/ $\vpi$ be a closed formula in\/ $\ls1k$, $k\le n+1$. 
Then\/ $\rL[\w]\mo\vva\vpi\w$ iff there exists a condition\/ 
$X\in\cX\cap\cg\w$, $X\fo\vpi$.
\ete

\bpf
We argue by induction on $k\le n + 1$. 
\rit{Starting with\/ $k=1$}, 
suppose that $\vpi$ is a $\ls11$ formula. 
By the $n$-Completeness of $\cX$ and the genericity of $\w$, 
there exists a condition $X\in\cX\cap\cg\w$, $X\fo\vpi$ or 
$X\fo\otr\vpi$.
Assume that $X\fo\vpi$. 
This claim can be naturally converted into a $\up12$ 
sentence with parameters in $\rL$, true in $\rL$. 
Then $\rL[\w]\mo\vva\vpi\w$ by the Shoenfield absoluteness.
Similarly, if $X\fo\otr\vpi$ (a $\up11$ sentence) then 
$\rL[\w]\mo\vva{\otr\vpi}\w$, by the same 
absoluteness argument.\vom

\rit{Step\/ $k\to k+1$}. 
Suppose that $\vpi:=\sus x\,\otr\psi(x)$ is a $\ls1{k+1}$ 
formula, $\psi(x)$ being $\ls1k$, and $k\le n$. 

\rit{Direction $\mpi$}. 
Assume that $\rL[\w]\mo\vva{\vpi}\w$, that is,  
$\rL[\w]\mo\vva{\otr\psi}\w(p)$, for a suitable real
$p\in\cN\cap\rL[\w]$. 
Then $p=\clo f(\w\dar\xi)$ for some $f\in\kont\xi$, 
$\xi\in\cpo$, 
by Corollary~\ref{1011}. 
Thus $\rL[\w]\mo\vva{\otr\psi(f)}\w$, and hence, by the 
inductive hypothesis, no condition $X\in\cX\cap\cg\w$ 
satisfies $X\fo\psi(f)$. 
We conclude by the $n$-Completeness that there is a 
condition $X\in\cX\cap\cg\w$ with $X\fo\otr\psi(f)$, 
and then $X\fo\vpi$ by \ref{fod4} of Definition~\ref{fod}.  

\rit{Direction $\imp$}. 
Assume that $X\fo\vpi$, that is, $X\fo\otr\psi(f)$ 
for some $f\in\kont\xi$, $\xi\in\cpo$, still 
by \ref{fod4} of Definition~\ref{fod}. 
Then no condition $X\in\cX\cap\cg\w$ 
satisfies $X\fo\psi(f)$. 
Hence, by the inductive hypothesis, 
$\rL[\w]\mo\neg\vva{\psi}\w(p)$, where  
$p=\clo f(\w\dar\xi)\in\cN\cap\rL[\w]$. 
We conclude that $\rL[\w]\mo\vva{\vpi}\w$, as required. 
\epf

Now we apply the truth theorem just proved, to show 
that the Completeness of a normal forcing implies 
the Odd-Expansion property, via 
the odd expansion theorem (Theorem~\ref{1121}).

\bte
[in $\rL$]
\lam{123}
Assume that $n\ge1$ and a normal forcing\/ 
$\cX\sq\pei$ is\/ $n$-complete and has the Fusion property. 
Then\/ $\cX$ has the\/ $n$-Odd-Expansion property 
of Definition~\ref{92}.
\ete

\bpf
Let\/ $\w$ be an\/ $\cX$-generic array over\/ $\rL$. 
Suppose that $\et\in\cpo$ and $\vpi(\cdot)$ is 
a $\ip1n$ formula,
with reals in $\rL[\w\dar\et]$ as parameters, and  
$\rL[\w]\mo\sus x\,\vpi(x)$. 
We have to find an odd expansion $\ta\in\cpo$ of $\et$, 
and some $q\in\rL[\w\dar\ta]$, 
such that $\rL[\w]\mo\vpi(q)$. 

If $p\in\cN\cap\rL[\w\dar\et]$ occurs in $\vpi$ then 
Corollary~\ref{1011} yields a code $f_p\in\kont\et$ such 
that $p=\clo{f_p}(\w\dar\et)$.  
Change each $p$ to $f_p$ in $\vpi(\cdot)$, 
and let $\psi(\cdot)$ be the $\xla$-formula obtained. 
Then $\vpi(\cdot)$ is identic to $\vva{\psi(\cdot)}\w$ 
and $\modd{\psi}=\et$. 

By Theorem~\ref{122}, there is a condition 
$X\in \cg \w\cap\cX$ satisfying $X\fo \sus x\,\psi(x)$.
Then by Theorem~\ref{1121} there is an odd expansion 
$\ta\in\cpo$ of\/ $\et$, and\/ $g\in\kont\ta$, such 
that\/ $X\fo\psi(g)$. 
Then $\rL[\w]\mo \vva{\psi(g)}\w$, that is, 
$\rL[\w]\mo \vpi(q)$, 
where $q=\clo g(v\dar\ta)\in\cN\cap\rL[\w\dar\ta]$, 
as required.
\epf

This theorem will allow us to replace the $\nn$-Odd-Expansion 
condition in Theorem~\ref{mt1b} by the $\nn$-Completeness of 
$\cX$ in $\rL$.

\bte
[in $\rL$]
\lam{mt1c}
If\/ $\nn\ge1$  
then there is a forcing\/ 
$\cX\in\RF$ with the Fusion, Structure, 
\dd\nn Definability, and\/ \dd\nn Completeness  
properties.
\ete

\bpf
[Theorems~\ref{mt1} and \ref{mt1b}  from Thm~\ref{mt1c}]
Apply Theorem~\ref{123}.
\epf

Thus Theorem~\ref{mt1c} implies Theorem~\ref{mt1}, 
the first main result of this paper. 
Chapters \ref{foax}--\ref{VII} below will contain the 
proof of Theorem~\ref{mt1c}, via  Theorem~\ref{mt1d}
as the next approximation,  
and thereby will 
accomplish the proof of Theorem~\ref{mt1}.



\sekt
[\ \ The construction of the final forcing begins. Rudiments]
{The construction of the final forcing begins} 
\las{sek6}

The purpose of Chapters \ref{sek6}--\ref{VII} is to 
define a normal forcing $\cX\in\rL$ satisfying 
requirements of Theorem~\ref{mt1c}. 
This will be a rather difficult task. 

As mentioned in the end of Section~\ref{ker}, in principle 
it suffices to first define an auxiliary \dd{\tud}kernel $\cK$ 
and then let $\cX=\noc{\cK\iex}$ by Lemmas \ref{22w} 
and \ref{k2r}. 
Unfortunately it does not seem to work that simple way. 
Instead, 
following \cite{vin79}, 
we'll make use of a kind of {\ubf limit} of an $\omi$-sequence 
of countable collections of iterated perfect sets, 
called {\ubf rudiments}. 
This construction realizes the idea of generalized \dd\tup 
iteration of Jensen's forcing somewhat differently than in 
\cite{jml19,gitpf,ww}, in particular, the CCC property will 
not be achieved.  

As for this chapter, 
Rudiments, rudiment hulls, and related notions are studied in 
Sections \ref{rud} and \ref{134*}.  
Then we introduce an important {\ubf refinement} 
relation $\ssq$ between rudiments. 
Basically, $\rU\ssq\rV$ will imply that the rudiment hull 
$\rh(\rU\cup\rV)$ has a rather transparent structure in 
terms of $\rU$ and $\rV$, at least \rit{locally}, \ie, 
in the context of projections $\rsq\i$.
We finally study {\ubf rudimentary sequences}. that is, 
transfinite sequences of rudiments increasing in the 
sense of $\ssq$ in Section~\ref{rs}.

{\ubf We argue in $\rL$ in this Chapter.}

\parf{Rudiments}
\las{rud}

Planning to maintain a construction of normal forcing 
notions in the form $\cX=\noc{\bigcup_{\al<\omi}\rP_\al}$, 
where each $\rP_\al$ is countable, we may note that the 
summands $\rP_\al$ cannot be normal forcing 
notions themselves, because each of conditions 
\ref{rfo3}, 
\ref{rfo5}, 
\ref{rfo6} of Section \ref{rfo} implies the uncountability 
of any normal forcing. 
Thus we have to somehow reduce the generality of those 
conditions. 
This is the content of this section.
We begin with two auxiliary notes.

First, suppose that $\et\sq\xi$ belong to $\cpo$. 
Say that $\et$ is a \rit{finite-type} in $\xi$, 
in symbol $\et\in\ft\xi$, 
\index{finite-type, $\ft\xi$}%
\index{zFTxi@$\ft\xi$}%
\kmar{ft xi}%
if $\et$ is obtained from sets of the form 
\bce
$\xi$ itself, \ \ \ 
$\ilq\i=\ens{\j\in\tup}{\j\sq\i}$, \ \ \ and \ \ \ 
$\xi\cap\tuq\al$,  
\ece
where 
$\al<\omi$, $\i\in\xi$,  and  
$\tuq\al=\ens{\i\in\tup}{\ran\i\sq\al}=\al\lom\bez\ans\La$, 
by a finite number of operations of 
set difference $\bez$ and 
(finite) $\cup$ and $\cap$. 
Clearly $\ft\xi$ is a countable or finite Boolean algebra.  

Second, if $\i\ekp\j$ belong to $\tup$, then there exists 
a canonical permutation $\pi_{\i\j}\in\per$ satisfying 
$\pi_{\i\j}(\i)=\j$ and $\pi_{\i\j}=\pi_{\i\j}\obr$, see 
Example~\ref{pi*ij}. 

\bdf
\lam{132}
\label{ruD}
Let $\al<\omi$.
A set $\rP$ is a {\em rudiment}
\index{rudiment, $\bfr\al$}%
%
{\em of width\/ $\al$},
in symbol $\rP\in\bfr\al$, 
\kmar{bfr al}%
\index{zRuda@$\bfr\al$}%
if $\rP$ 
satisfies the following conditions \ref{fr1}--\ref{nfr4}.
\ben
\stenu
\itlb{fr1}
$\pu\ne\rP\sq\pe{\tuq\al}$, where, we recall,  
$\tuq\al=\ens{\i\in\tup}{\ran\i\sq\al}$.


\vyk{
\pur
\itlb{fr3} 
If $\i\su\j$ belong to $\tuq\al$ then
$(\pro\rP\j)\rsdq\i=\pro\rP\i$, 
that is, 

(a) if $X\in \pro\rP\j$ then $X\rsq{\i}\in\pro\rP\i$, 
and 

(b) if $Y\in \pro\rP\i$ then there 
exists $X\in\pro\rP\j$ satisfying $Y=X\rsq{\i}$.
}

\itlb{nfr3}
If $\et\in\cpo,\,\et\sq\tuq\al$ is finite-type in 
$\tuq\al$, $X,Y\in\rP$,
and $Y\dar\et\sq X\dar\et$, then the 
set $X\cap (Y\dar\et\uai\al)$ belongs to $\rP$.

\itlb{fr2}
If $X\in \rP$, $Y\in\pe{\tuq\al}$, $Y\sq X$,  
$Y$ is clopen in $X$, then $Y\in \rP$.

\itlb{nfr4}
Invariance: if $\i,\j\in\tuq\al$, $\i\ekp\j$,
and $X\in\rP$,
then $\pi_{\i\j}\akt X\in\rP$. 
\een
If $\rP$ is such, and $\et\in\cpo,\,\et\sq\tuq\al$, 
then we let 
$\rP\dar\et=\ens{X\dar\et}{X\in\rP}$. 
In particular, if $\i\in\tuq\al$ then put 
\kmar{pro rP i}%
$\pro\rP\i=\ens{X\rsq\i}{X\in\rP}$. 
\index{zPi@$\pro\rP\i$}%
\edf

\vyk{\gol
A rudiment  $\rP$ is \rit{regular}, resp., 
\rit{semi-regular}, if in addition it 
satisfies resp.\ (A),(B):

\ben
\atc
\atc
\atc
\atc
\cenu
\itlb{fr5}
(A) If $X,Y\in \pro\rP\i$, $\i\in\tuq\al$, then $Y\cap X$ 
is clopen in $X$ and in $Y$, 

(B) If $X,Y\in\pro\rP\i$, $\i\in\tuq\al$, then $Y\cap X$ 
is clopen in $X$ and/or in $Y$. 

In both cases (a),(b), if $X\cap Y\ne\pu$ then 
$X\cap Y\in\pro\rP\i$ by \ref{fr2}.
\een
}%

Thus if $\cX\in\RF$ then 
$\cX\dar{\tuq\al}=\ens{X\dar{\tuq\al}}{X\in\cX}
\in\bfr\al$.%
\pagebreak[1]

The set 
$\pe{\tuq\al}$ belongs to $\bfr\al$ by Lemmas
\ref{less}, \ref{apro}, \ref{lin}. 
The set of all clopen sets $X\in\pe{\tuq\al}$
%
belongs to $\bfr\al$, too. 

The following lemma clarifies the connections between 
kernels, rudiments, and normal forcings.

\ble
\lam{r2n}
Assume that\/ $\rP\in\bfr\al$, $2\le\al<\omi$, 
$\cai\al\in\rP$, 
Then\/ $\Ker\rP=\sis{\pro\rP\i}{\i\in\tuq\al}$ 
\kmar{Ker rP}%
\index{kernel!$\Ker\rP$}%
\index{zKerP@$\Ker\rP$, kernel}%
is a strong $\tuq\al$-kernel, 
$\cX=\noc\rP\in\nf$, 
and\/ $\pro\cX\i=\pro\rP\i$ for all\/ $\i\in\tuq\al$.

Conversely, if\/ $\cK=\sis{\cK_\i}{\i\in\tuq\al}$ is 
an\/ $\tuq\al$-kernel, then the set\/ 
\bce
$\rP=\rP(\cK):=\ens{X\in\tuq\al}
{\kaz\i\in\tuq\al\,(X\rsq\i\in\cK_\i)}$
\index{rudiment!$\rP(\cK)$}%
\index{zPK@$\rP(\cK)$, rudiment}%
\ece
belongs to\/ $\bfr\al$, and\/ $\pro\rP\i=\cK_\i$ 
for all\/ $\i\in\tuq\al$. 
\ele
\bpf
Recall the notion of kernel in Section~\ref{ker}. 
Conditions \ref{ke1s}, \ref{ke2} of Section~\ref{ker}  
for $\Ker\rP$ are clear, 
and \ref{ke6} holds by 
\ref{nfr4} of Definition~\ref{ruD} for $\rP.$ 

To verify \ref{ke3} of Section~\ref{ker} 
for $\Ker\rP$, 
let $\j\su\i$ belong to $\xi=\tuq\al$,   
$X\in\pro\rP\i$, $Y\in\pro\rP\j$, and 
$Y\sq X\rsq\j$. 
Check $Z=X\cap(Y\uar^{\sq\i})\in\pro\rP\i$.
By definition, $Y=Y'\rsq\j$ and $X=X'\rsq\j$ 
for some $X',Y'\in\rP$. 
And we have $Y'\rsq\j=Y\sq X'\rsq\j$. 
Therefore the set 
\bce
$Z'=X'\cap(Y'\rsq\j\uar\xi)=X'\cap(Y\uar\xi)$ 
\ece
belongs to $\rP$ by \ref{nfr3}. 
Then
${Z'}\rsq\i=(X'\rsq\i)\cap(Y\uar^{\sq\i})
=X\cap (Y\uar^{\sq\i})=Z$, 
hence $Z\in\pro\rP\i$,
as required. 
 
To check \ref{ke4} of Section~\ref{ker}, 
assume that $\i\in\tuq\al$, 
$X\in\pro\rP\i$, 
a set $\pu\ne Y\sq X$ is clopen in\/ $X$, and prove that  
$Y\in\pro\rP\i$. 
We have $Y\in\pele\i$ by Lemma~\ref{lin}. 
By definition, $X=X'\rsq\i$ for some $X'\in\rP$. 
It follows by Lemma~\ref{apro} that the set 
$Y'=X'\cap (Y\uar\xi)$ belongs to $\pe{\xi}$, 
and $Y'$ is clopen in $X'$ by the choice of $W$. 
Therefore $Y'\in\rP$ by \ref{fr2} 
of Definition~\ref{ruD}. 
Hence $Y=Y'\rsq\i\in\pro\rP\i$, as required.  

Thus indeed $\rK=\Ker\rP$ is a strong $\tuq\al$-kernel. 
Then the expanded system\/ $\rK\iex$ is a strong 
$\tup$-kernel by Lemma~\ref{22w}. 
It follows by Lemma~\ref{k2r} that $\cZ=\noc{\rK\iex}$ 
is a normal forcing with 
$\pro{\cZ}\i={\rK\iex}_\i=\pro\rP\i$ 
for all $\i\in\tup$   
and accordingly 
$\pro{\cZ}\i=\rK_\i=\pro\rP\i$ for all $\i\in\tuq\al$.  
Therefore $\rP\sq\cZ$ by \ref{rfo5} of 
Section~\ref{rfo} for $\cZ$,  
hence $\cX\sq\cZ$ by the minimality of $\cX$.

We similarly get the inverse inclusion $\cZ\sq\cX$ 
by the minimality of $\cZ$. 
We conclude that $\cX=\cZ$, and hence the equality 
$\pro{\cX}\i=\pro\rP\i$ holds for all 
$\i\in\tuq\al$ by the above.  

The proof of the converse claim goes pretty similar to
the proof of Lemma~\ref{k2r}, and hence we leave 
the details for the reader.
\epf 

\vyk{
\ble
\lam{rudL}
Assume that\/ $\rP\in\bfr\al,\,2\le\al<\omi$. 
Then\/ $\Ker\rP=\sis{\pro\rP\i}{\i\in\tuq\al}$ 
is an\/ $\tuq\al$-kernel. 
Conversely, if\/ $\cK=\sis{\cK_\i}{\i\in\tuq\al}$ is 
an\/ $\tuq\al$-kernel, then\/ 
$\rP(\cK)=\ens{X\in\tuq\al}
{\kaz\i\in\tuq\al(X\rsq\i\in\cK_\i)}\in \bfr\al$. 
\ele

\bpf
Put $\xi=\tuq\al$.
To verify \ref{ke3} of Section~\ref{ker} 
for $\cK=\Ker\rP$, 
let $\j\su\i$ belong to $\xi$,   
$X\in\pro\rP\i$, $Y\in\pro\rP\j$, and 
$Y\sq X\rsq\j$. 
Check $Z=X\cap(Y\uar^{\sq\i})\in\pro\rP\i$.
By definition, $Y=Y'\rsq\j$ and $X=X'\rsq\j$ 
for some $X',Y'\in\rP$. 
And we have $Y'\rsq\j=Y\sq X'\rsq\j$. 
Therefore the set 
\bce
$Z'=X'\cap(Y'\rsq\j\uar\xi)=X'\cap(Y\uar\xi)$ 
\ece
belongs to $\rP$ by \ref{nfr3}. 
Then
${Z'}\rsq\i=(X'\rsq\i)\cap(Y\uar^{\sq\i})
=X\cap (Y\uar^{\sq\i})=Z$, 
hence $Z\in\pro\rP\i$,
as required. 
 
To check \ref{ke4} of Section~\ref{ker}, 
assume that $\i\in\tuq\al$, 
$X\in\pro\rP\i$, 
$\pu\ne Y\sq X$ is clopen in\/ $X$, and prove that  
$Y\in\pro\rP\i$. 
We have $Y\in\pele\i$ by Lemma~\ref{lin}. 
By definition, $X=X'\rsq\i$ for some $X'\in\rP$. 
It follows by Lemma~\ref{apro} that the set 
$Y'=X'\cap (Y\uar\xi)$ belongs to $\pe{\xi}$, 
and $Y'$ is clopen in $X'$ by the choice of $W$. 
It follows that $Y'\in\rP$ by \ref{fr2} 
of Definition~\ref{ruD}. 
Therefore $Y=Y'\rsq\i\in\pro\rP\i$, as required.  

\ref{ke1s}, \ref{ke2} for $\Ker\rP$ are clear, 
and \ref{ke6} holds by 
\ref{nfr4} of Definition~\ref{ruD} for $\rP.$ 

See Lemma~\ref{k2r} for the proof of the converse claim.
\epf
}

\vyk{
Thus rudiments$\to$kernels.
We also know that kernels$\to$normal forcings by 
Lemma~\ref{k2r}. 
The next corollary combines these two results. 
}

\vyk{
\bcor
\lam{r2r}
Assume that\/ $\rP\in\bfr\al,\,2\le\al<\omi$, 
and\/ $\cai\al\in\rP$. 
Let\/ $\cK=\Ker\rP$, and let\/ 
$\cK\iex$ be the expanded\/ $\tup$-kernel 
defined as in Lemma~\ref{22w}.  
Then the set\/ $\SC\rP:=\cX(\cK')$ is a normal 
forcing, and\/ $\pro{\SC\rP}\j=\pro\rP\j$ for all\/ 
$\j\in\tuq\al$. 
\ecor

\bpf
First of all, $\cK=\cK(\rP)$ is an\/ $\tuq\al$-kernel 
by Lemma~\ref{rudL}. 
Moreover, $\cK$ is a \rit{strong} kernel, that is, 
$\can{\ilq\j}\in\cK_\i=\pro\rP\i$ for all $\i\in\tuq\al$ 
because $\cai\al\in\rP$. 
Then the expanded system\/ $\cK'$ is a strong 
$\tup$-kernel by Lemma~\ref{22w}.  
It follows by Lemma~\ref{k2r} that $\cX=\cX(\cK')$ 
is a normal forcing. 
The equality $\pro{\cX}\j=\pro\rP\j$ for all 
$\j\in\tuq\al$ holds by construction.
\epf

}

\vyk{\gol
\bcor
\lam{lfutp}
Let\/ $2\le\al<\om$, $\rP\in\bfr\al$, 
$\cX=\xt\rP$, 
a map\/ $\phi:\om\to\za$ be\/ \dd{\tuq\al}admissible, 
sets\/ $\cD_n\sq\cX$ be\/ 
\dd{\psq}open dense in\/ $\cX,$ and\/ $X_0\in\cX.$ 

Then there is an increasing sequence\/ 
$\xi_0\sq\xi_1\sq\xi_2\sq\ldots\sq\tuq\al$ 
of finite sets\/ $\xi_n\in\cpo$, and a system\/ 
$\sis{X_u}{u\in \bse}$ of sets\/ $X_u$ 
satisfying\/ \ref{lfut1}, \ref{lfut2}, \ref{lfut3} 
of Theorem~\ref{lfut}, 
such that\/ $X_{u}\in \cX_n$ for all\/ $m$ and\/ 
$u\in2^n,$ and\/ $X_\La\sq X_0$.
\ecor

\bpf

\epf

\bdf
[\dd\rer extensions]
\lam{piext}
Let $2\le\al\le\omi$.
Any $\rP\in\bfr2$ admits a \dd\qer extensions  
\index{extension!pipar extension@\dd\qer extension $\xt\rP$}%
\index{pipar extension@\dd\qer extension $\xt\rP$}%
\index{zzPiparextension@\dd\qer extension $\xt\rP$}%
\index{zzP_@$\xt\rP$}%
$\pxt\rP\in\bfr\omi$ and $\pyt\rP\al\in\bfr\al$
defined 
\imaf{xt P}
as follows. 

Let $\i\in\tup\bez\tuq2$. 
There is unique $\j\in\tuq2$ satisfying $\i\ekp\j$; 
thus $\pro\rP\j$ is defined by \ref{fr1} above. 
Let $\pi'_{\j\i}\in\rer$ be the parity-preserving
permutation defined as above, in particular 
$\pi_{\j\i}(\j)=\i$. 
Put
\bce
$\pro{(\pxt\rP)}\i=\pi'_{\j\i}\pro\rP\j=
\ens{\pi'_{\j\i}X}{X\in\pro\rP\j}$.
\ece
%
Taking only those $\i$ which belong to $\tuq\al$, 
we get $\pyt\rP\al$.
\edf
 
\ble
[routine]
\lam{xt}
If\/ $\rP\in\bfr2,\,2\le\al<\omi,$  
then\/ $\pxt\rP\in\bfr\omi$ and 
$\pyt\rP\al\in\bfr\al$. 
Moreover\/ $\pxt\rP$ is\/ \dd\qer invariant, 
and\/ $\pyt\rP\al$ is\/ \dd{\qea\al}invariant. 

Accordingly the assembling hulls\/ 
$\xt{\pxt\rP}$ and\/ $\xt{\pyt\rP\al}$ are resp.\ 
\dd\qer invariant and\/ \dd{\qea\al}invariant.\qed
\ele 

}

\parf{Hulls, liftings and restrictions of 
rudiments}
\las{134*}

For any $\al<\omi$, if $\pu\ne\rU\sq\peq\al$ then 
there exists 
a least set $\rP\in\bfr\al$ with $\rU\sq\rP$. 
This $\rP$ will be denoted by $\rh(\rU)$, 
\kmar{rh rU}%
\index{rudimentary hull, $\rh(\rU)$}%
\index{zRHU@$\rh(\rU)$, rudimentary hull}%
the \rit{rudimentary hull} of $\rU$. 
Note that the number of finite-type sets 
$\et\sq\tuq\al$ is countable, 
and so is the number of clopen subsets. 
Therefore we have the following lemma:

\ble
\lam{rcL}
If\/ $\al<\omi$ and\/ $\pu\ne\rU\sq\peq{\al}$ 
is countable then\/ $\RC(\rU)$ is countable as well.\qed
\ele

Several next lemmas study   
\rit{liftings} of rudiments to bigger domains. 
Recall that if $\ga<\al<\omi$ and $\rP\sq\peq{\ga}$ 
then $\rP\uai\al=\ens{X\uai\al}{X\in\rP}$, where 
$X\uai\al\in\peq{\al}$ (lifting) 
is defined as in Section~\ref{prelim1}. 
If $\rP\in\bfr\ga$ then $\rP\uai\al$ is not a rudiment, 
but $\RC(\rP\uai\al)\in\bfr\al$, of course. 
It is not that easy to clearly describe the structure 
of $\RC(\rP\uai\al)$. 
Yet the next lemma at least claims that small 
projections do not change. 

\ble
\lam{134*0}
Assume that\/ $2\le\ga<\al<\omi$ 
and\/ $\rP\in\bfr\ga$. 
Let\/ $\rR=\RC(\rP\uai\al)$.
Then\/ $\rR\rsq\i=\rP\rsq\i$ for all\/ 
$\i\in\tuq\ga$.
\ele

\bpf
If $\i\in\tuq\al$ then 
let $\kn\i\in\tuq2$ be the only tuple in $\tuq2$ 
with $\i\ekp{\kn\i}$. 
Put $\rK_\i=\pi_{\i,\kn\i}\akt \pro\rP\i$. 
The system $\sis{\pro\rP\i}{\i\in\tuq\ga}$ is 
an\/ $\tuq\ga$-kernel 
(see the proof of Lemma~\ref{r2n}). 
It easily follows by \ref{nfr4} of Definition~\ref{132} 
that $\sis{\rK_\i}{\i\in\tuq\al}$ is an 
$\tuq\al$-kernel, and (*) $\rK_\i=\pro\rP\i$ for all 
$\i$ in the old domain $\tuq\ga$. 
Then 
$\rQ=\ens{X\in\tuq\al}
{\kaz\i\in\tuq\al(X\rsq\i\in\rK_\i)}\in\bfr\al$. 
Therefore $\rR\sq\rQ$. 
But $\pro\rQ\i=\rK_\i=\pro\rP\i$ for all $\i\in\tuq\ga$ 
by (*).
\epf

\ble
\lam{rres}
If\/ $\ga<\al<\omi$ and\/ $\rU\in\bfr\al$ 
then the set\/ 
$\rU\dai\ga=\ens{X\dai\ga}{X\in\rU}$ 
belongs to\/ $\bfr\ga$. 
\ele

\bpf
To check \ref{nfr3}  of Definition~\ref{132} 
for $\rU\dai\ga$, suppose that 
$X'=X\dai\ga$,  
$Y'=Y\dai\ga$,  
where $X,Y\in\rU$, and $\et\in\ft{\tuq\ga}$, 
$Y'\dar\et\sq X'\dar\et$. 
We have to prove that $Z'=X'\cap(Y'\dar\et\uai\ga)$ 
belongs to $\rU\dar\tuq\ga$. 
Note that $\et\in\ft{\tuq\al}$ as well because 
${\tuq\ga}$ itself belongs to $\ft{\tuq\al}$. 
It follows that $Z=X\cap(Y\dar\et\uai\al)$ 
belongs to $\rU$. 
However easily $Z'=Z\dai\ga$. 

Conditions \ref{fr2} and \ref{nfr4}  
are verified by similar routine arguments.
\epf

\bcor
\lam{rrec}
If\/ $\ga<\al<\omi$ and\/ 
$\cai{\ga}\in\cX\sq\peq{\ga}$, $\rP=\RC(\cX)$,  
then the sets\/ 
$\rQ'=\RC(\cX\uai\al)$ and\/ 
$\rQ=\RC(\rP\uai\al)$ coincide. 
\ecor

\bpf
Clearly $\rQ'\sq\rQ$. 
To prove the converse, 
note that $\rP'=\rQ'\dai\ga\in \bfr\ga$ 
by Lemma~\ref{rres}, and obviously $\cX\sq\rP'$. 
Therefore $\rP\sq\rP'$. 
On the other hand, $\rP'\uai\al\sq\rQ'$ because if 
$Y\in\rQ'$ and $X=Y\dai\ga\in\rQ'$ then 
$X\uai\al=\cai\al\cap Y\dai\ga\uai\al\in \rQ'$. 
(Note that $\cai\al\in\rQ'$ since $\cai\ga\in\cX$.) 
To conclude, 
$\rQ=\RC(\rP\uai\al)\sq \RC(\rP'\uai\al)\sq \RC(\rQ')
=\rQ'$. 
\epf

\ble
\lam{13.2.1}
Assume that\/ $\la<\omi$ is limit, $\rPi\ga\in\bfr\ga$ 
for all\/ $\ga<\la$, and\/ $\rPi\ga\uai\al \sq\rPi\al$ 
for all\/ $\ga<\al<\la$. 
Then\/ $\rP=\bigcup_{\ga<\la}(\rPi\ga\uai\la)\in\bfr\la$.
\ele

\bpf
$\rP\sq\peq{\la}$ holds by Lemma~\ref{apro}. 

We check \ref{fr2} of Definition~\ref{ruD}. 
Let $Y\in\peq{\la}$, $Y\sq X\in \rP$,  
$Y$ be clopen in $X$; prove $Y\in \rP$. 
By compactness, any clopen set is a finite union 
of basic clopen sets, hence there is $\ga<\la$ 
such that $X=X'\uai\la$ and 
$Y=Y'\uai\la$, where 
$X'=X\dai\ga\in\rPi\ga$ and $Y'=Y\dai\ga$. 
However $Y'\in\peq{\ga}$ by Lemma~\ref{less} and 
$Y'$ is clopen in $X'$ by  Lemma~\ref{darop}. 
Thus $Y'\in\rPi\ga$ by \ref{fr2} 
of Definition~\ref{ruD} for $\rPi\ga$. 
Therefore $Y=Y'\uai\la\in \rP$.

We check \ref{nfr3}. 
Assume that $\et\in\cpo$, $\et\sq\tuq\la$ is 
finite-type in $\tuq\la$, $X,Y\in\rP$,
and $Y\dar\et\sq X\dar\et$; prove that the  
set $Z=X\cap (Y\dar\et\uar {\tuq\la})$ belongs to $\rP$.
As above, there is $\ga<\la$ such that $X=X'\uai\la$ 
and $Y=Y'\uai\la$, where 
$X'=X\dai\ga$, $Y'=Y\dai\ga$, and $X',Y'\in\rPi\ga$. 
Further, $\et'=\et\cap\tuq\ga\in\cpo$ and $\et'$ is 
of finite-type in $\tuq\ga$, and clearly 
$Y'\dar\et'=Y\dar\et'\sq X'\dar\et'$. 
It follows by \ref{nfr3} for $\rPi\ga$ that the set 
$Z'=X'\cap (Y'\dar\et'\uai\ga)$ belongs to 
$\rPi\ga$. 
On the other hand, 
$Z\dai\ga=(X\dai\ga)\cap 
(Y\dar{\et'}\uai\ga)$ 
by Lemma~\ref{599}, 
so that $Z\dai\ga=Z'\in\rPi\ga$. 
Therefore $Z=Z'\uai\la\in \rP$.

\ref{nfr4}. 
Take $\i\ekp\j$ in $\tuq\la$,
and $X\in\rP$; show that $Y=\pi_{\i\j}\akt X\in\rP$. 
By construction, there is an index  $\ga<\la$ such that 
$\i,\j\in\tuq\ga$, and 
$X=X'\uai\la$, where $X'=X\dai\ga\in\rPi\ga$. 
Then $Y'=\pi_{\i\j}\akt X'\in \rPi\ga$
by \ref{nfr4} for $\rPi\ga$, and on the other hand  
easily $Y=Y'\uai\la\in\rP$, as required.
\epf

\parf{Refining rudimentary forcings}
\las{ref}

\bdf
[refinement]
\lam{134d}
Let $\rP,\rQ\in\bfr\al$, $\xi=\tuq\al$. 
Say that $\rQ$ is a {\em refinement\/} of $\rP$, 
\index{refinement $\rP\ssw\rQ$}%
\index{zPsqsQ@$\rP\ssq\rQ$}%
\index{zzzsqs@$\ssq$, refinement}%
in symbol $\rP\ssw\rQ$, 
\kmar{ssw}%
if the next three conditions hold:%
\ben
\atc
\atc
\atc
\atc
\stenu
{
\itlb{134i}
$\can\xi\in\rP$.
}%

\itlb{134ii}
\label{ref2x}%
If $\et\in\ft\xi$, $X\in\rP$,
$Y\in\rQ$, 
$Y\dar\et\sq X\dar\et$, then there is 
$Z\in\rQ$ such that 
$Z\sq X$ and $Z\dar\et=Y\dar\et$ --- in particular 
($\et=\pu$) if $X\in\rP$ then there is 
$Z\in\rQ$ such that $Z\sq X$.

\itlb{134iii}
\label{ref4}%
If $\i\in\xi$, $X\in\pro\rP\i$, $Y\in\pro\rQ\i$, 
then $X\cap Y$ is {\em clopen\/} in $Y$, hence if in
addition $X\cap Y\ne\pu$ then 
$X\cap Y\in\pro\rQ\i$ by \ref{fr2} of 
Definition~\ref{ruD}.\qed
\een
\eDf

\vyk{\gol 
We write $\rP\ssu\rQ$ (strict refinement)
\imaf{ssu}
\index{refinement strict $\rP\ssu\rQ$}%
\index{strict refinement $\rP\ssu\rQ$}%
\index{zPsqsuQ@$\rP\ssu\rQ$}%
if in addition the following holds:
\ben
%
\atc
\atc
\atc
\atc
\atc
\atc
\atc
\atc
\atc
\cenu
\itlb{ref4s}
If $\i\in\tuq\al$, $X\in\pro\rP\i$, $Y\in\pro\rQ\i$, 
then $X\ne Y$ 
--- then easily $X\cap Y$ 
is {\em meager} in $X$.\qed
\een
}

The transitivity of $\ssw$ does not 
necessarily hold.

\ble
\lam{133}  
Let\/ $\al<\omi$, 
$\rP\ssw\rQ$  belong to\/ $\bfr\al$,
$\j\su\i$ belong to\/ $\tuq\al.$ 
Then
\ben
\aenu
\itlb{ref1}
if $X\in\pro\rP\j$,
then there is $Y\in\pro\rQ\j$, $Y\sq X\,;$ 

\vyk{
\itlb{133b}
if $Z\in\pro\rQ\j$,
then there are $Y\in\pro\rQ\j$, 
$X\in\pro\rP\j$ with\/ 
$Y\sq X\cap Z\,;$ 
}

\itlb{ref2}
if\/ $X\in\pro\rP\i$,
$Y\in\pro\rQ\j$, 
$Y\sq X\rsq\j$, then there is $Z\in\pro\rQ\i$ such that 
$Z\sq X$ and $Z\rsq\j=Y;$ 

\itlb{zrc2}
if\/ $X\in\pro\rP\i$, 
$Y\in\pro\rQ\j$, $Y\sq X\rsq\j$, 
$W\in\pro\rQ\i$, 
and the set\/ $Z=X\cap(Y\rsuq\i)$ satisfies\/ $Z\sq W$, 
then\/ $Z\in\pro\rQ\i\,.$
\een
\ele
\bpf
\ref{ref1}
By definition, there exists $X'\in\rP$ with 
$X=X'\rsq\j$. 
By \ref{134ii} of Definition~\ref{134d} (with $\et=\pu$), 
there is  $Y'\in\rQ$, $Y'\sq X'$. 
Put $Y=Y'\rsq\j$.
\pagebreak[3]

\vyk{
\ref{133b}
By \ref{ref3} of Definition~\ref{134d}, $Y=Z\cap X\ne\pu$ 
for some $X\in \pro\rP\i$; 
$Y$ is clopen in $Z$ by \ref{ref4} 
of Definition~\ref{134d}. 
Thus $Y\in\pro\rQ\i$ by \ref{fr2} of Section \ref{rud}.
}

\ref{ref2}
There exist $X'\in\rP$, $Y'\in\rQ$ with 
$X=X'\rsq\i$, $Y=Y'\rsq\j$. 
Thus $Y'\rsq\j\sq X'\rsq\j$.
By \ref{134ii} of Definition~\ref{134d} 
(with $\et={\ilq\j}$), 
there is  $Z'\in\rQ$, $Z'\sq X'$, such that 
$Z'\rsq\j =Y'\rsq\j=Y$. 
Put $Z=Z'\rsq\i$.
 
\ref{zrc2} 
We have $Z\rsdq\j=Y\sq W\rsdq\j$, therefore  
$U=W\cap(Y\rsuq\i)\in\pro\rQ\i$   
as $\Ker\rQ$ is a kernel by Lemma~\ref{r2n}. 
Yet $Z= U\cap X$, hence
$Z$ is clopen in $U$ by \ref{ref4} of Definition~\ref{134d}. 
Thus $Z\in\pro\rQ\i$ by \ref{fr2} of Section \ref{rud}.
\epf

The next theorem deals with the set $\RC(\rP\cup\rQ)$ 
(the rudimentary hull) in case $\rP\ssw\rQ$. 
We expect that $\rQ$ is $\sq$-dense in $\RC(\rP\cup\rQ)$, 
in this case, but thus turns out to be too hard a problem. 
Still a result of this form holds in a local form as 
claim \ref{134a} of the next theorem shows.

\bte
\lam{134}  
Assume that $\rP\ssw\rQ$  belong to\/ $\bfr\al$ and\/ 
$\rR=\RC(\rP\cup\rQ)$. 
Then, for any\/ $\i\in\tuq\al$, 
$\pro\rQ\i$ is\/ $\sq$-open-dense in\/ $\pro\rR\i$, that is, 
\vyk{
$$
\kaz Z\in\pro\rR\i\,\sus Y\in\pro\rQ\i\,(Y\sq Z)
\qand 
\kaz Z\in\pro\rR\i\,\kaz Y\in\pro\rQ\i\,
(Z\sq Y\imp Z\in\pro\rQ\i).
$$
}%
%
\ben
\Renu
\itlb{134a}
$\kaz Z\in\pro\rR\i\,\sus Y\in\pro\rQ\i\,(Y\sq Z)$, \ \  and 

\itlb{134b}
$\kaz Z\in\pro\rR\i\,\kaz Y\in\pro\rQ\i\,
(Z\sq Y\imp Z\in\pro\rQ\i)$. 
\een
\vyk{\gol%
whereas 
$\pro\rP\j$ is\/ \dd\sq predense in\/ $\pro\rR\j$,
and moreover,  
\ben
\Renu
\atc\atc
\itlb{134c}
if\/  $Z\in\pro\rR\j$, 
then there are\/ $X\in\pro\rP\j$, $Y\in\pro\rQ\j$
with\/ $Y\sq Z\cap X.$
\een
}%
\ete

\bpf
Define sets $\rZ_\i\sq\pele\i$  by induction on $\lh\i$ 
as follows:
\ben
\Aenu
\itlb{lfr1} 
if $\lh\i=1$ then simply $\rZ_\i=\pro\rP\i\cup\pro\rQ\i\,;$ 

\itlb{lfr2} 
if $\lh\i=n+1\ge2$ and $\j=\i\res n$ then $\rZ_\i$ 
contains all $Z\in\pro\rQ\i$ and all sets 
$X\cap(Y\rsuq\i)$,
where $X\in\pro\rP{\i}$, $Y\in\rZ_\j$, $Y\sq X\rsq\j$. 
\een
Let\/ $\j\in\tuq\al$. 
We prove the following list of claims, one by one:
\ben
\nenu
\itlb{zrc0}
$\pro\rP\j\cup\pro\rQ\j\sq \qro\rZ\j\sq\pele\j\,;$ 

\itlb{zrc00}
if\/
$Z\in\qro\rZ\j$ then either\/
$Z\in\pro\rQ\j$ or\/ $Z\sq X$ for some\/ 
$X\in\pro\rP\j\,;$

\vyk
{\kra
\itlb{zrc1}
if\/
$X\in\qro\rZ\i$ 
then\/ $X\rsq\j\in\qro\rZ\j\,;$
\ref{ke2}
}

\vyk
{\kra
\itlb{zrc11}
if\/
$Z\in\qro\rZ\j$ then there is\/
$X\in\qro\rZ\i$ such that\/ $Z=X\rsq\i\,;$
\ref{ke1s}
}

\itlb{zrc3}
if\/
$Z\in\qro\rZ\j$, and\/ $Z\sq W\in\pro\rQ\j,$  
then\/ $Z\in\pro\rQ\j\,;$ 

\itlb{zrc4}
if\/
$\j\su\i$,
$Z\in\qro\rZ\i$, 
$W\in\qro\rZ\j$, $W\sq Z\rsq\j$, 
then\/ $P=Z\cap(W\rsuq\i)\in\qro\rZ\i\,;$  

\itlb{zrc7}
if\/ $X\in \qro\rZ\j$,
$\pu\ne Y\sq X,$ 
$Y$ is clopen in\/ $X,$ then\/ $Y\in\qro\rZ\j\,;$  

\itlb{AB} 
if $\j,\k\in\tuq\al$, $\k\ekp\j$,
and $X\in\qro\rZ\j$,
then $\pi_{\j\k}\akt X\in\qro\rZ\k\,;$ 

\vyk
{\kra
\itlb{zkern} 
the system\/ $\sis{\rZ_\i}{\i\in\tuq\al}$ is 
an\/ $\tuq\al$-kernel, 
as in Section~\ref{ker}$;$
}

\itlb{zrc6}\msur
$\pro\rQ\j$ is dense in\/ $\qro\rZ\j${\rm:}  
if\/
$Z\in\qro\rZ\j$ 
then there is\/ $X\in\pro\rQ\j,X\sq Z;$  
%

\itlb{ZZ}\msur
$\qro\rZ\i=\pro\rR\i\,.$ 
\een

\ref{zrc0} 
$\pro\rZ\j\sq\pele\j$ goes by induction on $\lh\j$, 
and the induction step via 
\ref{lfr2} above  
is carried out by Lemma~\ref{apro}. 
$\pro\rQ\j\sq \rZ_\j$ holds directly by the first option 
of \ref{lfr2}, whereas 
$\pro\rP\j\sq\pro\rZ\j$ is proved by induction using  
\ref{lfr2} and still Lemma~\ref{apro}. 
Claim  \ref{zrc00} 
are rather easy.  

\vyk
{\kra
\ref{zrc11} 
We can assume that $\lh\i=\lh\j+1$. 
By \ref{zrc00}, either $Z\in\pro\rQ\j$ --- 
and then a set $X$ required can be found even in 
$\pro\rQ\i$, --- 
or\/ $Z\sq Y$ for some $Y\in\pro\rP\j$. 
Take any $X'\in\pro\rP\i$ with $Y=X'\rsq\j$. 
Then $X=X'\cap(Y\usq\i)$ is as required.
}

\ref{zrc3} 
Argue by induction on $\lh\j$. 
If $\lh\j=1$ then use \ref{lfr1}  
and \ref{ref4} of Definition~\ref{134d}. 
Suppose that $\lh\j=n+1\ge 2$ and $\k=\j\res n$.
Then either $X\in\pro\rQ\j$ and we are done, or 
$Z=X\cap(Y\rsuq\j)$
where $X\in\pro\rP{\j}$, $Y\in\pro\rZ\k$, $Y\sq X\rsq\k$. 
It follows that $Y=Z\rsq\k\sq W\rsq\k\in\pro\rQ\k$.
Then $Y\in\pro\rQ\k$ by the inductive hypothesis.
Mow $Z\in\pro\rQ\j$ by Lemma~\ref{133}\ref{zrc2}.

\ref{zrc4} 
If $Z\in\pro\rQ\i$ then $Z'=Z\rsdq\j\in\pro\rQ\j$, hence 
$W\in\pro\rQ\j$ by \ref{zrc3}, and we are done. 
Consider the second case of \ref{lfr2}, that is, 
$\lh\i=n+1\ge2$, $\k=\i\res n$, and 
$Z=X\cap(Y\rsuq\i)$, where $X\in\pro\rP\i$, 
$Y\in\pro\rQ\k$, $Y\sq X\rsdq\k$. 
Then $W\sq Z\rsdq\j=Y\rsdq\j\in\rQ_\j$, hence 
$W\in\pro\rQ\j$ by \ref{zrc3}.
It follows that $U=Y\cap(W\rsuq\k)\in\pro\rQ\k$. 
Finally $P=X\cap(U\rsuq\i)\in\pro\rZ\i$. 

\ref{zrc7} 
Argue by induction. 
If $\lh\j=n+1\ge2$ and 
$Z=U\cap(Z'\rsuq\j)$, where $U\in\pro\rP\j$, 
$Z'\in\pro\rZ\k$, $\k=\j\res n$, $Z'\sq U\rsdq\k$, 
use Lemma~\ref{aclo} and then
use the inductive hypothesis. 

\ref{AB} 
A routine induction on \ref{lfr1}, \ref{lfr2},  
based on \ref{nfr4} of Definition~\ref{ruD}.

\vyk
{\kra
\ref{zkern} 
Collect \ref{zrc0}, \ref{zrc1}, \ref{zrc11}, 
\ref{zrc4}, \ref{zrc7}, \ref{AB}.
}

\vyk
{\kra
\ref{zrc5} 
We argue by induction on $\lh\i$, beginning with 
$\lh\i=\lh\j$, \ie, $\i=\j$, in which case $X=Y$ gives 
the result required. 
Now assume that $\lh\i=n+1>\lh\j$. 
If $Z\in\pro\rQ\i$ then there is nothing to prove. 

Suppose now that $Z=U\cap(Z'\rsuq\i)$, 
where $U\in\pro\rP\i$, 
$Z'\in\pro\rZ\k$, $\k=\i\res n$, and $Z'\sq U\rsdq\k$. 
By the inductive hypothesis there is a set 
$X'\in\pro\rQ\k$ such that $X'\sq Z'$ and $X'\rsdq\j=Y$. 
It remains to apply \ref{ref2}, getting a set 
$X\in\pro\rQ\i$ with $X\sq U$ and $X\rsdq\k=X'$.
}

\ref{zrc6} 
Argue by induction on $\lh\j$.  
If $\lh\j=1$, \ie, $Z\in\pro\rP\j\cup\pro\rQ\j$, 
then in case 
$Z\in\pro\rP\j$ apply Lemma~\ref{133}\ref{ref1}. 
Assume that $\lh\j=n+1\ge2$. 
If $Z\in\pro\rQ\j$ then there is nothing to prove. 
Suppose now that $Z=U\cap(Z'\rsuq\j)$, where $U\in\pro\rP\j$, 
$Z'\in\pro\rZ\k$, $\k=\j\res n$, $Z'\sq U\rsdq\k$. 
By the inductive hypothesis there is   
$X'\in\pro\rQ\k$ such that $X'\sq Z'$. 
Applying Lemma~\ref{133}\ref{ref2}, we get a set 
$X\in\pro\rQ\j$ with $X\sq U$ and $X\rsdq\k=X'$.

\ref{ZZ} 
\rit{Prove\/ $\sq$ by induction on $\lh\i$.} 
As case \ref{lfr1} is obvious, consider the step \ref{lfr2}. 
Thus suppose that $\lh\i=n+1\ge2$, $\j=\i\res n$, 
$Z=X\cap(Y\rsuq\i)\in\rZ_\i$, 
where $X\in\pro\rP{\i}$, $Y\in\rZ_\j$, $Y\sq X\rsq\j$, 
and in addition $Z\sq W\in\pro\rQ\i$. 
Then $Y\sq W\rsq\j\in \pro\rQ\j$, hence $Y\in\pro\rQ\j$ 
by the inductive hypothesis. 
Thus $Y=Y'\rsq\j$, $X=X'\rsq\i$, $X'\in\rP$, $Y'\in\rQ$, 
and $Y'\rsq\j\sq X'\rsq\j$. 
As $X',Y'\in\rR$, the set $Z'=X'\cap(Y'\rsq\j\uar\tuq\al)$ 
belongs to $\rR$ by \ref{nfr3} of Definition~\ref{ruD}. 
On the other hand, we have $Z'\rsq\i=Z$ by Lemma~\ref{599}.
Thus $Z\in\pro\rR\i$, as required.

\rit{To prove the direction\/ $\qs$,} 
consider the set $\rZ$ of all sets $X\in\pe{\tuq\al}$ 
satisfying $X\rsq\i\in\qro\rZ\i$ for all $\i\in\tuq\al$. 
Thus $\rP\cup\rQ\sq\rZ$ by \ref{zrc0}. 
We claim that $\rZ\in\bfr\al$. 

Indeed, if $Y\in\pe{\tuq\al}$, $Y\sq X\in \rZ$,    
$Y$ is clopen in $X$, then $Y\rsq\i$ is clopen in 
$X\rsq\i\in\qro\rZ\i$ 
for any $\i\in{\tuq\al}$ by Lemma~\ref{darop}, so that  
$X\rsq\i\in\qro\rZ\i$ by \ref{zrc7}, and we conclude that 
$Y\in \rZ$. 
Thus $\rZ$ satisfies \ref{fr2} of Definition~\ref{ruD}.

To check that $\rZ$ satisfies \ref{nfr3} of 
Definition~\ref{ruD}, assume that 
$\et\in\cpo,\,\et\sq\tuq\al$, 
$X,Y\in\rZ$, and $Y\dar\et\sq X\dar\et$.  
Prove that the set 
$Z=X\cap (Y\dar\et\uar {\tuq\al})$ belongs to $\rZ$. 
If $\i\in\et$ then $Z\rsq\i=Y\rsq\i\in\qro\rZ\i$. 
If $\i\in{\tuq\al}\bez\et$ and $\sg=\et\cap{\ilq\i}$ then 
$Z\rsq\i=X\rsq\i\cap (Y\dar\et)\rsuq\i $ 
by Lemma~\ref{599}, hence again $Z\rsq\i\in\qro\rZ\i$, 
as required. 

Now to check that $\rZ$ satisfies \ref{nfr4} of 
Definition~\ref{ruD}, make use of \ref{AB}.

To conclude, $\rZ\in\bfr\al$, and hence $\rR\sq\rZ$ 
and $\qro\rR\i\sq\pro\rZ\i$, as required.\vom 

\vyk{
\ref{zrc8}
If $Y\in\pro\rQ\j$ then a set $Z$ required exists even 
in $\pro\rQ\i$. 
Suppose that $Y\nin\pro\rQ\j$. 
Then $Y\sq P\in \pro\rP\j$ by \ref{zrc00}, and hence 
there is $X\in\pro\rP\i$ with $X\rsq\j=P$. 
It remains to refer to \ref{ref2}.\vom
}

\vyk{
If $\lh\j=1$ then by definition $W\in\pro\rP\j\cup\pro\rQ\j$,
and the result required follows from \ref{zrc0}. 
Thus suppose that $\lh\j>1$.
If still $W\in\pro\rQ\j$ then the same argument works.

By \ref{zrc00}, it remains to consider the case when  
$W\sq X\in\pro\rP{\j}$. 
Then by \ref{fr3}(b) for $\rP$ there exists a set 
$Z\in\pro\rP\i$ with $Z\rsdq\j=X$.
Note that $Z\in\pro\rZ\i$ by \ref{zrc0}.
Therefore the set $P=Z\cap(W\rsuq\i)$ belongs
to $\pro\rZ\i$ by \ref{zrc4} as well.
But $P\rsdq\j=W$ by construction.  
} 

Finally to prove claims 
\ref{134a}, \ref{134b} 
of the theorem, make use of \ref{ZZ}, and also of 
\ref{zrc6} and \ref{zrc3}. 
For instance, to check {\gol\ref{134a}}, 
note that $Z\in\qro\rZ\i$ by \ref{ZZ}, and hence 
there is $Y\in\pro\rQ\j$, $Y\sq Z$
by \ref{zrc6}.
\epf

\vyk{

\parf{The existence of refinements}
\las{exr}

It follows from the next theorem that countable 
rudiments admit countable refinements. 

\bte
[in $\rL$]
\lam{exrT}
Assume that\/ $\al<\omi$ and\/ $\rP\in\rud_\al$ 
is countable. 
Then there is a countable refinement\/ $\rQ\in\rud_\al$ 
of\/ $\rP$. 
\ete

\bpf
{\ubf We argue in $\rL$}. 
Let $\za=\tuq\al$, so that $\rP\sq\pe\za$. 
Fix any \dd\za admissible map $\phi:\om\to\za$ 
(see Section~\ref{fuz}).
We make use of the 
``refinement forcing'' $\Theta=\Theta_{\rP\phi}$, 
which consists of all functions $\vt$ such that: 
\ben
\nenu
\itlb{th1}
$\dom\vt\sq\om\ti\bse$ is finite, $\ran\vt\sq\rP$, 
and if 
$\ang{k,u}\in\dom\vt$ then $\ang{k,v}\in\dom\vt$ 
for any tuple $\vt\in\bse$ with $\lh v\le\lh u$; 

\itlb{th2}
if $\ang{k,u_0}\in\dom\vt$ and $m=\lh{u_0}$ then 
$\sis{\vt(k,u)}{u\in2^m}$ is a \dd\phi split system in 
the sense of Definition~\ref{splis}; 

\itlb{th3}
if $\ang{k,{v_0}}\in\dom\vt$ and $m+1=\lh{v_0}$ then 
$\sis{\vt(k,v)}{v\in2^{m+1}}$ is an expansion 
of $\sis{\vt(k,u)}{u\in2^{m}}$ in 
the sense of Definition~\ref{splis}. 
\een
The order by extension stipulates that 
if $\phi\sq\phi'$ then $\phi'$ is stronger 
in $\Ta$.

To set up a proper notion of genericity, note that all 
parameters $\omi\yi\rP\yi\za\yi\phi\yi\Theta$ involved 
belong to $\lomb$. 
Let $\cD$ be the (countable) 
collection of all sets $D\sq\lomb$, 
\dd\in definable over $\lomb$ with the five 
mentioned parameters. 

Now let $G\sq\Ta$ be a filter generic over $\cD$, 
that is, it intersects any set $D\in\cD$ dense in $\Ta$. 
We'll prove that $G$ leads to a refinement required. 
Namely, first of all let $g=\bigcup G$ be the union of 
all partial functions $\vt\in G$. 

\rit{We claim that $\dom g=\om\ti\bse.$} 
Indeed Lemma~\ref{pand} and \ref{fr2} 
of Section~\ref{rud} 
imply that each 
\dd\phi split system of height $m<\om$, 
of sets in $\rP$, 
has an expansion to a \dd\phi split system of 
height $m+1$, of sets in $\rP$. 
Therefore each set 
$D_{km}=\ens{\vt\in\Ta}{\ans{k}\ti2^m\sq\dom\vt}$ 
is dense in $\Ta$. 
As easily $D_{km}\in\cD$, we conclude that 
$G\cap D_{km}\ne\pu$ for all $k,m$. 
This implies the result required.

It follows that, for any $k$, $\sis{g(k,u)}{u\in\bse}$ 
is a \dd\phi fusion sequence of sets $g(k,u)\in\rP$, 
in the sense of Definition~\ref{fuzD}, and hence each 
set $X_k=\bigcap_m\bigcup_{u\in2^m}g(k,u)$ belongs to 
$\pe\za$. 
We now let 
$$
\bay{rcl}
\rQ_0 &=&
\bigcup_k\ens{Y\in\pe\za}
{Y\sq X_k\text{ is clopen in }X_k},\\[0.5ex]
\rQ &=&
\ens{\pi_{\i\j}\akt Y}
{Y\in\rQ_0\land \text{tuples }\i\ekp\j
\text{ belong to }\za=\tuq\al},
\eay
$$
and prove that $\rQ$ belongs to $\rud_\al$ 
and is a refinement of $\rP$. 

The next lemma presents the main ingredient of the proof.

\ble
[in $\rL$]
\lam{exrL}
If\/ $\i,\j$ belong to $\za=\tuq\al$, $k,l<\om$, 
and\/ $\ang{k,\i}\ne\ang{l,\j}$, then\/ 
$({X_k}\dir\i)\cap({X_l}\dir\j)=\pu$. 
{\rm(Recall that $X\dir\i=\ens{x(\i)}{x\in X}$.)}
\ele

Leaving the proof of the lemma to the next section, we 
proceed with the theorem.

\epf

\parf{The proof of the lemma}
\las{exr+}

\bpf[Lemma~\ref{exrL}]
It suffices to prove that the set 
$$
\Da^{\i\j}_{kl}=\ens{\vt\in\Ta}
{\sus m\,\sus n\,\kaz u\in2^m\,\kaz v\in2^n 
\big((\vt(k,u)\dir\i)\cap(\vt(l,v)\dir\j)=\pu\big)\!}
$$
is dense in $\Ta$.
To prove the density let $\vt\in\Ta$; we have to define 
$\vt'\in\Da=\Da^{\i\j}_{kl}$ with $\vt\sq\vt'$. 
By the density of the sets $D_{km}$ as above, we may assume 
that there exist $m,n<\om$ such that 
$\ans k\ti2^m\sq\dom\vt$, 
$\ans l\ti2^n\sq\dom\vt$, 
and $m,n$ are the largest with this property. 

{\ubf Case 1:} 
$k=l$, and accordingly $n=m$ and $\i\ne\j$. 
Thus the system $\sis{X_u}{u\in2^m}$ of sets 
$X_u=\vt(k,u)\in\rP$ is a \dd\phi split system. 
By Lemma~\ref{pand}, it admits an expansion by a
\dd\phi split system $\sis{Y_s}{s\in2^{m+1}}$ 
of relatively clopen sets 
$Y_{u\we e}\sq X_u$, $e=0,1$, hence all 
sets $Y_s$ belong to $\rP\in\rud_\al$. 
Now we are going to shrink the sets $Y_s$ as follows. 

{\it Step 1a.} 
Take any $s_0\in2^{m+1}.$ 
There is a set $Y\in\pe\za$, $Y\sq Y_{s_0}$, 
clopen in $Y_{s_0}$, with $(Y\dir\i)\cap(Y\dir\j)=\pu$, 
by Corollary~\ref{99z1}. 
Lemma~\ref{suz} 
(the clopen version, recall that $\rP$ is closed under 
clopen subsets in $\pe\za$) gives a 
\dd\phi split system $\sis{Y'_s}{s\in2^{m+1}}$ of 
sets $Y'_s\in\rP$, $Y'_s\sq Y_s$, with $Y'_{s_0}=Y.$ 
Do this one by one for all $s_0\in2^{m+1},$ 
and let $\sis{Z_s}{s\in2^{m+1}}$ be the final 
\dd\phi split system of 
sets $Z_s\in\rP$, $Z_s\sq Y_s$, with 
$(Z_s\dir\i)\cap(Z_s\dir\j)=\pu$ for all $s\in2^{m+1}.$ 

{\it Step 1b.} 
Now take any $s_0\ne s_1$ in $2^{m+1}.$ 
There are sets $U,V\in\pe\za$, $U\sq Z_{s_0}$, 
$V\sq Z_{s_1}$, 
clopen in resp.\ $Z_{s_0},Z_{s_0}$, with 
$(U\dir\i)\cap(V\dir\j)=\pu$, and still 
$U\dar\za_\phi[s_0,s_1]=V\dar\za_\phi[s_0,s_1]$, 
by Lemma~\ref{99z2}. 
Corollary~\ref{suz+} gives a 
\dd\phi split system $\sis{Z'_s}{s\in2^{m+1}}$ of 
sets $Z'_s\in\rP$, $Z'_s\sq Y_s$, 
with $Z'_{s_0}=U,$ $Z'_{s_1}=V.$ 
Do this one by one for all pairs of 
$s_0\ne s_1$ in $2^{m+1},$ 
and let $\sis{U_s}{s\in2^{m+1}}$ be the final 
\dd\phi split system of 
sets $U_s\in\rP$, $U_s\sq Z_s\sq Y_s$, with 
$(U_{s_0}\dir\i)\cap(U_{s_1}\dir\j)=\pu$ 
for all $s_0\ne s_1$ in $2^{m+1}.$  

Thus overall we have 
$(U_{s_0}\dir\i)\cap(U_{s_1}\dir\j)=\pu$ 
for all $s_0,s_1$ in $2^{m+1},$ 
no matter equal or different. 
It follows that the extended ``condition'' 
$\vt':=\vt\we \sis{\vt'(k,s)}{s\in2^{m+1}}$ 
belongs to $\Da^{\i\j}_{kk}$, where 
$\vt'(k,s):=U_s$ for all $s\in2^{m+1}.$ 

{\ubf Case 2:} 
$k\ne l$. 
The systems $\sis{X_u}{u\in2^m}$ and 
$\sis{Y_v}{v\in2^n}$ of sets 
$X_u=\vt(k,u)$ and $Y_v=\vt(l,v)$ in $\rP$ are 
\dd\phi split systems. 
By Lemma~\ref{pand}, they admits expansions by 
\dd\phi split systems 
$\sis{U_s}{s\in2^{m+1}}$, $\sis{V_t}{t\in2^{n+1}}$, 
of sets in $\rP$, as in Case 1. 
We are going to shrink the sets $U_s,V_t$ as follows. 

Take any $s_0\in 2^{m+1},$ $t_0\in 2^{n+1}.$ 
Lemma~\ref{99z2} ($\et=\pu$) yields 
sets $A,B\in\pe\za$, $A\sq U_{s_0}$, 
$B\sq V_{t_0}$, 
clopen in resp.\ $U_{s_0},V_{t_0}$, with 
$(A\dir\i)\cap(B\dir\j)=\pu$. 
Arguing as in Step 1a above, we independently get 
\dd\phi split systems 
$\sis{U'_s}{s\in2^{m+1}}$, $\sis{V'_t}{t\in2^{n+1}}$ 
of relatively clopen sets $U'_s\sq U_s$, $V'_t\sq V_t$ 
in $\rP$, 
such that $U_{s_0}=A$ and $V_{t_0}=B$, and hence 
$(U_{s_0}\dir\i)\cap(V_{t_0}\dir\j)=\pu$.  

Doing this consecutively for all pairs of 
$s_0\in 2^{m+1},$ $t_0\in 2^{n+1}$  
results in final \dd\phi split systems   
$\sis{A_s}{s\in2^{m+1}}$, $\sis{B_t}{t\in2^{n+1}}$ 
of sets $A_s,B_t\in\rP$, $A_s\sq U_s$, $B_t\sq V_t$, 
with 
$(A_{s}\dir\i)\cap(V_{t}\dir\j)=\pu$ 
for all $s\in 2^{m+1},$ $t\in 2^{n+1}.$  
It follows that the extended ``condition'' 
$\vt':=\vt\we \sis{\vt'(k,s)}{s\in2^{m+1}}
\we \sis{\vt'(l,t)}{t\in2^{n+1}}$ 
belongs to $\Da^{\i\j}_{kl}$, where 
$\vt'(k,s):=A_s$ and $\vt'(l,t):=B_t$,  
for all $s\in2^{m+1}, t\in2^{n+1}.$ 

This ends the proof of Lemma~\ref{exrL}.
\epf

}

\parf{Rudimentary sequences}
\las{rs}

The next definition introduces the notion of
transfinite sequences of rudiments, 
``$\ssw$-increasing'' in the sense that each term is a 
$\ssw$-successor of the rudimentary hull of the union 
of all previous terms, by condition~\ref{drs4} 
of Definition~\ref{drs} below. 
We use quotation marks because $\ssq$ is not
claimed to be a transitive relation.

{\ubf We still argue in $\rL$.}

\bdf
\lam{drs}
%
Let a \rit{rudimentary sequence} (or \ruds) 
of length $3\le \la\le\omi$ 
\kmar{ruds}
\index{rudimentary sequence}%
\index{rudimentary sequence@\ruds}%
\index{sequence!rudimentary}%
be any sequence
$\jqo=\sis{\rqu\al}{\al<\la}$, satisfying 
\ref{drs1},\ref{drs2},\ref{drs3},\ref{drs4} below: 
\ben
\Aenu
\itlb{drs1}
$\rqu0=\rqu{1}=\rqu{2}=
\ans{\text{all clopen sets }X\in\peq{2}}\in\bfr2$;

\itlb{drs2}
if $\nu<\la$ then  
$\rqu\nu\in\bfr{\nu}$ is at most countable;

\itlb{drs3}
if $\al<\nu<\la$, $\i\in\tuq\al$, 
and $X\in\pro{\rqu\nu}\i$ then  
$X\sqf\bigcup(\pro{\rqu\al}\i)$ 
in the sense of Section~\ref{frel}.
\een
For any such $\jqo$ we 
put $\uu\jqo=\bigcup_{\al<\la}\rqu\al$ and    
\kmar{noc jqo}
$\noc\jqo=\noc{\uu\jqo}$; thus 
 $\uu\jqo\sq\pei$ and 
$\noc\jqo\in\nf$ is a normal forcing.

If $\la<\omi$ strictly then we define 
\kmar{bsc jqo}%
\index{zzzbcs@$\bsc{}$}%
\index{zRHjqo@$\RC(\jqo):=\RC{(\bsc\jqo)}$}%
$\bsc\jqo=\bsc{}_{\al<\la}\rqu\al:=
\bigcup_{\al<\la}(\rqu\al\uar\tuq\la)$; 
thus $\bsc\jqo\sq\peq{\la}$, and then   
$\RC(\jqo):=\RC{(\bsc\jqo)}\in\bfr\la$. 
We add the last condition:
\ben
\Aenu
\atc\atc\atc
\itlb{drs4}
if $3\le\nu<\la$ then  
$\RC(\jqo\res\nu)\ssw\rqu\nu$ 
in the sense of
Definition~\ref{134d}; here 
$\RC(\jqo\res\nu)=\RC(\bsc(\jqo\res\nu))=
\RC(\bigcup_{\al<\nu}(\rqu\al\uar\tuq\nu))$.
\een
%
%
Let 
\index{rudimentary sequences, $\rs$}%
\index{rudimentary sequences, $\rs_\la$}%
\index{zRS@$\rs$}%
\index{zRSla@$\rs_\la$}%
$\rs_\la=$ all \ruds s of length $\la$, 
\kmar{rs}
$\rs=\bigcup_{\la<\omi}\rs_\la$.
\edf


\bte
[in $\rL$]
\lam{142}
Let\/
$\jqo=\sis{\rqu\al}{\al<\la}\in\rs_\la$, 
$3\le\la\le\omi$. 
Then\/
\ben
\renu
\itlb{142i}
$\rR=\RC(\jqo)\in\bfr{\la}$, $\cai\la\in\rR$,
and if\/ $\la<\omi$ then 
$\rR$ 
is countable$;$ 

\itlb{142j}
if\/ $\al<\la$ then$:$\/ \ 
{\rm(a)}   
the set\/
$\rpi\al=\RC(\jqo\res\al)\in\bfr{\al}$ 
is countable,\\[0.2ex] 
{\rm(b)} 
$\cai{\al}\in\rpi\al$, \\[0.3ex] 
{\rm(c)} 
$\rqm_\al\sq\rpi\al=\RC(\rqm_\al)$, where\/ 
$\rqm_\al=\bsc(\jqo\res\al)
=\bigcup_{\ga<\al}(\rqu\ga\uai\al)$,\\[0.6ex] 
{\rm(d)} 
$\kaz X\in\rpi\al\,\sus Y\in\rqu\al\,(Y\sq X);$

\itlb{142jb}  
if\/ $\ga<\al<\la$ then\/ 
$(\rpi\ga\cup\rqu\ga)\uai\al\sq\rpi\al;$ 

\itlb{142jc}  
if\/ $\la=\ga+1$ then\/ 
$\rR=\RC((\rqm_\ga\cup\rqu\ga)\uai\la)
=\RC((\rpi\ga\cup\rqu\ga)\uai\la);$  

\itlb{142jx}  
if\/ $\la=\ga+1$ and\/ $\j\in\tuq\ga$ then\/ 
$\pro\rR\j=\pro{(\RC(\rpi\ga\cup\rqu\ga))}\j;$

\itlb{142jd}  
if\/ $\la<\omi$ is a limit ordinal then\/ 
$\rR=\bigcup_{\al<\la}(\rpi\al\uai\la)$, 
and the set\/ ${(\uu\jqo)}\uai\la$ 
is\/ $\sq$-dense in\/ $\rR\,;$ 

\itlb{1421}
if\/ $\j\in\tuq\al$, $\al<\la<\omi$, 
then the set\/ 
$\bigcup_{\al\le\ba<\la}(\pro{\rqu{\ba}}\j)$ 
is\/ \dd\sq dense in\/ $\pro\rR\j;$ 

\itlb{1422}
if\/ $\j\in\tuq\al$, $2\le\al<\la<\omi$, then\/
$\pro{\rqu\al}\j$ is \dd\sq predense in\/
$\pro\rR\j\,;$

\itlb{1423}
if\/ $\la<\omi$, $\et\in\cpo$, $\et\sq\tuq\la$, and\/ 
$X\in\pe\et$, then\/ $X\in\cX:=\noc\jqo$ iff\/ 
$X\rsq\j\in\pro\rR\j$ for all\/ $\j\in\et$, where\/ 
$\rR=\RC(\jqo)$ by\/ \ref{142i}.

\itlb{1423x}
therefore, by\/ \ref{1423}, if\/ $\la<\omi$ and\/  
$\i\in\tuq\la$, 
then the sets\/ $\rR=\RC(\jqo)$ and\/ 
$\cX=\noc\jqo$ satisfy\/ $\pro\cX\j=\pro\rR\j$.
\een
\ete

\bpf
\ref{142i}, \ref{142j} are easy:  
$\cai\la\in\rR$ and \ref{142j}(b) hold 
by \ref{drs1} of Definition~\ref{drs}, 
\ref{142j}(d) holds by \ref{drs4} 
and \ref{134ii} of Section~\ref{ref} 
(the particular case).

\ref{142jb} 
We have 
$\rpi\ga\uai\al\sq\RC(\rqm_\ga\uai\al)$ 
by Corollary~\ref{rrec}, hence 
$$
(\rpi\ga\cup\rqu\ga)\uai\al
\,\sq\,
\RC(\rqm_\ga\uai\al)\cup (\rqu\ga\uai\al)
\,\sq\,
\RC(\rqm_\al)
\,=\,
\rpi\al, 
$$
as required. 
($\cai{\ga}\in\rpi\ga$ holds by  \ref{142j}(b).)

\ref{142jc}  
Let $\rU=\rqm_\ga$.
Then $\rU\sq\rpi\ga=\RC(\rU)$ and 
$$
\bay{rcccll}
\rR&=&\RC((\rU\cup\rqu\ga)\uai\la)
&\sq&
\RC((\rpi\ga\cup\rqu\ga)\uai\la)&\sq\\[0.5ex]
&&&\sq&
\RC(\RC(\rU\cup\rqu\ga)\uai\la)&,
\eay
$$
because 
$\rU\cup\rqu\ga\sq \rpi\ga\cup\rqu\ga\sq
\RC(\rU\cup\rqu\ga)$. 
On the other hand, by Corollary~\ref{rrec}, 
we have 
$\RC((\rU\cup\rqu\ga)\uai\la)=
\RC(\RC(\rU\cup\rqu\ga)\uai\la)$, 
so that both inclusions in the displayed formula are 
equalities, and we are done. 

\ref{142jx} 
We have 
$\rR
=\RC((\rpi\ga\cup\rqu\ga)\uai\la) 
=\RC(\RC(\rpi\ga\cup\rqu\ga)\uai\la)$, 
see the proof of \ref{142jc}. 
Therefore 
$\pro\rR\j=\pro{\RC(\rpi\ga\cup\rqu\ga)}\j$
by Lemma~\ref{134*0}.

\ref{142jd} 
As $\RC(\rqm_\al)=\rpi\al$ by \ref{142j}, 
$\rpi\al\uai\la\sq\RC(\rqm_\al\uai\la)\sq\rR$ 
by Corollary~\ref{rrec}, 
hence the set 
$\rR'=\bigcup_{\al<\la}(\rpi\al\uai\la)$ 
satisfies $\rR'\sq\rR$. 
Yet $\rR'=\RC(\rR')$ by Lemma~\ref{13.2.1} 
and \ref{142jb}.
Then, as $\rqm_\al\sq\rpi\al$, we have 
$$
\rR
=
\RC({\TS\bigcup}_{\al<\la}(\rqm_\al\uai\la))
\sq
\RC({\TS\bigcup}_{\al<\la}(\rpi\al\uai\la))
=\RC(\rR')=\rR',
$$
and clearly $\rR'\sq\rR$,
so that $\rR=\rR'=\bigcup_{\al<\la}(\rpi\al\uai\la)$, 
as required. 

To prove the density in \ref{142jd}, let $X\in\rR$. 
Then $X=Y\uai\la$, where $Y\in\rP_\al$ and $\al<\la$, 
by the above. 
However $\rP_\al=\RC(\jqo\res\al)\ssw\rqu\al$ by 
\ref{drs4} of Definition \ref{drs}. 
Therefore there iz $Y'\in\rqu\al$, $Y'\sq Y$, by 
\ref{134ii} of Definition \ref{134d}. 
It remains to take $X'=Y'\uai\la$.

\vyk{
To prove \ref{1421},  
let $\j\in\tuq\al$, $2\le\al<\la<\omi$. 
If $\la=3$,  then $\al=2$, $\j\in\tud$, and 
$\rR=\RC(\jqo)=\RC(\rqu2\uai3)$ by 
\ref{drs1} of Definition~\ref{drs}, hence 
$\pro{\rR}\j=\pro{\rqu2}\j$ by Lemma~\ref{134*0}, 
thus there is nothing to prove.
If $\la$ is limit   
and $Z\in\pro\rR\j$, then  
by \ref{142jd}, there is $\ba$, $\al<\ba<\la$, 
with $Z\in\pro{\rpi\ba}\j$. 
By \ref{142j}(d), there is 
$Y\in \pro{\rqu\ba}\j$ with $Y\sq Z$, as required.
}

The limit case in \ref{1421} easily follows from 
\ref{142jd}. 
Therefore suppose that $\la=\ga+1$ in \ref{1421}.  
Then $\al\le\ga$,    
$\rR=\RC((\rpi\ga\cup\rqu\ga)\uai\la)$ 
by \ref{142jc}, $\j\in\tuq\ga$. 
We convert this to 
$\rR=\RC(\RC(\rpi\ga\cup\rqu\ga)\uai\la)$ 
by Corollary~\ref{rrec}. 
Therefore 
$\pro\rR\j=\pro{(\RC(\rpi\ga\cup\rqu\ga))}\j$
by Lemma~\ref{134*0}. 
However $\pro{\rqu\ga}\j$ is dense in 
$\pro{(\RC(\rpi\ga\cup\rqu\ga))}\j$
by Theorem~\ref{134}\ref{134a}.

\ref{1422} 
Let $\la$ be limit and $X\in\pro\rR\la$. 
Then by \ref{1421} there is an ordinal $\ba$, 
$\al<\ba<\la$, and $Y\in\pro{\rqu\ba}\j$, such that 
$Y\sq X$. 
Then $Y\sqf\bigcup(\pro{\rqu\al}\j)$ by 
\ref{drs3} of Definition~\ref{drs}. 
We conclude that there is $Z\in \pro{\rqu\al}\j$ 
such that $Y\cap Z$ is not meager in $Y$. 
Therefore there us a set $\pu\ne U\sq Y\cap Z$ 
clopen in $Y$. 
Then $U\in \pro{\rqu\ba}\j$ by \ref{fr2} 
of Definition~\ref{132}, and we are done.

Now let $\la=\ga+1$ in \ref{1422}.
Suppose that 
$X\in\pro\rR\j$, 
where $\pro\rR\j=\pro{(\RC(\rpi\ga\cup\rqu\ga))}\j$ 
by \ref{142jx}. 
It follows from Theorem~\ref{134}\ref{134a} 
that there is a set $Y\in\pro{\rqu\ga}\j$ with $Y\sq X$. 
Then 
proceed as in the limit case. 

Finally check \ref{1423}. 
By definition the set 
$\cX=\noc\jqo$ satisfies the equality 
$\cX=\noc{\bigcup_{\al<\la}\rqu\al}$. 
As obviously $\rR=\RC(\jqo)\sq\noc\jqo$, we have 
$\cX=\noc{\rR}$ as well. 
It follows that $\pro\cX\j=\pro\rR\j$ 
for all $\j\in\tuq\la$ by Lemma~\ref{r2n}. 
It remains to refer to \ref{rfo5} of Section~\ref{rfo} 
for $\cX.$
\epf

\sekt[\ \,\ Specifying rudimentary sequences]
{Specifying rudimentary sequences}
\las{sek7}

The goal of this Chapter is to specify a list of conditions 
which imply that the normal forcing $\cX=\RC(\jqo)$, 
generated by a given \ruds\ $\jqo\in\rL$ of length $\omi$, 
satisfies Theorem~\ref{mt1c}.

We introduce {\ubf properties} \ref{C1}, \ref{C2}, 
\refn{C3}n, \ref{C4}, \ref{C5}, \refn{C6}n 
of an $\ssq$-increasing $\omi$-sequence 
$\jqo=\sis{\rqu\al}{\al<\omi}\in\rs_{\omi}$ 
of rudiments, which imply the following:
\bit
\item[$-$]
the set\/ $\rQ=\uu\jqo=\bigcup_{\al<\omi}\rqu\al$ is\/ 
dense in\/ $\cX=\noc\jqo$ ---  Lemma~\ref{143L}; 

\item[$-$]
$\cX$ has the Fusion property --- Lemma~\ref{144t}; 

\item[$-$]
$\cX$ is \dd ncomplete --- Lemma~\ref{145t}; 

\item[$-$]
$\cX$ has the Structure property --- Lemma~\ref{155t}; 

\item[$-$]
$\cX$ has the $n$-Definability property --- Lemma~\ref{154t};
\eit
and hence the limit 
forcing $\cX=\noc\jqo:=\noc{\uu\jqo}$ satisfies 
Theorem~\ref{mt1c}. 
The properties are summed up in the notion of 
{\ubf 1-5-$n$ extension}, which allows to transform 
the content of Theorem~\ref{mt1c} 
by Theorem~\ref{mt1d}. 

{\ubf We still argue in $\rL$ in this Chapter.}

\parf{Coding iterated perfect sets}
\las{cod}

Further study of \ruds s will involve  
a coding system of iterated perfect sets based on 
codes in 
$\hc$\,= all hereditarily countable sets.  

Clearly any set $X$ in some $\pe \xi$, $\xi\ne\pu$, 
is of cardinality
continuum, hence $X$ does not belong to $\hc$. 
This makes it difficult to evaluate the
complexity of different collections of sets $X$
of such kind.
To fix this problem, we make use of a 
coding by countable dense subsets.

\bdf
[codes]
\lam{cps}
If $\xi\in\cpo$ then let $\cpe\xi$ 
\index{codes!$\cpe\xi$}%
\index{zcIPSxi@$\cpe\xi$, codes}%
({\ubf c} from `codes') 
consist of all
\kmar{cpe xi\mns clo X}%
{\em at most countable} sets $A\sq\can\xi$ such that
the closure $\clo A$ in $\can\xi$ belongs to $\pe\xi$.
We put $\cpei=\bigcup_{\xi\in\cpo}\cpe\xi$; 
\index{codes!$\cpei$}%
\index{zcIPS@$\cpei$, codes}%
\kmar{cpei}
thus $\cpei\sq\hc$. 

If $\rA\sq\cpei$ then let $\cla\rA=\ens{\clo A}{A\in\rA}$. 
\kmar{cla rA}
\index{zA**@$\cla \rA$, set of closures}%
\index{operation!$\cla \rA$, set of closures}%
\index{zzz**@$\cla{\;}$, set of closures}%
\edf

In the trivial case $\xi=\pu$, the collection 
$\cpe\pu=\pe\pu$ contains 
the only one element $\bon=\ans\pu$, see Remark~\ref{pu}, 
and $\clo\bon=\bon$. 

\vyk{
Then let $\xnl\al=\bigcup_{\ba<\al}\xn\al$. 

We may note that $\lc\al\sq\ens{\lc\ba}{\ba<\al}$ in the 
G\"odel construction. 
It follows that $\rank{\lc\al}\le\al$. 
We conclude that $\dym{\xn\al}\sq\tuq\al$ for all $\al<\omi$.
}

\parf{Getting density }
\las{143}

This section is intended to define a condition which 
implies, for a given sequence 
$\jqo=\sis{\rqu\al}{\al<\omi}\in\rs_{\omi}$, that 
the set $\rQ=\uu\jqo=\bigcup_{\al<\omi}\rqu\al$ is 
\rit{\dd\psq dense} in $\cX=\noc\jqo:=\noc{\uu\jqo}$, 
\index{zNHq@$\noc\jqo$}%
that is, $\kaz X\in\cX\,\sus U\in\rQ\,(U\psq X)$.
This condition will be of \rit{step-wise form}, that is, 
in the form of a relation between each term $\rqu\al$ 
and the sequence $\jqo\res\al$ obtained before $\al$. 


{\ubf We continue to argue in $\rL$}. 
Under this assumption, the set $\hc$ of all 
\rit{hereditarily countable sets} satisfies $\hc=\lomi$, 
and hence $\hc$ is well-ordered by the canonical G\"odel 
relation $\lel$. 
\index{zzz-L@$\lel$, G\"odel \weo}%
\index{hereditarily countable, $\hc$}%
\index{zHC@$\hc$, hereditarily ctble}%
Thus $\hc=\ens{\lc\al}{\al<\omi}$ in 
$\rL$, where $\lc\al$ is the $\al$th element of $\hc$ 
via $\lel$.
Recall that $\bH\al=\ens{\lc\ga}{\ga<\al}$. 
See Section~\ref{cdia} on details.
%
We let 
\kmar{xpe al}
$$
\cxpe\al=\cpei\cap\bH\al 
\qand
\xpe\al=\ens{\clo B}{B\in\cxpe\al}
\,.
\kmar{cxpe al}
$$

To provide the density  property as above, we   
add some definitions based on the sequence of sets 
$\baS\al\sq\bH\al$, $\al<\omi$, 
satisfying Proposition~\ref{p62}. 

\ben
\Renu
\itlb{144A}
Let $\al<\omi$. 
If there is a unique triple of   
$M\in\cpei$ and $M',M''\in\hc$ such that 
$\ang{\om,M,M',M''}\in\baS\al$ 
\kmar{bM al\mns bmp \mns bmpp}%
\index{zMal@$\bM\al$}%
\index{zM'al@$\bmp\al$}%
\index{zM''al@$\bmpp\al$}%
%
then put $\bM\al=M,$  $\bmp\al=M',$  $\bmpp\al=M''.$  
Otherwise let $\bM\al=\bon=\ans{\pu}\in\pe\pu=\cpe\pu$ 
and  $\bmp\al=\bmpp\al=\pu$. 
Note that $\bM\al,\bmp\al,\bmpp\al\in\bH\al$ and 
$\bM\al\in\cxpe\al$.

\itlb{144B}
Let 
$\bPsi k\al=\ens{B\in\cxpe\al}{\ang{k,B}\in\baS\al}$ 
and $\bPsd k\al=\ens{\clo B}{B\in\bPsi k\al}$ 
\kmar{\,\mns\,\mns bPsi k al\mns bPsd k al}%
\index{zBkal@$\bPsi k\al$}%
\index{zBkal*@$\bPsd k\al$}%
for any $k$. 
Thus 
$\bPsi k\al\sq\cxpe\al,$ 
$\bPsd k\al\sq\xpe\al$ 
are countable.
\een


\ble
[in $\rL$]
\lam{144L}
If\/ $M\in\cpei$, $M',M''\in\hc$, 
and\/ $P_k\sq\cpei$, $\kaz k$, 
then the following set\/ $\bW$ is stationary 
in\/ $\omi:$
\bce
$\bW=\ens{\al}
{\bM\al=M \land\bmp\al=M'\land\bmpp\al=M''
\land \kaz k\,(P_k\cap\cxpe\al=\bPsi k\al)}.$ 
\ece
The sequences\/ 
$\sis{\ang{\bM\al,\bmp\al,\bmpp\al}}{\al<\omi}$   
and\/ $\sis{\bPsi k\al}{k<\om,\al<\omi}$ 
belong to\/ $\id\hc1$.
\ele

\bpf
Applying Proposition~\ref{p62} for the set 
\pagebreak[0]
$$
\baS{}=\ans{\ang{\om,M,M',M''}} \cup
\ens{\ang{k,B}}{k<\om\land B\in P_k}\,,
$$
we conclude that   
$D:=\ens{\al<\omi}{\baS{}\cap\bH\al=\baS\al}$
is stationary in $\omi$. 
On the other hand, the set $W'$ of all $\al<\omi$, 
such that
$$
\baS{}\cap\bH\al=\ans{\ang{\om,M,M',M''}} \cup
\ens{\ang{k,B}}{k<\om\land B\in P_k\cap\cxpe\al}\,,
$$
is a club. 
Thus $W'\cap D$ is still stationary. 
However $W'\cap D\sq \bW$ by construction.
To prove the definability claim apply 
Proposition~\ref{p62} yet again.
\epf 

Now we are sufficiently equipped to consider the 
density property.

\ble
[in $\rL$]
\lam{143L}
Assume that\/ $\jqo=\sis{\rqu\al}{\al<\omi}$ is a \ruds, 
satisfying the following condition$:$
{\rm
\ben
\genu
\itlb{C1}
\index{condition!\ref{C1}}%
\index{zP1@\ref{C1}, condition}%
for any\/  
$\la<\omi$, if\/ 
$\clo{\bM\la}\in \noc{\jqo\res\la}$ 
and $\dym{\clo{\bM\la}}\sq\tuq\la$, 
then there is\/ 
$Y\in\rqu\la$ satisfying\/ $Y\psq\clo{\bM\la}$.
\een
}
\noi
Then the set\/ $\rQ=\uu\jqo=\bigcup_{\al<\omi}\rqu\al$ is\/ 
\dd\psq dense in\/ $\noc\jqo$.
\ele

\bpf
Let 
$X\in\noc\jqo$. 
Then obviously 
$X\in \noc{\jqo\res\la}$ and 
$\dym{\clo{\bM\la}}\sq\tuq\la$ 
for all $\la$ larger than some $\la_0<\omi$. 
The set 
$\bW=\ens{\al}{\clo{\bM\al}=X}$ 
is stationary by Lemma~\ref{144L}, 
hence there is a limit ordinal  
$\la\in\bW\yi \la\ge\la_0$. 
Applying \ref{C1}, we complete the proof.
\epf

\parf{Getting Fusion}
\las{144}

The next lemma provides
another step-wise condition which implies the Fusion 
property as in Section~\ref{fup}.

\ble
[in $\rL$]
\lam{144t}
Assume that\/ $\jqo=\sis{\rqu\al}{\al<\omi}$ is a \ruds, 
satisfying both\/ \ref{C1} of Lemma~\ref{143L} 
and the following condition$:$
{\rm 
\ben
\genu
\atc
\itlb{C2}
\index{condition!\ref{C2}}%
\index{zP2@\ref{C2}, condition}%
for any 
{\ubf limit\/} 
$\al<\omi$, if\/ 
$\clo{\bM\al}\in \rqu{<\al}:=\bigcup_{\ga<\al}\rqu\ga$ 
then there is\/ 
$X\in\rqu\al$ satisfying\/ $X\psq\clo{\bbZ\al}$   
and\/ $X\sqd \bigcup\bPsd k\al$   
for all $k<\om$ such that\/ 
$\bPsd k\al\sq\rqu{<\al}$ 
and\/ $\bPsd k\al$ is dense in\/ $\rqu{<\al}$.
\een
}
\noi
Then the set\/ $\noc\jqo$ has the 
Fusion property of Section~\ref{fup}.
\ele

\bpf
We argue in $\rL$. 
Let $X_0\in \cX:=\noc\jqo$. 
Consider a sequence of dense sets $\cY_m\sq\cX.$ 
We have to find a set $Y\in\cX$ satisfying $Y\psq X_0$ 
and $Y\sqd\bigcup\cY_m$ for all $m$. 
Assume that $X_0\in \rqu{}:=\bigcup_{\al<\omi}\rqu\al$, 
by Lemma~\ref{143L}. 

We may \noo\ assume that each $\cY_m$ is in fact 
open-dense; 
then, still by Lemma~\ref{143L}, 
(*) each set 
$\cZ_m:=\cY_m\cap\rqu{}$ is open dense in 
$\rqu{}$. 
We let $P_m=\ens{B\in\cpei}{\clo B\in\cZ_m}$, so that  
$\cZ_m=\ens{\clo B}{B\in P_m}$, $\kaz m$. 
Pick a set $C\in\cpei$ satisfying $X_0=\clo {C}$. 
By Lemma \ref{144L}, the set\/ 
\pagebreak[1] 
$$
\bW=\ens{\al<\omi}
{\bbZ\al=C \land \kaz m\,(P_{\al m}=\bPsi m\al)}
$$ 
is stationary,   
where $P_{\al m}=P_m\cap\cxpe\al$. 
Let $\cZ_{\al m}=\ens{\clo B}{B\in P_{\al m}}$. 
Recall that $\rqu{<\al}:=\TS\bigcup_{\ga<\al}\rqu\ga$. 
Note that the set 
$$
E=\ens
{\al<\omi}
{X_0=\clo C\in \rqu{<\al}
\land
\kaz m\,(\cZ_{\al m}\text{ is open dense in }
\rqu{<\al})
} 
$$ 
is a club by (*) above. 
Thus there exists an ordinal $\al\in E\cap \bW.$ 

Then we have $\clo{\bbZ\al}=\clo C\in \rqu{<\al}$,  
and in addition each  $\bPsd m\al$ is dense in\/ $\rqu{<\al}$. 
Therefore by \ref{C2} of the lemma  there exists 
$X\in\rqu\al$ satisfying\/ $X\psq\clo{\bbZ\al}=X_0$   
and $X\sqd \bigcup\bPsd m\al$ for all $m<\om$. 
However $\bPsd m\al=\cZ_m\sq\cY_m$ by construction.
\epf

\vyk{
\bpf
\vyk{
Let $\lel$ be the canonical wellordering of $\rL$, as above.
\index{zz<L@$\lel$}%
\index{zzHC@$\hc$}%
Thus $\lel$ orders $\hc$ similarly to 
$\omi$, $\lel$ is a $\id\hc1$ relation, 
and $\lel$ has the  
\rit{goodness} property: the set of 
all \dd\lel initial segments 
$I_x(\lel)=\ens{y}{y\lel x}$, 
where $x\in\hc$, is $\id\hc1$.
}
The diamond 
principle $\Diamond_{\omi}$ is true in $\rL$  
by \cite[Thm 13.21]{jechmill}, 
hence there is a  $\id\hc1$ sequence of sets 
$S_\al\sq\al$, $\al<\omi$, such that 
\ben
\Aenu
\itlb{1gret}
if\/ 
$X\sq\hc$ then the set\/ $\ens{\al<\omi}{S_\al=X\cap\al}$ 
is stationary in\/ $\omi$,
\een
The $\id\hc1$-Definability property is achieved 
by taking the \dd\lel least possible $S_\al$ at each 
step $\al$ in the standard construction of $S_\al$.
We get the following two results \ref{2gret}, 
\ref{3gret} as easy corollaries. 

Let $A_\ga=\ens{C_\al}{\al\in S_\ga}$, where $C_\al$ is 
the \dd\al th element of $\hc=\lomi$ in the sense of G\"odel's 
well-ordering $\lel$. 
Then $\sis{A_\ga}{\ga<\omi}$ is still a $\id\hc1$ sequence, 
and we have the following:
\ben
\Aenu
\atc
\itlb{2gret}
if $d_\al\in\hc$ for all $\al<\omi$, and 
$D_\ga=\ens{d_\al}{\al<\ga}$,
then the set\/ $E=\ens{\ga}{A_\ga=D_\ga}$ 
is stationary in\/ $\omi$. \ 
\een
Indeed let $Q\sq\omi$ be a club. 
By \ref{1gret} for 
$X=\ens{\nu<\omi}{\sus\ba\,(d_\ba=C_\nu)}$, 
the set $F=\ens{\al<\omi}{S_\al=X\cap\al}$ 
is stationary.

and note that 
$K=\ens{\ka}{D_\ka=\ens{C_\al}{\al\in X\cap\ka}}$ 
is a club.)

Second, for any $\al$, if   
$A_\al=\sis{a_\da}{\da<\al}$, where each $a_\da$ itself is 
equal to an \dd\om sequence $\sis{a^n_\da}{n<\om}$, 
then let $B^n_\al=\sis{a^n_\da}{\da<\al}$ for all $n$. 
Otherwise let $B^n_\al=\pu$, $\kaz n$.
Then $\sid{B_\al^n}{n<\om}{\al<\omi}$ is a\/ $\id\hc1$ 
system of sets in $\hc$, such that 
\ben
\Aenu
\atc
\atc
\itlb{3gret}
if $Y^n_\al\in\hc$ for all $\al<\omi$ and $n<\om$, then 
the set\break  
$\ens{\mu}{\kaz n\,(B^n_\mu=\ens{Y^n_\al}{\al<\mu}}$ 
is stationary in $\omi$.
\een

Now things become somewhat more complex.  

Let $\mu<\omi$. 
We define $\nos z\mu=\bigcup B^0_\mu$. 
If $B^1_\mu\in\vmf$ and $\len{B^1_\mu}=\mu$ 
then let $\nos\vjpi\mu=B^1_\mu$; 
otherwise let $\nos\vjpi\mu$ be equal to the \dd\lel least 
\muq\ in $\vmf$ of length $\mu$. 
(Those exist 
by Corollary~\ref{xisc}\ref{xisc1}.) 
Finally we let $\nos D\mu=\bigcup B^2_{\mu+1}$. 

Let's show that the sequences of sets $\nos\vjpi\mu$,  
$\nos D\mu$, $\nos z\mu$ prove the theorem. 
Suppose that $\vjPi=\sis{\nor\jPi\nu}{\nu<\omi}\in\vmi$, 
$z\in\hc$, and $D\sq \mt\vjPi$.
Let $X^0_\al=z$, $X^1_\al=\ang{\al,\nor\jPi\al}$, 
$X^2_\al=D\cap\mt{\vjPi\res\al}$ for all $\al$.
The set 
$$
M=\ens{\mu<\omi}{B^n_\mu=\ens{X^n_\al}{\al<\mu}\:
\text{ for }\:n=0,1,2}
$$
is stationary by \ref{3gret}.
Assume that $\mu\in M$. 
Then $B^0_\mu=\ens{X^0_\al}{\al<\mu}=\ans z$, 
therefore $\nos z\mu=z$. 
Further 
$B^1_\mu=\ens{X^1_\al}{\al<\mu}
=\ens{\ang{\al,\nor\jPi\al}}{\al<\mu}=\vjPi\res\mu\in\vmf$, 
therefore $\nos\vjpi\mu=\vjPi\res\mu$.
Finally we have 
$\nos D\mu=\bigcup B^2_{\mu+1}=\bigcup_{\al\le\mu}{X^2_\al}
=D\cap\mt{\vjPi\res\mu}$, as required.
\epf
}

\parf{Getting Completeness}
\las{145}

Here we introduce another step-wise condition on a 
\ruds\ $\jqo$ which implies the $n$-Completeness 
property of Definition~\ref{121} for the according 
normal hull $\nh\jqo$.

\ble
[in $\rL$]
\lam{145t}
\sloppy 
Assume that\/ $n\ge2$ 
and\/ $\jqo=\sis{\rqu\al}{\al<\omi}$ is a\/ \ruds, 
satisfying\/ \ref{C1} of Lemma~\ref{143L} 
and the following condition$:$
{\rm
\ben
\genj  n
\atc
\atc
\itlb{C3}
\index{condition!\refn{C3}\nn}%
\index{zP3@\refn{C3}\nn, condition}%
if\/ $\nn\ge2$ then for any\/ 
$\la<\omi$, if\/  
$\clo{\bbZ\la}\in \rqu{<\la}:=\bigcup_{\ga<\la}\rqu\ga$, 
and\/ $\bmp\la$ is a closed formula\/ $\vpi$ in\/ 
$\bigcup_{k\le n}\ls1k$, 
then there is\/ 
$X\in\rqu\la$ satisfying\/ $X\psq\clo{\bbZ\la}$   
and either\/ $X \fo \vpi$ 
or\/ $X \fo \otr\vpi$.
\een
}
\noi
Then the set\/ $\nh\jqo$ is\/ $\nn$-complete.
\ele

We underline that 
condition \refn{C3}\nn\ is void in case $\nn=1$.

\bpf
We argue in $\rL$. 
Given $X_0\in \cX:=\noc\jqo$ and  a closed formula\/ $\vpi$ in\/ 
$\bigcup_{k\le n}\ls1k$, 
we have to find a set $Y\in\cX$ satisfying $Y\psq X_0$ 
and either\/ $X \fo \vpi$ or\/ $X \fo \otr\vpi$.  
We can \noo\ assume that 
$X_0\in \rqu{}:=\bigcup_{\al<\omi}\rqu\al$, 
by Lemma~\ref{143L}. 
Pick a set $C\in\cpei$ satisfying $X_0=\clo {C}$. 

The set\/ 
$ 
\bW=\ens{\al<\omi}
{\bM\al=C \land \bmp\al=\vpi}
$ 
is stationary by Lemma \ref{144L}, whereas  the set 
$ 
E=\ens
{\al<\omi}
{X_0=\clo C\in \rqu{<\al}
} 
$ 
is obviously a club. 
Thus there exists a limit ordinal $\la\in E\cap\bW.$ 
Then we have $\clo{\bbZ\la}=\clo C\in \rqu{<\la}$. 
Therefore by \refn{C3}n 
there exists a set 
$X\in\rqu\la$ satisfying\/ $X\psq\clo{\bbZ\la}=X_0$   
and either\/ $X \fo \vpi$ 
or\/ $X \fo \otr\vpi$, as required.
\epf

\vyk{

\sekt[\ \ Definability property]{Definability property}
\las{sek8}

In this Chapter, we deal with the $n$-Definability property 
of Definition~\ref{82d}, 
and our goal will be to define a group of three conditions 
which will imply that a normal forcing of the form 
$\cX=\noc\jqo$ satisfies that property. 
}

\parf{Uniform sets and the  Structure  property}
\las{155}

Our next goal will be to attack the  Structure   
property as in Definition~\ref{82d}.
We are going to define a  condition, 
related to uniform sets,  
which will imply that a normal forcing of the form 
$\cX=\noc\jqo$ satisfies that property. 

Recall that a set $X\in\pe\xi$ is \rit{uniform} 
(Section~\ref{174}), 
\index{set!uniform}%
\index{uniform}%
if for any pair of tuples $\i\su\j$ in $\xi$ 
and any $x,y\in X$, 
we have $x(\j)=y(\j)\imp x(\i)=y(\i)$.   

\ble
[in $\rL$]
\lam{155t}
Assume that\/ $\jqo=\sis{\rqu\al}{\al<\omi}$ is a\/ \ruds, 
satisfying both\/ \ref{C1} of Lemma~\ref{143L} 
and the following condition$:$
{\rm
\ben
\genu
\atc
\atc\atc
\itlb{C4}
\index{condition!\ref{C4}}%
\index{zP4@\ref{C4}, condition}%
for any\/  
$\la<\omi$, if\/ 
$\clo{\bbZ\la}\in \rqu{<\la}:=\bigcup_{\ga<\la}\rqu\ga$, 
then there is a uniform set\/ 
$X\in\rqu\la$, $X\psq\clo{\bbZ\la}$.
\een%
}
\noi
Then the set\/ $\cX=\nh\jqo$ satisfies the Structure   
property.
\ele

\bpf
Consider a pair of tuples $\i\su\j$ in $\tup$.
We claim that the set 
\bce
$C_{\i\j}=\ens{X\in\rqu{}}{X\text{ is uniform}\land
\i,\j\in\dym X}$ 
\ece
is $\psq$-dense in $\rqu{}=\bigcup_{\al<\omi}\rqu\al$. 
Indeed suppose that $Z\in\rqu{}$. 
The set 
$$
\bW=\ens{\al<\omi}{\i,\j\in\tuq\al\land\clo{\bM\al}=Z }
$$  
is stationary by Lemma~\ref{144L}. 
Therefore there is a limit $\la\in \bW$ with 
$Z=\clo{\bM\la}\in \rqu{}$ 
and $\i,\j\in\tuq\la=\dym Z$. 
Then \ref{C4} yields a set $X\in C_{\i\j}$, 
$X\sq Z$, as required. 

It follows by the density that there is a set 
$X\in C_{\i\j}$ such that 
$\i,\j\in\xi=\dym X$ and $\w\dar\xi\in\clo X$. 
Then $X$ is uniform, hence there is a continuous map 
$F:\can{}\to\can{}$ 
coded in $\rL$ such that $\clo f(x(\j))=x(\i)$ 
for all $x\in\clo X$. 
Then $\w(\i)=\clo f(\w(\j)\in\rL[\w(\j)]$, as required. 

On the other hand, if $\i\not\sq\j$, then $\i\nin\ilq\j$, 
and 
$\w(\i)\nin\rL[\w(\j)]$ follows from Corollary~\ref{632}.
\epf

\parf{Key formulas for the $n$-Definability property}
\las{151}

Approaching the $n$-Definability property, 
here formulas are introduced which will define 
sets in Definition~\ref{82d}.

Recall that if $X\in\pei$ and $\i\in\dym X$ then 
$X\dir\i=\ens{x(\i)}{x\in X}$, and if $\cX\sq\pei$ 
then $\cX\dir\i=\ens{X\dir\i}{X\in\cX\land \i\in\dym X}$.
Suppose, that, {\ubf in $\rL$}, 

\ben
\fenu
\itlb{151*}\msur
$\jqo=\sis{\rqu\al}{\al<\omi}\in\rL$ is a \ruds\  
and $\cX=\noc\jqo$ 
(as in Definition \ref{drs}), 
so that $\cX\in\nf$ is a normal forcing.
\een

The following formulas based on 
$\jqo=\sis{\rqu\al}{\al<\omi}\in\rL$ 
are considered.

\bde
\item[$\bev\jqo kx$:] 
\kmar{bev jqo kx}%
\index{formula!$\bev\jqo kx$}%
\index{zBeq@$\bev\jqo kx$, formula}%
$k<\om\land x\in\can{}\land 
\sus\j\in\tuq2\big(\lh\j=k \land \j\text{ is even\:}
\land$\\[0.3ex]
\vphantom.\hfill
$\kaz\al<\omi\,\sus Z\in\rqu\al\dir\j\,
(x\in\clo Z)
\big);$

\item[$\bod\jqo kx$:] 
\kmar{bod jqo kx}%
\index{formula!$\bod\jqo kx$}%
\index{zBoq@$\bod\jqo kx$, formula}%
$k<\om\land x\in\can{}\land 
\sus\j\in\tuq2\big(\lh\j=k \land \j\text{ is odd\:}
\land$\\[0.3ex]
\vphantom.\hfill
$\kaz\al<\omi\,\sus Z\in\rqu\al\dir\j\,
(x\in\clo Z)
\big).$
\ede

We'll prove that these formulas define 
the sets as in Definition~\ref{82d}  
in $\cX$-generic extensions of $\rL$ 
--- provided the basic \ruds\ $\jqo$ satisfies certain 
conditions. 
The next lemma proves this result in one direction.

\ble
\lam{151L}
Assume\/ \ref{151*} in\/ $\rL$ as above.  
Let\/ $\w\in\can\tup$ be a\/ $\cX$-generic array over\/ $\rL$,  
$\i\in\tup$, $k=\lh\i$, and\/ $x=\w(\i)$. 
Then\/ $\rL[x]\mo\bev\jqo kx$, resp.,  $\bod\jqo kx$,  
provided\/ $\i$ is resp.\  even, odd.
\ele

\bpf
Let $\j=\kn\i\in\tuq2$ (see Section~\ref{ker}), so that 
$\i\ekp\j$ (the parity-equivalence, Section~\ref{perm}), 
and 
$\pi_{\i\j}\in\per$ is parity-preserving.
We claim that 
\ben
\nenu
\itlb{151L1}
if $\al<\omi$ then $\pro{\rqu\al}\j$ is pre-dense in 
$\pro\cX\j$.
\een
As clearly 
$\cX=\bigcup_{\la<\omi}\cX_\la$, where 
$\cX_\la=\noc{\jqo\res\la}$, 
it suffices to check that 
\ben
\nenu
\atc
\itlb{151L2}
if $\al<\la<\omi$ and $\la$ is limit then 
$\pro{\rqu\al}\j$ is pre-dense in $\pro{\cX_\la}\j$.
\een
However $\pro{\cX_\la}\j=\pro{\rP_\la}\j$ by 
Theorem~\ref{142}\ref{1423x}, where 
$\rpi\la=\RC(\jqo\res\la)\in\bfr{\la}$. 
On the other hand, 
the set $\pro{\rqu\al}\j$ is pre-dense in 
$\pro{\rP_\la}\j$ by Theorem~\ref{142}\ref{1422}.
This implies \ref{151L2} and \ref{151L1}.

Now assume that $\al<\omi$ (in $\rL$), 
and let $\w'=\pi_{\i\j}\akt\w$. 
Then $\w'\in\can\tup$ is still $\cX$-generic over $\rL$ 
along with $\w$ since $\pi_{\i\j}\in\per$ and $\cX$ is 
necessarily $\per$-invariant.
It follows from  \ref{151L1} that 
$\pro {\w'}\j\in \clo P$ for some $P\in\pro{\rqu\al}\j$
by Lemma~\ref{gras}\ref{gras2}, and hence obviously 
${\w'}(\j)\in \clo Z$ for  
$Z=P\dir \j\in {\rqu\al}\dir\j$.

To conclude, the real $x=\w(\i)=\w'(\j)$ satisfies 
$\bev\jqo kx$, resp., $\bod\jqo kx$ in $\rL[x]$, 
in case $\i$ (and then $\j$ as well) is even, 
resp., odd.
\epf

\parf{The inverse of the lemma}
\las{152}

The condition \ref{C5} defined below will allow us   
to reverse Lemma~\ref{151L}. 
This condition involves a special notation. 
Recall definitions in Sections~\ref{173} and \ref{101}.

\bdf
[in $\rL$]
\lam{152d}
Let $\al<\omi$. 
If $\bmpp\al\in\komd$ and 
$\da_\al:=\modd{\bmpp\al}\sq\tuq\al$ 
then  define $\bbf\al\in\kond{\tuq\al}$ 
\kmar{bbf al}%
\index{zfal@$\bbf\al$}%
by $\bbf\al(x)=\bmpp\al(x\dar\da_\al)$ 
for all $x\in\rat{\tuq\al}$.
Otherwise define $\bbf\al\in\kond{\tuq\al}$ by 
$\bbf\al(x)=\om\ti0$ for all $x\in\rat{\tuq\al}$. 

In both cases define 
$\bbF\al=\clo{\bbf\al}\in\cfd_{\tuq\al}$. 
\kmar{bbF al}%
\index{zFal@$\bbF\al$}%

Let $\jqo=\sis{\rqu\al}{\al<\omi}$ be a \ruds. 
Define the following condition$:$
%
\ben
\genu
\atc\atc\atc\atc
\itlb{C5}
\index{condition!\ref{C5}}%
\index{zP5@\ref{C5}, condition}%
For any 
$\la<\omi$, if 
$\clo{\bM\la}\in\rqu{<\la}=\bigcup_{\ga<\la}\rqu\ga$ 
then there is a set $Y\in\rqu\la$, $Y\psq\clo{\bM\la}$, 
such that 
one of the two following claims holds:
\ben
\itlb{C51}
${\bbF\la}$ avoids every $E\in \rqu\al\dir\i$ on $Y$ 
for all $\i\in\tuq\la$;

\itlb{C52}
there is $\j\in\tuq\la$ such that ${\bbF\la}$ is an 
$\j$-axis map on $Y$ and ${\bbF\la}$ avoids each 
$E'\in \rqu\la\dir\i$ on $Y$ 
for all $\i\in\tuq\la$ with $\i\not\ekp\j$.
\qed
\een
\een
\eDf


\bte
\lam{153}
Assume that \ref{151*} of Section~\ref{151} holds, 
and\/ $\jqo$ satisfies\/ \ref{C1}, \ref{C2}, \ref{C5} 
in\/ $\rL$. 
Let\/ $\w$ be\/ $\cX$-generic over\/ $\rL$.  
Then 
$$
\gee {\rL[\w]}\w =
\ens{\ang{k,x}}{x\in\rL[\w]\land \rL[x]\mo\bev\jqo kx}.
%
$$
and the same for the `odd' case.
\ete

\bpf
The inclusions $\sq$ in both cases follow from 
Lemma~\ref{151L}. 
To establish the inverse inclusions,
let $k\ge1$, $x\in\rL[\w]\cap\can{}$, 
and  $\rL[x]\mo\bev\jqo kx$, so that there is an even 
tuple $\i\in\tuq2$ with $\lh\i=k$, satisfying   
\vyk{
$$
\kaz\al<\omi=\omil\,
\sus A\in\rqu\al\dir\i\,
(x\in\clo A).
\eqno(\dag)
$$
}
\busq
{555}
{
\kaz\al<\omi=\omil\,
\sus A\in\rqu\al\dir\i\,
(x\in\clo A).
}
We have to prove that $\ang{k,x}\in \gee {\rL[v]}\w$. 

By \ref{C2} and Lemma~\ref{144t}, the set $\cX=\nh\jqo\in\nf$ 
has the Fusion property. 
It follows, by Theorem~\ref{1022}\ref{1022ii} 
and Corollary~\ref{1011},  that 
$x=\clo f(\w\dar\sg)$ for some $\sg=\tuq{\al_0}\yi\al_0<\omi$, 
and $f\in\kond\sg$.  
We claim that the set $D_f=\bigcup_{\al_0<\la<\omi}D_{f\la}$ 
is $\psq$-dense in 
$\rqu{}=\bigcup_{\la<\omi}\rqu\la$, where
\vyk{
$$
\bay{rcl}
D_{f\la} &=&
\big\{Y\in\rqu\la: 
Y\text{ satisfies \ref{C51} or \ref{C52} of 
Definition \ref{152d}}  \\[0.3ex]
&& 
\phantom{\bigcup_{\al<\omi}\;\;\;\;}
\;\;\text{ for }f\text{ instead of }\bbF\la
\big\},
\eay
$$
}
$$
D_{f\la} =
\ens{Y\in\rqu\la}{ 
Y
\text{ satisfies \ref{C51} or \ref{C52} in Definiton~\ref{152d}}}
$$
Indeed suppose that $Z\in\rqu{}$. 
The set 
$\mathbf W=\ens{\la<\omi}
{\clo{\bM\la}=Z \land\bmpp\la=f}$  
is stationary by Lemma~\ref{144L}. 
Therefore there exists a limit ordinal $\la\in\mathbf W$ 
satisfying $\al_0<\la$, hence $\sg\sq\tuq\la$, 
$Z=\clo{\bM\la}\in \bigcup_{\ga<\la}\rqu\ga$, and 
$f=\bmpp\la$. 
Then \ref{C5} yields a set $Y\in D_f$, $Y\sq Z$, 
as required. 

By the density just proved, there exist $\la<\omi$ 
and $Y\in D_{f\la}$ satisfying $\w\res\tuq\la\in\clo Y$. 
(Note that $\modd Y=\tuq\la$ since $Y\in\rqu\la$.) 
We conclude from \eqref{555} and the choice of $f=\bmpp\la$ 
 that ${\bbF\la}$ does {\ubf not} avoid some 
 $E\in \rqu\la\dir{\i}$ on $Y.$ 
It follows that \ref{C51} definitely fails, 
and hence \ref{C52} holds 
for some $\j\in\tuq\la$ such that $\i\ekp\j$. 
In particular, ${\bbF\la}$ is 
a $\j$-axis map on $Y$, 
meaning that ${\bbF\la}(y\dar{\tuq\la})=y(\j)$ 
for all $y\in Y$, 
and hence $x={\bbF\la}(\w\dar{\tuq\la})=\w(\j)$.  
It remains to note that $\j$ is even and $\lh\j=k$ 
by the choice of $\i$, because $\i\ekp\j$.  
Thus $\ang{k,x}\in \gee {\rL[\w]}\w$, as required. 
\epf

\parf{Getting $n$-Definability}
\las{154}

Here we 
introduce another property, related to the definability of 
a \ruds\ as a whole, which will help us to reduce the formulas 
$\bev\jqo kx$, $\bod\jqo kx$ to $\ip1{n+1}$ as required by 
Definition~\ref{82d}, and thereby to fully establish the 
 $n$-Definability property of the ensuing normal forcing.

\bdf
[in $\rL$] 
\lam{154d}
Say that a sequence $\jba=\sis{\rba\al}{\al<\la}$ is a 
\rit{coded\/ \ruds}, if each $\rba\al\sq\cpei$ is at most countable 
\index{rudimentary sequence!coded}%
\index{sequence!rudimentary, coded}%
and the sets 
$\rqu\al=\cla{{\rba\al}}:=\ens{\clo A}{A\in\rba\al}$ 
form a \ruds\ 
$\jqo=\sis{\rqu\al}{\al<\la}$. 

We write $\jqo=\cla\jba$ in this case. 
\kmar{cla jba}%
\index{zzba**@$\cla\jba$}%
\index{symbol!$\cla{}$}%
\edf

\ble
[in $\rL$]
\lam{154t}
\sloppy
Let\/ $n\ge1$ 
and\/ $\jqo=\sis{\rqu\al}{\al<\omi}$ be a\/ \ruds, 
satisfying conditions\/ \ref{C1}, \ref{C2}, \ref{C4}, \ref{C5}, 
and the following condition$:$
{\rm
\ben
\genj n
\atc
\atc\atc
\atc\atc
\itlb{C6}
\index{condition!\refn{C6}\nn}%
\index{zP6@\refn{C6}\nn, condition}%
it is true in\/ $\rL$ that there is a coded\/ 
\ruds\ $\jba=\sis{\rba\al}{\al<\omi}$ for\/ $\jqo$, 
of the definability class\/ $\is\hc\nn$, such  that\/ 
$\jqo=\cla\jba$.
\een
}
\noi
Then\/ $\cX=\nh\jqo$ satisfies the\/ 
$n$-Definability property of Definition~\ref{82d}.
\ele

\bpf
We have to estimate the complexity of the  
relations $\rL[x]\mo\bev\jqo kx$ 
and $\rL[x]\mo\bod\jqo kx$ as in Theorem~\ref{153}. 

By \refn{C6}n, there exists a concrete parameter-free $\is{} n$ 
formula $\vpi(\cdot,\cdot)$ such that $\rQ=\rqu\al$ iff 
$\al,\rQ\in\lomi$ and $\lomi=(\hc)^\rL\mo\vpi(\al,\rQ)$. 
Let
$$
\bay{rcl}
\Psie(k,x)
\kmar{Psie(k,x)}%
\index{formula!$\Psie(k,x)$}%
\index{zzPhie(k,x)@$\Psie(k,x)$, formula}%
&\!:=\!&
 \kaz\al\,\kaz\rQ\,
\big[
\al,\rQ\in\rL\land \vpi(\al,\rQ)^\rL 
\imp \sus\j\in\tuq2\,\\[0.5ex]
&&
\hspace*{5ex}\big(
\lh\j=k\land{\j\text{ is even\,}}\land \sus A\in\rQ\dir\j\,
(x\in\clo A) 
\big)
\big],
\eay
$$
where $ \vpi(\al,\rQ)^\rL $ means the formal relativization of 
all unbounded quantifiers to $\rL$. 
(Compare to the formulas $\bev\jqo kx$ in Section \ref{151}.) 

Consider any $\cX$-generic array $\w\in\can\tup$ over\/ $\rL$, 
$k<\om$, and $x\in\rL[\w]\cap\can{}$. 
Recall that $\rL[\w]$ preserves $\omil$ 
by Theorem~\ref{1022}\ref{1022i}, 
and hence using $\omi=\omil=\omi^{\rL[\w]}$ 
does not lead to an ambiguity.
Theorem~\ref{153} implies that 
\setcounter{equation}0
\busq
{154a}
{
\ang{k,x}\in \gee {\rL[\w]}\w 
\leqv
\lomi[x]\mo\Psie( k,x).
}
Now assume that $\gM\sq\rL[\w]$ is a transitive class, 
closed under pairs, 
and\/ $\rL[x]\sq \gM$ for all\/ $x\in \gM$, 
by the Structure property as in 
Definition~\ref{82d}. 
Then   
\busq
{154b}
{\gee {}\w \cap\gM
\;=\;
\big\{\ang{k,x}\in\gM:\gM\mo\Psie( k,x)^{\lomi[x]}\big\} 
}
holds by \eqref{154a}, 
where the upper index $^{\lomi[x]}$ means the 
formal relativization of all unbounded quantifiers 
in $\Psie( k,x)$ to ${\lomi[x]}$.

Now note that $\vpi$ is $\is{} n$, and 
hence so is $\vpi(\al,\rQ)^\rL$ because 
``$x\in\rL$'' is $\is{}1$ by G\"odel. 
We conclude that $\Psie(k,x)$ is essentially a $\ip{}n$ formula. 
It follows that ${\gM\mo\Psie( k,x)^{\lomi[x]}}$ defines a 
$\ip{} n$ relation over $(\hc)^\gM$ since $y\in\lomi[x]$ is 
still a $\is{}1$ relation over $(\hc)^\gM$ by G\"odel. 
It follows by \eqref{154b} that $\gee {}\w \cap\gM$ is a 
$\ip\hc n$ set in $\gM$, hence a $\ip1{ n+1}$ set by 
Proposition~\ref{p60}, as required. 
The ``odd'' case is considered similarly.
\epf

\parf{Fourth form of the main theorem}
\las{156}

To summarize the results achieved above, we now formulate 
another form 
of Theorem \ref{mt1} in the introduction, 
that further develops the 
previous form given by Theorem~\ref{mt1c}.
This is based on the next definition, 
that gathers the step-wise properties 
\ref{C1}, \ref{C2}, \refn{C3}\nn, \ref{C4}, \ref{C5} 
in a single step-wise property. 

\bdf
\lam{15e}
Let $\la<\omi$, $n\ge1$. 
Say that a term $\rqu\la$ is a \rit{1-5-$n$ extension} of a 
\index{rudimentary sequence@\ruds!1-5-$n$ extension}%
\index{extension!1-5-$n$ extension}%
\index{1-5-$n$ extension}%
\ruds\ $\jqo=\sis{\rqu\ga}{\ga<\la}$ 
if the following \ref{C*},\ref{C1+},\ref{C1++} hold:
\ben
\Aenu
\itlb{C*}
the extended sequence 
$\jqo\we\rqu\la=\sis{\rqu\ga}{\ga\le\la}$ 
is still a \ruds; 
%
\itlb{C1+}
as in \ref{C1},  
if 
$\clo{\bM\la}\in \noc{\jqo}$ 
and $\dym{\clo{\bM\la}}\sq\tuq\la$ 
then there is\/ 
$Y\in\rqu\la\yi Y\psq\clo{\bM\la};$ 

\itlb{C1++}
if 
$\clo{\bM\la}\in \rqu{<\la}:=\bigcup_{\ga<\la}\rqu\ga$ 
then there is a set $Y\in\rqu\la$ satisfying\/ 
$Y\psq\clo{\bM\la}$  
and  the following conditions \ref{C2+}--\ref{C5+}:
%
%
\atm
\atm
\Cnenu
\itlb{C2+} 
as in \ref{C2}, if $\la$ is limit then 
$Y\sqd \bigcup\bPsd k\la$ holds for all $k<\om$ such that  
$\bPsd k\la\sq\rqu{<\la}$ and\/ $\bPsd k\la$ is dense 
in\/ $\rqu{<\la}$;

\itlb{C3+}
as in \refn{C3}n,
if $n\ge2$ and $\bmp\la$ is a closed formula $\vpi$ in 
$\bigcup_{k\le n}\ls1k$ then  $Y \fo \vpi$ or $Y\fo \otr\vpi$ 
--- {\em void in case\/ $n=1$};

\itlb{C4+}\msur
as in \ref{C4},
$Y$ is a uniform set;

\itlb{C5+}
as in \ref{C5} of Definition \ref{152d}, 

either (a) 
${\bbF\la}$ avoids every $E\in\rqu\la\dir\i$ on $Y$ 
for all $\i\in\tuq\la$, 

or (b)
there is $\j\in\tuq\la$ such that ${\bbF\la}$ is an 
$\j$-axis map on $Y$ but ${\bbF\la}$ avoids each 
$E'\in\rqu\la\dir\i$ on $Y$ 
for all $\i\in\tuq\la$ satisfying $\i\not\ekp\j$.\qed
\een
\eDf


\bte
[in $\rL$]
\lam{mt1d}
Assume that\/ $\nn\ge1$. 
Then there is a\/ \ruds\/ $\jqo=\sis{\rqu\al}{\al<\omi}$ 
satisfying the global definability condition\/ \refn{C6}\nn\ 
and such that, for any ordinal\/ $\la<\omi$, 
$\rqu\la$ is a\/  a  1-5-$\nn$ extension of\/ $\jqo\res\la$.
\ete

\vyk{
\bte
[in $\rL$]
\lam{mt1c}
Assume that\/ $\nn\ge1$. 
Then there is a \ruds\/ $\jqo=\sis{\rqu\al}{\al<\omi}$ 
satisfying conditions\/ \ref{C1}, \ref{C2}, \refn{C3}\nn, \ref{C4}, \ref{C5}, 
and\/ \refn{C6}\nn.
\ete
}

\bpf
[Theorem~\ref{mt1} from Theorem~\ref{mt1d}] 
Let $\jqo$ be such a \ruds\ as in Theorem~\ref{mt1d}.
Consider the associated normal forcing  
$\cX=\noc\jqo\in\RF$. 

Lemma~\ref{144t} implies that $\cX$ 
has the Fusion property. 

Lemma~\ref{155t} implies that $\cX$ 
has the Structure  property. 

Lemma~\ref{145t} implies that the set\/ $\cX=\nh\jqo$ is\/ 
$\nn$-complete.  

Finally, $\cX$ satisfies the $\nn$-Definability property by Lemma~\ref{154t}. 

To conclude, the set $\cX$ is as Theorem~\ref{mt1c} requires. 
 
But Theorem~\ref{mt1c} implies Theorem~\ref{mt1}, 
see Section~\ref{98}.\vom

\epF{Thms~\ref{mt1c} and \ref{mt1} from Thm~\ref{mt1d}}


Thus Theorem~\ref{mt1d} implies Theorem~\ref{mt1}, 
the first main result of this paper. 
Chapters \ref{VI} and \ref{VII} below will contain the 
proof of  Theorem~\ref{mt1d},  
and thereby will accomplish the proof of Theorem~\ref{mt1}.

\sekt[\ \ The existence of 1-5-${\nn}$ extensions] 
{The existence of 1-5-${\nn}$ extensions}
\las{VI}

Working towards the proof of Theorem~\ref{mt1d}, 
the goal of this Chapter will be the existence of 1-5-$n$ 
extensions of \ruds s of countable length.

\parf{The existence theorem and basic notation}
\las{181}

\bte
[in $\rL$]
\lam{18t}
Let\/ $\la<\omi$ and\/ ${\nn}\ge1$. 
Then every\/ \ruds\/ $\jqo=\sis{\rqu\al}{\al<\la}$ admits 
a\/ 1-5-${\nn}$ extension\/ $\rqu\la$. 
\ete
\pagebreak[2]%

{\ubf Notation, in $\rL$.}
We fix $\la,{\nn},\jqo,\rqu\al$ as in the theorem. 
Put 
\index{notation Chapter~\ref{VI}!$\nn$}%
\index{zn@$\nn$, notation Chapter~\ref{VI}}%
\index{notation Chapter~\ref{VI}!$\la$}%
\index{zzla@$\la$, notation Chapter~\ref{VI}}%
\index{notation Chapter~\ref{VI}!$\jqo$}%
\index{zzkopp@$\jqo$, notation Chapter~\ref{VI}}%
\index{notation Chapter~\ref{VI}!$\rqu\al$}%
\index{zQal@$\rqu\al$, notation Chapter~\ref{VI}}%
\index{notation Chapter~\ref{VI}!$\rql\la$}%
\index{zQ<al@$\rql\al$, notation Chapter~\ref{VI}}%
\index{notation Chapter~\ref{VI}!$\bta$}%
\index{zztab@$\bta$, notation Chapter~\ref{VI}}%
\index{notation Chapter~\ref{VI}!$\rU_\la$}%
\index{zUla@$\rU_\la$, notation Chapter~\ref{VI}}%
\index{notation Chapter~\ref{VI}!$\cX_\la$}%
\index{zXla@$\cX_\la$, notation Chapter~\ref{VI}}%
\bce
$\rql\la=\bigcup_{\al<\la}\rqu\al$, \ 
$\bta=\tuq\la$, \ 
$\rU_\la=\rh(\rql\la\uar\bta)$, \ 
$\cX_\la=\nh{\rql\la}$. 
\ece

\bre
\lam{18r} 
$\rU_\la\in\bfr\la$ is a countable rudiment, 
$\cX_\la\in\NF$ is a normal forcing,  
$\rql\la\uar\bta\sq\rU_\la$, $\can\bta\in\rU_\la$. 
In addition, 
$\rU_\la\sq\cX_\la$, 
and $\cX_\la\rsq\i=\rU_\la\rsq\i$ for all $\i\in\bta$ 
by Lemma~\ref{r2n}.  
\ere

We'll use the sets 
$\bM\la\in\cxpe\la$; 
$\bmp\la, \bmpp\la\in\bH\la$;  
$\bPsi k\la\sq\cxpe\la$ and
$\bPsd k\la\sq\xpe\la$ (both countable sets);
defined in \ref{144A},\ref{144B} of Section~\ref{143}. 
\ben
\nenu
\itlb{181a}
If $\clo{\bM\la}\in \rql\la$ then put 
$\tX=\clo{\bM\la}\uar\bta$, otherwise let 
\index{notation Chapter~\ref{VI}!$\tX$}%
\index{zXover@$\tX$, notation Chapter~\ref{VI}}%
\index{set!Xt@$\tX$}%
$\tX=\can\bta$, so that $\tX\in\rU$ in both cases.

\itlb{181b}
If $\bmp\la$ is a closed formula in 
$\bigcup_{k\le {\nn}}\ls1k$, 
then let $\bvpi_\la$ be that formula, otherwise let 
\index{notation Chapter~\ref{VI}!$\bvpi_\la$}%
\index{zzFila@$\bvpi_\la$, notation Chapter~\ref{VI}}%
$\bvpi_\la$ be say $0=0$.

\itlb{181c}
Use $\bmpp\la$ to define 
$\bbf\la\in\kond\bta$ and $\bbF\la\in\cfd_{\bta}$ as in 
\index{notation Chapter~\ref{VI}!$\bbf\la$}%
\index{zfla@$\bbf\la$, notation Chapter~\ref{VI}}%
\index{notation Chapter~\ref{VI}!$\bbF\la$}%
\index{zFla@$\bbF\la$, notation Chapter~\ref{VI}}%
Definition~\ref{152d}. 
\een

On the basis of this notation, our proof of Theorem~\ref{18t} 
will proceed as follows.  
We define the notion of {\ubf generic iterated perfect sets}, 
and prove 
the existence lemma and some properties of such sets in Section~\ref{182}. 
Then we {\ubf pick a generic set $Y_0\sq\tX$} 
in Section~\ref{183} and   
shrink it to a set $\tY\sq Y_0$ 
satisfying some {\ubf conditions related to} 
\ref{C2+}, 
\ref{C3+}, 
\ref{C4+}, 
\ref{C5+} 
of Definition~\ref{15e} above. 
The next step is {\ubf the lifting theorem} 
of Section~\ref{184}; it says roughly 
that any generic set in $\pel\i$ can be extended to a 
generic set in $\pele\i$.  
This theorem allows us to define a 
{\ubf rudiment $\rP\sq\pe\bta$}  
in Section~\ref{187}, 
of all sets 
$X\in\pele\bta$ 
{\ubf whose all projections $X\rsq\i$ are generic} 
(but not necessarily $X$ itself is such). 
This rudiment contains $\tY$ and refines $\rU_\la$. 
After a short but necessary work related to 
condition \ref{C1+} of Definition~\ref{15e}, 
we then take a suitable countable {\ubf sub-rudiment} 
of $\rP$ to be the 
layer $\rQ_\la$ for Theorem~\ref{18t}.


\parf{Generic perfect sets}
\las{182}  

We continue to argue in $\rL$. 
Consider the set $\homb=\lomb$,
\kmar{homb,lomb}
and define the following countable sets:
$$
\pC=
\index{notation Chapter~\ref{VI}!$\pC$}%
\index{zCfrak@$\pC$, notation Chapter~\ref{VI}}%
\bta\cup\rU_\la\cup
\ans{\la,\bta,\jqo,\omi,\hc}\sq \homb=\lomb;
$$
$$
\pD=\ans{\text{all sets $X\sq\homb$ $\in$-definable over 
$\homb$ with parameters in $\pC$}}.
\index{notation Chapter~\ref{VI}!$\pD$}%
\index{zDfrak@$\pD$, notation Chapter~\ref{VI}}%
$$

\bre
\lam{inC}
Such sets as $\omi,\hc,\pei,\cpei$, as well as 
many sets related to $\jqo$ this or another way, 
like $\rql\la,\rU_\la,\cX_\la, 
\sis{\bPsi k\la}{k<\om},\sis{\bPsd k\la}{k<\om},$ \etc\ 
belong to $\pD\cap\homb$, 
and can be used as parameters to define sets in $\pD$.
\ere

\bdf
\lam{Dgen}
Assume that $\et\in\cpo,\,\et\sq\bta$. 
A set $X\in\pe\et$ is \dd\pD\rit{generic} iff 
\index{set!Dgeneric@$\pD$-generic}%
\index{D-generic@$\pD$-generic set}%
\index{generic!$\pD$-generic set}%
$X\sqf\bigcup D$ holds 
for any set $D\in\pD$, $D\sq\rU_\la\dar\et$, dense in 
$\rU_\la\dar\et$.

Recall that 
$\rU_\la\dar\et=\ens{Y\dar\et}{Y\in\rU_\la}$. 
See Section~\ref{frel} on $\sqf,\sqd$.
\edf

\ble
\lam{egen}
If\/ $U\in\rU_\la$  then there is a\/ 
\dd\pD generic set\/ $X\in\pe\bta$, $X\sq U.$  
\ele
\bpf
Fix  any \dd{\bta}admissible map  $\phi:\om\onto\bta$. 
The next claim is a consequence of property \ref{nfr3}   
of the rudiment $\rU_\la$, the density, and 
Corollary~\ref{suz} applied consecutively enough many times:
\ben
\nenu
\itlb{1810}
If $m<\om$ and a set $D\in\pD$, $D\sq \rU_\la$, 
is dense in $\rU_\la$  
then any $\phi$-split system $\sis{X_u}{u\in 2^m}$  
of sets $X_u\in\rU_\la$ 
admits a narrowing $\sis{X'_{u}}{u\in 2^{m}}$    
in $\rU_\la$ such that $X'_u\in D$ for all $u\in 2^m.$ 
\een
Using \ref{1810} and the countability of $\pD$, 
we get a fusion sequence $\sis{X_u}{u\in 2^{<\om}}$  
of sets in $\rU_\la$, 
such that $X_\La\sq U$, and, for each  $D\in\pD$ 
dense in $\rU_\la$, there is $m<\om$  
with $X_u\in D$ for all $u\in 2^m.$  
Then $X=\bigcap_m\bigcup_{u\in2^m}X_u\in\pe\bta$, 
$X\sq U$, and $X\sqd\bigcup D$ for each set $D\in\pD$, 
$D\sq \rU_\la$, dense in $\rU_\la$. 
\epf

The next theorem provides some basic properties of 
\kmar{\ \mns\ \mns182T}
\dd\pD generic sets.

\bte
\label{182T}
\ben
\renu
\itlb{182a}
If\/ $X\in\pe\bta$ is\/ \dd\pD generic and\/ 
$\et\in\ft\bta$ 
{\em(an initial segment of finite type, Section~\ref{rud})} 
then\/ $X\dar\et$ is\/ \dd\pD generic as well$;$

\itlb{182b}
moreover, if, in\/ \ref{182a}, 
$D\in\pD$, $D\sq\rU_\la\dar\et$, 
$D$ is pre-dense in\/ $\rU_\la\dar\et$, then\/ 
$X\sqd\bigcup D;$

\itlb{182c}
if\/ $\al<\la$, $\i\in\tuq\al$, 
$X\in\pele\i$ is\/ \dd\pD generic 
then\/ $X\sqd\bigcup(\rqu\al\rsq\i);$ 

\itlb{182d}
if\/ $\et\in\ft\bta$, $U\in \rU_\la\dar\et$, and\/ 
$X \in\pe\et$ is\/ \dd\pD generic then\/ 
$X\cap U$ is clopen in\/ $X;$

\itlb{182e}
if\/ $\i\ekp\j$ belong to\/ $\bta$ and\/ 
$X\in\pele\i$ 
is\/ \dd{\pD}generic then so is 
$\pi_{\i\j}\akt X.$
\een 
\ete

\bpf
\ref{182a}
Assume that $D\in\pD$, $D\sq\rU_\la\dar\et$, is dense in 
$\rU_\la\dar\et$; prove that $X\dar\et\sqd\bigcup D$. 
It follows from property \ref{nfr3}   
of the rudiment $\rU_\la$ that the set 
$D'=\ens{U\in\rU_\la}{U\dar\et\in D}$ 
is dense in $\rU_\la$. 
Moreover $D'$ belongs to $\pD$ because so do $D$ and 
$\et\in\ft\bta$. 
(Not necessarily true for an arbitrary 
$\et\in\cpo,\et\sq\bta$.) 
Thus $X\sqf\bigcup D'$ by the genericity, hence 
$X\dar\et\sqf\bigcup D$.

\ref{182b}
Apply \ref{182a} for the dense set 
$D_1=\ens{V\in\rU\dar\et}{\sus U\in D\,(V\sq U)}$.

\ref{182c}
We know that $\rqu\al\rsq\i$ is predense in 
$\rU_\la\rsq\i$ by Theorem~\ref{142}\ref{1422}. 
It remains to apply \ref{182b} with $\et=[{\sq}\i]$.

\ref{182d}
Recall that $\rU_\la$ is a rudiment, hence it satisfies 
\ref{fr2} of Section~\ref{rud}. 
It easily follows that $\rU_\la\dar\et$ satisfies 
\ref{fr2} as well: 
if $\pu\ne Z\sq Y\in \rU_\la\dar\et$, $Z\in\pe\et$, 
and $Z$ is clopen in $Y$ then $Z\in\rU_\la\dar\et$. 
(Indeed if $Y=U\dar\et$, $U\in\rU_\la$, then 
$U'=U\cap (Z\uar\bta)\in\pe\bta$ by Lemma~\ref{apro}, 
and $U'$ is clopen in $U$ by the choice of $Z$ --- 
thus $U'\in\rU_\la$. 
But $Z=U'\dar\et$.) 
We conclude that the set $D$ of all 
$Y\in\rU_\la\dar\et$, satisfying 
$Y\sq U$ or $Y\cap U=\pu$, is dense in $U_\la\dar\et$. 
We conclude that $X\sqf\bigcup D$ by the genericity, 
in other words, $X\sq\bigcup D'$, where $D'\sq D$ is 
finite. 
Thus $D'=D'_1\cup D'_2$, where 
$D'_1=\ens{Y\in D'}{Y\sq U}$ and  
$D'_2=\ens{Y\in D'}{Y\cap U=\pu}$. 
Thus $X\sq Y_1\cup Y_2$, where 
$Y_e=\bigcup D'_e$ are two disjoint closed sets. 
Finally, $X\cap U=X\cap Y_1=X\bez Y_2$, which 
implies the result required.

\ref{182e}
This is clear as $\pi_{\i\j}\in\pD$.
\epf

\vyk{
Assume that $\et\in\cpo,\,\et\sq\bta$. 
A set $X\in\pe\et$ is \dd\pD\rit{generic} iff 
$X\sqd\bigcup D$ holds 
for any set $D\in\pD$, $D\sq\rU\dar\et$, dense in 
$\rU\dar\et$.
}

\parf{The choice of \boldmath $\tY$}
\las{183}  

Using Lemma~\ref{egen}, fix a \dd\pD generic set 
$Y_0\in\pe\bta$, $Y_0\sq \tX$. 
Using consecutively 
Lemma~\ref{111}\ref{111c1},  
Lemma~\ref{174L},  
Theorem~\ref{176L},  
and Theorem~\ref{175t},  
we obtain a set $\tY\in\pe\bta$, $\tY\sq Y_0\sq\tX$, 
\index{notation Chapter~\ref{VI}!$\tY$}%
\index{zYover@$\tY$, notation Chapter~\ref{VI}}%
\index{set!Yt@$\tY$}%
satisfying the following \ref{1811} -- \ref{1814}:
\ben 
\qenu
\atc
\atc
\itlb{1811}
$\tY\fo \bvpi_\la$ or\/ $\tY\fo \otr\bvpi_\la$;

\itlb{1812}
$\tY$ is uniform;

\itlb{1813}
either (a) ${\bbF\la}$ avoids $\tY\dir\i$ on $\tY$ 
for all $\i\in\bta$, 
or (b) 
${\bbF\la}$ is a\/ $\bj$-axis map on\/ $\tY$ 
for some\/ $\bj\in\bta$, and ${\bbF\la}$ avoids\/  
\index{notation Chapter~\ref{VI}!$\bj$}%
\index{zjb@$\bj$, notation Chapter~\ref{VI}}%
$\tY\dir\i$  on\/ $\tY$ 
for all $\i\in\bta\bez\ans\bj$;  

\itlb{1814}
the image\/ $S=\ima{{\bbF\la}\!}\tY$ is\/ $U$-avoidable 
on\/ $\i$ for all\/  $\i\in\bta$, $U\in\rU_\la\rsq i$.
\een

\bre
\lam{1815r}
The set $\tY\sq\tX$ is \dd\pD generic along 
with $\tX$, and hence 
\ben 
\qenu
\atc
\itlb{1815}
if $\la$ is limit, $k<\om$, 
$\bPsd k\la\sq\rqu{<\la}$, and 
$\bPsd k\la$ is dense in\/ $\rqu{<\la}$, 
then $\tY\sqf \bigcup\bPsd k\la$.  
\een
This needs some work. 
By the density assumption, the derived set 
$\Phi_{\la k}:= \bPsd k\la\uar\bta$ is dense in 
$\rU'=\rqu{<\la}\uar\bta$.
Howevere $\rU'$ itself is dense in $\rU_\la=\rh(\rU')$ 
by Theorem~\ref{142}\ref{142jd} --- here we use that $\la$ is limit. 
Thus $\Phi_{\la k}$ is dense in $\rU_\la$. 
It follows that $\tY\sq\tX\sqf\bigcup\Phi_{\la k}$, 
by the \dd\pD genericity. 
($\Phi_{\la k}\in\pD$ holds since 
$\bPsd k\la\in\pD$.)  
This implies $\tY\sqf \bigcup\bPsd k\la$ as well.
\ere

\bre
\lam{stra}
A certain oddity in the numbering above is 
caused by the fact that we want to indicate a connection 
with the numbering of items in Definition~\ref{15e}. 
Thus say \ref{1811} corresponds to condition \ref{C3+} 
in \ref{15e}, \etc\ 
In addition, \ref{1814} will assist \ref{1813} in getting 
to \ref{C5+} in \ref{15e}, whereas \ref{C1+} will be 
considered in Section~\ref{188} below by means not 
related to $\tY$.
\ere

\bre
\lam{bda}
Coming back to \ref{1813}, we may note that $\bj$ is 
unique in case (b) by \ref{1761} in the proof of 
Theorem~\ref{176L}. 
Moreover (a) and (b) are incompatible. 
(If (b) holds then take $\i=\bj$ in (a), 
getting a contradiction.)  
This allows us to define $\bda=\bta$ in case (a) of 
\index{notation Chapter~\ref{VI}!$\bda$}%
\index{zzdab@$\bda$, notation Chapter~\ref{VI}}%
\ref{1813}, and $\bda=\ens{\i\in\bta}{\i\not\ekp\bj}$ 
in case (b).
\ere

Let $\et\in\cpo,\et\sq\bta$.
Say that $Z\in\pe\et$ is a \rit{\dd\bda set} iff it 
is similar to $\tY$ in the sense that 
${\bbF\la}$ avoids $Z\dir\i$  on\/ $\tY$ 
for all $\i\in\bda\cap\et$.

\ble
\lam{bdaL}
$\tY$ is a\/ \dd\pD generic\/ \dd\bda set.
\ele
\bpf
$\tY$ is \dd\pD generic since $Y_0$ is such and 
$\tY\sq Y_0$. 
$\tY$ is a \dd\bda set by \ref{1813}.
\epf

\parf{Lifting theorem}
\las{184}  

Our further major goal will be to include $\tY$ in a 
suitable rudiment, by Corollary~\ref{1871} below. 
The following is the key technical result.

\bte
[in $\rL$]
\lam{184t}
Let\/ $\i\in\bta$, $U\in\rU_\la\rsq\i$, $X\in\pel\i$ 
be a\/ \dd\pD generic\/ \dd\bda set, and\/ $X\sq U\rsl\i$. 
Then there exists a\/ \dd\pD generic\/ \dd\bda set\/ 
$X'\in\pele\i$, 
such that\/ $X'\sq U$, $X'\rsl\i=X$. 
\ete

\bpf
This is a rather long argument. 
We fix $\i,U,X$ during the course of the proof. 
We can assume, by \ref{1814}, that
\ben
\fenu
\itlb{184*}
${\bbF\la}$ avoids $U\dir\i$ on $\tY$. 
\een

Let \rit{an atom} be any set of the form 
\index{set!atom}%
\index{atom}%
$V=W\cap(P\usq\i),$ where 
$\pu\ne P\sq X$ is clopen in $X$  
(then $P\in\pel\i$), 
$W\in\rU_\la\rsq\i$, $W\sq U$, and 
$P\sq W\rsl\i$. 
Let $\rQ$ = all finite non-empty 
unions of atoms. 
We claim that 
\ben
\Aenu
\itlb{184A} 
If $Q\in\rQ$ then $Q\rsl\i\sq X$ and $Q\rsl\i$ 
is clopen in $X$ 
(as a finite union of relatively clopen sets);

\itlb{184B} 
$\rQ\sq\pele\i$;

\itlb{184C} 
if $\pu\ne Q'\sq Q\in \rQ$, $Q'$ is clopen in $Q$, 
then $Q'\in\rQ$.

\itlb{184D} 
if $Q,Q'\in \rQ$, $\et\in\cpo$, $\et\sq[{\su}\i]$, 
$Q\dar\et\sq Q'\dar \et$, 
then the set $Q''=Q'\cap(Q\dar\et\usq\i)$ 
belongs to $\rQ$.
\een

To prove \ref{184B}, assume that 
$Q=V_1\cup\ldots\cup V_n\in \rQ$, where each 
$V_e=W_e\cap(P_e\usq\i)$ is an atom, 
so that  
$\pu\ne P_e\sq X$ is clopen in $X$ 
(then $P_e\in\pel\i$ is\/ \dd\pD generic),  
$W_e\in\rU_\la\rsq\i$, $W_e\sq U$, and 
$P_e\sq W_e\rsl\i$. 
We have $V_e\in\in\pel\i$ by Lemma~\ref{apro}, 
and obviously $V_e\rsl\i=P_e$.

Let $e=1,\dots,n$. 
Coming back to Section~\ref{closu}, put 
$\cT_e(x)=\der(\cs {V_e,}x\i)$ for 
all $x\in P_e=V_e\rsl\i$, so that 
$\cT_e:P_e\to\Der$ is continuous by Lemma~\ref{hh+}. 
Define the extended map $\cT'_e:X\to\Der$ by 
$\cT'_e(x):=\cT_e(x)$ for $y\in P_e$ and 
$\cT'_e(x):=\pu$ for $x\in X\bez P_e$. 
Then $\cT'_e$ is continuous since $P_e$ 
is clopen in $X$.  

We conclude that 
$\cT(x):=\cT_1(x)\cup\ldots\cup\cT_n(x):X\to\Der$ 
is continuous. 
It follows by Lemma~\ref{hh} that the set
\bce
$Q'=\ens{z\in\can{\ilq\i}}
{z\rsl\i\in P=P_1\cup\ldots\cup P_n\land z(i)\in[\cT(x)]}$
\ece
belongs to $\pele\i$. 
On the other hand easily $Q'=Q$.

It suffices to prove \ref{184C} in case when 
$Q=W\cap(P\usq\i)$ is an atom, so that  
$\pu\ne P\sq X$ is clopen in $X,$ 
$W\in\rU_\la\rsq\i$, $W\sq U$, and 
$P\sq W\rsl\i$. 
By Lemma~\ref{aclo} we have 
$Q'=W'\cap(P'\usq\i)$, 
where $W'\sq W$ and $P'\sq P$ are relatively clopen 
and still $P\sq W\rsl\i$. 
Thus $Q'$ is an atom as well.

To prove \ref{184D} note that the sets 
$Q\dar\et$, $Q'\dar \et$ are clopen in $U\dar\et$ 
by Lemma~\ref{darop}. 
Thus $Q''$ is clopen in $Q'$. 
It remains to refer to \ref{184C}. 

\ble
\lam{185}
Let\/ $\cY\in\pD$, $\cY\sq\rU_\la\rsq\i$, $\cY$ is dense 
in\/ $\rU_\la\rsq\i$, and\/ $Q\in\rQ$. 
Then there is\/ $Q'\in \rQ$, $Q'\sq Q$, such that\/ 
$Q'\rsl\i=Q\rsl\i$ and\/ $Q'\sqf\bigcup\cY$.
\ele

\bpf[Lemma]
We \noo\ assume that 
$Q=W\cap(P\usq\i)$ is an atom, where 
$\pu\ne P\sq X$ is clopen in $X,$  
$W\in\rU_\la\rsq\i$, $W\sq U$,  
$P\sq W\rsl\i$. 
Then
\ben
\nenu
\itlb{185a}
$Q\rsl\i=P\sq X$ and $Q\rsl\i$ is clopen in $X$;

\itlb{185b}
$Q\rsl\i=P\sq X$ is a \dd\pD generic \dd\bda set 
(because such is $X$). 
\een
We claim that the set\pagebreak[1]
$$
\cY_1\,=\,\ens{A\rsl\i}{A\in\cY\land A\sq W}\,\cup\,
\ens{Z\in\rU_\la\rsl\i}{Z\cap W\rsl\i=\pu}
$$
is dense in $\rU_\la\rsl\i$. 
Indeed let $S\in\rU_\la\rsl\i$; we have 
to find $Z\in\cY_1$, $Z\sq S$. 

{\ubf Case 1:} 
$S\sq W\rsl\i$. 
Then the set $W'=W\cap(S\usq\i)$ belongs to 
$\rU_\la\rsq\i$ as $\rU_\la$ is a rudiment.  
Thus, by the density of $\cY$, there is a set 
$A\in\cY$, $A\sq W'$. 
Then $Z=A\rsl\i\in\cY_1$ is as required.

{\ubf Case 2:} $S'=S\bez (W\rsl\i)\ne\pu$. 
Then there is a set $\pu\ne Z\sq S'$ clopen in $S$. 
As $\rU_\la$ is a rudiment, we have 
$Z\in\rU_\la\rsl\i$. 
Thus $Z\in\cY_1$, as required. 

The density of $\cY_1$ is established. 
As obviously $\cY_1\in\pD$, it follows that 
$P\sqf\bigcup\cY_1$ by \ref{185b}, hence 
$P\sq Z_1\cup\ldots\cup Z_m$, $Z_e\in\cY_1$, $\kaz e$. 
By the choice of $P$, we can \noo\ assume that each 
$Z_e$ belongs to the first part of $\cY_1$, \ie, 
$Z_e=A_e\rsl\i\in\rU_\la\rsl\i$, 
where $A_e\in \cY$, $A_e\sq W$. 
Let $P_e=P\cap Z_e$. 

As $P$ is \dd\pD generic, 
each $P_e$ is clopen in $P$ by 
Theorem~\ref{182T}\ref{182d}, 
and hence clopen in $X$ by \ref{185a}. 
It follows that each $V_e=A_e\cap (P_e\usq\i)$ 
is an atom (or $\pu$). 
Therefore $Q'=V_1\cup\ldots\cup V_n\in\rQ$, 
$Q'\sqf\bigcup\cY$ 
(as each $A_e$ belongs to $\cY$), 
and $Q'\rsq\i=P_1\cup\ldots\cup P_n=P=Q\rsq\i$,
as required.

\epF{Lemma~\ref{185}}

To proceed with another lemma, we fix a 
$\ilq\i$-admissible function $\phi\in\pD$, 
$\phi:\om\to\ilq\i$ 
(meaning that if $\j\sq\i$ then $\phi(k)=\j$ 
for infinitely many $k$).

\ble
\lam{186}
Let\/ $n<\om$, and\/ $\sis{Y_s}{s\in2^n}$ be a 
system of sets\/ $Y_s\in\rQ$, satisfying\/ \ref{prct} 
of Definition~\ref{splis} with\/ $\za=\ilq\i$. 
Let\/ $\cY\in\pD$, $\cY\sq\rU_\la\rsq\i$, 
$\cY$ be dense in\/ $\rU_\la\rsq\i$. 
Then there is a system\/ $\sis{Q_s}{s\in2^n}$ of 
sets\/ $Q_s\in\rQ$, $Q_s\sq Y_s$,   
satisfying\/ \ref{prct} and\/ 
$Q_s\rsl\i=Y_s\rsl\i$, $Q_s\sqf\bigcup\cY$ 
for all\/ $s\in2^n.$  
\ele

\bpf
If $s\in 2^n$ then, by Lemma~\ref{185}, 
pick a set $Q_s\in \rQ$, $Q_s\sq Y_s$, such that\/ 
$Q_s\rsl\i=Y_s\rsl\i$ and\/ $Q_s\sqf\bigcup\cY$. 
The system $\sis{Q_s}{s\in2^n}$ still 
satisfies \ref{prct} (with $\za=\ilq\i$) 
because if $s\ne t$ belong to 
$2^n$ then $\za_\phi[s,t]\sq\ile\i$, 
hence 
$Q_s\dar\za_\phi[s,t]=Y_s\dar\za_\phi[s,t]=
Y_t\dar\za_\phi[s,t]=Q_t\dar\za_\phi[s,t]$.
\epF{Lemma~\ref{186}}

{\ubf Finalization.} 
Now we are able to accomplish the proof of 
Theorem~\ref{184t}. 
We define a \dd\phi fusion sequence 
$\sis{Q_u}{u\in\bse}$ of sets $Q_u\in\rQ$ 
(still with $\za=\ilq\i$ in Definition~\ref{fuzD}) 
satisfying 
\ben
\nenu
\itlb{184f1} 
$Q_\La= U\cap (X\usq{\i})$ --- this is even an atom 
by the choice of $U,X$ in Theorem~\ref{184t};

\itlb{184f2} 
if $\cY\in\pD$, $\cY\sq\rU_\la\rsq\i$, 
$\cY$ is dense in\/ $\rU_\la\rsq\i$, then there is 
$m<\om$ such that $Q_u\sqf\bigcup\cY$ 
for all\/ $u\in 2^m$;

\itlb{184f3} 
if $m<\om$ then 
$\bigcup_{u\in2^m}Q_u\rsl\i=Q_\La\rsl\i=X$.
\een
Namely suppose that a layer $\sis{Q_u}{u\in2^m}$ 
has been defined so that both \ref{prct},\ref{bprct} 
of Definition~\ref{splis}, and \ref{184f3}, hold for 
this $m$. 
Let $Y_{u\we e}=\spli{(Q_u)}\i e$ for all $u\in2^m$ 
and $e=0,1$, where $\i=\phi(m)$, so that 
$\sis{Y_s}{s\in2^{m+1}}$ is a clopen expansion 
of $\sis{Q_u}{u\in2^m}$ by Lemma~\ref{pand}. 
Each $Y_s$ belongs to $\rQ$ by \ref{184C} above. 
Lemma~\ref{186} yields a system 
$\sis{Q_s}{s\in2^{m+1}}$ of sets $Q_s\in\rQ$, 
$Q_s\sq Y_s$, satisfying \ref{prct} and  
$Q_s\rsl\i=Y_s\rsl\i$, $Q_s\sqf\bigcup\cY$ 
for all\/ $s\in2^{m+1},$  as required.

Having an \ref{184f1}-\ref{184f2}-\ref{184f3}  
fusion sequence in $\rQ$, we define 
$X'=\bigcap_m\bigcup_{u\in2^m}Q_u$. 
Then $X'\in\pele\i$ by Theorem~\ref{fut}, 
$X'\sq Q_\La=U$ by construction, 
$X'\rsl\i=X$ by \ref{184f3}, 
$X'$ is \dd\pD generic by \ref{184f2}. 

Further, ${\bbF\la}$ avoids $X'\dir\i$ on $\tY$ 
by \ref{184*} and because $X'\sq U$. 
Moreover, if $\j\su\i$ and $\j\in\bda$ then 
${\bbF\la}$ avoids $X'\dir\j$ on $\tY$ since 
$X'\rsl\i=X$ and $X$ is a \dd\bda set. 
Thus overall ${\bbF\la}$ avoids $X'\dir\j$ on $\tY$ 
for every $\j\in\ilq\i\cap\bda$, and hence 
$X'$ is a \dd\bda set, as required.
\epF{Theorem~\ref{184t}}

\parf{Consequences of the lifting theorem}
\las{187}  

Consider the system\/ $\cK=\sis{\cK_\i}{\i\in\bta}$ 
\index{notation Chapter~\ref{VI}!$\cK=\sis{\cK_\i}{\i\in\bta}$, system}%
\index{zKi@$\cK=\sis{\cK_\i}{\i\in\bta}$, notation Chapter~\ref{VI}}%
of sets 
$$
\cK_\i=\ens{X\in\pele\i}
{X\text{\rm\ is a \dd\pD generic \dd\bda set}}.
$$

\bcor
\label{187c}
\ben
\renu
\itlb{187c1}
\kmar{\ \mns 187c}
Let\/ $\j\su\i$ belong to\/ $\bta$, 
$U\in\rU_\la\rsq\i$, $X\in\cK_\j$, $X\sq U\rsq\j$. 
Then there is a set\/ $X'\in\cK_\i$, 
$X'\sq U$, such that\/ $X'\rsq\j=X;$ 

\itlb{187c2}
in particular, with\/ $U=\cai\i$, if\/ $X\in\cK_\j$ 
then there is a set\/ $X'\in\cK_\i$ such that\/ 
$X'\rsq\j=X;$ 

\itlb{187c3}
the system\/ $\cK=\sis{\cK_\i}{\i\in\bta}$
is a\/ $\bta$-kernel.
\een
\ecor
\bpf
\ref{187c1} is an immediate corollary of 
Theorem~\ref{184t} 
(applied by induction on $\lh\i-\lh\j$), 
with \ref{187c2} being a particular case 
of \ref{187c1}.

To prove \ref{187c3}, note that 
\ref{187c2} implies 
\ref{ke1s} of Section~\ref{ker} for $\cK$. 
Condition \ref{ke2} in Section~\ref{ker} is obvious, 
whereas \ref{ke3}, \ref{ke4} hold  
because the property of being 
a \dd\pD generic \dd\bda set 
is transferred to all smaller sets still in $\pei$. 
(Note that $Z$ in \ref{ke3} and $Y$ in \ref{ke4} 
belong to $\pele\i$ by Lemma~\ref{apro}, resp., 
Lemma~\ref{lin}.) 
Finally \ref{ke6} holds because all notions related 
to the property of being a \dd\pD generic \dd\bda set 
are invariand under the action of $\pi_{\i\j}$ 
because $\pi_{\i\j}\in\pD$.
\epf


Following Section~\ref{rud}, we consider the rudiment 
\bce
$\rP=\rP(\cK):=\ens{X\in\bta=\tuq\la}
{\kaz\i\in\bta\,(X\rsq\i\in\cK_\i)}\in\bfr\la$.%
\ece
\index{notation Chapter~\ref{VI}!$\rP=\rP(\cK)$, rudiment}%
\index{zPPK@$\rP=\rP(\cK)$, notation Chapter~\ref{VI}}%

\bcor
\label{1871}
\ben
\renu
\itlb{1871a}
$\rP\in\rud_\la$  
\kmar{\ \mns1871}
and\/ $\pro\rP\i=\cK_\i$ for all\/ $\i\in\bta;$ 

\itlb{1871b}
$\rP$ is a refinement of\/ $\rU_\la:$ 
$\rU_\la\ssw\rP$ in the sense of Section~\ref{ref}$;$

\itlb{1871c}
 $\tY\in\rP$. 
\een
\ecor

\bpf
\ref{1871a}
holds by Lemma~\ref{r2n}.%
\pagebreak[1]

\ref{1871b}
We have to check \ref{134i}, \ref{134ii}, \ref{134iii} 
of Section~\ref{ref}. 

Of them, \ref{134i} (\ie, $\can\bta\in\rU_\la$) holds by 
Theorem~\ref{142}\ref{142i}.

To prove \ref{134ii}, assume that 
$\et\in\ft\xi$, $U\in\rU_\la$, $Y\in\rP$, 
$Y\dar\et\sq U\dar\et$, and the goal is to find 
$Z\in\rP$ satisfying $Z\sq U$ and $Z\dar\et=Y\dar\et$. 
For that purpose, we define a system of sets 
$X_\i\in\cK_\i$, $\i\in\bta$, such that 
\ben
\aenu
\itlb{186a}
$X_\i=Y\rsl\i$ for all $\i\in\et$;

\itlb{186b}
$X_\i\sq U\rsq\i$ for all $\i$; 

\itlb{186c}
if $\j\su\i$, $\lh\i=\lh\j+1$, then 
$X_\i\rsl\i=X_\j$. 
\een
The construction goes on as follows. 
Assume that $\j\su\i$ in $\bta$, $\lh\i=\lh\j+1$, 
$\i\nin\et$, and a set $X_\j\in\cK_\j$, 
$X_\j\sq U\rsq\j=U\rsl\i$ 
has been defined. 
Use Corollary~\ref{187c}\ref{187c1} to get a set 
$X_\i\in\cK_\i$, $X_\i\sq U\rsq\i$, with 
$X_\i\rsl\i=X_\j$. 

After the construction of sets $X_\i\in\cK_\i$ 
satisfying \ref{186a},\ref{186b},\ref{186c} is 
accomplished, the set 
$Z=\ens{x\in\can\bta}{\kaz\i\in\bta\,(x\dar\i\in X_\i)}$ 
is as required for \ref{134ii}. 

To prove \ref{134iii}, assume that 
$\i\in\bta$, $U\in\pro{\rU_\la}\i$, $Y\in\pro\rP\i$. 
Then $U\cap Y$ is clopen in $Y$ by 
Theorem~\ref{182T}\ref{182d}, as required.

\ref{1871c}
As $\tY$ is \dd\pD generic by 
Lemma~\ref{bdaL}, we conclude that 
each $\tY\rsq\i$ is \dd\pD generic 
as well by Theorem~\ref{182T}\ref{182a}. 
And $\tY\rsq\i$ is a \dd\bda set since such is $\tY$ 
itself still by Lemma~\ref{bdaL}.
\epf

\parf{The construction of a sub-rudiment}
\las{188}  

We know that the set $\tY$ chosen in Section~\ref{183} belongs 
to $\rP$ by Corollary~\ref{1871}\ref{1871c}. 
Here we define another special set $\tyj\in\rP$, related rather to 
condition \ref{C1+} of Definition~\ref{15e}, and then define 
a set $\rP'$ required, in the form of  
a countable sub-rudiment of $\rP$ containing both $\tY,\tyj$. 
In some similarity to \ref{181a} of Section~\ref{181}, 
we first define $\txj$ as follows:
\ben 
\qenu
\itlb{1818}
if $\clo{\bM\la}\in \noc{\jqo\res\la}$ and 
$\dym{\clo{\bM\la}}\sq\bta=\tuq\la$  
then let $\txj=\clo{\bM\la}\uar\bta$, otherwise let 
$\txj=\can\bta$, so 
$\txj\in \noc{\jqo\res\la}\yi \dym{\txj}=\bta$ in both cases.
\een

\bcor
[of Lemma~\ref{r2n}]
\lam{188a}
If\/ $\i\in\bta$ then\/ $\txj\rsq\i\in\rU_\la\rsq\i$. \qed
\ecor

Note that $\txj$ as a whole is not assumed to belong 
to $\rU_\la$!

\bcor
[of Corollaries~\ref{188a} and \ref{187c}\ref{187c1}]
\lam{188b}
There is a system of sets\/ $Y_\i\in\cK_\i\yi \i\in\bta$ such that\/ 
$Y_\i\sq\txj\rsq\i$ and if\/ $\j\su\i$ then\/ $Y_\j=Y_i\rsl\j$. \qed
\ecor

Recall that $\cK_\i$ and $\rP$ were defined in Section~\ref{187}.

\bcor
[of Corollary~\ref{pe5}]
\lam{188c}
There is a set\/ $\tyj\in\pe\bta$ such that\/ 
\index{notation Chapter~\ref{VI}!$\tyj$}%
\index{zY1over@$\tyj$, notation Chapter~\ref{VI}}%
\index{set!Y1t@$\tyj$}%
$\tyj\rsl\i=Y_\i$ for all\/ $\i\in\bta$.  
Note that then\/ $\tyj\in\rP$ as\/ $Y_\i\in\cK_\i\yi\kaz\i$.\qed
\ecor

To conclude, we  have got a set $\tyj\in\rP$ satisfying $\tyj\sq\txj$ 
(because $\tyj\rsl\i=Y_\i\sq\txj\rsq\i$). 
Recall that $\tY\in\rP$ and $\rU_\la\ssq\rP,$ by Corollary~\ref{1871}.

\ble
\lam{188d}
There is a {\bfit countable} sub-rudiment\/ $\rP'\sq\rP$ 
still containing\/ $\tY,\tyj$ and satisfying\/ $\rU_\la\ssq\rP'$.
\ele
\bpf
A routine ``elementary substructure'' argument.
\epf

\ble
\lam{188e}
Such a\/ $\rP'$ is a\/ 1-5-$\nn$ extension of\/ 
$\jqo=\sis{\rqu\al}{\al<\la}$.
\ele
\bpf
Basically we have to check \ref{C*}, \ref{C1+}, and 
\ref{C1++} (including \ref{C2+}----\ref{C5+}) of 
Definition~\ref{15e} for $\rQ_\la:=\rP'$. 

\ref{C1+}
Suppose that
$\clo{\bM\la}\in \noc{\jqo\res\la}$ and 
$\dym{\clo{\bM\la}}\sq\bta=\tuq\la$.  
Thus $\txj=\clo{\bM\la}\uar\bta$ by \ref{1818} above. 
However $\tyj\in\rP'$ and $\tyj\sq\txj$ by construction, 
and this completes the proof of \ref{C1+}.

\ref{C2+} --\ref{C5+}. 
In accordance to Definition~\ref{15e}, we assume that 
$\clo{\bM\la}\in \rqu{<\la}:=\bigcup_{\ga<\la}\rqu\ga$ 
--- then $\tX=\clo{\bM\la}\uar\bta$ by \ref{181a} of 
Section~\ref{181} ---   
and the goal is to find  a set $Y\in\rqu\la$ satisfying both 
$Y\psq\clo{\bM\la}$  
and  each of \ref{C2+},\ref{C3+},\ref{C4+},\ref{C5+}. 
Let's chech that the set $\tY$ defined in Section~\ref{183}
is as required. 
First of all, note that $\tY\in\rP'$ and 
$\tY\sq\tX\psq\clo{\bM\la}$ by construction. 
It remains to check \ref{C2+} --\ref{C5+} 
of Definition~\ref{15e} for $\tY$.

\ref{C2+} 
Suppose that if $\la$ is limit, $k<\om$, 
$\bPsd k\la\sq\rqu{<\la}$, and $\bPsd k\la$ is dense 
in\/ $\rqu{<\la}$.
Then $\tY\sqd \bigcup\bPsd k\la$ holds  by \ref{1815}, 
as required.

\ref{C3+} and \ref{C4+} are immediate corollaries 
of \ref{1811}, \ref{1812}.

\ref{C5+}
This is not so straightforward. 
First of all we claim that 
\ben
\fenu
\itlb{188*}
if $Z\in\rP$ and $\i\in\bda$ then ${\bbF\la}$ avoids 
$Z\dir\i$ on $\tY.$ 
\een
Indeed $Z'=Z\rsq\i\in\cK_\i$ by 
Corollary~\ref{1871}\ref{1871a}, meaning that $Z'$ 
is a \dd\bda set. 
It follows that ${\bbF\la}$ avoids $Z\dir\i=Z'\dir\i$  on $\tY$ 
because $\i\in\bda$, as required. 

{\ubf Case 1:} 
(a) of \ref{1813} in Section~\ref{183}. 
Then $\bda=\bta$, and hence ${\bbF\la}$ avoids 
$Z\dir\i$ on $\tY$ for all $\i\in\bta=\tuq\la$ by 
\ref{188*}. 
Thus we have \ref{C5+}(a) of Definition~\ref{15e}. 

{\ubf Case 2:} (b) of \ref{1813} in Section~\ref{183}. 
Then accordingly $\bda=\ens{\i\in\bta}{\i\not\ekp\bj}$ 
(see Remark~\ref{bda}) for some $\bj\in\bta$ as in 
\ref{1813}(b). 
In other words, ${\bbF\la}$ is a 
$\bj$-axis map on $\tY$, and ${\bbF\la}$ avoids\/  
$Z\dir\i$  on $\tY$ 
for all $\i\in\bta$, $\i\not\ekp\bj$,  by \ref{188*}. 
Thus we have \ref{C5+}(b) of Definition~\ref{15e}, 
as required.  
\epF{Lemma~\ref{188e}}

\qeD{Theorem~\ref{18t}}

\sekt[\ \ The final forcing construction]
{The final forcing construction}
\las{VII}

Theorem~\ref{18t} obviously allows to define, in $\rL$, 
a \ruds\ $\jqo=\sis{\rqu\la}{\la<\omi}$ of length $\omi$, 
such that each term $\rqu\la$ is a 
{\ubf 1-5-${\nn}$ extension} 
of the subsequence $\sis{\rqu\al}{\al<\la}$, for a given 
$\nn\ge1$ of Theorem~\ref{mt1d}. 
Our next and the final step in the proof of 
Theorems \ref{mt1d}--
\ref{mt1b}--\ref{mt1} 
will be to maintain such a construction so that the 
{\ubf global definability condition} \refn{C6}\nn\ also holds. 

The content of this Chapter will mainly be the assessment 
of the {\ubf complexity} of different constructions related to the 
concept of 1-5-${\nn}$ extension. 
To evaluate complexity in terms of definability over $\hc$, 
we'll have to code various uncountable objects considered by 
sets in $\hc$, and evaluate the complexity of the 
{\ubf coding of some principal relations}. 
This line is concluded by Theorem~\ref{199t} which asserts 
that the notion of 1-5-${\nn}$ extension is 
{\ubf essentially $\id\hc\nn$ in the codes}. 
Therefore the {\ubf ``G\"odel-least''} choice of a (code of a)
1-5-${\nn}$ extension at each step still results is 
a \ruds\ satisfying \refn{C6}\nn. 
This will complete 
{\ubf the proof of Theorem~\ref{mt1d}} and thereby 
{\ubf Theorem~\ref{mt1}}, our first main result,  
in Section~\ref{200}.

{\ubf We argue in $\rL$ in this Chapter.}

\vyk{
\parf{Rank-bounded $\is{}1$ definability}
\las{191}

It is well-known that the class $\is\hc1$ of \paf\ $\is{}1$ definability   
over $\hc=\lomi$ 
is substantially wider than the class  $\id\hc1=\is\hc1\cap\ip\hc1$. 
Nevertheless there are meaningful subclasses of $\is\hc1$ included 
in $\id\hc1$, for instance, the class $\is11$ of the effective projective 
hierarchy. 
The goal of this section is to define a bigger subclass of $\is\hc1$, 
still included in $\id\hc1$.
 
Recall that $\is{}1$ formulas are \dd\in formulas of the form 
$\sus y\,\vpi$, 
where $\vpi$ is a \rit{bounded} formula, \ie, it is required that every 
quantifier in $\vpi$ is \rit{bounded}, 
that is, has the form $\sus a\in b$ or $\kaz a\in b$.  
The next definition intends to somehow restrict the domain of the 
unbounded quantifier of a $\is{}1$ formula. 

Along with the ordinary supremum $\tsup X$, we let $\ssup X$ be 
\kmar{ssup}
\rit{the strict supremum}, \ie, the least (in the order considered) 
element strictly bigger than all $x\in X$ (if exists). 
Recall that $\rank x\in\Ord$ is the set-theoretic rank of a set 
$x$, defined so that 
$\rank x = \ssup_{y\in x}\rank y= \tsup_{y\in x}(\rank y+1)$. 
In particular, $\rank\pu=0$, $\rank\om=\om$, $\rank x<\omi$ for 
all $x\in\hc$. 

Let $\hcr\al=\ens{x\in\hc}{\rank x<\al}$; 
\kmar{hcr al}
this can be viewed as the 
von Neumann class $\rV_\al$ \rit{within $\hc$}. 
In particular $\pws\om\sq\hcr{\om+1}$.

\bdf
\lam{191d}
Let $\susr\al x\,\Phi(x)$ mean 
$\sus x \big(\rank x<\al\land\Phi(x)\big)$, 
and let $\kazr\al x\,\Phi(x)$ be understood similarly. 
\kmar{susr al x}
Say that a set $X\sq\hc$ is \rit{$\isr 1$ over $\hc$} if there exists a 
\kmar{isr 1}
\paf\ $\is{}1$ formula $\sus y\,\Phi(x,y)$ ($\Phi$ being bounded), 
such that $X=\ens{x\in\hc}{\hc\mo \susr{\rank x+\om}y\,\Phi(x,y)}$. 
The superscript {\rm\tt rb} means \rit{rank-bounded}.
\edf

\bte
\lam{191t}
Every set\/ $\isr 1$ over $\hc$ is\/ $\id\hc1$.
\ete

The plan of the proof is to present a given set $Y\sq\hc$, 
$\isr 1$ over $\hc$, in the form $Y=\ima SZ$, where 
$Z\sq{\cN}$ is a $\id12$ set (hence $\id\hc1$) 
whereas $S$ is a map still $\isr 1$ over $\hc$. 
We make use of a rather well-known technique of coding 
sets by \wf\ relations. 

If $f\in{\cN}$ then let 
$\rep f=\ens{\ang{k,l}}{f(2^k\cdot3^l)=1}$,  
\kmar{rep f\mns abc f}
$\abc{f}=\dom {\rep f}\cup\ran {\rep f}$,  
\bce
$
p_f(t)=\ens{k<\om} 
{\sus k_0=k\rep f k_1 \rep f\dots\rep f k_m=t}
$ \ 
for all $t\in \abc{f}$ 
\ece
(all hereditary $\rep f$-predecessors of $t$).
We let ${f \ser t}$ be the characteristic function 
\kmar{f ser t}
of the set 
$\ens{2^k\cdot 3^l}{k\rep f l\land l\in\ans t\cup p_f(t)}$. 
If there is unique $t\in\abc{{\rep f}}$ satisfying 
$f={f \ser t}$ 
(or equivalently, $\abc{f}=p_f(t)$) 
then we put $t=\vrh{f}$. 
\kmar{vrh}

Let $\chc$ (\rit{codes} of sets in $\hc$) 
\kmar{chc}
contain all $f\in{\cN}$ satisfying 
\ref{191a},\ref{191b},\ref{191c}, where \ref{191a} is 
\rit{the extensionality}  
and \ref{191c} \rit{the well-foundedness} conditions.
\ben
\Renu
\itlb{191a}
if $k,k'\in \abc{f}$ and $\kaz l\,({l\rep f k}\eqv{l\rep f k'})$ 
then $k=k'$;

\itlb{191b}
$\vrh f$ exists; 

\itlb{191c}
if $\pu\ne u\sq\abc{f}$ then there is $k\in u$ such that 
$\kaz l\in u\,\neg\,(l\rep f k)$.
\een
If $f\in\chc$ then using \ref{191b} we can define 
$S_f(k)=\ens{S_f(l)}{l\rep f k}\in\hc$ by induction, and 
then $S(f)=S_f(\vrh f)\in\hc$. 
If $f,g\in\chc$ then:
\bit
\item
say that $h\in{\cN}$ \rit{translates\/ $f$ to\/ $g$} 
if $\ima h{\abc f}=\abc g$ and the equivalence 
${k\rep f l}\eqv{h(k)\rep g h(l)}$ holds for all 
$k,l\in\abc f$; 

\item
define $f \bra g$ 
\kmar{bra} 
iff there is $h\in\bI$ which translates $f$ to $g$; 

\item
define $f \bel g$ 
\kmar{bel} 
iff there is $m\in\abc g$ such that 
$m\rep g \vrh g$ and  $k\bra {g \ser m}$.
\eit

\ble
[standard]
\label{191L}
\ben
\renu
\itlb{191L1}
\sloppy
$S:\chc\onto\hc$ is a\/  $\id\hc1$ map, 
whereas\/  
$\chc\sq{\cN}$ is a\/ $\ip11$ set$;$

\itlb{191L2}
if\/ $f,g\in\chc$ then\/ 
${f\bra g}\eqv {S(f)=S(g)}$,  
${f\bel g}\eqv {S(f)\in S(g)};$ 

\itlb{191L3}
if\/ $f\in\chc$ and\/ $k\in\abc f$ then\/ 
$f\ser k\in\chc$ and\/ $S_f(k)=S(f\ser k);$ 

\itlb{191L4}
%
if\/ $f\in\chc$ then\/ 
$S(f)=\ens{S_f(k)}{k\rep f\vrh f}
.$ 
\qed
\een
\ele

The following is the key lemma in the proof of 
Theorem~\ref{191t}.

\ble
\label{191L'}
There is a\/ $\is11$ formula\/ $\Ta(f,g)$ such that
\ben
\renu
\itlb{191L'1}
if\/ $f\in{\cN}$ satisfies \ref{191a} and\/ \ref{191b},  
$g\in\chc$, and\/ $\Ta(f,g)$, then\/ 
$f\in\chc;$

\itlb{191L'2}
if\/ $f,g\in\chc$ then\/ 
$\Ta(f,g)\eqv{\rank S(f)<\rank S(g)+\om}.$
\een
\ele
\bpf
{\ubf Step 1}. 
Let $\Phi(f,g)$ be the $\is11$ formula saying:
\begin{quote}
there exists $h\in{\cN}$ such that $\ima h{\abc f}\sq g$, 
$k\rep f l\imp h(k)\rep g h(l)$, and if $k\in\abc f$ then 
$h(k)=\ssup_{\rep g}\ens{h(l)}{l\rep f k}$, and the last 
equality means that the $\ssup$ exists and is equal to $h(k)$.
\end{quote}
Standard arguments show that 1) if $f\in{\cN}$ 
satisfies \ref{191a} and\/ \ref{191b},  
$g\in\chc$, and\/ $\Phi(f,g)$, then\/ 
$f\in\chc$ as well, and 2) if both $f,g\in\chc$ then 
$\Phi(f,g)$ iff $\rank{S(f)}\le\rank{S(g)}$.\vom

{\ubf Step 2}. 
For any $g\in\chc$ we are going to define $\xg\in\chc$ 
such that $\rank{S(\xg)}=\rank{S(g)}+\om$. 
Let $g\in{\cN}$ satisfy \ref{191b} above, so that 
$t=\vrh g$ exists. 
The values of $\xg$ are defined as follows:\vom 
\ben
\aenu
\itlb{191*a}
if $k\rep g l$ then $2k\rep{\xg} 2l$, and let $\wt=2t$,

\itlb{191*b}
$\wt\rep{\xg} (2l+3)$, and 
$(2l+2m+3)\rep{\xg} (2l+2m+5)$,  
$(2l+2m+3)\rep{\xg} (2l+1)$ for all $m<\om$, 

\itlb{191*c}
if $K\rep{\xg} L$ then this is only via 
either \ref{191*a} or \ref{191*b}.
\een
Then clearly $g\mto\xg$ is an arithmetic map, and if 
$g\in\chc$ then $\xg\in\chc$ 
and $\rank{S(\xg)}=\rank{S(g)}+\om$. \vom

{\ubf Step 3}. 
To finalize the proof of the lemma, we let $\Ta(f,g)$ be the 
formula ``$f,g$ satisfy \ref{191b} above, and $\Phi(f,\xg)$''.
Then $\Ta$ is as required by the results at Steps 1,2.
\epF{Lemma~\ref{191L'}}

\bpf[Theorem~\ref{191t}]
Let 
$Y=\ens{y\in\hc}{\hc\mo\susr{\rank y+\om}x\,\gvpi(x,y)},$  
where $\gvpi$ is a \paf\ bounded $\in$-formula. 

If $\psi$ is a bounded $\in$-formula then define 
$\wpsi$ by induction as follows:
\ben 
\Aenu
\itlb{hat1}
$\widehat{x=y}$ is $x\bra y$ and 
$\widehat{x\in y}$ is $x\bel y$;

\itlb{hat2}
$\widehat{A\land B}$ is $\widehat{A}\land\widehat B$, 
and the same for $\lor,\neg,\imp,\eqv$;  

\itlb{hat3}
if $\psi$ is $\sus a\in b\,\chi(a,b)$ then $\wpsi$  
is 
$\sus m\,(m\rep b\vrh b\land \widehat\chi(b\ser m,b)$;  

\itlb{hat4}
if $\psi$ is $\kaz a\in b\,\chi(a,b)$ then $\wpsi$ 
is 
$\kaz m\,(m\rep b\vrh b\imp \widehat\chi(b\ser m,b)$.  
\een
Then, by construction and Lemma~\ref{191L}, we have 
for any bounded $\psi$:
\ben
\fenu
\itlb{191*}
 $\wpsi$ is a $\is0\iy$ (arithmetic) formula;

\itlb{191**}
if $f,g\in\chc$ then \ 
$\wpsi(f,g)\eqv \big(\hc\mo\psi (S(f),S(g))\big)$.
\een
It follows by \ref{191**} and Lemma~\ref{191L'} that 
\addtocounter{equation}{-2}
\busq
{191E3}
{%
\left.
\bay{rcl}
Y&=&\ens{y\in\hc}
{\sus g\in\chc\,\big(y=S(g)\land g\in G\big)}\\[0.5ex]
&=&\ens{y\in\hc}{\kaz g\in\chc\,\big(y=S(g)\imp g\in G\big)},
\eay
\right\}
}
where 
\busq
{191E4}
{G=\ens{g\in\bI}
{\underbrace{g\in\chc}_{\ip11}
\;\land\;\sus f\in\bI 
\underbrace{\big(\Ta(f,g)\land\wgvpi(f,g)\big)}_
{\is11\text{ by \ref{191*} as $\Ta$ is }\is11}}
}%
is a $\id12$ set, hence a $\id\hc1$ set. 
It follows then by the first line of \eqref{191E3} 
that $Y$ is $\is\hc1$ since 
$\chc$ is $\ip11$ (therefore $\id\hc1$), $G$ is $\id\hc1$ by 
\eqref{191E4}, and $S$ is $\id\hc1$ by Lemma~\ref{191L}. 
To see that $Y$ is $\ip\hc1$ use the 2nd line of 
\eqref{191E3}.\vtm

\epF{Theorem~\ref{191t}}
}

\vyk{
Let $\comi$ (codes of ordinals) = 
all $f\in\chc$ satisfying the following: 
\ben
\Renu
\atc
\atc
\atc
\itlb{191d}
the relation $\rep f$ strictly linearly orders the set $\abc f$.
\een

The following definitions, intended to express a relation 
between the ranks of $S(f)$ and $S(g)$ in terms of relations 
within $\chc$, are somewhat more tricky. 
Let $f,d\in\bI$. 
We write $f\rank d$ iff the following holds:
\ben
\nenu
\itlb{191r1}
$f$ and $d$ satisfy \ref{191a} and \ref{191c} above and 
$d$ also satisfies \ref{191d};


\itlb{191r2}
there is $h\in\bI$ such that $\ima h{\abc f}=\abc d$ 
and in addition
\ben
\itlb{191r3a}
if $k\rep f l$ then $h(k)\rep d h(l)$, so $h$ is a homomorphism; 

\itlb{191r3b}
if $m\in\abc f$ then $h(m)=\ssup_{\rep d} \ens{h(k)}{k\rep f m}$ 
(also meaning that the $\ssup$ actually exists).
\een
\een
We finally write $f\rra g$ in case  $f,g\in\bI$ 
satisfy \ref{191a} and \ref{191c} above, and 
there is $d\in\bI$ and $k,m\in\abc d$ such that 
$d$ satisfies \ref{191d}, 
$f\rank{(d\ser k)}$, 
$g\rank{(d\ser m)}$, 
and the order type of $\ens{l}{m\rep d l}$ under 
$\rep d$ is $\om+1$.

}


\parf{Some simple definability claims}
\las{193}

{\ubf We continue to argue in $\rL$.} 
As usual, $\pfs X=\ens{Y\sq X}{Y\text{ is finite}}$. 
\kmar{pfs}
\index{finite power set $\pfs X$}%
\index{power set!finite, $\pfs X$}%
\index{zPfin@$\pfs X$, finite power set}%
To countably code the topology of spaces $\can\xi$, 
put $U^\xi(\i,k,e)=\ens{x\in\can\xi}{x(\i)(k)=e}$ 
for all $\i\in\xi\in\cpo$,  $k<\om$, $e=0,1$.
If $u\sq\xi\ti\om\ti2$ is finite and {\em consistent\/} 
(that is, for no $\i,k$ both $\ang{\i,k,0}$
and $\ang{\i,k,1}$ belong to $u$)
then put
$U^\xi_u=\bigcap_{\ang{\xi,k,e}\in u}U^\xi(\i,k,e)$
(a basic clopen cube in $\can\xi$).
Finally, if 
$$
b\in \cb\xi:=\pfs{\pfs{\xi\ti\om\ti2}}
\kmar{cb xi}%
\index{zcCOxi@$\cb\xi$, codes of clopen sets}%
\index{codes!clopen sets, $\cb\xi$}%
$$
is {\em consistent}, in the sense that each
$u\in b$ is such, then put
$\okr\xi b=\bigcup_{u\in b}U^\xi_u$, 
\kmar{okr xi b}
\index{zcCObxi@$\okr\xi b$, clopen set}%
\index{clopen set!coded, $\cb\xi$}%
an arbitrary clopen subset of $\can\xi$.
($\cb{}$ from \rit{codes of ClOpen (sets)}.)

If $\xi\in\cpo$ then let 
$\ct\xi=\ens{X\sq\can\xi}{X\text{ is finite or countable}}$.
\kmar{ct xi}%
\index{zctblxi@$\ct\xi$, countable sets}%
\index{countable sets!$\ct\xi$}%

If $X,Y\in\ct\xi$ then let 
$X\caps Y=(X\cap\clo Y)\cup (Y\cap\clo X)$; 
\kmar{caps}
\index{zzzcaps@$\caps$, coded intersection}%
\index{intercection!coded, $\caps$}%
then clearly $X\caps Y\in\ct\xi$ 
and $\clo{(X\caps Y)}=\clo X\cap\clo Y$. 

If $\et\sq\xi$ belong to $\cpo$ and $Y\in\ct\et$ then 
let $Y\uars\xi$ consist of all points $x\in\can\xi$ 
\kmar{uars}
\index{zzzuars@$\uars$, coded lifting}%
\index{lifting!coded, $\uars$}%
such that $y=x\dar\et\in Y$ and the set 
$\ens{\ang{\i,k}}{\i\in\xi\bez\et\land x(\i)(k)=1}$ 
is finite. 
Thus $Y\uars\xi\in\ct\xi$ provided $Y\in\ct\et$, whereas 
$Y\uar\xi$ is not necessarily countable, of course, 
but still $\clo{(Y\uars\xi)}=Y\uar\xi$.

\ben
\penu
\itlb{1930}
the sets $\omi$, $\om$, $\bI$, $\cpo$, 
$\ct{}=\bigcup_{\xi\in\cpo}\ct\xi$, 
$\ens{\ang{\xi,X}}{\xi\in\cpo\land X\in\ct\xi}$ are 
$\id\hc1$ (as subsets of $\hc$);

\itlb{1931}
the maps $\xi\mto\cb\xi$ and $\xi,X,b\mto X\cap \okr \xi b$ 
belong to $\id\hc1$;

\itlb{1932}
the set  $\ens{\ang{X,Y}}{X,Y\in\ct{}\land \clo X\sq\clo Y}$ is 
$\id\hc1$;

\itlb{1933}
the map $\ang{\i,\j,X}\mto \pi_{\i\j}\akt X$ 
(Section~\ref{perm}) belongs to $\id\hc1$;

\itlb{1934}
the maps $\ang{X,Y}\mto X\caps Y$ 
and $\ang{\xi,Y}\mto Y\uars\xi$ belong to $\id\hc1$;

\itlb{1935}
the maps $\xi\in\cpo\mto \ft\xi$ 
(subsets of finite type, Section~\ref{rud}) 
and $\al\mto\tuq\al:=\al\lom\bez\ans\La$ 
(Section~\ref{prelim1})
belong to $\id\hc1$.
\een

The proof of \ref{1930}--\ref{1935} is based on one 
common principle. 
Let $\zct$ be the theory of Zermelo $\Z$ 
\kmar{zct}%
\index{theory!$\zct$}%
\index{zyAleph@$\zct$, theory}%
sans the Power Set axiom, 
plus the axiom saying that every set $x$ is at most 
countable. 
An $\in$-formula $\vpi(x,y,\dots)$ is 
\index{absolute!ZAabsolute@$\zct$-absolute formula}%
\index{formula!ZAabsolute@$\zct$-absolute}%
\index{YAabsolute@$\zct$-absolute formula}%
$\zct$-\rit{absolute}, if for any transitive model 
$\gM\in\hc$, $\gM\mo\zct$, and any $x,y,\ldots\in\gM$, 
the equivalence 
$(\hc\mo\vpi(x,y,\dots))\eqv (\gM\mo\vpi(x,y,\dots))$ 
holds.

\bte
\lam{193ta}
If\/ $\vpi(x,y,\dots)$ is a\/ $\zct$-absolute\/ 
$\in$-formula then the set\/ 
$X=\ens{\ang{x,y,\ldots}}{\hc\mo\vpi(x,y,\dots)}$ 
is of the definability class\/ $\id\hc1$.
\ete
\bpf
The relation $\ang{x,y,\ldots}\in X$ is equivalent to 
each of the two formulas
$$
\bay{l}
\sus\gM\in\hc\,\big(\gM\mo\zct\land\gM\text{ is transitive}
\land\gM\mo\vpi(x,y,\dots)\big),\\[0,5ex]
\kaz\gM\in\hc\,\big(\gM\mo\zct\land\gM\text{ is transitive}
\imp\gM\mo\vpi(x,y,\dots)\big).
\eay
$$
The first formula provides $X\in\is\hc1$, 
the second one gives $X\in\ip\hc1$.
\epf

Now to prove \ref{1930}--\ref{1935} it suffices to check
that some natural formulas, which define the sets and relations 
mentioned in \ref{1930}--\ref{1935}, are $\zct$-absolute. 
This is entirely routine, except perhaps for the relation 
$\clo X\sq\clo Y$, which we have to rewrite as follows. 
If $X,Y\in\ct\xi$ for one and the same $\xi\in\cpo$ 
then we let $\xi(X,Y)=\xi$, otherwise keep  $\xi(X,Y)$ 
undefined. 
Now, $\clo X\sq\clo Y$ is equivalent to the following 
formula, easily shown to be $\zct$-absolute:
$$
\xi=\xi(X,Y)\text{ is defined and }
\kaz b\in\cb\xi
(X\cap \okr\xi b\ne\pu\imp Y\cap \okr\xi b\ne\pu).
$$

\parf{Definability of iterated perfect sets}
\las{192}

Recall that $\cpe\xi=\ens{X\in\ct\xi}{\clo X\in\pe\xi}$ 
and $\cpei=\bigcup_{\xi\in\cpo}{\cpe \xi}$.
\index{codes!$\cpe\xi$}%
\index{zcIPSxi@$\cpe\xi$, codes}%
\index{codes!$\cpei$}%
\index{zcIPS@$\cpei$, codes}%
\kmar{cpei}

\bte
\lam{192t}
$\cpei$ and\/
$\ens{\ang{\xi,A}}{\xi\in\cpo\land A\in\cpe\xi}$
are\/ $\id\hc1$ sets.
\ete
\bpf
We use the notation of Section~\ref{193}.
Let $\Psi(\xi,A)$ say the following:
\ben
\nenu
\itlb{P1}
$\xi\in\cpo$ and 
$A\sq\can\xi$, and

\itlb{P2}
there is a set
$C\sq\can\xi$ and a bijection $h:C\onto A$ such
that:
\ben
\itlb{P3}
$C$ is topologically dense in $\can\xi$; 

\itlb{P4}
if $b_1\in \cb\xi$ and $\okr\xi {b_1}\cap A\ne\pu$
then there is $b\in \cb\xi$ such that the
image $\ima h {(C\cap \okr\xi b)}$ is equal to
$\okr\xi {b_1}\cap A$; 

\itlb{P5}
if $\i\in\xi$ and $x,y\in Z$
then $x\rsd{\sq\i}=y\rsd{\sq\i}$ iff
$h(x)\rsd{\sq\i}=h(y)\rsd{\sq\i}$.
\een
\een

We assert that (*) 
\big($\hc\mo\Psi(\xi,A)$\big) iff
\big($\xi\in\cpo$ and $A\in\cpe\xi$\big).

In the nontrivial direction, assume that $\xi,A\in\hc$
and $\Psi(\xi,A)$ is true in $\hc$.
Then $\xi\in\cpo$ by \ref{P1}, thus it remains to prove
that $\clo A\in\pe\xi$.

Let, by \ref{P2}, a set
$C\sq\can\xi$ and a bijection $h:C\onto A$ satisfy
\ref{P3}, \ref{P4}, \ref{P5} in $\hc$, so that in fact
$C\in\ct\xi$ is dense in $\can\xi$ by \ref{P3}.
In particular, $\clo C=\can\xi$.
Let $H=\clo h$ be the topological closure of $H$
in $\can\xi\ti\can\xi$.

It easily follows from \ref{P4}
(and the compactness of the spaces considered)
that $H$ is a homeomorphism
from $\clo C=\can\xi$ onto $\clo A$.
Finally, \ref{P5} implies that $H$ is \prok,
hence $\clo A\in\pe\xi$, as required.
This ends the proof of (*). 

It remains to prove that $\Psi$ defines a $\id\hc1$ 
relation. 
This looks somewhat doubtful 
(in spite of the rather obvious $\zct$-absoluteness of 
\ref{P1}, \ref{P3}, \ref{P4}, \ref{P5} and 
Theorem~\ref{193ta}), 
because the $\exists$ 
quantifier in \ref{P2} does not seem to be 
replaceable by a $\forall$ quantifier. 
Yet we can apply the following trick. 

Recall that $\tuo=\om\lom\bez\ans\La\in\cpo$. 
Clearly each $\xi\in\cpo$ can be embedded in $\tuo$ 
via a map $\pi\in\Ga_\xi$, where $\Ga_\xi$ consists of 
all \dd\su preserving and length-preserving 
injections $\pi:\xi\to\tuo$. 
Thus
$$
\bay{rcl}
\Psi(\xi,A)
&\eqv& 
\sus\pi\in\Ga_\xi\,\sus\xi'\,\sus A'
\big(
\xi'=\pi\akt\xi\land A'=\pi\akt A\land\Psi(\xi',A')
\big)\\[1mm]
&\eqv& 
\kaz\pi\in\Ga_\xi\,\kaz\xi'\,\kaz A'
\big(
\xi'=\pi\akt\xi\land A'=\pi\akt A\imp\Psi(\xi',A')
\big).
\eay
$$
On the other hand, 
if it is assumed that $\xi'\sq\tuo$ and 
$A'\in\ct{\xi'}$, then the formula 
$\Psi(\xi',A')$ is convertible to an equivalent 
$\is11$ form by a suitable coding of $\xi',A'$ by reals, 
and hence $\Psi$ defines a $\id\hc1$ relation in this 
particular domain by Proposition~\ref{p60}. 
It follows that the first line of the double equivalence 
above provides a $\is\hc1$ definition of the relation 
defined by $\Psi$, whereas the second line provides 
its $\ip\hc1$ definition, as required.
\epf

\parf{Definability of rudiments}
\las{194}

We come back to Definition~\ref{154d}. 

Given any set $\rB\sq\cpei$ 
(so that $\rB$ consists of codes of sets in $\pei$), 
we let $\cla{\rB}:=\ens{\clo A}{A\in\rB}$; 
\kmar{cla B}%
\index{zB*@$\cla B$}%
\index{symbol!$\cla{}$}%
\index{zzz**@$\cla{\;}$, set of closures}%
thus $\cla\rB\sq\pei$. 
Let $\al<\omi$, $\xi=\tuq\al$.
Say that $\rB\sq\cpe\xi$  is 
\rit{a coded rudiment of width $\al$}, in symbol  
\kmar{crud al}
\index{rudiment!coded, $\crud_\al$}%
\index{zcRuda@$\crud_\al$}%
$\rB\in\crud_\al$, if $\cla\rB\in\rud_\al$. 
To evaluate the complexity of $\crud_\al$ in the next theorem,  
we define several related notions. 
If $\al<\omi$, $\xi=\tuq\al$, $\rB\sq\cpe\xi$ then  
let $\rB^+=\rB^+_1\cup\rB^+_2\cup\rB^+_3$ 
\index{zB+@$\rB^+=\rB^+_1\cup\rB^+_2\cup\rB^+_3$}%
be the union of the three following sets: 
$$
\bay{rcl}
\rB^+_1 &=&
\ens{X\caps ((Y\dar\et)\uars\xi)}
{X,Y\in\rB\land \et\in\ft\xi
\land \clo{(Y\dar\et)}\sq\clo{(X\dar\et)}};\\[0.5ex]
\rB^+_2 &=&
\ens{X\cap \okr\xi b}
{X\in\rB\land b\in\cb\xi\land X\cap  \okr\xi b\in\cpe\xi};\\[0.5ex]
\rB^+_3 &=&
\ens{\pi_{\i\j}\akt X}
{X\in\rB\land \i,\j\in\xi\land\i\ekp\j}. 
\eay
$$
We also define $\crh\rB=\bigcup_n\rB_n$ 
(the \rit{coded} rudimentary hull), where 
\kmar{crh rB}%
\index{coded rudimentary hull, $\crh\rB$}%
\index{rudimentary hull!coded, $\crh\rB$}%
\index{zRHB@$\crh\rB$, coded rudimentary hull}%
$\rB_0=\rB$ and $\rB_{n+1}=(\rB_n)^+$, $\kaz n$. 
Then: (1) $\crh\rB\in \crud_\al$,\vom

(2) $\cla{(\crh\rB)}=\rh{(\cla\rB)}$ 
(rudimentary hull, Section~\ref{134*}),\vom 

(3) 
$\rB\in\crud_\al$ iff $\cla\rB=\cla{(\crh\rB)}$.

\bte
\lam{194t}
The following sets belong to\/ $\id\hc1:$
\ben
\renu
\itlb{1943}
$W_1=\ens{\ang{\rB,\crh\rB}}
{\sus\al<\omi\,(\rB\sq\peq\al)\land\rB\text{\rm\ is countable}};$

\itlb{1944}
$W_2=\ens{\ang{\al,\rB}}
{\al<\omi\land\rB\sq\peq\al\land\rB\in\crud_\al}.$
\een
\ete
\bpf
\ref{1943}
For any $\rB$, if there is an ordinal $\al$ such that 
$\rB\sq\peq\al$ then let $\al(\rB):=\al$. 
Then $\ang{\rB,\rB'}\in W_1$ iff $\Phi_1(\rB,\rB')$ holds 
in $\hc$, where 
$$
\Phi_1(\rB,\rB')\;:=\;
\big(\rB'=\crh\rB \land \al(\rB)=\al\text{ exists}
\land \rB,\rB'\sq\peq\al\big).
$$
In this formula, the two first summands are 
$\zct$-absolute, hence $\id\hc1$ by Theorem~\ref{193ta}, 
whereas the rightmost summand is $\id\hc1$ 
by Theorem~\ref{192t}.

\ref{1944}
Quite similarly, $\ang{\al,\rB}\in W_2$ iff 
$\Phi_2(\rB,\rB')$ holds in $\hc$, where 
$$
\Phi_2(\al,\rB)\;:=\;
\big(\al(\rB)=\al\land  \rB\sq\peq\al \land 
\rB=\crh\rB\big), 
$$
and then replace $\rB=\crh\rB$  
by $\ang{\rB,\rB}\in W_1$ and refer to \ref{1943}.
\epf

\parf{Definability of rudimentary sequences}
\las{195}

Recall that a sequence $\jba=\sis{\rba\al}{\al<\la}$ is a 
\rit{coded\/ \ruds\ of length\/ $\la$}, or a \rit{\cruds},  
\index{coded rudimentary sequence}%
\index{rudimentary sequence!coded, @\cruds}%
\index{sequence!rudimentary, coded}%
if each $\rba\al\in\crud_\al$ is countable 
and the sets 
$\rqu\al=\cla{\rba\al}:=\ens{\clo A}{A\in\rba\al}\in\rud_\al$ 
form a \ruds\ $\jqo=\cla\jba:=\sis{\rqu\al}{\al<\la}$. 
\kmar{cla jba}%
\index{zzba**@$\cla\jba$}%

\bte
\lam{195t}
The following set belongs to\/ $\id\hc1:$ 
\bce
$W=\ens{\ang{\al,\jba}}
{\al<\omi\land\jba \text{\rm\ is a coded\/ \ruds\ of length }\al}$.
\ece
\ete

\bpf
Conditions \ref{drs1}, \ref{drs2}, \ref{drs3} of 
Definition~\ref{drs} find their $\id\hc1$ forms by 
different results above. 
In particular, as far as \ref{drs3} is concerned,  
make use of \ref{1932} in Section~\ref{193}. 
Recall the remaining condition \ref{drs4}:
\ben
\Aenu
\atc\atc\atc
\itlb{drs4x}
if $3\le\nu<\la$ then  
$\RC(\jqo\res\nu)\ssw\rqu\nu$ 
in the sense of
Definition~\ref{134d}; here 
$\RC(\jqo\res\nu)=\RC(\bsc(\jqo\res\nu))=
\RC(\bigcup_{\al<\nu}(\rqu\al\uar\tuq\nu))$ 
and $\rqu\al=\jqo(\al)$.
\een
In terms of a coded\/ \ruds\ 
$\jba=\sis{\rba\al}{\al<\la}$, 
it takes the form:
\ben
\cAenu
\atc\atc\atc
\itlb{drs4c}
if $3\le\nu<\la$ then  
$\rba{<\nu}\sswc \rba\nu$ --- 
where $\rB\sswc\rB'$ means 
$\cla{\rB}\ssq\cla{\rB'}$ 
provided $\rB,\rB'\sq\cpe\nu$, 
$\rba{<\nu}=\crh{\ens{A\uars{\tuq\nu}}
{A\in\bigcup_{\al<\nu}\rba\al}}$, whereas
$\crh\rB$ and $\uars$ are defined in Sections~\ref{194}, 
resp., \ref{193}.
\een
Thus it remains to prove that $\rB\sswc\rB'$ is a 
$\zct$-absolute, hence a $\id\hc1$ relation 
by Theorem~\ref{193ta}. 
To check this, we return to Definition~\ref{134d}. 
In terms of $\rP=\cla\rB$ and $\rQ=\cla{\rB'}$, 
conditions \ref{134i}, \ref{134ii}, \ref{134iii} there 
take the form:
\ben
\atc
\atc
\atc
\atc
\cstenu
\itlb{c134i}
There is $A\in\rB$ dense in $\can\xi$, 
so that $\clo A=\can\xi$.

\itlb{c134ii}
If $\et\in\ft\xi$, $A\in\rB$, $B\in\rB'$, 
$\clo{(B\dar\et)}\sq \clo{(A\dar\et)}$, then there is 
$C\in\rB'$ such that 
$\clo C\sq \clo A$ and $\clo{(C\dar\et)}=\clo{(B\dar\et)}$.

\itlb{c134iii}
If $\i\in\xi$, $A\in\pro\rB\i$, $B\in\pro{\rB'}\i$, 
then $\clo A\cap \clo B$ is {\em clopen\/} in $\clo B$. 
\een

That \ref{c134i} is $\zct$-absolute, is pretty clear. 

See the end of Section~\ref{193} regarding the conversion 
of formulas like $\clo C\sq \clo A$ in \ref{c134ii} to a 
$\zct$-absolute form. 

Finally, $\clo A\cap \clo B=\clo{A\caps B}$. 
Then the clopenness of $\clo{A\caps B}$ in $\clo B$ is 
equivalent to the following  $\zct$-absolute formula:
$$
\sus b\in \cb\xi\,
(\clo{A\caps B}\cap \okr\xi {b}=\clo B). 
$$
Thus  
\ref{c134i}+\ref{c134ii}+\ref{c134iii}, as a whole,
is $\zct$-absolute, and $\id\hc1$, as required.
\epf

\parf{Definability claims related to continuous functions}
\las{196}

Recall the notions $\rat\xi$, $\kont\xi$,  
$\kond\xi$, 
$\komt=\bigcup_{\xi\in\cpo}\kont\xi$, 
and $\komd=\bigcup_{\xi\in\cpo}\kond\xi$, 
related to codes of continuous functions $\can\xi\to{\cN}$ 
and $\can\xi\to\dn=\can{},$ 
$\xi\in\cpo$, and 
defined in $\rL$ in Section~\ref{101}. 

See Sections~\ref{173},\ref{176},\ref{175} on axis maps 
and avoidance.

\bte
[in $\rL$]
\lam{196t}
The following sets belong to\/ $\id\hc1:$
\ben
\renu
\itlb{196t1}
$\ens{\ang{\xi,f}}{\xi\in\cpo\land f\in\kont\xi}$ and\/ 
$\ens{\ang{\xi,f}}{\xi\in\cpo\land f\in\kond\xi};$

\itlb{196t2}
$\{\ang{\xi,A,f,\i}:
\xi\in\cpo\land f\in\kond\xi\land A\in\cpe\xi\;\land$\\[0.4ex]
$\,$\hspace*{0.4\textwidth}
$\land\;\clo f\text{\rm\ is an $\i$-axis map on }\clo A\};$

\itlb{196t3}
$\{\ang{\xi,A,f,\rU}:
\xi\in\cpo\land f\in\kond\xi\land A\in\cpe\xi\land 
\rU\in\hc\text{\rm\ consists of}$\\[0.4ex]
$\,$\hspace*{0.02\textwidth}
$\text{\rm{}countable subsets of }\can{}
\land\clo f\text{\rm\ avoids $\clo E$ on }\clo A 
\text{\rm\ for any }E\in\rU\};$

\itlb{196t4}
$\ens{\ang{\xi,A}}{\xi\in\cpo\land A\in\cpe\xi\land
\clo A
\text{\rm\ is uniform as in Section~\ref{174}}}$.
\een
\ete

\bpf
\ref{196t1} 
Let $f:\rat\xi\to{\cN}$. 
Then $f\in\kont\xi$ iff for any $m,k<\om$ 
there exists $b\in \cb\xi$ (a code of a clopen set in $\can\xi$)
such that for all $x\in\rat\xi$ the equivalence  
$x\in\okr\xi {b}\eqv f(x)(m)=k$ holds. 
This yields a $\zct$-absolute definition, 
and hence the class $\id\hc1$, for the first set. 

\ref{196t2} 
Let $f\in\kond\xi$. 
Then  $\clo f$ is an $\i$-axis map on $\clo A$ iff 
for all $b\in \cb\xi$, $k<\om$, and $e=0,1$ the 
following holds: 
\pagebreak[0]
$$
\kaz x\in\okr\xi {b}\cap\rat\xi 
\big(x(\i)(k)=e\land f(x)(k)=1-e\big) 
\;\imp\;
A\cap \okr\xi {b}=\pu, 
$$
and this is a $\zct$-absolute formula. 

\ref{196t3} 
According to the compactness of the spaces considered, 
if a continuous map $\clo f$ avoids $\clo E$ on $\clo A$ 
then there exist clopen supersets $X\supseteq\clo A$ 
and $Y\supseteq\clo E$ such that 
$\clo f$ avoids $Y$ on $X$. 
We conclude that the relation 
``$\clo f$ avoids $\clo E$ on $\clo A$'' 
is equivalent to the following $\zct$-absolute formula: 
$$
\sus b,c\in \cb\xi\,
\big(A\sq\okr\xi {b}\land E\sq\okr\xi c
\land\kaz x\in\okr\xi {b}\cap\rat\xi\, 
(f(x)\nin\okr\xi{c})\big) . 
$$

\ref{196t4} 
For $\clo A$ to be uniform it's necessary that $A$ itself 
is uniform, \ie, if $\i\su\j$ belong to $\xi=\dym A$ and 
$x,y\in A$ satisfy $x(\j)=y(\j)$ then $x(\i)=y(\i)$ as well. 
In other words, there is a map 
$h_{\i\j}:A\dir\j\to A\dir \i$ satisfying 
$x(\i)=h_{\i\j}(x(\j))$ for all $x\in A$. 
Thus the condition that (*) every closure $\clo{h_{\i\j}}$ 
in the according space $\can{}\ti\can{}$ remains a map, 
is necessary and sufficient for $\clo A$ to be uniform. 
On the other hand, (*) is $\zct$-absolute by an argument 
similar to used in the proof of \ref{196t1}.
We leave the details to the reader.
\epf

\parf{Definability of the forcing approximation}
\las{198}

{\ubf Still arguing in $\rL$,} 
now we come back to the notion of forcing approximation 
$\fo$ introduced by Definition~\ref{fod}. 
The goal of the next theorem 
is to evaluate the complexity of the sets 
$$
\bay{rcl}
\kmar{For us1n}%
\index{set!forS1n@$\For{\us1n}$}%
\index{zForS@$\For{\us1n}$}%
\For{\us1n} &\!\!=\!\!& \ens{\ang{X,\vpi}}
{X\in\cpei\land \vpi\text{ a closed $\ls1n$ formula} \land 
\clo X\fo\vpi};\\[0.5ex]
\For{\up1n} &\!\!=\!\!& \ens{\ang{X,\vpi}}
{X\in\cpei\land \vpi\text{ a closed $\lp1n$ formula} \land 
\clo X\fo\vpi}.
\kmar{For up1n}%
\index{set!forP1n@$\For{\up1n}$}%
\index{zForP@$\For{\up1n}$}%
\eay
$$

\bte
\lam{198t}
The set\/ $\For{\up11}$ belongs to\/ $\id\hc1$.

The set\/ $\For{\us11}$ belongs to\/ $\ip\hc1$.

If\/ $n\ge1$ then\/ $\For{\us1{n+1}}$ belongs to\/ $\is\hc n,$ 
$\For{\up1{n+1}}$ belongs to\/ $\ip\hc n.$ 
\ete

\bpf
{\it Case\/ $\ip11$.} 
Assume that $X\in\cpei$, $\vpi$ is a closed $\lp11$ formula, 
$\xi=\dym X\cup\dym\vpi\in\cpo$. 
Using the same trick as in the end of Section~\ref{192},
note that $\xi$ can be embedded in $\tuo$ 
via a map $\pi\in\Ga_\xi$, where $\Ga_\xi$ consists of 
all \dd\su preserving and length-preserving 
injections $\pi:\xi\to\tuo$. 
Then $\clo X\fo\vpi$ is equivalent to each of the two formulas:
$$
\bay{c}
\sus\pi\in\Ga_\xi\,\sus\xi'\,\sus A'\,\sus\vpi'\,
\big(
\xi'=\pi\akt\xi\land X'=\pi\akt X\land\vpi'=\pi\akt\vpi
\land {\clo {X'}}\fo\vpi'
\big),\\[1mm]
\kaz\pi\in\Ga_\xi\,\kaz\xi'\,\kaz A'\,\kaz\vpi'\,
\big(
\xi'=\pi\akt\xi\land X'=\pi\akt X\land\vpi'=\pi\akt\vpi
\imp {\clo {X'}}\fo\vpi'
\big).
\eay
$$
On the other hand, 
if it is assumed that $\xi'\sq\tuo$ 
then ``$\clo {X'}\fo\vpi'$'' is essentially a 
$\ip11$ relation via a suitable coding of 
$\vpi',X'$ by reals, by \ref{fod1} of Definition \ref{fod},
and hence we have a $\id\hc1$ relation in this 
particular domain by Proposition~\ref{p60}. 
It follows that the first line of the double equivalence 
above provides a $\is\hc1$ definition of the relation 
``$\clo {X'}\fo\vpi'$'', whereas the second line provides 
its $\ip\hc1$ definition, as required.

{\it Case\/ $\is11$.} 
Essentially the same argument, but if $\vpi$ is a $\us11$ 
formula then \ref{fod1} of Definition \ref{fod} yields 
a $\ip12$ relation, hence $\ip\hc1$ relation. 

{\it Inductive step\/ $\ip1n\to\is1{n+1}$, $n\ge1$.} 
By \ref{fod4} of Definition \ref{fod}, 
$\For{\us1{n+1}}$ consists of all pairs 
$\ang{X,\sus x\,\vpi(x)}$, where $\vpi$ is a $\lp1n$ 
formula and there is $f\in\komt$ satisfying 
$\ang{X,\vpi(f)}\in\For {\up1{n}}$. 
Thus if $\For{\up1{n}}$ belongs to\/ $\ip\hc {n-1}$ 
or at worst $\id\hc {n}$
then $\For{\us1{n+1}}$ belongs to\/ $\is\hc n$.

{\it Inductive step\/ $\is1{n+1}\to\ip1{n+1}$, $n\ge1$.} 
By \ref{fod5} of Definition \ref{fod}, 
$\For{\up1{n+1}}$ consists of all pairs 
$\ang{X,\otr\vpi}$, where $X\in\cpei$, 
$\vpi$ is a closed $\ls1{n+1}$ formula, 
and there is no $Y\in\cpei$ satisfying $\clo Y\psq\clo X$ 
and $\ang{Y,\vpi}\in\For {\us1{n+1}}$. 
Thus if $\For{\us1{n+1}}$ belongs to\/ $\ip\hc {n}$
then $\For{\up1{n+1}}$ belongs to\/ $\ip\hc n$.
\epf

\parf{Definability of being an 1-5-$\nn$ extension}
\las{199}

Here we collect all the previous results of this chapter to prove the following main definability theorem. 
If $\nn\ge 1$ then let $\prdo\nn$ be the set 
of all pairs\/ $\ang{\jba,\rba\la}$, 
where\/ $\jba=\sis{\rba\al}{\al<\la}$ is a 
coded\/ \ruds\ of length some\/ $\la<\omi$,  
$\rba\la\in\crud_\la$, and the set\/ 
$\rqu\la=\cla{\rba\la}:=\ens{\clo A}{A\in\rba\la}$ 
is an 1-5-$\nn$ extension of the \ruds\

\bte
[in $\rL$]
\lam{199t}
Let\/ $\nn\ge1$. 
Let\/ $\prdo\nn$ be the set 
of all pairs\/ $\ang{\jba,\rba\la}$, 
where\/ $\jba=\sis{\rba\al}{\al<\la}$ is a 
coded\/ \ruds\ of length some\/ $\la<\omi$,  
$\rba\la\in\crud_\la$, and the set\/ 
$\rqu\la=\cla{\rba\la}:=\ens{\clo A}{A\in\rba\la}$ 
is an 1-5-$\nn$ extension of the\/ \ruds\ 
$\jqo=\cla\jba:=\sis{\rqu\al}{\al<\la}$, 
where\/ 
$\rqu\al=\cla{\rba\al}:=
\ens{\clo A}{A\in\rba\al},\,\kaz\al$.

Then\/ $\prdo\nn$ belongs to\/ $\id\hc\nn$. 
\ete

\bpf 
We have to evaluate coded forms of 
conditions \ref{C*}, \ref{C1+}, \ref{C1++}
(including \ref{C2+}--\ref{C5+} in the last one) 
as in Definition~\ref{15e}. 

\bit
\item[\ref{C*}]
\rit{The extended sequence\/ 
$\jba\we\rba\la$ 
is a\/ \cruds\ {\rm(of length\/ $\la+1$)}.} 
\eit

This condition is $\id\hc1$ by Theorem~\ref{195t}. 

\bit
\item[\ref{C1+}]
\rit{If 
$\clo{\bM\la}\in \noc{\cla\jba}\yi\dym{\clo{\bM\la}}\sq\tuq\la$ 
then\/ 
$\sus A\in\rba\la\,(\clo A\psq\clo{\bM\la}).$}
\eit

This needs some bit of work.
Recall that the map $\al\mto{\bM\al}$ is $\id\hc1$ by 
Lemma~\ref{144L}. 
The relation $\clo A\sq\clo B$ is $\id\hc1$ by 
\ref{1932} in Section~\ref{193}. 
Thus the 2nd and 3rd subformulas in \ref{C1+} 
define $\id\hc1$ relations. 
Let's focus on the 1st subformula 
$\clo{\bM\la}\in \noc{\cla\jba}$. 
Here $\noc{\cla\jba}=\noc{\jqo}:=\noc{\rqu{<\la}}$, 
where 
$$
\TS
\rqu{<\la}=\bigcup_{\al<\la}\rqu\al=\cla{{\rB_{<\la}}} 
\quad\text{and}\quad
\rB_{<\la}=\bigcup_{\al<\la}\rba\al\,,
$$ 
and $\noc{\cdot}$ is the normal hull, Definition~\ref{noc}. 

To eliminate the operation $\noc{\cdot}$ of indefinite 
complexity, we define $\rU=\rh(\rqu{<\la}\uar\tuq\la)$ 
(the rudimentary hull, Section~\ref{rud}), so that 
$\rU\in\rud_\la$ is countable. 
At the level of codes, we put 
$\rA = {\rB_{<\la}}\uars{\tuq\la}$ 
(see Section~\ref{193} on $\uars$), 
so that $\rA\sq\cpe\la$ is countable and 
$\cla{{\rA}} = {\rQ}\uar{\tuq\la}$. 

We further define $\rC=\crh{\rA}$ 
(the coded rudimentary hull, Section~\ref{194}), 
hence $\rC\in\crud_\la$ and 
$\rU=\cla\rC:=\ens{\clo C}{C\in\rC}$. 

Now suppose that ${\bM\la}\in\cpei$ and 
$\xi=\dym {\bM\la}\sq\tuq\la$. 
We are going to define the relation $\clo{\bM\la}\in \cX$, 
where $\cX=\noc{\cla\jba}$,  
in terms of the above notation, so that it becomes $\id\hc1$. 
First of all, $\clo {\bM\la}\in \cX$ iff 
$\clo {\bM\la}\rsq\i\in\cX\rsq\i$ for 
all $\i\in\xi$, by  \ref{rfo5} of Section~\ref{rfo}. 
On the other hand, $\cX\rsq\i=\rU\rsq\i$ 
by Lemma~\ref{r2n}. 
Thus (*) $\clo {\bM\la}\in \cX$ iff 
$\clo {\bM\la}\rsq\i\in\rU\rsq\i$ 
for all $\i\in\xi=\dym {\bM\la}$. 

On the other hand, 
the relation $\clo {\bM\la}\rsq\i\in\rU\rsq\i$ 
is equivalent to 
\bce
$\sus C\in\rC\,\big(\clo {(\bM\la\rsq\i)}=\clo {(C\rsq\i)}\big)$.
\ece
This allows to rewrite (*) as follows: 
$$
\clo{\bM\la}\in \noc{\cla\jba}
\;\eqv\;
\kaz \i\in\dym{\bM\la}\,
\sus C\in\rC\,\big(\clo {(\bM\la\rsq\i)}=\clo {(C\rsq\i)}\big),
\eqno(\mdag)
$$
where $\rC=\crh{\rA}=\crh{({\rB_{<\la}}\uars{\tuq\la})}$.  
Finally note that the right-hand side of (\mdag) contains 
only $\id\hc1$ relations and operations 
by \ref{1932} and \ref{1934} in Section~\ref{193} 
and Theorem~\ref{194t}.  
We conclude that ``$\clo{\bM\la}\in \noc{\cla\jba}$'' is a 
$\id\hc1$ relation, 
and hence so is \ref{C1+} as a whole 
(with $\la$, $\rB_\la$, $\jba$ as arguments). 

\bit
\item[\ref{C1++}]
{%
$\clo{\bM\la}\in \rqu{<\la}
\imp 
\sus Y\in\rqu\la\,(Y\psq\clo{\bM\la}\land 
\text{\rm\ref{C2+}--\ref{C5+}})$, \ \ or equivalently, 

$\sus B\in\rB_{<\la}\,(\clo{\bM\la}=\clo B)
\imp 
\sus A\in\rB_\la\,(\clo A\psq\clo{\bM\la}\land 
\text{\rm\ref{C2+}--\ref{C5+}})$.}
\eit
Temporarily leaving \ref{C2+}--\ref{C5+} aside in the 2nd 
line of \ref{C1+} here,  note that the subrelations  
${\sus B\in\rB_{<\la}\,(\clo{\bM\la}=\clo B)}$ and  
$\sus A\in\rB_\la\,(\clo A\psq\clo{\bM\la})$ 
are $\id\hc1$ by \ref{1932} 
in Section~\ref{193}. 
Now consider \ref{C2+}--\ref{C5+} one by one, 
assuming that $\clo{\bM\la}\in \rQ_{<\la}$, or equivalently, 
that some $B\in\rB_{<\la}$ satisfies 
$\clo{\bM\la}=\clo B$.

\bit
\item[\ref{C2+}]
{%
If $\la$ is limit, $k<\om$, $\bPsd k\la\sq\rqu{<\la}$, 
and\/ $\bPsd k\la$ is dense in\/ $\rqu{<\la}$, 
then $\clo A\sqd \bigcup\bPsd k\la$. 
}
\eit
Here we recall that $\la,k\mto \bPsd k\la$ is a $\id\hc1$ 
map by Lemma~\ref{144L}. 
Then replace the subformula $\bPsd k\la\sq\rqu{<\la}$ by 
$\kaz A\in \bPsd k\la\,\sus B\in\rB_{<\la}\, 
(A\in\cpei\land \clo A=\clo B)$ --- 
which defines a $\id\hc1$ relation by \ref{1932} 
in Section~\ref{193}.  
Similar routine $\id\hc1$ replacements apply also for 
the subformulas ``$\bPsd k\la$ is dense in\/ $\rqu{<\la}$''
and $\clo A\sqd \bigcup\bPsd k\la$, with an extra reference 
to \ref{1934} in Section~\ref{193}.  
After that, we conclude that  \ref{C2+} is a $\id\hc1$ 
relation.

\bit
\item[\ref{C3+}]
{%
If\/ $\nn\ge2$ and $\bmp\la$ is a closed formula\/ 
$\vpi$ in\/ $\bigcup_{k\le \nn}\ls1k$ then\/  
$\clo A \fo \vpi$ or\/ $\clo A\fo \otr\vpi$ 
--- {\em this is void in case $\nn=1$}. 
}
\eit
Use Theorem~\ref{198t} to see that \ref{C3+} 
is a $\id\hc\nn$ condition. 

\bit
\item[\ref{C4+}]
{%
$\clo A$ is a uniform set.%
} --- 
\it 
Still a $\id\hc1$ condition by Theorem~\ref{196t}\ref{196t4}.
\eit
%

\bit
\item[\ref{C5+}]
{%
Either (a) 
$\bbF\la:=\clo{\bbf\la}$ avoids $\clo E$ on $\clo A$ 
for all $\i\in\tuq\la$ and $E\in\rB_\la\dir\i$, 

or (b)
there is $\j\in\tuq\la$ such that $\clo{\bbf\la}$ is an 
$\j$-axis map on $\clo A$ but $\clo{\bbf\la}$ 
avoids $\clo E$ on $\clo A$ 
for all $E\in\rB_\la\dir\i$ and 
$\i\in\tuq\la$, $\i\not\ekp\j$.
}
\eit

Theorem~\ref{196t} (different items) implies that \ref{C5+} 
is $\id\hc1$, too.

This completes the proof of Theorem~\ref{199t}: 
all components of the definition of $\prdo\nn$ are $\id\hc1$ 
except for \ref{C3+} which is $\id\hc\nn$.
\epf

\parf{The final forcing construction}
\las{200}

\bpf[Theorem~\ref{mt1d}, finalization, in $\rL$] 
Let\/ $\nn\ge1$. 
Theorem~\ref{18t} implies that for any 
coded \ruds\ $\jba'$ of length $\la=\dom{\jba'}<\omi$ 
there exists a coded rudiment $\rba\la\in\crud_\la$ 
satisfying $\ang{\jba',\rba\la}\in \prdo\nn$. 
Let $\rba\la(\jba')$ be the $\lel$-\rit{minimal} of 
such coded rudiments $\rba\la\in\crud_\la$. 

Define a coded \ruds\ $\jba=\sis{\rba\la}{\la<\omi}$ 
\kmar{jba}
\index{zzba@$\jba$, final coded \ruds}%
so that $\rba\la=\rba\la(\jba\res\la)$ for all $\la<\omi$. 
Then, by Theorem~\ref{199t}, $\jba$ belongs to $\id\hc\nn$, 
because it is known that iterated constructions, by taking the 
$\lel$-minimal choice in the domain bounded by a $\id\hc\nn$ 
relation, lead to $\id\hc\nn$ final results 
(say by Proposition~\ref{p61}\ref{p613}). 
It follows that the according \ruds\ 
$\jqo=\sis{\rqu\la}{\la<\omi}$, 
where $\rqu\la=\cla{(\rba\la)}$, $\kaz\la$, 
satisfies the global definability condition \refn{C6}\nn\ 
via $\jba$.

On the other hand, 
each $\rqu\la$ is a\/  a  1-5-$\nn$ extension 
of\/ $\jqo\res\la$,  
because $\ang{\jba\res\la,\rba\la}\in \prdo\nn$ 
by construction. 

Thus the sequence $\jqo$ witnesses Theorem~\ref{mt1d}. 
\epf

\bpf[Theorem~\ref{mt1}, finalization] 
It remains to recall that Theorem~\ref{mt1d} implies 
Theorem~\ref{mt1}, see Section~\ref{156}. 
\epF{Theorem~\ref{mt1}}

\sekt
[\ \ Proof of the second main theorem] 
{Proof of the second main theorem}
\las{b}

Here we prove Theorem~\ref{mt2}. 
The model 
$\caL\w\cap\dn$ 
defined in Section~\ref{b1} will be a set 
in an $\cX$-generic extension $\rL[\w]$, where 
$\cX$ is given by Theorems~\ref{mt1b} and \ref{mt1c}. 
Here $\caL\w$ is defined in such a way that 
\ben
\Renu
\itlb{be1}
if $\al<\omi$ 
is an odd ordinal then the real 
$x_{\al-1}=\w(\ang{\al{-}1})\in\dn$ 
{\ubf does not} belong to $\caL\w$, 
but on the other hand 

\itlb{be2}
the real $x_{\al-1}$ is {\ubf definable over}  
$\caL\w$ by a suitable $\is1{n+2}$ formula, 
with $x_{\al}=\w(\ang{\al})\in\cam$ as the only 
parameter, by means of the structure of reals 
of the form $\w(\i)$, where $\i\in\tup$ is even and 
$\i(0)=\al$ --- which we put in $\caL\w$.
\een
Using \ref{be1} and \ref{be2}, we prove in 
Section~\ref{it1} that 
$\xCA{\fs1{\nn+2}}$ fails in $\caL\w$.
Then using the $\nn$-Odd-Expansion property of $\cX$ 
we show that 
$\caL\w$ is an elementary submodel of $\rL[\w]$ 
\poo\ all $\is1{\nn+1}$ formulas with reals in 
$\caL\w$ as parameters, 
and infer $\xCA{\fs1{\nn+1}}$ in $\caL\w$ in 
Section~\ref{it2}.
We finally establish the \paf\ $\xAC{\is1\iy}$ in $\caL\w$ in 
Section~\ref{it3} by permutations-related arguments.

\parf{The model}
\las{b1}

If $\w\in\can\tup$ is an $\tup$-array 
of reals then let $\tux\w$ consist of all tuples $\i\in\tup$ 
such that 
\ben
\fenu
\itlb{b1*} 
the ordinal $\al=\i(0)$ is odd, 
hence $\al{-}1$ is well-defined, and:  
if $1\le k<\lh\i$ and 
$\i({k})$ is even then $\w(\ang{\al{-}1})(k)=0$. 
\een
We put $\Omv =\ens{\xi\in\cpo}{\xi\sq\tux\w}$  
and 
$ 
\caL\w=\bigcup_{\xi\in \Omv}\rL[\w\dar\xi]
$.

Quite obviously, $\caL\w$ is not necessarily a model of $\zf$.

\ble
\label{b2}
\ben
\renu
\itlb{b22}
If\/ $\et\osq\xi$ belong to\/ $\cpo$  
then\/ $\et\in\Omv\imp \xi\in\Omv$.

\itlb{b23}
If\/ $\al<\omil$ is odd and\/ $k\ge1$ then TFAE$:$
$1)$ there is an even tuple\/ $\i\in\tux\w$ with\/ 
$\i(0)=\al$ and\/ $\lh\i=k+1$, and\/  
$2)$ $\w(\ang{\al{-}1})(k)=0$.

\itlb{b24}
If\/ $\i=\ang\al\in\tup$ then\/ $\i\in\tux\w$ iff\/ 
$\al$ is odd. 
\qed
\een
\ele

\ble
\lam{b3}
Let\/ $\cX\in\RF$ be a normal forcing in\/ $\rL$, 
and\/ $\w\in\can\tup$ be\/ \dd\cX generic over\/ $\rL$. 
Then\/ $\tux\w\yi\Omv\in\rL[\w]$ 
{\rm(not necessarily\/ $\in\rL$)}
and$:$
\ben
\renu
\itlb{b31}
if\/ $\i\in\tup$ then\/ $\w(\i)\in\caL\w$ iff\/ $\i\in\tux\w$$;$

\itlb{b32}
if\/ $\i=\ang\al\in\tup$ then\/ $\w(\i)\in\caL\w$ iff\/ $\al$ 
is odd.
\een
\ele 

\bpf
\ref{b31}
If $\i\in\tux\w$ then obviously $\ilq\i\in\Omv$ and we are done. 
To prove the converse suppose that $\w(\i)\in\caL\w$, hence 
$\w(\i)\in\rL[\w\dar\xi]$ for some $\xi\in\Omv$.
Then $\i\in\xi$ by Corollary~\ref{632}, hence $\i\in\tux\w$.

To prove \ref{b32} use \ref{b31} and Lemma~\ref{b2}\ref{b24}.
\epf

\bte
\lam{bt}
Assume that\/ $\nn\ge1$ and\/ $\cX\in\RF$ is a normal 
forcing 
as in Theorem~\ref{mt1b}, \ie, 
$\cX$ has the Fusion, Structure, $\nn$-Odd-Expansion, and 
$\nn$-Definability properties in\/ $\rL$.
Let\/ $\w\in\can\tup$ be\/ \dd\cX generic over\/ $\rL$. 
Then$:$
\ben
\renu
\itlb{bt1}
$\xCA{\fs1{\nn+2}}$ (with parameters) 
fails in\/ $\stk\om{\caL\w\cap\dn}$.

\itlb{bt2}
$\xCA{\fs1{\nn+1}}$ (with parameters) 
holds in\/ $\stk\om{\caL\w\cap\dn}$.

\itlb{bt3}
$\xAC{\is1{\iy}}$ and\/ $\xCA{\is1{\iy}}$ 
(\paf) 
hold in\/ $\stk\om{\caL\w\cap\dn}$.
\een
\ete 

Reals $x\in\caL\w\cap\dn$ are identified with sets 
$\ens{k}{x(k)=0}$, so that we view $\caL\w\cap\dn$ as a subset of 
$\pws\om$ in the context of this theorem. 

Quite obviously Theorem~\ref{bt} implies Theorem~\ref{mt2}. 

The proof of Theorem~\ref{bt} goes on below in this Chapter.

\parf{Item 1: violation of Comprehension at the level $\nn+2$}
\las{it1}

\bpf
[item \ref{bt1} of Theorem~\ref{bt}]
By 
the $\nn$-Definability property of $\cX$ as in 
Definition~\ref{82d}, the set 
$E=\gee{}\w\cap\caL\w$ 
is $\ip1{\nn+1}$ over $\caL\w$, where 
$$
\gee \gM\w
= 
\ens{\ang{k,\w(\i)}}
{k\ge1\land \i\in\tup\text{ is even\,}\land
\lh\i=k}.
$$
Here it is not claimed that $E\in\caL\w$. 
What {\ubf is} asserted is that there is a \paf\ 
$\ip1{\nn+1}$ formula $\vpi(k,x)$ such that 
$$
E=\ens{\ang{k,x}}{x\in\caL\w\land \caL\w\mo\vpi(k,x)}.
\eqno(1)
$$
Now we claim that, for any $k\ge1$,  
$$
\w(\ang0)(k)=0 
\eqv 
\sus x (\ang{k{+}1,x}\in E\land \w(\ang1)\in\rL[x]). 
\eqno(2)
$$

From left to right, let $\w(\ang0)(k)=0$. 
By Lemma~\ref{b2}\ref{b23} ($\al=1$), there is 
an even tuple $\i\in\tup[\w]$ with $\i(0)=1$ and\/ $\lh\i=k+1$. 
Let $x=\w(\i)$. 
By definition, $\ang{k{+}1,x}\in E$. 
Moreover $\w(\ang1)\in\rL[x]$ by 
the Structure property,
since $\ang1\sq\i$ 
by construction. 
Thus the right-hand side of (2) holds. 

From left to right, suppose that the right-hand side of (2) 
holds, and this is witnessed by some $x$.  
Then $x=\w(\i)$, where $\i\in\tup$ is even and $\lh\i=k{+}1$, 
and, as $\ang{k{+}1,x}\in E\sq\caL\w$, 
we have $x\in\caL\w$, and hence $\i\in\tup[\w]$ by \ref{b31}. 
Moreover, as $\w(\ang1)\in\rL[x]$, we have $\ang1\sq\i$ 
by the Structure property, 
hence $\i(0)=1$. 
To conclude, $\i\in\tux\w$ is even, $\lh\i=k{+}1$, 
$\i(0)=1$.  
This implies $\w(\ang0)(k)=0 $ by Lemma~\ref{b2}\ref{b23} 
($\al=1$), as required. 

Combining (1) and (2), it is clear now that $\w(\ang0)$ is 
definable over $\caL\w$ by a $\is1{\nn+2}$ formula 
(note the quantifier $\sus x$ in (2)!), 
with $\w(\ang1)\in\caL\w$ as the only parameter.
However $\w(\ang0)\nin\caL\w$ by Lemma~\ref{b3}\ref{b32}.
\epf

\parf{Item 2: verification of Comprehension at the level $\nn+1$}
\las{it2}

\bpf
[item \ref{bt2} of Theorem~\ref{bt}] 
The first step is the following claim, motivated by the 
 the $\nn$-Odd-Expansion property of $\cX$ and 
Lemma~\ref{b2}\ref{b22}:

\ben
\nenu
\itlb{it21}\msur
$\caL\w$ is an elementary submodel of $\rL[\w]$ 
\poo\ all $\is1{\nn+1}$ formulas with reals in 
$\caL\w$ as parameters.
\een

Now let $\vpi(p,k)$ be a $\is1{\nn+1}$ formula with 
some $p\in\caL\w\cap\dn$ as the only parameter. 
We are going to prove that the set 
$X=\ens{k}{\caL\w\mo\vpi(p,k)}$ 
belongs to $\caL\w.$  
By definition, $p\in\rL[v\dar\et]$ for some $\et\in\Omv$. 
Let 
\bce
%
$\Om_\et=\ens{\xi\in\cpo}{\et\osq\xi}$, \ \ all odd expansions 
of $\et$ in $\cpo$, 
\ece
and $\cbL\et\w=\bigcup_{\xi\in \Om_\et}\rL[\w\dar\xi]$. 
Note that $\Om_\et\sq\Omv$ by Lemma~\ref{b2}\ref{b22}, 
and $\Om_\et$ obviously satisfies the same property, that is, 
if\/ $\et\osq\xi$ belong to\/ $\cpo$  then\/ 
$\et\in\Om_\et\imp \xi\in\Om_\et$. 
Therefore, similarly to \ref{it21}, we obtain:

\ben
\nenu
\atc
\itlb{it22}\msur
$\cbL\et\w$ is an elementary submodel of $\rL[\w]$ --- 
and hence of $\caL\w$ as well by \ref{it21} --- 
\poo\ all $\is1{\nn+1}$ formulas with reals in 
$\cbL\et\w$ as parameters. 

\itlb{it23}
Hence in particular $X=\ens{k}{\cbL\et\w\mo\vpi(p,k)}$. 
\een

Note finally that unlike $\Omv$ the set $\Om_\et$ belongs 
to $\rL$, and is closed under countable unions. 
It follows that 
$\cbL\et\w\cap\dn=\rL[\w\dar\tux\w]\cap\dn$, 
hence the set $\cbL\et\w\cap\dn$ satisfies the full 
schema of $\CA$.  
It follows that $X\in \cbL\et\w\sq\caL\w$ 
by \ref{it23}, as required. 
\epf

\parf{Item 3: verification of the \paf\ Choice}
\las{it3}

\bpf
[item \ref{bt3} of Theorem~\ref{bt}]
This will be rather similar to the proof of Theorem~\ref{65} 
in the version of its last claim. 

To begin with, consider the subgroup $\Gav\in\rL$ of 
the group $\per$ of parity-preserving permutations $\pi$ 
of $\tup$ (Section~\ref{perm}) which consists of all 
$\pi\in\per$ such that, for each odd $\al$, if 
$\pi(\ang\al)=\ang\ga$ (also odd!) then 
$\pi(\ang{\al{-}1})=\ang{\ga{-}1}$. 

\ble
\lam{it31}
Let\/ $\w\in\can\tup$ be\/ \dd\cX generic over\/ $\rL$, 
and\/ $\pi\in\Gav$. 
Then\/ 
\begin{multicols}{2}
\ben
\renu
\itlb{it311}
$\pi\akt\w$ is \dd\cX generic over\/ $\rL$,
\itlb{it312}
$\tux{\pi\akt\w}=\pi\akt{\tux\w}$, 
\itlb{it313}
$\Omx{\pi\akt\w}=\pi\akt{\Omv}$, 
\itlb{it314}
$\caL\w=\caL{\pi\akt\w}$.
\een
\end{multicols}
\ele
\bpf[lemma]
\ref{it312}
Let $\w'=\pi\akt\w$, 
$\i\in\tup$, $\al=\i(0)$, $\j=\pi\akt\i$, $\al'=\j(0)$, 
so that $\ang{\al'}=\pi\akt\ang\al$. 
If $\al$ is even then so is $\al'$ (as $\pi$ is parity-preserving), 
and we have $\i\nin\tux\w$,   $\j\nin\tux{\w'}$. 
Thus suppose that $\al$ is odd. 

Then $\al'$ is odd too, and the even ordinals 
$\ga=\al{-}1$, $\ga'=\al'{-}1$ are defined and satisfy 
$\ga'=\pi\akt\ga$ since $\pi\in\Gav$, and moreover 
(I) $\w'(\ga')=\w(\ga)$. 
It remains to note that (II) if $1\le k<\lh\i=\lg\j$ then the 
ordinals $\i(k)$ and $\j(k)$ are both even or both odd. 
We conclude from (I),(II) that condition \ref{b1*} of 
Section~\ref{b1} holds for $\i,\w$ and $\j,\w'$ simultaneously, 
as required.

This completes the proof of \ref{it312}. 
The other two equalities \ref{it313}, \ref{it314} 
are easy corollaries.
\epF{lemma}

To begin the proof of the theorem, 
fix a \paf\ $\is1\iy$ formula $\vpi(k,x)$, and assume that 
(*) $\caL\w\mo\kaz k\,\sus x\,\vpi(k,x)$. 
By necessity, the arguments somewhat change \poo\ the 
proof of Theorem~\ref{65}. 
First of all, for any $\al\in\Ord$ and suitable set $z$, 
$\gor\al z$ will denote the $\al$th element of $\rL[z]$ 
\kmar{gor al z}
in the sense of the G\"odel well-ordering of $\rL[z]$. 
Then it follows from (*) that, in $\rL$, 
there exist sequences of conditions $X_k\in\cX$, 
ordinals $\al_k$, 
and sets $\xi_k\in\Omv$, satisfying 
\vyk{
In other words, for any $k$ there is a real 
$x_k\in{\cN}\cap \caL\w$ satisfying 
$\caL\w\mo\vpi(k,x_k)$,  
and hence there is $\xi_k\in\Omv$ such that 
$x_k\in\rL[\w\res\xi_k]$. 
In other words, 
\ben
\nenu
\itlb{x65*1}
$\caL\w\mo \sus x\in\rL[\w\dar\xi_k]\,
\vpi(k,x)$.
\een
}%
\ben
\nenu
\itlb{x65*1}
$X_k\fox\cX\,\big(\caL\pv\mo 
\vpi(k,\gor{\al_k}{\pv\dar\xi_k})\big)$ --- for all $k<\om$.
\een

Now assume to the contrary that 
$\caL\w\mo \neg \,\sus f\,\kaz k\,\vpi(k,f(k))$, and hence 
there exists a condition $X\in\cX$, satisfying 
\ben
\nenu
\atc
\itlb{x65*2}
$X\fox\cX\,\big( 
\caL\pv\mo \neg\, \sus f\,\kaz k\,\vpi(k,f(k))\big)$.  
\een

Let $\ta=\dym X$, $\ta_k=\dym{X_k}$. 
Arguing in $\rL$, we get   a sequence of permutations 
$\pi_k\in\Gav$ 
by induction, 
satisfying  $\vt_k\cap\vt_j=\vt_k\cap\ta=\pu$ 
whenever $k\ne j$, where 
$\vt_k=\pi_k\akt\ta_k\in\cpo$. 
Let $Y_k=\pi_k\akt X_k$, thus $Y_k\in\cX_{\vt_k}$. 
Let $\sg_k={\pi_k}\akt \xi_k$; $\sg_k\in\Omv$ 
by Lemma~\ref{it31}. 
Then \ref{x65*1} implies by Theorem~\ref{621}:
\ben
\nenu
\atc
\atc
\itlb{x65*3}
$Y_k\,\fox\cX \,\big(\caL{\pi_k\akt\pv}\mo 
\vpi(k,\gor{\al_k}{(\pi_k\akt\pv)\dar\xi_k})\big)$, 
\ \ 
$\kaz k<\om$, 
\een
Here $\caL{\pi_k\akt\pv}$ can be replaced by 
just $\caL\pv$ by Lemma~\ref{it31}\ref{it314}, 
whereas $(\pi_k\akt\pv)\dar\xi_k$ can be replaced  
by $\pi_k\akt(\pv\dar\sg_k)$. 
This implies
\ben
\nenu
\atc
\atc
\atc
\itlb{x65*4}
$Y_k\,\fox\cX \,
\big(\caL{\pv}\mo 
\vpi(k,\gor{\al_k}{\pi_k\akt(\pv\dar\sg_k)})\big)$, 
\ \ 
$\kaz k$. 
\een

Now let $\vt=\bigcup_k\vt_k$.  
Then the set $Y=\bigcap_k(Y_k\uar\vt)$ belongs to 
$\cX\dar\vt$ 
by Lemma~\ref{rfoL} 
(\poo\ Lemma~\ref{pe4}). 
As obviously $Y\psq Y_k$, \ref{x65*4} implies:
\ben
\nenu
\atc
\atc
\atc
\atc
\itlb{x65*5}
$Y\fox\cX\,
\big(\caL{\pv}\mo 
\vpi(k,\gor{\al_k}{\pi_k\akt(\pv\dar\sg_k)})\big)$. 
\een

Now follows the key step. 
The set $\sg=\bigcup_k\sg_k$ belongs to $\Omv$ 
because so does each $\sg_k=\pi_k\akt\xi_k$.
The term  
$\gor{\al_k}{\pi_k\akt(\pv\dar\sg_k)}$ 
in \ref{x65*5}, as a function of $k$ and $\pv\dar\sg$, 
is defined in  $\rL[\pv\dar\sg]$ by an absolute formula 
with parameters $k\mto\al_k$,  $k\mto\pi_k$, 
$k\mto\sg_k$ 
(all three maps belong to $\rL$ by construction). 
Therefore the map 
$f(k)=\gor{\al_k}{\pi_k\akt(\pv\dar\sg_k)}$ 
is forced by $Y$ to belong to $\rL[\pv\dar\sg]$. 
We conclude that 
\ben
\nenu
\atc
\atc
\atc
\atc
\atc
\itlb{x65*6}
$Y\fox\cX\,\sus f\,\kaz k<\om\,
\big(\caL{\pv}\mo \vpi(k,f(k))\big)$. 
\een
Thus conditions $Y$ and $X$ force contradictory 
statements by \ref{x65*2}. 
Yet  $\dym Y\cap\dym X=\sg\cap\ta=\pu$ by 
construction, which implies that $Y$ and $X$ are   
compatible in $\cX$. 
This is a contradiction.
\epF{item \ref{bt3} of Thm~\ref{bt}}

\qeD{Theorem~\ref{bt} and Theorem~\ref{mt2}}

\sekt
[\ \ Final remarks and questions]
{Final remarks and questions}
\las{frq}

In this final chapter, we begin with an explanation in 
Section~\ref{why} as of the principal necessity of the 
{\ubf separate treatment of even and odd tuples} 
in the proof of Theorem~\ref{mt1}.
Then we outline some {\ubf further applications} 
of our methods in Section~\ref{furth},  
discuss the possibility to prove the results like 
our Theorems~\ref{mt1} and \ref{mt2} on the basis 
of second order arithmetic in Section~\ref{70},
and finish with conclusive remarks and a 
commented list of problems in Section~\ref{rq}.

\parf{Why the even/odd distinction?}
\las{why}

One may ask whether a simpler version of the 
construction, which merges the even and odd tuples under 
common treatment, will not work. 
For instance, redefine $\Oma'$ to be the set of all $\xi\in\cpo$ 
in $\rL$ such that $\sus m\,\kaz\i\in\xi\,(\lh\i\le m)$.
(Compare with the actual definition of $\Oma$ 
in Definition~\ref{71}.)  
The following argument, presented rather tentatively, 
shows that this does not go towards the proof of 
Theorem~\ref{mt1}\ref{mt11}. 

To begin with, 
consider the whole $\pei\in\RF$ as the forcing notion. 
It has the Fusion property by Theorem~\ref{pfp} and is 
\dd\nn complete for every $\nn$. 

Let $\w\in\can\tup$ be an array \dd\pei generic over $\rL$. 
It is known from the studies of generalized iterated Sacks 
extensions (see \eg\ \cite{kl69}) that if $m\ge1$ then 
the \dd\rL degrees of the reals of the form $x=\w(\i)$, 
where $\i\in\tup$, $\lh\i=m$, can be described by 
the following \paf\ $\is14$ formula in $\rL[\w]$: 
\bde
\item[$D(m,x):=$]
there is exactly $m$ different \dd\rL degrees of reals 
strictly below 
$x$, and those degrees are linearly ordered by the relation 
$a\lec b$ iff $a\in\rL[b]$.
\ede
In other words, we have, for all $m\ge1$, 
\ben
\nenu
\itlb{why1}
$\rL[\w]\mo\kaz m\,\kaz x\in\cN 
\big(D(m,x)\eqv 
{\sus \i\in\tup\,
(\lh\i=m\land x\approx \w(\i))} 
\big)$,
\een
where $a\approx b$ iff $a\lec b\land b\lec a$, 
and $\cN=\bn.$   

\bte
\lam{whyt}
If\/ $\cX\in\RF$ is a forcing with the Fusion 
and \dd{5}Completeness properties, and\/ 
$\w\in\can\tup$ is an array \dd\cX generic over $\rL$, 
then\/ \ref{why1} holds in\/ $\rL[\w]$. 
\ete

Before the proof starts, we may note that, 
by rather standard arguments \eg\ in \cite{kl69}, 
\ref{why1} implies that $\xAC{\is14}$ {\ubf fails} 
in $\rL(\W{\Oma'}[\w])$, and hence 
$\rL(\W{\Oma'}[\w])$ is {\ubf not} a model for 
Theorem~\ref{mt1}\ref{mt11} 
for any $\nn\ge6$. 
Thus to complete the task in this Section, it suffices 
to prove the theorem.
 

\bpf
[Theorem~\ref{whyt}] 
First of all, Corollary~\ref{1011} implies that \ref{why1} 
is equivalent to the following claim: 
if $f\in\komt$ in $\rL$ and $m\ge1$ then
\ben
\nenu
\atc
\itlb{why2}
$\rL[\w]\mo
\big(D(m,\vva f\w)\eqv 
\underbrace{\sus \i\in\tup\,
(\lh\i=m\land \vva f\w\approx \w(\i))}_{\Phi} 
\big)$. 
\een
(See Section~\ref{103} on the  valuation $\vva\vpi\w$.)
Thus we have to prove \ref{why2} for all \dd\cX generic $\w$ 
{\ubf provided it holds for all \dd\pei generic $\w$}. 
We fix $f\in\komt$ and $m\ge1$ in the course of our 
arguments, and let $\xi=\modd f$, so that 
$\xi\in\cpo$ in $\rL$ and $f\in \kont\xi$.

\vyk{
We would like to infer that \ref{why1} and hence the 
subsequent $\xAC{\is14}$-failure claim, 
are also forced by any forcing $\cX\in\RF$ with the Fusion 
and \dd{\nn}Completeness properties, for all $\nn$ 
big enough.
Unfortunately, the quantifier $\sus\i$ and the related 
occurrence of $\w(\i)$ do not allow us to immediately 
apply Theorem~\ref{122}, since \ref{why1} is just not 
a formula of that type. 

To circumvent this obstacle, we need some extra work. 
}

{\ubf Arguing in $\rL$}, 
fix 
a surjection (not necessarily $1{-}1$) 
$s:\om\onto\xi_m=\ens{\i\in\xi}{\lh\i=m},$  and 
define a map $H:\can\xi\to\cam$ such that 
$(H(u))_k=u(s(k))$ for all $u\in\can\xi$ 
and $k$. 
(Recall that $(x)_k(l)=x(2^k(2l+1)-1)$, $\kaz l$.) 
Thus 
$H$ is a continuous map 
$\can\xi\to\cam$, and hence  
the restriction 
$h=H\res\rat\xi$ 
belong to $\kont\xi$ and satisfies $H=\clo{h}$.
(See Section~\ref{101} on the notation involved.)

Now consider the formulas $\Phi$ as in \ref{why2} and 
$\Psi:= \sus k\,(f\approx (h)_k)$,  
an $\xla$-formula with $f,h$ as 
the only parameters in $\kont\xi$.
See Section~\ref{103} on $\xla$-formulas. 

\ble
\lam{why3}
Under the assumptions of the theorem, 
$\rL[\w]\mo\big(\Phi\eqv \vva\Psi\w\big)$.
\ele
\bpf[\rm sketch]
Note that $\vva\Psi\w$ is essentially 
$\sus \i\in\xi\,(\lh\i=m\land \vva f\w\approx \w(\i))$ 
by the definition of $H$, and this is $\Phi$. 
\epf
%
\vyk{
Thus the whole right-hand side of the equivalence is
$$
\kaz\xi\in\cpo\,
\big[
x\in\rL[\w\dar\xi]\imp 
\sus \i\in\xi\,\big(\lh\i=k\land x\approx \w(\i)\big)
\big].
$$
It remains to note that any $x\in\dn\cap\rL[\w]$ satisfies 
$x\in\rL[\w\dar\xi]$ for some $\xi\in\cpo$ as $\w$ is 
$\omil$-preserving.
}

\bcor
\lam{why4}
If\/ $\w\in\can\tup$ is\/ \dd\pei generic over\/ $\rL$  
then\/ 
$\vva{\big(D(m,f)\eqv \Psi\big)}\w$ 
holds in\/ $\rL[\w]$. 
In other words, $\pei$ forces\/ 
$\vva{\big(D(m,f)\eqv \Psi\big)}\w$.
\ecor
\bpf[\rm sketch]
Apply the assumption that \ref{why2} holds for such a $\w$.
\epf

Now recall that $\lec$ is a $\is12$ relation by 
Addison~\cite{add2,add1}, 
hence such is $\approx$ as well. 
Therefore $\Psi$ is a $\ls12$ \dd\xla formula, and hence 
the equivalence $D(m,f)\eqv \Psi$ in brackets is essentially 
a $\ls15$ \dd\xla formula 
with $f,h\in \kont\xi$ as 
the only parameters. 
(Recall that $D$ is $\is14$.)  
We conclude (omitting details) that  
$\vva{\big(D(m,f)\eqv \Psi\big)}\w$ 
is also forced by $\cX$ --- 
by Theorem~\ref{122} and the assumption that $\cX$ 
is \dd5complete. 
It follows by Lemma~\ref{why3} that $\cX$ forces 
$D(m,\vva f\w)\eqv \Phi$, 
which is \ref{why2}, as required.
\epF{Theorem}

\vyk{
This failure is naturally expressed by a conjunction $\gA$ 
of a $\ip16$ formula and a $\ip16$ formula. 
(Note that some usual laws of transformation of analytic 
formulas do not work in non-$\AC$ environment;  
this is why we are getting to level $6$ here.) 
We conclude that $\pei$ {\ubf forces} $\gA$, in other words, 
forces the failure of $\xAC{\is14}$ in 

Now let $\nn>6$  and consider a normal forcing 
$\cX\in\RF$ with the Fusion, 
Structure, \dd\nn Definability, and\/ \dd\nn Completeness  
properties, as in Theorem~\ref{mt1c}. 
}

\parf{Some further results}
\las{furth}

\noindent{\ubf Applications to Separation problem.}\vom

The Separation problem was known in the early years of 
descriptive set theory. 
The following two statements are considered:
\bde
\item[\Sep{\fs1n}:] 
any two disjoint $\fs1n$ sets 
(in the same Polish space) are separated by a $\fd1n$ 
set. 

\item[\Sep{\fp1n}:]
any two disjoint $\fp11$ sets 
are separated by a $\fd1n$ 
set. 
\ede
Luzin~\cite{lus:ea} established \Sep{\fs11}. 
Novikov~\cite{nov1931} proved 
\nSep{\fp11}, so that there exist two disjoint $\fp11$ sets 
not separated by a $\fd11$ set. 
Then Novikov~\cite{nov1935} demonstrated that, 
on the second projective level, on the contrary, 
we have  \Sep{\fp12} but \nSep{\fs12}. 
See monographs of Kechris \cite{Kdst} and 
Moschovakis \cite{mDST} on the modern treatment of 
Separation.

As for the higher projective classes, Addison~\cite{add1,add2} 
proved that, assuming the axiom of constructibility $\rV=\rL$, 
if $n\ge3$ then we have  \Sep{\fp1n} but \nSep{\fs1n}, 
that is, similar to level $n=2$. 
As usual in such cases, a problem has been raised 
of building generic models in which, on the contrary, 
we have \nSep{\fp1n} and/or \Sep{\fs1n},  see  \eg\ 
an early survey \cite{matsur} by Mathias. 

Working on the first part of this problem, 
Harrington sketched a model for \nSep{\fp1n}, 
based on the almost-disjoint forcing of \cite{jsad},  
in his handwritten notes \cite{h74} 
(never published, but rather known to set theorists). 
We succeeded to implement Harrington's idea on the basis 
of product Jensen forcing in \cite{kl49}. 
The next theorem outlines another model for \nSep{\fp1n},  
based on the technique developed in this paper.

For $e=0,2,4$ let 
$A_e=\ens{\la+6k+e}
{\la<\omi\text{ limit}\land k<\om}$. 

Let $J$ consist of all tuples $\i\in\tup$ such that 
1) if $\lh\i\ge2$ and $\i(0)\in A_0$ then $\i(1)$ is odd, 
and 2) if $\lh\i\ge3$ and $\i(0)\in A_2$ then $\i(2)$ is odd. 

\bte
\lam{nonsep}
Assume that\/ $\nn\ge1$ and, in\/ $\rL$, $\cX\in\nf$ is a 
normal forcing satisfying four conditions of 
Theorem~\ref{mt1b} for this\/ $\nn$. 
Let\/ $\w\in\can\tup$ be an array\/ \dd\cX generic over\/ $\rL$. 
Then\/ $\rL[\w\dar J]\mo\ZFC$ is a model of\/  
\nSep{\fp1{\nn+2}}. 

To be more exact, it holds in\/ $\rL[\w\dar J]$ that the sets
$$
H_0 =\ens{\w(\ang\al)}{\al\in A_0}
\qand\/
H_2 =\ens{\w(\ang\al)}{\al\in A_2}
$$
are disjoint\/ $\ip1{\nn+2}$ sets non-separable by a\/ 
$\fd1{\nn+2}$ set.
\ete

The proof will appear elsewhere. 
To prove that $H_0,H_2$ belong to $\ip1{\nn+2}$ 
in $\rL[\w\dar J]$, we use sets $P_1,P_3$ from the proof 
of Theorem~\ref{82}, show that 
$$
\bay{rcl}
H_0 &=&\ens{x}
{\ang{1,x}\in P_1 \land \kaz y\,
(\ang{3,y}\in P_1\imp \ang{x,y}\nin P_3)};\\[0.8ex]
H_2 &=&\ens{x}
{\ang{1,x}\in P_1 \land \kaz y\,
(\ang{2,y}\in P_1\imp \ang{x,y}\nin P_3)};
\eay
$$
and apply Lemma~\ref{83}.
The proof of non-separability is more involved.  

The Reduction principle (Kuratowski~\cite{kursep}) 
for a class $K$ is as follows:
\bde
\item[\Red{K}:] 
any two sets $X,Y$ in $K$ contain subsets 
$X'\sq X$, $Y'\sq Y$, still in $K$, such that 
$X'\cap Y'=\pu$ and $X'\cup Y'=X\cup Y$. 
\ede
It is known that \Red K implies Separation \Sep{\dop K} 
for the complementary class $\dop K$ 
(containing all complements of sets in $K$), and 
accordingly \Red{\fp11} and \Red{\fs12} hold 
whereas \Red{\fs11} and \Red{\fp12} fail, and under 
$\rV=\rL$ \Red{\fs1n} hold and \Red{\fp1n} fail for 
any $n\ge3$. 
(See \cite{Kdst} for a full account of related results.) 

\bvo
\lam{sepred} 
Let $n\ge3$. 
Define a model in which \Sep{\fs1n}, or stronger, 
\Red{\fp1n} holds. 
Define a model in which \Sep{\fp1n} holds but 
\Red{\fp1n} fails. 
\evo

See a new interesting approach in a preprint \cite{hoff} 
on these questions. 
\vtm

\noindent
{\ubf Applications to the Uniform Projection problem.}\vom

By definition, a set $X$ in the Baire space ${\cN}=\bn$ 
belongs to $\fs1{n+1}$ 
iff it is equal to the projection 
$\dom P=\ens{x}{\sus y\,P(x,y)}$ of a ``planar'' 
$\fp1n$ set $P\sq{\cN}\ti{\cN},$  so that in breef 
$\fs1{n+1}$ = projections of $\fp1n$.
In particular, this is true for $n=0$; by definition, 
$\fp10\:$= all closed sets $P\sq{\cN}\ti{\cN}$. 

The picture drastically changes if we consider only 
\rit{uniform} sets $P\sq{\cN}\ti{\cN},$ 
\index{set!uniform}%
\ie, those satisfying 
${P(x,y)\land P(x,z)}\imp {y=z}$. 
It was established by Luzin \cite{lus:ea} that 
projections of uniform $\fp10$ sets, and even 
uniform $\fd11$ (that is, Borel) sets in ${\cN}\ti{\cN}$ 
are $\fd11$, 
which is a proper subclass of $\fs11$, 
and the other way around, every $\fd11$ set in ${\cN}$ 
is the projection of a uniform $\fp10$ set. 
On the contrary, the Novikov -- Kondo uniformization 
theorem \cite{luno,kond} asserts that every $\fp11$ 
set $P\sq{\cN}\ti{\cN}$ contains a uniform $\fp11$ 
subset $Q\sq P$ with $\dom Q=\dom P$, and hence 
\bce
$\fs12$ = projections of $\fp11$ = projections of 
uniform $\fp11$.
\ece
(See Luzin~\cite{lhad2} or Moschovakis~\cite[4F]{mDST} 
on uniformization of $\fs11$ sets.)

Even before the the Novikov -- Kondo uniformization, Luzin 
raised several problems in \cite[pp.\ 274-276]{lbook}, 
the general content of which was a comparison of the class 
$\ufp1n$ of projections of uniform 
$\fp1n$ sets $P\sq{\cN}\ti{\cN}$ with classes 
$\fs1{n+1}$,  $\fd1{n+1}$,  $\fs1{n}$.  
The following two theorems show that, for each $n\ge2$, 
both $\fs1{n+1}=\ufp1n$ and 
$\fs1n\not\sq \ufp1n$ 
(or $\fd1{n+1}\not\sq \ufp1n$ 
in case $n=2$) 
are statements consistent with $\zfc$. 

\bte
\lam{up1}
Assuming\/ $\rV=\rL$, we have\/ 
$\fs1{n+1}=\ufp1n$ 
for all\/ $n\ge 2$, albeit there is no uniformization 
theorem for\/ $\fp1n$ similar to the 
Novikov -- Kondo uniformization theorem for\/ $\fp11$.
\ete

\bte
\lam{up2}
Under the assumptions of Theorem~\ref{nonsep}, 
we have\/ 
$\fs1{\nn+2}\not\sq \ufp1{\nn+2}$
and\/ 
$\fd1{\nn+2}\not\sq \ufp1{\nn+1}$ 
in appropriate submodels of\/ $\rL[w]$.
\ete

The proofs will appear elsewhere. 
We may note that the 2nd non-inclusion of the last theorem 
follows from the 1st non-inclusion for $\nn{-}1$ instead of $\nn$. 
Yet this reduction leaves aside the case $\nn=1$ in the 
2nd non-inclusion because the 1st one is false for $\nn=0$.

\parf{Working on the basis of the consistency of $\pad$}
\las{70}

The main results of this paper, 
Theorems \ref{mt1} and \ref{mt2}, can be naturally 
viewed as formal consistency results related to 
certain subsystems of second order Peano arithmetic $\pad$ 
and obtained by means of forcing technique 
and other tools of $\zfc$ which go way beyond $\pad$ itself.
Therefore it is usually a tempting problem in such cases to 
reproduce the consistency results obtained on the basis 
of $\cons\pad$, the formal consistency of $\pad$.  

Such a reproduction of another result, the consistency 
of the assertion $\wo_{n}\land \neg\wo_{n-1}$, 
based of the consistency of $\pad$, where 
\bde
\item[\ \ \ \ $\wo_n$:]
\it
there is a wellordering of the reals of class $\id1{n}$,
\ede
has been recently achieved, for any given $n\ge3$, 
by adapting the proof of the consistency of 
$\wo_n\land \neg\wo_{n-1}$ with $\zfc$ 
in an earlier paper \cite{kl67}.

The adaptation of this $\zfc$-based proof to $\pad$ was 
carried out in \cite{kl73}.   
There we utilize  $\zfcm$, a subtheory of $\zfc$ obtained by 
removing the Power Set axiom and some changes in other 
axioms, as a proxy theory.  
(See \eg\ \cite{gitPWS} for a comprehensive account 
of $\zfcm$.) 
The advantage of $\zfcm$ is that this theory is 
equiconsistent with $\pad$, while it 
is still  a rather forcing-friendly theory. 
The equiconsistency of $\zfcm$ and $\pad$ is considered to 
be a well-known result, although, as far as we know, no 
complete proof has ever been published. 
A sketch given in \cite{kl73} involves  
some results of  \cite{aptm,
Kr} 
and other earlier papers.

On the other hand, $\zfcm$ allows to adapt many typical 
forcing notions related to reals, in the form of 
\rit{pre-tame class forcings}, based on appropriate coding 
of the ``continual'' forcing conditions by real-like objects, 
and the general class forcing theory 
set up in \cite{sdffs,AF,AG}. 
Such an adaptation contains a lot of routine 
(but nevertheless time and space consuming) work. 
In addition, regarding the $\zfcm$-adapted proof  in 
\cite{kl73}, there are two non-routine issues. 
Firstly, this is getting rid of countable transitive models, 
of theories similar to $\zfcm$, 
in evaluation of the definability level of some 
constructions, as in Theorem~\ref{193ta} above. 
Secondly, circumventing the use of diamond, 
which is definitely not 
a $\zfcm$ result in its common formulation and proof. 
Note that the requirement of cardinal-preservation  
of the forcing notion considered in the $\zfc$ setting 
is a {\it conditio sine qua non\/} for such an adaptation, 
because generic collapse of cardinals   
is definitely beyond the formal consistency of $\pad$.

Anyway, we were able to overcome these difficulties in 
\cite{kl73} and prove the consistency of 
$\wo_{n}\land \neg\wo_{n-1}$ 
(for any given $n\ge3$) with $\pad$, based on the 
consistency of $\pad$ itself (equivalently,  of $\zfcm$). 
Metamathematically, this means that  $\cons\pad$ implies 
$\cons(\pad+\wo_n+\neg\wo_{n-1})$. 

The methods developed in \cite{kl73} 
(and in \cite{kl75} with respect to another problem) 
are also applicable 
to the main results of this article 
(Theorems \ref{mt1} and \ref{mt2}). 
Adapting their proofs, we are able to establish the following 
form of our main results: 

\bte
[1st main theorem for $\pad$]
\lam{mt1'}
Assume that\/ $\nn\ge1$. 
Then\/ $\cons\pad$ implies the consistency of the 
following theories$:$  
\ben
\nenu 
\itlb{mt1'1} 
$\pad + \xDC{\fp1\nn}+ \neg\xAC{\ip1{\nn+1}}\,;$ 

\itlb{mt1'2} 
$\pad +\xAC\od+\xDC{\ip1{\nn+1}}+\neg\xAC{\fp1{\nn+1}}\,;$ 

\itlb{mt1'3} 
$\pad +\AC+\xDC{\fp1{\nn}}+\neg\xDC{\ip1{\nn+1}}\,;$ 

\itlb{mt1'4} 
$\pad +\AC+\xDC{\ip1{\nn+1}}+\neg\xDC{\fp1{\nn+1}}\,.$ 
\qed
\een

\ete

\bte
[2nd main theorem for $\pad$]
\lam{mt2'}
Assume that\/ $\nn\ge1$. 
Then\/ $\cons\pad$ implies\/ 
$\cons
(\pao + \xAC{\is1\iy} + \xCA{\fs1{\nn+1}}+ 
\neg\xCA{\fs1{\nn+2})}$. 
\qed
\ete

The details will appear elsewhere.

Identifying theories with their deductive 
closures, we may present the concluding statement  
of Theorem \ref{mt2'} as follows:
\bce
$\pao+\xAC{\is1\iy}+\xCA{\fs1{\nn+1}}
\;\;\sneq\;\; 
\pao+\xAC{\is1\iy}+\xCA{\fs1{\nn+2}}$.
\ece
Studies on subsystems of $\pad$ have discovered many 
cases in which $S\sneq S'$ holds for a given pair 
of subsystems $S,S'$, see \eg\ \cite{simp}. 
And it is a rather typical case that such a strict 
extension is established by demonstrating that $S'$ 
proves the consistency of $S$.
One may ask whether this is the case for the 
result in the displayed line above. 
The answer is in the negative: 
namely 
\bce
\rit{the theories\/ $\pao+\xAC{\is1\iy}$ 
and the full\/ $\pad$ are equiconsistent\/} 
\ece
by a result in \cite[Lemma 3.1.7]{HFuse81}, 
also mentioned in \cite{schindt}. 
This equiconsistency result also follows from 
a somewhat sharper 
theorem in \cite[1.5]{Schm}.

\parf{Remarks and questions}
\las{rq}

In this study, the technique of countable-support generalized 
iterations of Jensen forcing, combined with the~method of 
definable generic forcing notions,  
was employed to the construction of models of $\zf$ and $\pad$ 
with different effects related to the Choice and Comprehension 
axioms. 
The main results obtained show that the strength of a 
Choice or Comprehension principle naturally depends on the 
next three factors in essential way:\vom

1) the type of the principle considered: 
$\CA$, $\AC$, or $\DC$;\vom

2) the level considered in the projective hierarchy, \vom

3) admission or non-admission of parameters.\vom   

\noi
These results (Theorems \ref{mt1} and \ref{mt2}) 
are significant strengthening of previously known results 
in this area, including our own earlier results in 
\cite{kl69,kl71},
especially with regard to the transfer of 
ensuing independence results to an arbitrary 
level of the projective hierarchy.
These are new results and valuable improvements upon 
much of known independence results in this area.  
The technique developed in this paper may lead 
to further progress in studies of different aspects 
of the projective hierarchy.

This theorem 
continues our series of recent research such as
\bit
\vyk{
\item[$-$]
a $\ip1\nn$ \dd\Eo equivalence class containing no 
$\od$ elements, while 
every countable \dd{\is1{\nn}}set of reals 
contains only 
$\od$ reals \cite{kl34}, 
}%

\item[$-$]
a $\varPi^1_\nn$ real singleton $\ans a$ such that 
$a$ codes a cofinal map 
$f:\omega\to\omega_1^\rL$,
while every $\varSigma^1_\nn$
set $X\sq\om$ is constructible and hence cannot code 
a cofinal map $\omega\to\omega_1^\rL$, \cite{kl36},

\item[$-$]
a non-\ROD-uniformizable $\ip1{\nn}$ set with  
countable cross-sections, 
while all $\fs1\nn$ sets with countable  
cross-sections are $\fd1{\nn+1}$-uniformizable 
\cite{kl38},

\item[$-$]
a model of $\ZFC$, in~which the separation theorem 
fails for a given projective class $\fp1n$, 
is defined in~\cite{kl49};

\item[$-$]
a model of $\ZFC$, in~which the full basis theorem holds 
in the absence of analytically definable 
well-orderings of the reals, is defined in~\cite{kl29}.
\eit

These results also bring us closer to solving the following 
extremely important problem by S.\:D.\:Friedman   
{\cite[P. 209]{sdffs},  
\cite[P. 602]{sdf_ccf}}: 
assuming the consistency of an inaccessible 
cardinal, find a model for a given $n$ in which all 
$\fs1n$ sets of reals are Lebesgue measurable 
and have the Baire and perfect set properties, 
but there is a $\fd1{n+1}$ 
\weo\ of the reals. 

From our study, it is concluded that the technique of 
\rit{definable generic} inductive constructions  
of forcing notions in $\rL$, developed for 
Jensen-type generalized  forcing iterations, 
succeeds to solve important descriptive set theoretic problems. 

We present several questions related to possible 
extensions of the results achieved in this paper, 
that arise from our study.

\bvo
\lam{24} 
Recall that $\xDC \od \eqv\xDC \ROD $ by 
Lemma~\ref{23}\ref{236}.
Is the ordinal-definable 
$\xDC \od $ equivalent to the full $\DC$ in $\zf$? 
\evo

\bvo
\lam{pro1} 
Still about the Dependent Choices principle.   
Three different forms of this axiom were introduced 
by Definition~\ref{02}: 
$\xDC K$, $\xDCm K$, $\xDCs K$. 
Lemma~\ref{23} contains several results on the 
relationship of these forms of $\DC$ to each other. 
But still some questions remain unresolved. 
For instance, consider the implications 
$\xDCs K\imp \xDC K\imp \xDCm K$
in Lemma~\ref{23}\ref{231}. 
The first implication is actually an equivalence for 
appropriate classes $K$ by Lemma~\ref{23}\ref{234}.
What about $\xDC K\imp \xDCm K$, the second one?
Can we split it by suitable models, 
provided $K=\ip1n$ or $\fp1n$?
\evo

\bvo
\lam{pro-} 
Does the implication  
$\xDCm{\ip1{n+1}}\imp\xDCm{\fp1{n}}$ 
hold, similarly to \ref{235} of Lemma~\ref{23}?
\evo

\bvo
[Communicated by Ali Enayat]
\lam{84ae}
A natural question is whether the main results of this paper  
also hold for second order set theory 
(the Kelley-Morse theory of classes). 
This may involve a generalization of the Sacks forcing 
to uncountable cardinals, as in \cite{fg,kanam}, 
as well as the new models of set theory recently defined 
by Fuchs~\cite{fu}, on the basis of further development of 
the methods of \rit{class forcing} introduced by 
S.\,D.\,Friedman~\cite{sdffs}. 
\evo

\bvo
\lam{84vik}
Another natural question is whether the main results 
of this paper (Theorems \ref{mt1} and \ref{mt2}) can 
be achieved on the basis of the 
{\em finite-support\/} generalized iterations of the Jensen 
forcing, developed in \cite{kl27,jml19}. 
Unlike the {\em countable-support\/} approach, used 
in this paper, these 
iterations preserve CCC but generally do not allow to use 
the splitting/fusion technique.
\evo

Now we return to the result on consistency of hypothesis 
$\wo_{n}\land \neg\wo_{n-1}$, 
discussed in Section~\ref{70}.
The generic model used to prove this consistency claim in 
\cite{kl67} definitely satisfies the continuum 
hypothesis $2^{\alo}=\ali$. 
The problem of obtaining models of \ZFC\  
in which $2^{\alo}>\ali$  
and there is a projective \weo\ of the 
real line, has been  
known since the early years of modern set theory. 
See, \eg, problem 
3214 in an early survey \cite{matsur} by Mathias. 
Harrington \cite{hlong} solved it by 
getting a generic model of $\ZFC$, 
in which $2^{\alo}>\ali$  
and there is a $\id13$ \weo\ of the continuum. 
This model involves various forcing notions 
like the almost-disjoint 
forcing \cite{jsad} and a forcing notion 
by Jensen and Johnsbr\aa ten \cite{jj}.

\bvo
\lam{pro2}
Prove the consistency of  $\wo_{n}\land \neg\wo_{n-1}$ 
by a model satisfying the requirement 
that the negation $2^{\alo}>\ali$ of the 
continuum hypothesis holds. 
\evo

Finally, a very recent paper \cite{ww} presents another 
study of interrelations between various forms of Choice from 
somewhat different point of view. 
In particular Theorem in \cite[page 5]{ww} claims a 
model of 
\bce
$\zf \:+\: 
\DC(\dR,\fp1n)\:+\:
\neg\AC(\dR,\text{unif\,}\ip1{n+1})
\:+\:\neg\AC(\dR,\text{\ubf Ctbl})$ 
\ece
for any $n\ge1$, where:
\bde
\item[$\DC(\dR,\fp1n)$] 
asserts that if $\pu\ne X\sq{\cN}$ is a $\fp1n$ set and 
\kmar{DC(RK)}%
\index{axiom!$\DC(\dR,\fp1n)$}%
\index{zDCRP@$\DC(\dR,\fp1n)$}%
$P\sq X\ti X$ is a $\fp1n$ relation with $\dom P=X$, 
then there is a chain $\sis{x_k}{k<\om}$ of reals 
$x_i\in X$ satisfying $x_k\mathrel{P}x_{k+1}$ 
for all $k$ --- this is equivalent to our $\xDC{\fp1n}$ by 
Lemma~\ref{23}\ref{234};

\item[$\AC(\dR,\text{\rm unif\,}\ip1{n+1})$] 
asserts that if $\pu\ne X_k\sq{\cN}$ are sets in $\ip1{n+1}$ 
and the set $\ens{k\we x}{k<\om\land x\in X_k}$ belongs 
to $\ip1{n+1}$ as well --- equivalent to our 
$\xAC{\ip1{n+1}}$ 
as in Definition~\ref{02};

\item[$\AC(\dR,\text{\ubf Ctbl})$] 
asserts that any family of countable or finite 
sets $\pu\ne X_k\sq{\cN}$ admits a choice function --- note 
that in $\zf$ the union $\bigcup_kX_k$ is not necessarily 
countable, and the set 
$\widehat X=\ens{\ang{k,x}}{k<\om\land x\in X_k}$
is not necessarily even $\ROD$, in this case under $\zf$.
\ede 

\vyk{
\bvo
\lam{prow}
Study the relations between the $\DC$ versions as in 
Definition~\ref{02} and $\DC(\dR,\fp1n)$ as above. 
\evo
}

\bvo
\lam{proww}
Find out whether axiom $\AC(\dR,\text{\ubf Ctbl})$ 
as above 
is fulfilled in the models that are built to prove our 
Theorem~\ref{mt1}. 
\evo

It should be noted that, 
when dealing with $\AC(\dR,\text{\ubf Ctbl})$ in 
the choice-less environment of $\zf$, 
the behavior of countable sets can be different from what 
one is accustomed with in ordinary mathematics, see \eg\ 
\cite{milLB,milDF}.


\back
The authors are thankful to Ali Enayat, Gunter Fuchs,  
Victoria Gitman, and Kameryn Williams, 
for their enlightening comments that made it possible to 
accomplish this research, and separately Ali Enayat for  
references to \cite{HFuse81,schindt,Schm} in 
matters of Theorem~\ref{mt2}. 
\vyk{
The authors are thankful to 
the anonymous referee for their 
comments and suggestions, which significantly 
contributed to improving the quality of the publication.
}
\eack


\np

\renek{\refname}
{{References}\addcontentsline{toc}{section}{References}}
\small

\bibliographystyle{plain} 

\bibliography{74d.bib,74dkl.bib}

\index{zzz @\ \hspace*{15ex}{\ubf blank line intentional}}%
\index{zz @\ \hspace*{15ex}{\ubf Greek Index}}%

\renek{\indexname}
{{Index}\addcontentsline{toc}{section}{Index}}

\printindex


\end{document}